\pdfoutput=1
\documentclass[10pt,reqno]{amsart}

\numberwithin{equation}{section}

\usepackage[margin=1in]{geometry}
\usepackage{cancel}
\usepackage[tt=false]{libertine}
\usepackage{mathtools}
\usepackage{amssymb}
\usepackage{mathrsfs}
\usepackage[varbb]{newpxmath}

\let\savedbigtimes\bigtimes
\let\bigtimes\relax
\usepackage{mathabx} 
\let\bigtimes\savedbigtimes

\usepackage{textcomp}
\usepackage{enumerate}
\usepackage{bm}

\usepackage[usenames,dvipsnames]{xcolor}
\usepackage[colorlinks=true,
	linkcolor=green!50!blue!50!black,
	citecolor=PineGreen,
	urlcolor=RedViolet]{hyperref}
\hypersetup{bookmarksopen=true}
\usepackage[footnotesize,
	justification=centering]{caption}
\usepackage[font=footnotesize]{subcaption}
\usepackage{graphicx}

\newtheorem{thm}{Theorem}[section]
\newtheorem{lem}[thm]{Lemma}
\newtheorem{ppn}[thm]{Proposition}
\newtheorem{cor}[thm]{Corollary}

\theoremstyle{definition}
\newtheorem{dfn}[thm]{Definition}
\newtheorem{rmk}[thm]{Remark}
\newtheorem*{rmk*}{Remark}

\newtheorem{ass}{Assumption}

\definecolor{Green}{HTML}{239B56}
\definecolor{Blue}{HTML}{2471A3}
\definecolor{Orange}{HTML}{D35400}
\definecolor{DarkBlue}{HTML}{1B4F72}
\definecolor{DarkGray}{HTML}{1C2833}

\newcommand{\beq}{\begin{equation}}
\newcommand{\eeq}{\end{equation}}

\newcommand{\<}{\langle}
\renewcommand{\>}{\rangle}

\newcommand{\f}{\frac}
\newcommand{\set}[1]{\{#1\}}
\newcommand{\Ind}[1]{\mathbf{1}\{#1\}}
\newcommand{\I}{\mathbf{1}}
\newcommand{\pd}{\partial}

\DeclareMathOperator{\Cov}{Cov}
\DeclareMathOperator{\Var}{Var}
\DeclareMathOperator{\Th}{th}
\DeclareMathOperator{\Ch}{ch}
\DeclareMathOperator{\Hess}{Hess}
\DeclareMathOperator{\diag}{diag}
\DeclareMathOperator{\tr}{tr}
\DeclareMathOperator{\sgn}{sgn}

\newcommand{\E}{\mathbb{E}}
\newcommand{\R}{\mathbb{R}}
\renewcommand{\P}{\mathbb{P}}
\newcommand{\Q}{\mathbb{Q}}

\newcommand{\gm}{\gamma}
\newcommand{\lm}{\lambda}
\newcommand{\Gm}{\Gamma}

\newcommand{\bB}{\bm{B}}
\newcommand{\bG}{\bm{G}}
\newcommand{\bR}{\bm{R}}
\newcommand{\bC}{\bm{C}}
\newcommand{\bZ}{\bm{Z}}
\newcommand{\bGam}{\bm{\Gamma}}
\newcommand{\bLam}{\bm{\Lambda}}
\newcommand{\Gprime}{\bm{G}'}

\newcommand{\bn}{\mathbf{n}}
\newcommand{\bh}{\mathbf{h}}
\newcommand{\bW}{\bm{W}}
\newcommand{\bV}{\bm{V}}
\newcommand{\bc}{\mathbf{c}}
\newcommand{\zroM}{\mathbf{0}}
\newcommand{\oneM}{\mathbf{1}}
\newcommand{\bx}{\mathbf{x}}
\newcommand{\eM}{\mathbf{e}}
\newcommand{\bX}{\bm{X}}
\newcommand{\bY}{\bm{Y}}

\newcommand{\bxi}{\bm{\xi}}
\newcommand{\bze}{\bm{\zeta}}
\newcommand{\bsi}{\bm{\sigma}}
\newcommand{\bmag}{\mathbf{m}}
\newcommand{\bH}{\mathbf{H}}
\newcommand{\br}{\mathbf{r}}
\newcommand{\zroN}{\mathbf{0}}
\newcommand{\oneN}{\mathbf{1}}
\newcommand{\by}{\mathbf{y}}
\newcommand{\eN}{\mathbf{e}}
\newcommand{\bu}{\mathbf{u}}
\newcommand{\bv}{\mathbf{v}}
\newcommand{\bw}{\mathbf{w}}

\newcommand{\bg}{\mathbf{g}}

\newcommand{\st}{\textsf{t}}

\newcommand{\ROW}{\textup{\textsf{R}}}
\newcommand{\PRF}{\textup{\textsf{P}}}
\newcommand{\COL}{\textup{\textsf{C}}}
\newcommand{\SAT}{\textup{\textsf{S}}}
\newcommand{\ADM}{\textup{\textsf{A}}}
\newcommand{\BBB}{\textup{\textsf{B}}}

\DeclareMathOperator{\spn}{span}
\DeclareMathOperator{\vol}{vol}
\DeclareMathOperator{\proj}{proj}

\newcommand{\bP}{\mathbf{P}}
\newcommand{\bQ}{\mathbf{Q}}

\newcommand{\cA}{\mathcal{A}}

\newcommand{\cL}{\mathcal{L}}
\newcommand{\cP}{\mathcal{P}}

\newcommand{\FF}{\mathscr{F}}
\newcommand{\Fprime}{\mathscr{F}'}

\title[Gardner formula for Ising perceptron models at small densities]{Gardner formula for Ising perceptron models at small densities}

\newcommand{\hpi}{\hat{\pi}}
\newcommand{\cpi}{\acute{\pi}}
\newcommand{\vpi}{\varpi}

\newcommand{\bemph}{\textbf}

\newcommand{\RS}{\textup{\textsf{RS}}}
\newcommand{\AT}{\textup{\textsf{AT}}}
\newcommand{\annealed}{\textup{\textsf{ann}}}

\newcommand{\mbH}{\mathbb{H}}

\newcommand{\bE}{\bm{E}}
\newcommand{\ap}{\mathbf{p}}

\newcommand{\tcA}{\tilde{\mathcal{A}}}
\newcommand{\tcK}{\tilde{\mathcal{K}}}
\newcommand{\tcL}{\tilde{\mathcal{L}}}
\newcommand{\tbX}{\tilde{\bm{X}}}

\newcommand{\ii}{\mathfrak{i}}

\newcommand{\CSTZRO}{c_0}
\newcommand{\CSTONE}{c_1}
\newcommand{\csttwo}{c_2}
\newcommand{\BIGCTWO}{C_2}

\newcommand{\ERR}{\textup{\textsf{ERR}}}

\newcommand{\PL}{\textup{\textsf{PL}}}

\newcommand{\dvpi}{\dot{\vpi}}
\newcommand{\dpi}{\dot{\pi}}

\newcommand{\starpi}{\pi_*}
\newcommand{\stardpi}{\dpi_*}
\newcommand{\starvpi}{\vpi_*}
\newcommand{\stardvpi}{\dvpi_*}

\newcommand{\sss}{\mathfrak{s}}
\newcommand{\bareps}{\bar{\epsilon}}
\newcommand{\ETA}{\eta} 

\newcommand{\bbx}{\bar{\bx}}
\newcommand{\bby}{\bar{\by}}

\author[E.\ Bolthausen]{$^\star$Erwin Bolthausen}
\author[S.\ Nakajima]{$^\circ$Shuta Nakajima}
\author[N.\ Sun]{$^\dagger$Nike Sun}
\author[C.\ Xu]{$^\ddagger$Changji Xu}
\thanks{\raggedright$^\star$Institute of Mathematics, University of Zurich. $^\circ$Department of Mathematics and Computer Science, University of Basel. $^\dagger$Department of Mathematics, Massachusetts Institute of Technology. $^\ddagger$Center for Mathematical Sciences and Applications, Harvard University.}

\begin{document}

\date{\today}

\begin{abstract} We consider the Ising perceptron model with $N$ spins and $M=N\alpha$ patterns, with a general activation function $U$ that is bounded above. For $U$ bounded away from zero or $U(x)=\Ind{x\ge\kappa}$, it was shown by Talagrand \cite{MR1782273,MR3024566} that for small densities $\alpha$, the free energy of the model converges as $N\to\infty$ to the replica symmetric formula conjectured in the physics literature \cite{krauthmezard1989} (see also \cite{gd1988optimal}). We give a new proof of this result, which covers the more general class of all functions $U$ that are bounded above and satisfy a certain variance bound. The proof uses the (first and second) moment method conditional on the approximate message passing iterates of the model. In order to deduce our main theorem, we also prove a new concentration result for the perceptron model in the case where $U$ is not bounded away from zero.
\end{abstract}

\maketitle

\setcounter{tocdepth}{1}
\tableofcontents

\section{Introduction}
\label{s:intro}

\subsection{Overview}
We study a class of generalized \bemph{Ising perceptron} models, defined as follows. Let $\bG$ be an $M\times N$ matrix with i.i.d.\ standard gaussian entries. Let $U:\R\to[0,\infty)$ be a bounded measurable function
(the \bemph{activation function})
and denote $u\equiv \log U : \R\to [-\infty, \infty)$. The associated \bemph{Ising perceptron partition function} is
	\beq\label{e:Z}
	\bZ 
	\equiv\bZ(\bG)
	\equiv \sum_J
	\exp\bigg\{
	\sum_{a\le M} 
	u\bigg(
	 \f{(\eM_a)^\st \bG J}{N^{1/2}} \bigg)
	\bigg\}
	\,,
	\eeq
where the sum goes over $J\in\set{-1,+1}^N$. The $J_i$ are called the \bemph{spins}, while the vectors $\bg^a$ are called the \bemph{patterns}.
A special case is the \bemph{half-space intersection} model defined by the function $U(x)=\Ind{x\ge\kappa}$, where $\kappa\in\R$ is a fixed parameter. In this paper we develop a method to compute the asymptotic free energy of the generalized model \eqref{e:Z} with $M=\alpha N$ for small $\alpha$, and $N\to\infty$. Note that if $U$ is scaled by any factor $c$, then the partition function \eqref{e:Z} is scaled by $c^M$ --- therefore, since we assume $U$ is bounded, we may as well assume that $U$ maps into $[0,1]$. More precisely, we work throughout under the following:

\begin{ass}\label{a:bdd}
The function $U$ is a measurable mapping from $\R$ into $[0,1]$. Moreover, with $\E_\xi$ denoting expectation over the law of a standard gaussian random variable $\xi$, we have
	\beq\label{e:not.annealed}
	\E_\xi[\xi U(\xi)]
	= \int z U(z)\,\varphi(z)\,dz
	\ne0\,,
	\eeq
where $\varphi$ denotes the standard gaussian density.
\end{ass}

\begin{ass}\label{a:Lip} 
Writing $\E_{\xi,\xi'}$ for expectation over i.i.d.\ standard gaussians $\xi,\xi'$, the quantity
	\[(K_2)'(U)
	\equiv\max\bigg\{1,\sup\bigg\{
	\f{\E_{\xi,\xi'}[ (\xi-\xi')^2
		U(x+c\xi)U(x+c\xi')]}
		{\E_{\xi,\xi'}[U(x+c\xi)U(x+c\xi')]}
	: x\in\R,\f25 \le c \le \f73 \bigg\}\bigg\}\,.\]
is finite. This assumption implies that the quantity
	\[K_2(U)\equiv\max\bigg\{1,
	\sup\bigg\{
	\f{\E_{\xi,\xi'}[ (\xi-\xi')^2
		U(x+c\xi)U(x+c\xi')]}
		{\E_{\xi,\xi'}[U(x+c\xi)U(x+c\xi')]}
	: x\in\R,
	\f12 \le c \le 2 \bigg\}
	\bigg\}
	\]
is also finite, and indeed $K_2(U)\le (K_2)'(U)$. (The bound on  
$(K_2)'(U)$ further ensures that $K_2(U_\ETA)$ is bounded, where $U_\ETA$ is a smoothed approximation of $U$;
see Lemma~\ref{l:smoothed.Ktwo}.)
\end{ass}

See Remark~\ref{r:assumptions} below for more discussion on the above assumptions --- in particular, we will explain that assumption~\eqref{e:not.annealed} only rules out an easier case of the problem. To state our main results we introduce some further notation. As above, let $\xi$ denote an independent standard gaussian random variable, and let $\E_\xi$ denote expectation over the law $\xi$. Given $q\in[0,1)$ let
	\beq\label{e:L}
	L_q(x)
	\equiv \log \E_\xi
		U\Big( x + (1-q)^{1/2}\xi\Big)
	\equiv \log \int 
	U\Big( x + (1-q)^{1/2}z\Big)
	\varphi(z)\,dz\,,
	\eeq
where $\varphi$ denotes the standard gaussian density as above. Let
	\beq\label{e:F}
	F_q(x)
	= (L_q)'(x)
	= \f1{(1-q)^{1/2}} 
	\f{\E_\xi [\xi U( x + (1-q)^{1/2}\xi) ]}
	{\E_\xi U( x + (1-q)^{1/2}\xi) }
	\,.\eeq
We will sometimes abbreviate $L\equiv L_q$ and $F\equiv F_q$.

\begin{ppn}[\hyperlink{proof:p.fp}{proved in Section~\ref{s:technical}}]\label{p:fp}
If $U$ satisfies Assumption~\ref{a:bdd}, then there exists a positive constant $\alpha(U)>0$ such that for all $0<\alpha\le\alpha(U)$ there exists a unique pair $(q,\psi) \in [0,1/25]\times[0,\infty)$ satisfying
	\beq\label{e:fp}
	\begin{pmatrix} q\\ \psi
	\end{pmatrix}
	= \begin{pmatrix}
	\bar{q}(\psi)\\
	\alpha \bar{r}(q)
	\end{pmatrix}
	\equiv \begin{pmatrix}
	\E[  \Th(\psi^{1/2}Z)^2] \\ 
	\alpha \E[ F_q(q^{1/2}Z)^2 ]
	\end{pmatrix}\,,
	\eeq
where $F_q$ is defined by \eqref{e:F} (and depends on $U$). Moreover we can take 
	\beq\label{e:alpha.U}
	\alpha(U)\equiv \f1{e^{10}
	\cdot \CSTONE 
		\cdot C_1(U)^6
		\cdot K_2(U)^4}\,,
	\eeq
where $\CSTONE$ is an absolute constant 
characterized by Lemma~\ref{l:dpsi.dq} and 
Corollary~\ref{c:rs.estimate},
and $C_1(U)$ is a finite constant depending only on $U$
which is characterized by Lemma~\ref{l:poly}.
The solution $(q,\psi)$ of \eqref{e:fp} satisfies 
	\beq\label{e:fp.bounds}\f{(\E_\xi[\xi U(\xi)])^2}{2}
	\le \f{q}{\alpha}
	\le \f{\psi}{\alpha}
	\le 3 \cdot C_1(U)^2\eeq
for all $0\le \alpha\le\alpha(U)$.
\end{ppn}

For any $U$ and $\alpha$ such that 
\eqref{e:fp} has a unique solution $(q,\psi)\in [0,1)\times[0,\infty)$,
the \bemph{replica symmetric formula} for the free energy of the corresponding perceptron model
\eqref{e:Z} is given by
	\beq\label{e:rs}
	\RS \equiv \RS(\alpha;U)
	= -\f{\psi(1-q)}{2}
	+\E\bigg\{ \log 2\Ch(\psi^{1/2}Z)
	+\alpha  L_q(q^{1/2}Z)\bigg\}\,,\eeq
where the expectation is over an independent standard gaussian random variable $Z$. In this paper we show:

\begin{thm}[\bemph{main theorem}]\label{t:main}
If $U$ satisfies Assumptions~\ref{a:bdd} and \ref{a:Lip}, then there exists a positive constant $\alpha=\alpha'(U)>0$ such that,
if $\bG$ is an $M\times N$ matrix with i.i.d.\ standard gaussian entries and $M/N\to\alpha$ with $0\le \alpha\le\alpha'(U)$, then
for the (generalized) Ising perceptron model \eqref{e:Z} we have
	\[\lim_{N\to\infty} \f1N \log \bZ(\bG)
	=\RS(\alpha;U)\,,\]
where the limit is in probability. Moreover we can take 
	\beq\label{e:alpha.p.U}
	\alpha'(U)
	\equiv 
	\f1{e^{16}
	\cdot \CSTONE 
		\cdot C_1(U)^6
		\cdot (K_2)'(U)^4}
	\le \f{\alpha(U)}{e^6}\,,
	\eeq
for $\alpha(U)$ as defined by \eqref{e:alpha.U}.
\end{thm}

For $U$ bounded away from zero, as well as for $U(x)=\Ind{x\ge\kappa}$, the result of Theorem~\ref{t:main} was previously shown by Talagrand \cite{MR1782273,MR3024566}. Our proof is very different from Talagrand's, and uses the idea of ``conditioning on the AMP iteration,'' as previously introduced by \cite{ding2018capacity,MR4015008} (see also \cite{alaoui2020algorithmic, fan2021replica, brennecke2021note}). By contrast, Talagrand's proof uses an interpolation approach, which seemingly necessitates more conditions on $U$. Our result extends to the more general class of functions $U$ satisfying Assumptions~\ref{a:bdd} and \ref{a:Lip}. See \S\ref{sss:lit.ising} below for further discussion and comparison.

\begin{rmk}\label{r:assumptions}
We make some further comments on our assumptions:
\begin{enumerate}[1.]
\item From our perspective, Assumption~\ref{a:bdd} is relatively mild. It may be possible to relax 
the condition $U\le1$ to accommodate functions $U(x)$ that do not grow too quickly in $|x|$, but we will not pursue this here. Next, if the condition \eqref{e:not.annealed} fails --- meaning that $\E_\xi[\xi U(\xi)]=0$ --- then the fixed-point equation \eqref{e:fp} is solved by $q=\psi=0$, and the replica symmetric free energy
\eqref{e:rs} reduces to the \bemph{annealed free energy}
	\beq\label{e:annealed}
	\annealed(\alpha;U)
	=\f1N \log \E\bZ(\bG)
	=\log2+\alpha\log \E U(\xi)\,.
	\eeq
In this case, it is known that the limiting free energy can be obtained by a direct first and second moment method approach, without the need of a conditioning scheme. This is done for the case of symmetric $U$ by \cite{MR3983947}, and the argument of that paper can be extended to cover the case $\E_\xi[\xi U(\xi)]=0 $. Moreover it is expected that this case may be more tractable to analyze for finer properties of the solution space, following \cite{perkins2021frozen,abbe2021proof} (further discussed in \S\ref{sss:lit.other} below).

\item We view Assumption~\ref{a:Lip} as the slightly more restrictive condition, although we will show (by straightforward arguments) that it holds if $U$ is bounded away from zero, compactly supported, or logconcave (see Proposition~\ref{p:logc} below). Moreover, Assumption~\ref{a:Lip} is essentially necessary to ensure that the function $F_q$ in \eqref{e:F} is Lipschitz --- this is by an easy calculation, which we give in Lemma~\ref{l:F.Lipschitz}. This allows us to use existing results on AMP and state evolution  (\cite{MR2810285,MR3147441}; see \S\ref{ss:intro.AMP} and \S\ref{ss:review.AMP}) which all require the message-passing functions to be Lipschitz. 
\end{enumerate}
\emph{Assumption~\ref{a:bdd} holds throughout this paper, even if not explicitly stated.} However, we will point out explicitly each place where Assumption~\ref{a:Lip} is used.
\end{rmk}

\begin{ppn}[proved in \S\ref{ss:logc}]\label{p:logc} Suppose $U$ satisfies Assumption~\ref{a:bdd}. If in addition $U$ is bounded away from zero, compactly supported, or logconcave, then $U$ also satisfies Assumption~\ref{a:Lip}.
\end{ppn}

\subsection{Background and related work}
\label{ss:literature}

In this subsection we give some background on the perceptron model, and survey the related work. \bemph{Some high-level discussion of 
key ideas in this paper is given in
 \S\ref{sss:lit.ising}--\ref{sss:lit.AMP}.}

\newcommand{\asat}{\alpha_\textup{c}}

The perceptron problem originates from a toy model of a single-layer neural network, as follows. Suppose we have $N+1$ input nodes, labelled $0\le j\le N$. Likewise we have $N+1$ output nodes, labelled $0\le i\le N$. For all $i\ne j$, between the $j$-th input node and the $i$-th output node there is an edge weight $J_{i,j}$, to be determined. It will be convenient to fix $J_{i,i}\equiv0$ for all $i$. The system is given $M$ input ``patterns'' $g^1,\ldots,g^M$, which are  vectors in $\R^{N+1}$. We then say that the system \bemph{memorizes} the pattern $g^a$ if
	\beq\label{e:mem.sign}
	\sgn \bigg( \sum_{j=0}^N
	J_{i,j} (g^a)_j\bigg)
	= \sgn\Big( (g^a)_i\Big)
	\eeq
for all $0\le i\le N$. One can then ask, given $M=N\alpha$ i.i.d.\ random patterns, whether there exists a choice of edge weights $J$ such that the system memorizes all $M$ patterns. The \bemph{storage capacity} $\asat$ of the model is the supremum of all $\alpha=M/N$ for which memorization of all $M$ given patterns is possible with probability $1-o_N(1)$. Models of this type have been considered at least since the mid-20th century
(e.g.\ \cite{mcculloch1943logical, hebb1949organization, little1974existence, MR652033}).

One can consider the constraint \eqref{e:mem.sign} separately for each $0\le i\le N$, and by symmetry it suffices to understand the case $i=0$.
Recall that $J_{0,0}\equiv0$, and denote $J_i\equiv J_{i,0}$ for $1\le i\le N$. Denote $g_{a,j}\equiv (g^a)_0(g^a)_j$ for all $1\le a\le M$ and $1\le j\le N$, and note that the $g_{a,i}$ are i.i.d.\ standard gaussian random variables. Thus \eqref{e:mem.sign} is equivalent to
	\[
	\f1{N^{1/2}}\sum_{j=0}^N g_{a,j} J_j \ge \kappa
	\]
for $\kappa=0$. Of course, one can then generalize the model by taking a non-zero parameter $\kappa$: taking $\kappa<0$ weakens the original constraint \eqref{e:mem.sign}, while taking $\kappa>0$ gives a more restrictive constraint than \eqref{e:mem.sign}. This is equivalent to the model \eqref{e:Z} with $U(x)=\Ind{x\ge\kappa}$.
The two most commonly studied variants of the model are the \bemph{Ising perceptron} where $J_i\in\set{-1,+1}$ (as in this paper), and the \bemph{spherical perceptron} where $J=(J_i)_{i\le N}$ is restricted to the sphere of radius $N^{1/2}$. 

\subsubsection{Non-rigorous results from statistical physics}

In the physics literature, the spherical perceptron with $U(x)=\Ind{x\ge\kappa}$ for $\kappa\ge0$  was analyzed in a series of celebrated works of Gardner and Derrida \cite{gardner1987maximum,gardner1988space,gd1988optimal,gd1989three}, using the non-rigorous replica method. This method also applies to the Ising perceptron with  $U(x)=\Ind{x\ge\kappa}$ for any $\kappa\in\R$, but the original Gardner--Derrida analysis contained an error leading to incorrect predictions. A corrected replica calculation for the Ising model was first given by Krauth and M\'ezard \cite{krauthmezard1989}. The same results were rederived using the cavity method by M\'ezard \cite{mezard1989space}. (While the replica and cavity methods are both non-rigorous, the cavity method may be generally considered to yield more transparent derivations.)

The Gardner--Derrida and Krauth--M\'ezard results cover the \bemph{replica symmetric} regime, where the system is expected to exhibit some form of correlation decay. The spherical perceptron 
with $U(x)=\Ind{x\ge\kappa}$ (also called the \emph{positive spherical perceptron}) is expected to be replica symmetric only for $\kappa\ge0$, whereas the Ising perceptron
with $U(x)=\Ind{x\ge\kappa}$ 
is expected to be replica symmetric 
for all $\kappa\in\R$. More recently there have been several works in the physics literature investigating the \emph{negative spherical perceptron} and its potential consequences in statistical applications, e.g.\ \cite{franz2016simplest,franz2017universality}.

\subsubsection{Rigorous results on the spherical perceptron}
\label{sss:lit.spherical}

The mathematical literature contains numerous very strong results on the spherical perceptron for $U(x)=\Ind{x\ge\kappa}$, especially for $\kappa\ge0$ (conjecturally the replica symmetric regime). For $\kappa=0$, the storage capacity $\asat=2$ was known since the 1960s \cite{MR146858,cover1965geometrical}. For general $\kappa\ge0$, the storage capacity $\asat(\kappa)$ was proved by a short and elegant argument \cite{stojnic2013another}, using convex duality together with Gordon's gaussian minimax comparison inequality \cite{MR800188,MR950977}. However, perhaps the most striking result for this model is that of Shcherbina and Tirozzi \cite{MR1964377}, proving the Gardner free energy formula for the spherical perceptron for
 all $\kappa\ge0$ and all $\alpha$ up to $\asat(\kappa)$. The proof of \cite{MR1964377} makes crucial use of the classical Brunn--Minkowski inequality for volumes of bodies in euclidean space \cite{lusternik1935brunn,MR82697}. The main result of \cite{MR1964377} was reproved by Talagrand (\cite[Ch.~3]{MR2731561} and \cite[Ch.~8]{MR3024566}) with a perhaps slightly simpler argument, using instead the functional Brunn--Minkowski (Pr\'ekopa--Leindler) inequality \cite{MR315079,MR0430188,MR404557}. This inequality implies concentration of Lipschitz functionals under strongly logconcave measures \cite{MR1097258}, which can be used to deduce concentration of overlaps and cavity equations (see e.g.\ \cite[Thm.~3.1.11]{MR2731561}).\footnote{In this work we have also used the result of \cite{MR1097258} (restated in Theorem~\ref{t:maurey}), but only to prove Proposition~\ref{p:logc} which is not required for the main theorem.} As noted by \cite{MR1964377} and \cite[\S3.4]{MR2731561}, similar concentration results can also be obtained using instead the Brascamp--Lieb inequality \cite{MR0450480}; and indeed this idea appears in earlier work on the Hopfield model \cite{MR1601727}. Thus, all existing results on the positive spherical perceptron (excluding the case $\kappa=0$) use powerful tools from \bemph{convex geometry}.\footnote{The Pr\'ekopa--Leindler inequality generalizes the Brunn--Minkowski inequality, and also can be used to deduce the Brascamp--Lieb inequality \cite{MR1800062}. For more on the relations among these inequalities we refer to the survey \cite{MR1898210}.}

\subsubsection{Rigorous results on the Ising perceptron}
\label{sss:lit.ising}

The mathematical literature on the Ising perceptron is far less advanced than for the spherical perceptron. For the half-space model, the free energy was computed heuristically by \cite{krauthmezard1989}; their method applies also to the more general model \eqref{e:Z}. One consequence of the \cite{krauthmezard1989} calculation is an explicit prediction $\alpha_\star$ for the storage capacity $\asat$ for the model $U(x)=\Ind{x\ge\kappa}$ --- for $\kappa=0$, the conjectured threshold $\alpha_\star$ is approximately $0.83$.

In the rigorous literature, most existing results concern the half-space model $U(x)=\Ind{x\ge0}$.\footnote{The existing results for $U(x)=\Ind{x\ge0}$ can likely be extended to cover
$U(x)=\Ind{x\ge\kappa}$ for any $\kappa\in\R$.} For this model,
it was shown by \cite{MR1629627,MR1716771} that there is a small absolute constant $\epsilon>0$ such that the transition must occur between $\epsilon$ and $1-\epsilon$: that is, the partition function \eqref{e:Z} is non-zero with high probability for $\alpha\le\epsilon$, and zero with high probability for $\alpha\ge1-\epsilon$. A more recent work \cite{ding2018capacity} (further discussed below) uses some of the methods of this paper to show, under a certain variational hypothesis, that the partition function is non-zero with \emph{non-negligible} probability for $\alpha<\alpha_\star$, where $\alpha_\star$ is the conjectured threshold from \cite{krauthmezard1989}. A more recent work \cite{MR4317708} confirms that the model indeed has a \emph{sharp} threshold.\footnote{To be precise, the result of \cite{ding2018capacity} is with gaussian noise $\bG$ (as in this paper), while the other results \cite{MR1629627,MR1716771,MR4317708} are for the Bernoulli noise model where $g_{a,i}$ are i.i.d.\ symmetric random signs. It is reasonable to expect that the result of \cite{MR1629627,MR1716771,MR4317708} can be transferred to the Bernoulli noise model.}

For the situation where we have a more general function $U$ in \eqref{e:Z},
Talagrand \cite{MR1782273} (see also \cite[Ch.~2]{MR2731561}) proves that the limiting free energy is given by the replica symmetric formula \eqref{e:rs}, for small enough $\alpha$, under the assumption that the function $u\equiv\log U$ is \bemph{uniformly bounded}. This corresponds to the case of our main result Theorem~\ref{t:main} where $u$ is bounded, which we \hyperlink{proof:t.main.lbd.bdd}{prove at the end of Section~\ref{s:second.mmt}.} Even for bounded $u$, the two proofs are very different: \cite{MR1782273} uses an interpolation method to derive replica symmetric equations, while this paper uses first and second moments conditional on the AMP iteration. We remark also that the argument of \cite{MR1782273} seemingly needs to go through a smoothed approximation of $u$, while our proof for bounded $u$ requires no smoothing.

In comparison with previous work of Talagrand, the main new result of this work is that the limiting free energy is given by the replica symmetric formula \eqref{e:rs}, for small enough $\alpha$, for all $U$ satisfying Assumptions~\ref{a:bdd} and \ref{a:Lip}. A special case of this result, for the half-space model $U(x)=\Ind{x\ge\kappa}$, was previously obtained in \cite[Ch.~9]{MR2731561} (with partial results appearing in a previous work \cite{MR1680236}).\footnote{The function $U(x)=\Ind{x\ge\kappa}$
satisfies the hypothesis of Theorem~\ref{t:main}: it clearly satisfies Assumption~\ref{a:bdd},
and one can check that it satisfies Assumption~\ref{a:Lip} either by direct calculation or by applying Proposition~\ref{p:logc}.} Talagrand's proof for the half-space model relies crucially on an estimate \cite[Thm.\ 8.2.4]{MR3024566} which says roughly that if $(u_i)_{i\le n}$ is a near-isotropic gaussian process, then the fraction of indices $i$ where $u_i \ge \kappa$ cannot be too small. The proof of this estimate uses a gaussian comparison inequality (see \cite[Lem.~1.3.1]{MR2731561} and
\cite[Propn.~8.2.2]{MR3024566}), and does not extend for instance to the event $u_i\in E$ where $E$ is a bounded measurable subset of $\R$. In this paper we prove an analogous (weaker) estimate for general $E$ by different methods (Proposition~\ref{p:gaus.poly}), and use this in the proof of Theorem~\ref{t:main} in the case of unbounded $u$.

\subsubsection{TAP, AMP, and conditioning}
\label{sss:lit.AMP}

The main idea in the proof of Theorem~\ref{t:main},
which we discuss further in \S\ref{ss:intro.AMP} below,
is to \bemph{compute (first and second) moments of the partition function \eqref{e:Z} conditional on the AMP filtration}. The motivation originates from the \bemph{TAP (Thouless--Anderson--Palmer) framework}, which were introduced for the classical Sherrington--Kirkpatrick model \cite{sherrington1975solvable} by \cite{thouless1977solution} (and further investigated by \cite{AT1978stability,plefka1982convergence}).  For the model \eqref{e:Z}, the TAP equations read
	\begin{align}
	\label{e:perceptron.TAP.col}
	\bmag &
	\equiv \Th(\bH)
	= \Th\bigg(
		\f{\bG^\st\bn}{N^{1/2}}
		-\beta\bmag
	\bigg)\,,\\
	\bn
	&\equiv F(\bh)
	= F\bigg(
		\f{\bG\bmag}{N^{1/2}} 
		- \acute{\beta}\bn
	\bigg)\,,
	\label{e:perceptron.TAP.row}
	\end{align}
where the functions $\Th$ and $F$ are applied coordinatewise,
$\bmag\equiv \Th(\bH)$ is a vector in $\R^N$,
and $\bn\equiv F(\bh)$ is a vector in $\R^M$.
For the model \eqref{e:Z} at small $\alpha$, it is conjectured that the TAP equations
\eqref{e:perceptron.TAP.col} and \eqref{e:perceptron.TAP.row}
have a unique solution $(\bmag^\star,\bn^\star)$, such that $\bmag^\star$ approximates the mean value of a random configuration $J$ sampled from the Gibbs measure
	\beq\label{e:Gibbs}
	\mu(J)
	\equiv
	\f1{\bZ(\bG)}
	\prod_{a\le M} U\bigg(\f{(\eM_a)^\st \bG J}{N^{1/2}} \bigg)\,.
	\eeq
Meanwhile, the vector $\bn^\star$ describes the distribution of the vector $\bG J/N^{1/2}$ where $J$ is sampled from $\mu$; see \cite{mezard1989space}.
It is further expected that $N^{-1}\log \bZ$ concentrates very well around a \bemph{TAP free energy} $\Phi(\bmag^\star,\bn^\star)$, which in turn concentrates around the replica symmetric value \eqref{e:rs}. The TAP equations and TAP free energy  can be derived as a dense limit of the belief propagation equations and Bethe free energy; see \cite{mezard2017mean}. For more recent work on the TAP framework in a variety of settings, we refer to \cite{chen2018generalized,MR4207445,MR4203332,arous2021shattering,MR4261707}.

As we commented in Remark~\ref{r:assumptions} above, if we have $\E_\xi[\xi U(\xi)]=0$ (i.e.\ if the assumption~\eqref{e:not.annealed} does not hold), then the (unconditional) second moment method can be used to
analyze the partition function $\bZ$ from \eqref{e:Z}, following \cite{MR3983947}. If $\E_\xi[\xi U(\xi)]\ne0$, however, it is well known that the unconditional second moment method does not say anything about the random variable $\bZ$, at any positive $\alpha=M/N$. Since the TAP fixed point $(\bmag^\star,\bn^\star)$ is described by a relatively simple set of equations \eqref{e:perceptron.TAP.col} and \eqref{e:perceptron.TAP.row}, and is conjectured to carry a great deal of information about the random measure \eqref{e:Gibbs}, it is natural to consider the second moment method \bemph{conditional on the TAP solution} $(\bmag^\star,\bn^\star)$. The problem with this approach is that it is not in fact known that the equations \eqref{e:perceptron.TAP.col} and \eqref{e:perceptron.TAP.row} have a unique solution. A way around this is to use the \bemph{AMP (approximate message passing) iteration}, which constructs approximate solutions of the TAP equations \cite{MR2810285,MR3147441}.

The idea of \bemph{conditioning on the AMP iteration} was introduced by \cite{ding2018capacity,MR4015008} and has been developed in subsequent works \cite{alaoui2020algorithmic, fan2021replica, brennecke2021note}. Of these prior works, \cite{MR4015008} and \cite{brennecke2021note} concern the classical Sherrington--Kirkpatrick (SK) model
with a gaussian coupling matrix
(i.e., the Hamiltonian is a scalar multiple of $J^\st \bG J$ where $\bG$ is an $N\times N$ matrix with i.i.d.\ random gaussian entries). The work \cite{fan2021replica} concerns (more general) SK models with random orthogonally invariant coupling matrices, and uses a simplified ``memory-free'' AMP iteration that was developed and analyzed by \cite{opper2001adaptive,opper2016theory,cakmak2019memory,fan2020approximate}. The works \cite{ding2018capacity} and \cite{alaoui2020algorithmic} concern the perceptron model, but only use the AMP conditioning method to prove lower bounds. In the current work, we show that the AMP conditioning method gives sharp upper and lower bounds for the generalized perceptron \eqref{e:Z} at small $\alpha$.

\subsubsection{Other related work}
\label{sss:lit.other}

As noted above, in the special case that $U$ satisfies $\E_\xi[\xi U(\xi)]=0$ (i.e.\ if assumption~\eqref{e:not.annealed} does not hold), the model \eqref{e:Z} is mathematically much more tractable, and can be analyzed by an (unconditional) second moment method. The condition $\E_\xi[\xi U(\xi)]=0$ holds for instance if $U$ is a bounded \bemph{symmetric} function. The second moment analysis was done for the cases $U(x)=\Ind{|x|\le\kappa}$ and $U(x)=\Ind{|x|\ge\kappa}$ in \cite{MR3983947}. For the model $U(x)=\Ind{|x|\le\kappa}$, much finer structural results (on the typical geometry of the solution space) were obtained by \cite{perkins2021frozen,abbe2021proof}. These results were inspired in part by questions raised in the physics literature about the algorithmic accessibility of CSP solutions (see e.g.\ \cite{baldassi2016unreasonable,budzynski2019biased}). Finally, for the perceptron model in statistical settings, there is an extensive literature
which we will not describe here; we refer the reader for instance to \cite{MR3939767,mzz2021negper} and many references therein.

\subsection{AMP iteration}
\label{ss:intro.AMP}
Our convention throughout is that if $f:\R\to\R$
and $\mathbf{z}\equiv (z_j)_j$ is any vector, then
	\beq\label{e:componentwise}
	f(\mathbf{z})
	\equiv (f(z_j))_j
	\eeq
denotes the vector of the same length
which results from applying $f$ componentwise to $\mathbf{z}$.
Recall $F\equiv F_q$ from \eqref{e:F}. Let $\bmag^{(0)}=\zroN\in\R^N$, $\bn^{(0)}=\zroM\in\R^M$, $\bmag^{(1)}=q^{1/2}\oneN\in\R^N$, $\bn^{(1)}=(\psi/\alpha)^{1/2}\oneM\in\R^M$.
 The \bemph{approximate message passing (AMP)} iteration for the perceptron model is given by 
(cf.\ \eqref{e:general.AMP.row}
and \eqref{e:general.AMP.col})
	\begin{align}
	\label{e:perceptron.AMP.col}
	\bmag^{(t+1)} &
	\equiv \Th(\bH^{(t+1)})
	= \Th\bigg(
		\f{\bG^\st\bn^{(t)}}{N^{1/2}}
		-\beta\bmag^{(t-1)}
	\bigg)\,,\\
	\bn^{(t+1)} 
	&\equiv F(\bh^{(t+1)})
	= F\bigg(
		\f{\bG\bmag^{(t)}}{N^{1/2}} 
		- \acute{\beta}\bn^{(t-1)}
	\bigg)\,,
	\label{e:perceptron.AMP.row}
	\end{align}
where $\beta\bmag^{(t-1)}$
and $\acute{\beta}\bn^{(t-1)}$ are the \bemph{Onsager correction terms}, whose coefficients are defined by
	\beq\label{e:perceptron.onsager}
	\begin{pmatrix}\beta \\ \acute{\beta}\end{pmatrix}
	= \begin{pmatrix}
	\alpha \E F'( q^{1/2}Z) \\
	\E\Th'(\psi^{1/2}Z)
	\end{pmatrix}\,.
	\eeq
We remark that since $\Th'(x)=1-(\Th x)^2$, it follows using \eqref{e:fp} that $\acute{\beta}=1-q$. Recall from the discussion of \S\ref{sss:lit.AMP}
that the main idea in the proof of Theorem~\ref{t:main} is to \bemph{compute (first and second) moments of the partition function \eqref{e:Z} conditional on the AMP filtration}
	\beq\label{e:tap.condition}
	\FF\equiv \FF(t)\equiv\sigma\bigg(
	\Big(\bG\bmag^{(s)},\bn^{(s+1)} : s\le t\Big),
	\Big( \bG^\st\bn^{(\ell)},\bmag^{(\ell+1)}
		: \ell \le t-1 \Big)
	\bigg)
	\eeq
in the limit $t\to\infty$. The computation relies on existing results on the asymptotic behavior of AMP in the large-$N$ limit \cite{MR2810285,MR3147441}
(see also \cite{donoho2009message,MR3311445,MR3876443,MR4079177}). In \S\ref{ss:review.AMP} we review the relevant results
from \cite{MR2810285,MR3147441} that are used in our proofs. The results from our conditional method of moments calculation are summarized as follows:

\begin{thm}[conditional first moment]\label{t:ubd}
If $U$ satisfies Assumptions~\ref{a:bdd} and \ref{a:Lip}, then there exists a positive constant $\alpha(U)>0$ such that,
if $\bG$ is an $M\times N$ matrix with i.i.d.\ standard gaussian entries and $M/N\to\alpha$, 
and $\FF(t)$ is the AMP filtration defined by \eqref{e:tap.condition}, then for all $0\le \alpha\le\alpha(U)$ we have
	\[
	\E\Big( \bZ \,\Big|\,\FF(t)\Big)
	\le \exp \bigg\{ N\Big(
		 \RS(\alpha;U) + o_t(1)
		 \Big)\bigg\}
	\]
with high probability (i.e., with probability $1-o_N(1)$).
\end{thm}

Theorem~\ref{t:ubd} implies the upper bound 
in Theorem~\ref{t:main} by standard arguments, using Markov's inequality. The \hyperlink{proof:t.main.ubd}{proof of the upper bound in Theorem~\ref{t:main}} is therefore given at the end of Section~\ref{s:analysis.first.mmt}, after the 
\hyperlink{proof:t.ubd}{proof of Theorem~\ref{t:ubd}}.

\begin{thm}[conditional second moment]\label{t:lbd} 
Suppose $U$ satisfies Assumptions~\ref{a:bdd} and \ref{a:Lip}, and $0\le \alpha \le\alpha(U)$
as defined by \eqref{e:alpha.U}.
If $\bG$ is an $M\times N$ matrix with i.i.d.\ standard gaussian entries and $M/N\to\alpha$, 
we can construct a random variable $\bar{\bZ}\le\bZ$ such that
	\beq\label{e:lbd.restricted.first.mmt}
	\E\Big( \bar{\bZ}(\bG) \,\Big|\,\FF(t)\Big)
	\ge \exp \bigg\{ N\Big(
		 \RS(\alpha;U) - o_t(1)
		 \Big)\bigg\}\eeq
with high probability,
and for which we have the second moment estimate
	\beq\label{e:ubd.restricted.second.mmt}
	\E\Big( \bar{\bZ}(\bG)^2 \,\Big|\,\FF(t)\Big)
	\le \exp \bigg\{2 N\Big(
		 \RS(\alpha;U) + o_t(1)
		 \Big)\bigg\}\,,\eeq
also with high probability.
\end{thm}

In the \emph{bounded case} $\|u\|_\infty<\infty$
(recall $u\equiv\log U$),
 Theorem~\ref{t:lbd}
implies the lower bound in  Theorem~\ref{t:main} by standard arguments, using the Azuma--Hoeffding martingale inequality.
The \hyperlink{proof:t.main.lbd.bdd}{proof of the lower bound in Theorem~\ref{t:main} in the bounded case}
is given at the end of Section~\ref{s:second.mmt},
after the \hyperlink{proof:t.lbd}{proof of Theorem~\ref{t:lbd}}. In the more general setting where $u$ may be unbounded, 
the proof of Theorem~\ref{t:main} requires further estimates, as we outline in the next subsection.
 
\subsection{Concentration results for unbounded case}
\label{ss:intro.unbdd}

Assumption~\ref{a:bdd} implies that we must have
	\beq\label{e:positive.measure.set}
	1\ge U(x) > \delta' \Ind{x\in E(U)}
	\eeq
where $\delta'$ is a positive constant, and $E(U)$ is a subset of the real line of positive Lebesgue measure (which we denote $|E(U)|$). Moreover we can assume without loss that $E(U)$ is bounded, i.e., $E(U)\subseteq[-E_{\max}(U),E_{\max}(U)]$ for some finite $E_{\max}(U)$. Following \cite[\S8.3]{MR3024566}, define the truncated logarithm
	\[\log_A(x)\equiv\max\Big\{-A,\log x\Big\}\,.\]
The following is an adaptation of \cite[Propn.~9.2.6]{MR3024566} (see also \cite[Propn.~8.3.6]{MR3024566}):

\begin{ppn}\label{p:poly.conc}
Suppose $U$ satisfies Assumption~\ref{a:bdd}, and let $\delta'$ and $E(U)$ be as above. Then for $\tau=\exp(-12)$ we have
	\[
	\P\bigg(\f1N\Big|\log_{N\tau}
	\bigg(\f{\bZ}{2^N}\bigg)
	-\E\log_{N\tau}\bigg(\f{\bZ}{2^N}\bigg)\Big|
	\ge \f{(\log N)^2}{N^{1/2}}
	\bigg)
	\le \f1{N^2}\]
for all $N$ large enough (depending on $|E(U)|$, $E_{\max}(U)$, and $\delta'$).
\end{ppn}

Next let $\ETA$ be a small positive constant, and consider the smoothed function
	\beq\label{e:U.ETA}
	U_\ETA(x)
	\equiv (U*\varphi_\ETA)(x)
	= \int U(x+\ETA z)\varphi(z)\,dz
	= \E_\xi U(x+\ETA\xi)\,.
	\eeq
Let $\bZ(\ETA)$ denote the perceptron partition function with $U_\ETA$ in place of $U$:
	\beq\label{e:tilde.Z}
	\bZ(\ETA)
	\equiv \sum_J \prod_{a\le M}
		U_\ETA\bigg( \f{(\bg^a,J)}{N^{1/2}}\bigg)\,.
	\eeq
Note that $U_\ETA$ satisfies Assumption~\ref{a:bdd}: it is a smooth mapping from $\R$ into $[0,1]$ for any $\ETA>0$,
and condition \eqref{e:not.annealed} holds for $\ETA$ small enough.
We will show (see Lemma~\ref{l:smoothed.Ktwo}) that $K_2(U_\ETA)$ can be  bounded in terms of $(K_2)'(U)$. We then have the following approximation result:

\begin{ppn}\label{p:Z.eps.approx.Z}
Suppose $U$ satisfies Assumption~\ref{a:bdd}, and let $\delta'$ and $E(U)$ be as above. Then we have
	\[
	\limsup_{N\to\infty}\f1N \bigg| \E\bigg[
	\log_{N\tau}\bigg(\f{\bZ(\ETA)}{2^N}\bigg)
		-\log_{N\tau}\bigg(\f{\bZ}{2^N}\bigg)
	\bigg] \bigg|
	\le o_\ETA(1)
	\]
for $\tau=\exp(-12)$.
\end{ppn}

Propositions~\ref{p:poly.conc} and \ref{p:Z.eps.approx.Z} are proved in Section~\ref{s:conc.log.Z}. The proofs rely on a bound for near-isotropic gaussian processes, Proposition~\ref{p:gaus.poly}, which we mentioned in \S\ref{sss:lit.ising} above. Finally, we have the following:

\begin{ppn}\label{p:RS.U.eps}
If $U$ satisfies Assumption~\ref{a:bdd}, then we have
	\[
	\lim_{\ETA\downarrow0} \RS(\alpha;U_\ETA)
	= \RS(\alpha;U)
	\]
for all $0\le \alpha\le\alpha'(U)$ (as defined by \eqref{e:alpha.p.U}).
\end{ppn}

\begin{ppn}\label{p:conc.smoothed}
Suppose $U$ satisfies Assumption~\ref{a:bdd}, and let $\bZ(\ETA)$ be as in \eqref{e:tilde.Z}. Then we have
	\[
	\P\bigg( \Big|\log \bZ(\ETA)
		-\E\log\bZ(\ETA)\Big| \ge Nx\bigg)
	\le 32N \cdot
	\exp\bigg\{ -\f{Nx^2}{  32 \BIGCTWO C_1(U;\ETA)^2 }	\bigg\}
	\]
for all $0\le x\le 5 (\BIGCTWO)^{1/2} C_1(U;\ETA)$.
\end{ppn}

The \hyperlink{proof:p.RS.U.eps}{proof of Proposition~\ref{p:RS.U.eps}} is given in Section~\ref{s:technical}, while the
\hyperlink{proof:p.conc.smoothed}{proof of Proposition~\ref{p:conc.smoothed}}
is given in Section~\ref{s:conc.log.Z}. Then Propositions~\ref{p:poly.conc}, \ref{p:Z.eps.approx.Z}, \ref{p:RS.U.eps}, and \ref{p:conc.smoothed} can be combined to finish the proof of Theorem~\ref{t:main} in the unbounded case $\|u\|_\infty=\infty$.
\hypertarget{sketch:proof.t.main}{The argument goes roughly as follows:}
by Propositions~\ref{p:poly.conc} and \ref{p:Z.eps.approx.Z}, with high probability
	\[
	\f1N\log_{N\tau} \f{\bZ}{2^N} - o_N(1)
	= \E\log_{N\tau} \f{\bZ}{2^N}
	=
	\f1N \E\log_{N\tau} \f{\bZ(\ETA)}{2^N}
	+o_\ETA(1)\,.
	\]
By applying Theorem~\ref{t:lbd} to $U_\ETA$, and combining with Proposition~\ref{p:RS.U.eps} and Proposition~\ref{p:conc.smoothed}, we obtain
	\[
	\f1N \E\log_{N\tau} \f{\bZ(\ETA)}{2^N}-o_N(1)
	= \RS(\alpha;U_\ETA)-\log2
	=\RS(\alpha;U)-\log2 + o_\ETA(1)\,.
	\]
For $0<\alpha\le\alpha(U)$, the above is $\ge -\tau/2$ by straightforward estimates (Corollary~\ref{c:rs.estimate}). It follows that
	\[
	-\f{\tau}{2}
	\le \RS(\alpha;U)-\log2
	= o_N(1) +	\f1N\log_{N\tau} \f{\bZ}{2^N}
	= o_N(1) + \f1N\log \f{\bZ}{2^N}
	\]
with high probability, as desired. 
At the end of Section~\ref{s:conc.log.Z} we give the \hyperlink{proof:t.main.conclusion}{conclusion of the  proof of Theorem~\ref{t:main}}, where the above sketch is made precise.

\subsection*{Organization}
The remaining sections of the paper are organized as follows:
\begin{itemize}
\item In Section~\ref{s:first.mmt} we give a preliminary expression
(see Theorem~\ref{t:ubd.slice})
 for the first moment of the perceptron partition function conditional on $\FF(t)$.

\item In Section~\ref{s:technical} we collect some basic technical results, including basic consequences of Assumptions~\ref{a:bdd} and \ref{a:Lip}. We also give the proofs of Propostions \ref{p:fp}, \ref{p:logc}, and \ref{p:RS.U.eps}.

\item In Section~\ref{s:analysis.first.mmt} we analyze the conditional first moment calculations from 
Section~\ref{s:first.mmt} and complete 
the \hyperlink{proof:t.ubd}{proof of Theorem~\ref{t:ubd}}.
This leads to 
the \hyperlink{proof:t.main.ubd}{upper bound in Theorem~\ref{t:main}}, presented at the end of the section.

\item In Section~\ref{s:second.mmt} we prove Theorem~\ref{t:lbd}, which bounds
the first and second moments of the (truncated) perceptron partition function conditional on $\FF(t)$.
From this we deduce the
\hyperlink{proof:t.main.lbd.bdd}{lower bound in Theorem~\ref{t:main} for the case $\|u\|_\infty<\infty$}.

\item In Section~\ref{s:local.clt} we prove a local central limit theorem (Proposition~\ref{p:density.bound}) which is required for the calculations of Sections \ref{s:first.mmt}--\ref{s:second.mmt}.
\item In Section~\ref{s:conc.log.Z} we
 prove Propositions~\ref{p:poly.conc}, \ref{p:Z.eps.approx.Z}, and \ref{p:conc.smoothed}; and use these to \hyperlink{proof:t.main.conclusion}{conclude the proof of Theorem~\ref{t:main}}.
\item Lastly, in Appendix~\ref{s:AMP} we prove a gaussian resampling identity
(Lemma~\ref{l:resampling}) which is used in the conditional moment calculations of 
Sections \ref{s:first.mmt}--\ref{s:second.mmt}. We also give a heuristic review of the state evolution limit of AMP, which was rigorously established in earlier works \cite{MR2810285,MR3147441}.
\end{itemize}

\subsection*{Acknowledgements}
We are grateful to Andrew Lawrie, Joe Neeman, Elchanan Mossel, and Ofer Zeitouni for many helpful conversations. Research of S.N.\ is supported by SNSF grant 176918. Research of N.S.\ is supported by
NSF CAREER grant DMS-1940092 and NSF-Simons grant DMS-2031883.

\section{First moment conditional on AMP}
\label{s:first.mmt}

We consider the perceptron model \eqref{e:Z} with an independent copy $\Gprime$ of the disorder matrix $\bG$ --- this is clearly equivalent (in law) to the original model. The (random) weight of the configuration $J$ is
	\beq\label{e:weight}
	\SAT 
	\equiv \SAT_J(\Gprime)
	\equiv \exp \bigg\{ \bigg(\oneM,
	u\bigg(\f{\Gprime J}{N^{1/2}}
		\bigg)\bigg)\bigg\}\,,
	\eeq
where
$u\equiv\log U:\R\to[-\infty,0]$ is applied componentwise according to the convention \eqref{e:componentwise}. As in \eqref{e:Z}, the corresponding perceptron partition function is
	\beq\label{e:Z.G.prime}
	\bZ(\Gprime)
	\equiv
	\sum_J
	\SAT_J(\Gprime)\,.
	\eeq
Let $\bmag^{(s)}$ and $\bn^{(\ell)}$ be generated from the AMP iteration \eqref{e:perceptron.AMP.col} and \eqref{e:perceptron.AMP.row} with $\Gprime$ in place of $\bG$ (and with the same initial values
for $\bmag^{(0)}$, $\bn^{(0)}$, $\bmag^{(1)}$, $\bn^{(1)}$ as before). Then, similarly as in \eqref{e:tap.condition}, let
	\beq\label{e:tap.condition.prime}
	\Fprime(t)\equiv\sigma\bigg(
	\Big(\Gprime\bmag^{(s)},\bn^{(s+1)} : s\le t\Big),
	\Big( (\Gprime)^\st\bn^{(\ell)},\bmag^{(\ell+1)}
		: \ell \le t-1 \Big)
	\bigg)\,.
	\eeq
We emphasize that $\Fprime(t)$ in \eqref{e:tap.condition.prime} is defined with respect to $\Gprime$ while $\FF(t)$ in \eqref{e:tap.condition} was defined with respect to $\bG$. This section is organized as follows:
\begin{itemize}
\item In \S\ref{ss:review.AMP} we give a brief review of known results
\cite{MR2810285,MR3147441} on the state evolution limit of AMP.
\item In \S\ref{ss:first.moment.slice}
we decompose $\bZ(\Gprime)$
into two parts (see \eqref{e:Z.outside}):
one part $\bZ_\circ(\Gprime)$ roughly captures the contribution of configurations $J\in\set{-1,+1}^N$
which lie close to $\bmag^{(t)}$ in some sense
(see \eqref{e:N.circ}),
while $\bZ_\bullet(\Gprime)$ is the remainder of the partition function.
We then state the main result of this section,
Theorem~\ref{t:ubd.slice},
which gives the conditional first moment upper bound
for $\bZ_\circ(\Gprime)$.

\item In \S\ref{ss:single.config}
we state and prove Proposition~\ref{p:first.mmt}, which gives a conditional first moment upper bound
for a single configuration	$J\in\set{-1,+1}^N$.

\item In \S\ref{ss:entropy} we complete the  \hyperlink{proof:t.ubd.slice}{proof of 
	Theorem~\ref{t:ubd.slice}}.
We also supply some large deviations bounds,
Lemmas~\ref{l:AH.pi} and \ref{l:AH.vpi},
which will be used later to bound $\bZ_\bullet(\Gprime)$
(see Corollary~\ref{c:AH.bound} in \S\ref{ss:AH.conclusion}).
\end{itemize}
The bound from Theorem~\ref{t:ubd.slice} will be analyzed in Section~\ref{s:analysis.first.mmt}
to conclude the 
\hyperlink{proof:t.ubd}{proof of Theorem~\ref{t:ubd}}. Throughout this section, $U$ satisfies Assumption~\ref{a:bdd} and \ref{a:Lip}. 

\subsection{Review of AMP state evolution}\label{ss:review.AMP}

In this subsection we review the main results on approximate message passing (as introduced in \S\ref{ss:intro.AMP}) that will be used in our proofs. What follows is primarily based on \cite{MR2810285,MR3147441}. A more detailed review
(with heuristic derivations) is given in Section~\ref{s:AMP}.

\begin{dfn}[state evolution recursions]
\label{d:state} Let $(q,\psi)$ be as given by Proposition~\ref{p:fp}, and abbreviate $F\equiv F_q$. Let
	\beq\label{e:intro.lm.gm.one}
	\rho_1\equiv\lm_1 \equiv \bigg(\f1q\bigg)^{1/2}
		\E \Th(\psi^{1/2} Z) =0\,,\quad
	\mu_1\equiv\gm_1\equiv
	\bigg(\f{\alpha}{\psi}\bigg)^{1/2}
		\E F(q^{1/2} Z)\eeq
(cf.\ \eqref{e:product.with.first.iterate}). 
Next let $\xi,\xi'$ be independent standard gaussian random variables, and for $s\ge1$ let
	\begin{align}\nonumber
	\rho_{s+1} \equiv
	\rho(\mu_s) &\equiv 
	 \f1q 
	\E\bigg[ 
	\Th\bigg(\psi^{1/2}
	\Big\{
	\mu_s \xi + [1-(\mu_s)^2]^{1/2} \xi'
	\Big\}
	\bigg)
	\Th(\psi^{1/2}\xi)
	\bigg] \\
	\mu_{s+1}
	\equiv\mu(\rho_s)
	&\equiv
	\f{\alpha}{\psi}
	\E \bigg[
	F\bigg( q^{1/2}\Big\{
		\rho_s \xi + [1-(\rho_s)^2]^{1/2} \xi'
	\Big\} \bigg)
	F(q^{1/2}\xi)\bigg]
	\label{e:def.rho.mu}
	\end{align}
(cf.\ \eqref{e:rho.mu.two} and \eqref{e:rho.mu.three}).
Supposing that $\gm_1,\ldots,\gamma_{s-1}$ and $\lm_1,\ldots,\lm_{s-1}$ have been defined, we let
	\beq\label{e:intro.lm.gm.rec}
	\lm_s
	=\f{\rho_s-\Lambda_{s-1}}{(1-\Lambda_{s-1})^{1/2}}\,,
	\quad
	\gm_s=\f{\mu_s-\Gamma_{s-1}}
		{(1-\Gamma_{s-1})^{1/2}}
	\eeq
(cf.\ \eqref{e:lm.gm.recursion}), where we have used the abbreviations
	\beq\label{e:sum.of.squares}
	\Gamma_{s-1}\equiv \sum_{\ell\le s-1}
		(\gm_\ell)^2\,,\quad
	\Lambda_{s-1}\equiv \sum_{\ell\le s-1}
		(\lm_\ell)^2\,.\eeq
The above recursions are standard in the AMP literature, so we defer the explanations to Section~\ref{s:AMP}.  We will confirm in Lemma~\ref{l:AMP.well.def} that the recursions result in well-defined quantities for all $s\ge1$.\end{dfn}

We now explain how the constants given in Definition~\ref{d:state} describe the large-$N$ behavior of the AMP iteration.
To this end, we define the (deterministic) matrices
	\begin{align}\label{e:Gamma.matrix}
	\bGam
	&\equiv\begin{pmatrix}
	1 & \\
	\gm_1 & (1-\Gamma_1)^{1/2} \\
	\gm_1 & \gm_2 & (1-\Gamma_2)^{1/2} \\
	\vdots & & & \ddots \\
	\gm_1 & \gm_2 & \cdots && (1-\Gamma_{t-2})^{1/2}
	\end{pmatrix} \in\R^{(t-1)\times(t-1)} \,,\\
	\bLam
	&\equiv \begin{pmatrix}
	1 & \\
	\lm_1 & (1-\Lambda_1)^{1/2} & \\
	\lm_1 & \lm_2 & (1-\Lambda_2)^{1/2} \\
	\vdots &&& \ddots \\
	\lm_1 & \lm_2 & \cdots &  & (1-\Lambda_{t-1})^{1/2}
	\end{pmatrix}
	\in\R^{t\times t}\,.
	\label{e:Lambda.matrix}
	\end{align}
It will follow from Lemma~\ref{l:AMP.well.def} below that in our setting we will have $\Gamma_s\in[0,1)$ and $\Lambda_s\in[0,1)$ for all $s\ge0$, which implies that both $\bGam$ and $\bLam$ are non-singular matrices.  As in \eqref{e:tap.condition.prime}, let $\bmag^{(s)}$ and $\bn^{(\ell)}$ be generated from the AMP iteration \eqref{e:perceptron.AMP.col} and \eqref{e:perceptron.AMP.row} with $\Gprime$ in place of $\bG$. Recall that $\bmag^{(s)}\equiv\Th(\bH^{(s)})$ and $\bn^{(s)}\equiv F(\bh^{(s)})$, where $F=F_q$ is given by \eqref{e:F}. We define vectors $\by^{(s)}$
and $\bx^{(s)}$ by setting 
\begin{align}\label{e:gs.H.y}
	\f{\bH[t-1]}{\psi^{1/2}}
	\equiv \f1{\psi^{1/2}}
	\begin{pmatrix}
	(\bH^{(2)})^\st \\
	\vdots\\
	(\bH^{(t)})^\st \end{pmatrix}
	&\equiv \bGam\begin{pmatrix}
	(\by^{(1)})^\st \\ \vdots \\ (\by^{(t-1)})^\st
	\end{pmatrix}
	\equiv\bGam\by[t-1]
	\in\R^{(t-1)\times N} \,,\\
	\f{\bh[t]}{q^{1/2}}
	\equiv \f1{q^{1/2}}
	\begin{pmatrix}
	(\bh^{(2)})^\st \\
	\vdots\\
	(\bh^{(t+1)})^\st \end{pmatrix}
	&\equiv
	\bLam\begin{pmatrix}
	(\bx^{(1)})^\st \\ \vdots \\ (\bx^{(t)})^\st
	\end{pmatrix}
	\equiv \bLam\bx[t]
	\in\R^{t\times M}\,,
	\label{e:gs.h.x}
	\end{align}
for $\bGam$ and $\bLam$ as in \eqref{e:Gamma.matrix} and \eqref{e:Lambda.matrix}. Then the $\bx^{(s)}$ ``behave like'' i.i.d.\ standard gaussian vectors in $\R^M$, while the $\by^{(s)}$
``behave like'' i.i.d.\ standard gaussian vectors in $\R^N$. For an intuitive explanation we refer to the heuristic derivation of \eqref{e:amp.decomp.H.y} and \eqref{e:amp.decomp.h.x} given in 
in Section~\ref{s:AMP}. The formal version is given by the next definition and lemma: 

\begin{dfn}[pseudo-Lipschitz functions]
\label{d:PL}
Following \cite{MR2810285}, we say that a function $f:\R^\ell\to\R$ (where $\ell$ is any positive integer) is  \bemph{pseudo-Lipschitz of order $k$} if
there exists a constant $L>0$ such that
	\[
	\|f(x)-f(y)\| \le L 
		\bigg(1+\|x\|^{k-1}+\|y\|^{k-1}\bigg)
		\|x-y\|
	\]
for all $x,y\in\R^\ell$. We say for short that $f$ is a $\PL(k)$ function.
\end{dfn}

\begin{lem}[{\cite[Lem.~1]{MR2810285}}]
\label{l:AMP.PL} Suppose $U$ satisfies Assumptions~\ref{a:bdd} and \ref{a:Lip}. In particular, this guarantees that the function $F_q$ of \eqref{e:F} is Lipschitz (see Lemma~\ref{l:F.Lipschitz}). Let $\bG$ be an $M\times N$ matrix with i.i.d.\ standard gaussian entries, such that $M/N=\alpha$. Assume $0\le \alpha\le \alpha(U)$, and let $(q,\psi)$ be the solution given by Proposition~\ref{p:fp}. Then let $\bmag^{(s)}\equiv\Th(\bH^{(s)})$ and $\bn^{(\ell)}\equiv F_q(\bh^{(\ell)})$ be  generated from the AMP iteration \eqref{e:perceptron.AMP.col} and \eqref{e:perceptron.AMP.row}, with the same initial values
for $\bmag^{(0)}$, $\bn^{(0)}$, $\bmag^{(1)}$, $\bn^{(1)}$ as before. If $f:\R^{t-1}\to\R$ is a
$\PL(k)$ function, then
	\[
	\f1N\sum_{i\le N}
		f\Big((\bH[t-1]\eN_i)^\st\Big)
	\stackrel{N\to\infty}{\longrightarrow}
	\E f(\psi^{1/2}\bGam \bxi )
	\]
where $\bxi$ here denotes a standard gaussian vector in $\R^{t-1}$, and the convergence holds in probability as $N\to\infty$ for any fixed $t$. Likewise, if $f:\R^t\to\R$
is a $\PL(k)$ function, then
	\[
	\f1M\sum_{a\le M}
		f\Big((\bh[t]\eM_a)^\st\Big)
	\stackrel{N\to\infty}{\longrightarrow} 
	\E f(q^{1/2}\bLam\bxi)
	\]
where $\bxi$ here denotes a standard gaussian vector in $\R^t$.
\end{lem}

We remark that the results of \cite{MR2810285} are for a more general setting where the AMP iteration starts from a random initialization with bounded moments up to order $2k-2$; the result then holds for any $f$ which is $\PL(k)$. In this paper we start from an initialization with bounded moments of all finite orders, so in Lemma~\ref{l:AMP.PL} we can take $f$ to be in $\PL(k)$ for any finite $k$. We now present a few applications of Lemma~\ref{l:AMP.PL} which illustrate how some of the recursions from Definition~\ref{d:state} naturally arise. First, it follows from Lemma~\ref{l:AMP.PL} and the definition \eqref{e:def.rho.mu} that
	\[
	\f{(\bmag^{(r)},\bmag^{(s)})}{Nq}
	=\f{(\Th(\bH^{(r)}),\Th(\bH^{(s)}))}{Nq}
	\simeq
	\rho((\bGam\bGam^\st)_{r-1,s-1})\,.
	\]
In the above and throughout this paper, we write $f\simeq g$ to indicate that $f-g$ converges to zero in probability as $N\to\infty$. In the case $r=s$ we have
	\[
	(\bGam\bGam^\st)_{r-1,r-1}
	\stackrel{\eqref{e:Gamma.matrix}}{=} 
	\sum_{\ell\le r-2} (\gamma_\ell)^2
		+ (1-\Gamma_{r-2})
	\stackrel{\eqref{e:sum.of.squares}}{=}
	\Gamma_{r-2} +(1-\Gamma_{r-2})
	=1\,.
	\]
If $r\ne s$, we can suppose without loss that $r<s$, in which case
	\[
	(\bGam\bGam^\st)_{r-1,s-1}
	\stackrel{\eqref{e:Gamma.matrix}}{=} 
		\sum_{\ell\le r-2} (\gamma_\ell)^2
		+ \gamma_{r-1}(1-\Gamma_{r-2})^{1/2}
	\stackrel{\eqref{e:sum.of.squares}}{=}
		\Gamma_{r-2} + \gamma_{r-1} (1-\Gamma_{r-2})^{1/2}
	\stackrel{\eqref{e:intro.lm.gm.rec}}{=}
	\mu_{r-1}\,.
	\]
It follows that 
$\|\bmag^{(r)}\|^2 \simeq Nq$ for all $r$, and for $r<s$ we have
	\beq\label{e:m.scalar.products}
	\f{(\bmag^{(r)},\bmag^{(s)})}{Nq}
	\simeq
	\rho(\mu_{r-1})
	\stackrel{\eqref{e:def.rho.mu}}{=} \rho_r
	\stackrel{\eqref{e:intro.lm.gm.rec}}{=}
	\Lambda_{r-1} + \lm_r(1-\Gamma_{r-1})^{1/2}
	\stackrel{\eqref{e:Gamma.matrix}}{=} 
	(\bLam\bLam^\st)_{r,s}\,.
	\eeq
A similar calculation gives that
$\|\bn^{(r)}\|^2\simeq N\psi$ for all $r$, and
for $r<s$ we have
	\beq\label{e:n.scalar.products}
	\f{(\bn^{(r)},\bn^{(s)})}{N\psi}
	\simeq
	\mu(\rho_{r-1})
	=\mu_r = (\bGam\bGam^\st)_{r,s}\,.
	\eeq
Let
 $\br^{(s)}$ be the Gram--Schmidt orthogonalization of the vectors $\bmag^{(s)}$ for $s\ge1$: thus $\br^{(1)}=\bmag^{(1)}/\|\bmag^{(1)}\|=\oneN/N^{1/2}$,
	\[
	\br^{(2)}
	= \f{\bmag^{(2)} - (\bmag^{(2)},\br^{(1)})\br^{(1)}}
		{\|\bmag^{(2)} - (\bmag^{(2)},\br^{(1)})\br^{(1)}\|}\,,
	\]
and so on. The $\br^{(s)}$ form an orthonormal set in $N$-dimensional space (assuming the number of iterations is much smaller than the dimension). Likewise, let $\bc^{(s)}$ be the Gram--Schmidt orthogonalization of the vectors $\bn^{(s)}$ for $s\ge1$; these form an orthonormal set in $M$-dimensional space. Let $\bLam_N,\bGam_N$ be the (random) matrices such that
	\begin{align}\label{e:gs.m.r.EXACT}
	\f{\bmag[t]}{(Nq)^{1/2}}
	\equiv\f1{(Nq)^{1/2}}
	\begin{pmatrix}
	(\bmag^{(1)})^\st \\ \vdots \\ (\bmag^{(t)})^\st
	\end{pmatrix}
	&=\bLam_N
	\begin{pmatrix}
	(\br^{(1)})^\st \\ \vdots \\ (\br^{(t)})^\st
	\end{pmatrix}
	\equiv \bLam_N\br[t]\in\R^{t\times N}
	\,,\\
	\f{\bn[t-1]}{(N\psi)^{1/2}}
	\equiv\f1{(N\psi)^{1/2}}
	\begin{pmatrix} (\bn^{(1)})^\st\\
		\vdots \\ (\bn^{(t-1)})^\st
		\end{pmatrix}
	&=\bGam_N \begin{pmatrix}
	(\bc^{(1)})^\st\\
	\vdots \\
	(\bc^{(t-1)})^\st \end{pmatrix}
	\equiv \bGam_N\bc[t-1]
	\in\R^{(t-1)\times M}
	\,.
	\label{e:gs.n.c.EXACT}
	\end{align}
It can be deduced from \eqref{e:m.scalar.products} and \eqref{e:n.scalar.products} that
	\beq\label{e:Lm.Gm.N.conv}
	\begin{pmatrix}\bLam_N\\\bGam_N\end{pmatrix}
	\simeq\begin{pmatrix}\bLam \\ \bGam
		\end{pmatrix}\,.
	\eeq
(This means $\bLam_N-\bLam$ and $\bGam_N-\bGam$ converge entrywise to zero, in probability,
as $N\to\infty$.) Since $\br[t]$ and $\bc[t-1]$ have orthonormal rows, the above implies
	\begin{align*}
	\f{\bmag[t]\bmag[t]^\st}{Nq}
	\stackrel{\eqref{e:gs.m.r.EXACT}}{=}
	\bLam_N\br[t]\br[t]^\st(\bLam_N)^\st
	= \bLam_N(\bLam_N)^\st 
	&\simeq \f{\bh[t]\bh[t]^\st}{N\alpha q}
	\in\R^{t\times t}
	\,,\\
	\f{\bn[t-1]\bn[t-1]^\st}{N\psi}
	\stackrel{\eqref{e:gs.n.c.EXACT}}{=}
	\bGam_N\bc[t-1]\bc[t-1]^\st (\bGam_N)^\st
	= \bGam_N(\bGam_N)^\st
	&\simeq \f{\bH[t]\bH[t]^\st}{N\psi}
	\in\R^{(t-1)\times(t-1)}\,,
	\end{align*}
where the approximations on the right-hand side use 
Lemma~\ref{l:AMP.PL}. This above of course consistent with the previous calculations
\eqref{e:m.scalar.products} and
 \eqref{e:n.scalar.products} (cf.\ \cite[eq.\ (3.18) and (3.19)]{MR2810285}).

A further consequence
of Lemma~\ref{l:AMP.PL}
is that for all $k,\ell\ge1$ we have
	\begin{align*}
	\f{(\bmag^{(k+1)},\by^{(\ell)})}{Nq^{1/2}}
	&=
	\f{(\Th(\bH^{(k+1)}),\by^{(\ell)})}{Nq^{1/2}}
	\stackrel{\eqref{e:gs.H.y}}{=}
	\f1{Nq^{1/2}}
	\bigg(
	\Th\bigg(
	\psi^{1/2}\bigg\{
	\sum_{\ell'\le t-1}
	\Gamma_{k,\ell'} \by^{(\ell')}
	\bigg\} \bigg), \by^{(\ell)}
	\bigg) \\
	&\simeq
	\f{\Gamma_{k,\ell}}{q^{1/2}}
	 \E\Big[ Z\Th(\psi^{1/2} Z)\Big]
	=\f{\Gamma_{k,\ell}}{q^{1/2}}
	\psi^{1/2}
	\E\Big[ \Th'(\psi^{1/2} Z)\Big]
	\stackrel{\eqref{e:fp}}{=}
	\f{\Gamma_{k,\ell}}{q^{1/2}}
	\psi^{1/2} (1-q)
	\,,\end{align*}
having used the gaussian integration by parts identity.
Recall also that
$\bmag^{(1)}=q^{1/2}\oneN$, so Lemma~\ref{l:AMP.PL} also implies
	\[\f{(\bmag^{(1)},\by^{(\ell)})}{Nq^{1/2}}
	\simeq \E\xi = 0\]
for all $\ell\le t-1$, where $\xi$ is a standard gaussian random variable. The above calculations can be summarized as
	\beq\label{e:m.dot.y}
	\bigg\|\f{\by[t-1]\bmag[t]^\st}{Nq^{1/2}}
	- \begin{pmatrix}
	\mathbf{0} & \displaystyle
	\f{\psi^{1/2}}{q^{1/2}}(1-q) \bGam^\st
	\end{pmatrix}
		\bigg\|_\infty
	\le \ERR_{t,1} \simeq0\,,
	\eeq
where $\mathbf{0}$ denotes the zero vector in $t-1$ dimensions, and $\ERR_{t,1}$ is an an $\FF(t)$-measurable random variable that converges to zero in probability as $N\to\infty$ 
(cf.\ \cite[eq.\ (3.20) and (3.21)]{MR2810285}).
This concludes our review of the required results on the state evolution of AMP, and we turn next to the conditional moment calculations. We introduce some notation which will be used later in the paper:

\begin{rmk}[bounds on $\bLam_N$ and $\bGam_N$]
\label{r:singval.bound}
Since $\bLam$ and $\bGam$ are both non-singular (this will be verified in Lemma~\ref{l:AMP.well.def} below), we can define a large finite constant $\varsigma_t$ such that we have the bound
	\beq\label{e:varsigma.t}
	\max\bigg\{ \|\bLam_N\|_\infty,
		\|(\bLam_N)^{-1}\|_\infty,
	\|\bGam_N\|_\infty, \|(\bGam_N)^{-1}\|_\infty
	\bigg\} \le \bigg(\f{\varsigma_t}{t}\bigg)^{1/2}
	\eeq
with high probability.  In the above, and throughout this paper, $\|\cdot\|_\infty$ denotes the entrywise maximum absolute value of a vector or matrix. On the other hand, we write $\|u\|$ for the euclidean norm of a vector $u$, and $\|A\|$ for the spectral norm a matrix $A$. It follows from \eqref{e:varsigma.t} that we also have
	\[
	\max\bigg\{ \|\bLam_N\|,
		\|(\bLam_N)^{-1}\|,
	\|\bGam_N\|, \|(\bGam_N)^{-1}\|
	\bigg\} \le (\varsigma_t)^{1/2}\]
with high probability. 
\end{rmk}

The \hyperlink{proof:p.AT}{proof of the following proposition} is deferred to \S\ref{ss:logc}. It amounts to checking that an Almeida--Thouless (AT) condition
(\cite{AT1978stability}; see Lemma~\ref{l:AT.conclusion})
 is satisfied.

\begin{ppn}
\label{p:AT}
Suppose $U$ satisfies Assumptions~\ref{a:bdd} and \ref{a:Lip}. 
For $0<\alpha\le\alpha(U)$ as defined by \eqref{e:alpha.U}, the state evolution recursions
from Definition~\ref{d:state}
result in $\Gamma_t\to1$ and $\Lambda_t\to1$
as $t\to\infty$.
\end{ppn}

\subsection{Positions of configurations relative to AMP iterates}
\label{ss:first.moment.slice}

We now define parameters $\pi(J)$ and $\vpi(J)$
which summarize the position of configurations $J\in\set{-1,+1}^N$ relative to the vectors $\br^{(s)}$ and $\by^{(\ell)}$ from \eqref{e:gs.m.r.EXACT} and \eqref{e:gs.H.y}.

\begin{dfn}[parameters $\pi$ and $\vpi$]\label{d:pi.vpi}
Let $\Fprime(t)$ be as in \eqref{e:tap.condition.prime}. For $J\in\set{-1,+1}^N$, define
	\begin{align}\label{e:def.pi}
	\pi(J)
	&\equiv \f{\br[t]J}{N^{1/2}}
	= \bigg(\f{(\br^{(s)},J)}{N^{1/2}} \bigg)_{s\le t}
	\in\R^t\,,\\
	\vpi(J)
	&\equiv \f{\by[t-1]J}{N}
	= \bigg(\f{(\by^{(\ell)},J)}{N}\bigg)_{\ell\le t-1}
	\in\R^{t-1}\,.
	\label{e:def.vpi}
	\end{align}
Note that for any given $J\in\set{-1,+1}^N$, its parameters $\pi(J)$ and $\vpi(J)$ are measurable with respect to $\Fprime(t)$.
\end{dfn}

Recall that the vectors $\br^{(s)}$ and $\bmag^{(s)}$ ($1\le s\le t$) are linearly related by \eqref{e:gs.m.r.EXACT}, while the vectors $\by^{(\ell+1)}$ and $\bH^{(\ell)}$  ($1\le\ell\le t-1$) are linearly related by \eqref{e:gs.H.y}. For part of our calculation it is more convenient to work with $\bmag^{(s)}$ and $\bH^{(\ell+1)}$ rather than with $\br^{(s)}$ and $\by^{(\ell)}$. For this reason we also define the following parameters:

\begin{dfn}[parameters $\hpi$ and $\delta$]
\label{d:J.params}
Given $\Fprime(t)$ as in \eqref{e:tap.condition.prime},
and given any $J\in\set{-1,+1}^N$,
we decompose $J$ as $J=J'+J''$ where $J'$ is the orthogonal projection of $J$ onto the span of the vectors $\bmag^{(s)}$, $1\le s\le t$. We let $\hpi_s$ for $1\le s\le t$ be the coefficients such that
	\beq\label{e:def.hpi}
	J' = \sum_{s\le t}\hpi_s
		\f{\bmag^{(s)}}{q^{1/2}}
	=\f{\bmag[t]^\st\hpi}{q^{1/2}}
	\,.
	\eeq
Next let $\bv\equiv J''/\|J''\|$, and let
$\delta\in\R^{t-1}$ be defined by
	\beq\label{e:defn.delta}
	\bGam_N(\bGam_N)^\st\delta
	=\f{\bH[t-1]\bv}{(N\psi)^{1/2}}\,.\eeq
Note that for any given $J\in\set{-1,+1}^N$, its parameters $\hpi(J)$ and $\delta(J)$ are measurable with respect to $\Fprime(t)$.
\end{dfn}

The parameters 
$(\pi,\vpi)$ of
Definition~\ref{d:pi.vpi}
are related as follows to the parameters $(\hpi,\delta)$ of Definition~\ref{d:J.params}:

\begin{lem}[change of basis]
\label{l:change.of.basis.EXACT}
Given $\Fprime(t)$ as in \eqref{e:tap.condition.prime},
suppose $J\in\set{-1,+1}^N$ has parameters 
$\pi(J)$, $\vpi(J)$, $\hpi(J)$, $\delta(J)$ as in 
Definitions~\ref{d:pi.vpi} and \ref{d:J.params}.
Then we have $\pi(J) = (\bLam_N)^\st\hpi(J)$, and
 	\[
	\vpi(J) = \f{\by[t-1]\bmag[t]^ \st}{N q^{1/2}}
		\hpi(J)
	+\Big(1-\|\pi(J)\|^2\Big)^{1/2}
	\bGam^{-1}\bGam_N(\bGam_N)^\st
	\delta(J)\,.\]

\begin{proof}
For convenience we will often abbreviate $\pi\equiv\pi(J)$, etc. The expression \eqref{e:def.hpi} can be rewritten as
	\[
	\f{J'}{N^{1/2}}
	\stackrel{\eqref{e:def.hpi}}{=}
	\f{\bmag[t]^\st\hpi}{(Nq)^{1/2}}
	\stackrel{\eqref{e:gs.m.r.EXACT}}{=}
	\br[t]^\st (\bLam_N)^\st\hpi
	\,,\]
so by comparing with \eqref{e:def.pi} we see that
$\pi(J)=(\bLam_N)^\st\hpi(J)$. Next
we have
	\beq\label{e:H.dot.J.prime.EXACT}
	\f{\bH[t-1]J'}{N\psi^{1/2}}
	\stackrel{\eqref{e:def.hpi}}{=}
	\f{\bH[t-1]\bmag[t]^\st\hpi }{N(\psi q)^{1/2}}
	\stackrel{\eqref{e:gs.H.y}}{=}
	\f{\bGam\by[t-1]\bmag[t]^ \st\hpi}{N q^{1/2}}\,.
	\eeq
It is clear from \eqref{e:def.pi}
that $\|J'\|/N^{1/2}=\|\pi\|$, and since $\bv\equiv J''/\|J''\|$, it follows that
	\beq\label{e:H.dot.J.pprime}
	\f{\bH[t-1]J''}{N\psi^{1/2}}
	=\f{\|J''\|}{N^{1/2}}
	\cdot \f{\bH[t-1]\bv}{(N\psi)^{1/2}}
	= \Big(1-\|\pi\|^2\Big)^{1/2}
		\f{\bH[t-1]\bv}{(N\psi)^{1/2}}
	\stackrel{\eqref{e:defn.delta}}{=}
	\Big(1-\|\pi\|^2\Big)^{1/2}
		\bGam_N(\bGam_N)^\st\delta\,.
	\eeq
Combining \eqref{e:gs.H.y}, \eqref{e:H.dot.J.prime.EXACT},
and \eqref{e:H.dot.J.pprime} gives
	\begin{align*}
	\vpi(J)
	&\stackrel{\eqref{e:def.vpi}}{=}
	\f{\by[t-1]J}{N}
	\stackrel{\eqref{e:gs.H.y}}{=}
	\f{\bGam^{-1}\bH[t-1]J}{N\psi^{1/2}}
	\stackrel{\eqref{e:H.dot.J.pprime}}{=}
	\f{\bGam^{-1}\bH[t-1]J'}{N\psi^{1/2}}
	+\Big(1-\|\pi\|^2\Big)^{1/2}
	\bGam^{-1}\bGam_N(\bGam_N)^\st \\
	&\stackrel{\eqref{e:H.dot.J.prime.EXACT}}{=}
	\f{\by[t-1]\bmag[t]^ \st}{N q^{1/2}}\hpi
	+\Big(1-\|\pi\|^2\Big)^{1/2}
	\bGam^{-1}\bGam_N(\bGam_N)^\st
	\delta\,.
	\end{align*}
This concludes the proof.
\end{proof}
\end{lem}

\begin{lem}[approximate change of basis]
\label{l:change.of.basis.APPROX}
Given $\Fprime(t)$ as in \eqref{e:tap.condition.prime},
suppose again that $J\in\set{-1,+1}^N$ has parameters 
$\pi(J)$, $\vpi(J)$, $\hpi(J)$, $\delta(J)$ as in 
Definitions~\ref{d:pi.vpi} and \ref{d:J.params}.
Define also $\dpi(J)\equiv\bLam^\st\hpi(J)$ and
	\beq\label{e:vpi.exact}
	\dvpi(J)\equiv
	(\bGam_N)^\st
	\bigg\{
	\f{\psi^{1/2}}{q^{1/2}}(1-q)\cpi(J)
	+ \Big(1-\|\pi(J)\|^2\Big)^{1/2}\delta(J)
	\bigg\}\,,
	\eeq
where $\cpi\equiv\cpi(J)\equiv(\hpi_2,\ldots,\hpi_t)\in\R^{t-1}$. Then
	\[
	\max\bigg\{
	\Big\|\pi(J)-\dpi(J)\Big\|_\infty+
	\Big\|\vpi(J)-\dvpi(J)\Big\|_\infty
	: J\in\set{-1,+1}^N
	\bigg\}
	\le\ERR_{t,2}\,,\]
where $\ERR_{t,2}$ is an an $\Fprime(t)$-measurable random variable that converges to zero in probability as $N\to\infty$.

\begin{proof} 
It follows trivially from the definition \eqref{d:pi.vpi} and the Cauchy--Schwarz inequality that
	\[
	\|\pi(J)\|_\infty
	\le \max\bigg\{
		\f{\|\br^{(s)}\| \cdot \|J\|}{N^{1/2}}
		: s\le t\bigg\}
	= 1\,,
	\]
where we emphasize that the bound clearly holds uniformly over all $J\in\set{-1,+1}^N$. Therefore
	\[\Big\|\pi(J)-\dpi(J)\Big\|_\infty
	\le \sup\bigg\{
		\bigg\| \bigg( (\bLam_N)^{-1}(\bLam_N-\bLam)
		\bigg)^\st u\bigg\|_\infty
		: \|u\|_\infty\le1\bigg\}\,.\]
The right-hand side above is $\Fprime(t)$-measurable and 
does not depend on $J$, and it follows from \eqref{e:Lm.Gm.N.conv} that it tends to zero in probability as $N\to\infty$.
Next, to compare $\vpi(J)$ with $\dvpi(J)$, we note that $\vpi(J)-\dvpi(J)$ can be expressed as $\textup{I}(J)+\textup{II}(J)$ where
	\begin{align*}
	\textup{I}(J)
	&\equiv\bigg\{
	\f{\by[t-1]\bmag[t]^ \st}{N q^{1/2}}
	- \begin{pmatrix}
	\mathbf{0} & \displaystyle
	\f{\psi^{1/2}}{q^{1/2}}(1-q)
		(\bGam_N)^\st
	\end{pmatrix}
	\bigg\}
	(\bLam_N)^{-1}\pi(J)
	\,,\\
\textup{II}(J)
	&\equiv 
	\Big(1-\|\pi(J)\|^2\Big)^{1/2}
	\bGam^{-1}(\bGam_N-\bGam)(\bGam_N)^\st
	\delta(J)\,.	
	\end{align*}
Since $\|\pi(J)\|_\infty\le1$ as noted above, it follows using \eqref{e:Lm.Gm.N.conv} and \eqref{e:m.dot.y} that
$\|\textup{I}(J)\|_\infty$ can be bounded
 uniformly over $J$
 by an $\Fprime(t)$-measurable quantity that tends to zero in probability as $N\to\infty$. Next we note that \eqref{e:defn.delta} combined with the Cauchy--Schwarz inequality gives, for all $J$,
 	\[
	\Big\|\bGam_N(\bGam_N)^\st\delta(J)\Big\|_\infty
	\le
	\max \bigg\{
		\f{\|\bH^{(\ell)}\|}{N\psi^{1/2}} : \ell \le t-1\bigg\}\,.
	\]
The right-hand side above is $\Fprime(t)$-measurable,
and it can be deduced from Lemma~\ref{l:AMP.PL} that it converges in probability to $1$ as $N\to\infty$. It follows by combining with \eqref{e:Lm.Gm.N.conv}
that $\|\textup{II}(J)\|_\infty$ can also be bounded  uniformly over $J$ by an $\Fprime(t)$-measurable quantity that tends to zero in probability as $N\to\infty$. 
This proves the claim.\end{proof}
\end{lem}

We next use the AMP iteration to define a convenient change of measure on the discrete cube:

\begin{dfn}[change of measure]
\label{d:cube.P.Q}
Let $\bP$ denote the uniform probability measure on $\set{-1,+1}^N$, and let $\bQ$ be the probability measure on the same space which is given by
	\[
	\f{d\bQ}{d\bP}
	= \prod_{i\le N}
	\f{\exp((\bH^{(t)})_i J_i)}
		{ \Ch (\bH^{(t)})_i}
	= \f{\exp\{(\bH^{(t)},J)\}}
		{\exp\{(\oneN,\log\Ch\bH^{(t)})\}}\,.\]
If $J$ is sampled from the measure $\bQ$, its expected value is exactly $\Th(\bH^{(t)})=\bmag^{(t)}$. We now compute the expected values under $\bQ$ of the parameters from Definition~\ref{d:pi.vpi}. First we note that
	\beq\label{e:def.pi.star}
	\stardpi \equiv
	\f{\br[t]\bmag^{(t)}}{N^{1/2}}
	= \f{\br[t]\bmag[t]^\st \hat{e}_t}{N^{1/2}}
	\stackrel{\eqref{e:gs.m.r.EXACT}}{=} 
	q^{1/2} \br[t]\br[t]^\st
		(\bLam_N)^\st\hat{e}_t
	= q^{1/2} (\bLam_N)^\st\hat{e}_t\,,
	\eeq
where $\hat{e}_s$ denotes the $s$-th standard basis vector in $\R^t$. Let us define also $\starpi\equiv  q^{1/2} \bLam^\st\hat{e}_t$, and note that $\starpi\simeq\stardpi$ by \eqref{e:Lm.Gm.N.conv}.
Next we note that 
	\beq\label{e:def.vpi.star}
	\stardvpi
	\equiv\f{\by[t-1]\bmag^{(t)}}{N}
	\stackrel{\eqref{e:m.dot.y}}{\simeq}
	\psi^{1/2}(1-q)(\bGam^\st\acute{e}_{t-1})
	\equiv \starvpi
	\in\R^{t-1}\,,
	\eeq
where $\acute{e}_\ell$ denotes the $\ell$-th standard basis vector in $\R^{t-1}$.
\end{dfn}

Recalling \eqref{e:Z} and \eqref{e:weight}, we now define
	\beq\label{e:N.circ}
	\bm{N}_\circ
	\equiv \bigg\{
	(\pi,\vpi)
	: \max\Big\{
	\|\pi(J)-\starpi\|,
	\|\vpi(J)-\starvpi\|
	\Big\}
	\le 16\cdot C_1(U)
	\alpha^{1/2}
	\bigg\}\,,\eeq
where the constant $C_1(U)$ 
comes from Lemma~\ref{l:poly} below. We also let
	\beq\label{e:H.circ}
	\mbH_\circ \equiv
	\bigg\{J\in\set{-1,+1}^N
	: (\pi(J),\vpi(J))\in \bm{N}_\circ\bigg\}\,,
	\eeq
and we let $\mbH_\bullet\equiv
\set{-1,+1}^N\setminus\mbH_\circ$.
We then decompose
$\bZ(\Gprime)=\bZ_\circ(\Gprime)+\bZ_\bullet(\Gprime)$ where
	\beq\label{e:Z.outside}
	\bZ_\circ(\Gprime)
	\equiv \sum_{J\in\mbH_\circ}
		\SAT_J(\Gprime)\,,\quad
	\bZ_\bullet(\Gprime)
	\equiv \sum_{J\in\mbH_\bullet}
		\SAT_J(\Gprime)\,,
	\eeq
The main result of this section is as follows:

\begin{thm}\label{t:ubd.slice}
Suppose $U$ satisfies Assumptions~\ref{a:bdd} and \ref{a:Lip}, and let $\Fprime(t)$ be as in \eqref{e:tap.condition.prime}. Given $\bareps\in\R$, define
	\[
	\bX(\pi,\vpi)
	\equiv
	\bx[t]^\st\starpi
	+\bigg\{ \bx[t]^\st(\pi-\starpi)
	+ N^{1/2}\bareps\bc[t-1]^\st
		(\vpi-\starvpi)
		\bigg\}
	\in\R^M\,,
	\]
for $\starpi$ and $\starvpi$ as in Definition~\ref{d:cube.P.Q}. (The parameter $\bareps$ will be fixed later
in \eqref{e:eps.U}.)
Then define
	\[\Psi(\pi,\vpi) 
	\equiv
	\f{\|\vpi
	-\bareps(\vpi-\starvpi)\|^2}
		{2(1-\|\pi\|^2)}
	- \f{(\starvpi,\vpi)}{1-q}
	+\f1N \sum_{a\le M} L_{\|\pi\|^2}
	(\bX_a(\pi,\vpi))\,.\]
If $\bQ$ is the measure on $\set{-1,+1}^N$ 
from Definition~\ref{d:cube.P.Q}, then we have
	\[
	\f{\E(\bZ_\circ(\Gprime) \,|\,\Fprime(t))}
		{\exp\{(\oneN,\log(2 \Ch(\bH^{(t)})))\}}
	\le 
	\sum_{J\in\mbH_\circ}
	\bQ(J) \exp\bigg\{N\bigg[
		\Psi(\pi(J),\vpi(J))
	+ \ERR_{t,3}
	\bigg]\bigg\}\,,
	\]
where $\ERR_{t,3}$ is an an $\Fprime(t)$-measurable random variable that converges to zero in probability as $N\to\infty$.
\end{thm}

The \hyperlink{proof:t.ubd.slice}{proof of Theorem~\ref{t:ubd.slice}} is given in \S\ref{ss:entropy}.

\subsection{First moment for a single configuration}
\label{ss:single.config} The main result of this subsection is the following:

\begin{ppn}\label{p:first.mmt}
Suppose $U$ satisfies Assumption~\ref{a:bdd} and \ref{a:Lip}, and let $\Fprime(t)$ be as in \eqref{e:tap.condition.prime}. Define
	\[\cA(\pi,\dpi,\dvpi,\theta)
	\equiv \f{\|\dvpi-\theta\|^2}{2(1-\|\pi\|^2)}
	+\f1N\sum_{a\le M}L_{\|\pi\|^2}
		\bigg(
		\bx[t]^\st\dpi
		+ N^{1/2}\bc[t-1]^\st\theta
		\bigg)\,,\]
where the function $L$ is defined by \eqref{e:L}.
Recall $\SAT_J(\Gprime)$ from \eqref{e:weight}.
There exists a finite constant $\wp_{t,1}$ such that for any large finite constant $\theta_{\max}$, it holds with probability $1-o_N(1)$ that
	\[
	\f1N
	\log \E\Big(\SAT_J(\Gprime)\,\Big|\,\Fprime(t)\Big)
	\le
	\inf\bigg\{
	\cA\Big(\pi(J),
		\dpi(J),\dvpi(J),\theta\Big) 
	: \|\theta\| \le \theta_{\max}\bigg\}
	+ \f{\wp_{t,1}}{N}
	\]
uniformly over all $J\in\set{-1,+1}^N$ with $\|\pi(J)\|\le 4/5$.
\end{ppn}

The \hyperlink{proof:p.first.mmt}{proof of Proposition~\ref{p:first.mmt}} is given at the end of this subsection.

\begin{dfn}[row and column subspaces]\label{d:row.col}
Given $\Fprime(t)$ as in \eqref{e:tap.condition.prime}, define the linear subspaces
	\begin{align*}
	V_\ROW \equiv V_\ROW(t)
		&\equiv \spn\bigg\{
	\eM_a(\bmag^{(s)})^\st
	: 1\le a\le M, 1\le s\le t\bigg\}\,,\\
	V_\COL \equiv V_\COL(t-1)
		&\equiv \spn\bigg\{
	\bn^{(\ell)}(\eN_i)^\st
	: 1\le i\le N, 1\le \ell\le t-1\bigg\}\,.
	\end{align*}
Let $V_{\ROW\COL}\equiv V_\ROW+V_\COL$.
Let $\proj_\ROW$ denote the orthogonal projection
onto $V_\ROW$, and define analogously $\proj_\COL$
and $\proj_{\ROW\COL}$.
Note that $(\Gprime)_{\ROW\COL}\equiv\proj_{\ROW\COL}(\Gprime)$
is measurable with respect to $\Fprime(t)$.
\end{dfn}

\begin{dfn}[row and column events]
\label{d:R.C.events}
We now let $\bG$ be an independent copy of $\Gprime$, and define
	\begin{align}\label{e:row}
	\ROW
	&\equiv\Big\{\proj_\ROW(\bG)=(\Gprime)_\ROW\Big\}
	=\bigg\{
	\f{\bG\bmag^{(s)}}{N^{1/2}} 
		= \bh^{(s+1)}+
		\acute{\beta}\bn^{(s-1)}
		\textup{ for all $1\le s\le t$}
	\bigg\}\,,\\
	\COL
	&\equiv
	\Big\{\proj_\COL(\bG)=(\Gprime)_\COL\Big\}
	=\bigg\{
	\f{\bG^\st\bn^{(\ell)}}{N^{1/2}}
	=\bH^{(\ell+1)}+\beta\bmag^{(\ell-1)}
	\textup{ for all $1\le \ell\le t-1$}
	\bigg\}\,.
	\label{e:col}
	\end{align}
We shall refer to $\ROW$ as the \bemph{row event} (since it constrains the rows of the matrix $\bG$). Likewise we shall refer to $\COL$ as the \bemph{column event}.
\end{dfn}
Our calculation is based on the following resampling principle
(\hyperlink{proof:l.resampling}{proved in \S\ref{ss:resampling}}):

\begin{lem}[resampling]\label{l:resampling}
Let $\Fprime(t)$ be as in \eqref{e:tap.condition.prime}. If $f:\R^{M\times N}\to\R$ is any bounded measurable function, then
	\[
	\E\Big( f(\Gprime)\,\Big|\,\Fprime(t)\Big)
	= \E\bigg( f(\bG)
	\,\bigg|\, \ROW,\COL,  (\Gprime)_{\ROW\COL}
	\bigg)
	\]
where $\bG$ denotes an independent copy of $\Gprime$;
and the events $\ROW$ and $\COL$ are defined by \eqref{e:row} and \eqref{e:col}.\end{lem}

\begin{dfn}[configuration-dependent subspaces]
\label{d:prf.adm.subspaces}
Given $\Fprime(t)$ as in \eqref{e:tap.condition.prime},  and $J\in\set{-1,+1}^N$, recall from Definition~\ref{d:J.params} that we decompose
$J=J'+J''$, and let $\bv\equiv J''/\|J''\|$.
We then define the linear subspaces
	\begin{align*}
	V_\PRF
	&\equiv\spn\bigg\{
	\eM_a\bv^\st : 1\le a\le M
	\bigg\}\,,\\
	V_\ADM
	&\equiv
	\spn\bigg\{
	\bn^{(\ell)}\bv^\st: 
		1\le \ell\le t-1
	\bigg\}
	\end{align*}
Note that $V_\ADM$ is a subspace of $V_\PRF$,
and is also a subspace of $V_\COL$.
Let $\proj_\ADM$ denote the orthogonal projection onto $V_\ADM$, and note that $(\Gprime)_\ADM\equiv\proj_\ADM(\Gprime)$ is measurable with respect to $\Fprime(t)$.
\end{dfn}

\begin{dfn}[admissibility event]
\label{d:ADM.event}
As in Definition~\ref{d:R.C.events}, let $\bG$ be an independent copy of $\Gprime$, and define
	\beq\label{e:adm}
	\ADM
	\equiv
	\Big\{\proj_\ADM(\bG)=(\Gprime)_\ADM\Big\}
	\stackrel{\eqref{e:col}}{=} 
	\bigg\{
	\f{\bn[t-1]
	\bG\bv}{N \psi^{1/2}}
	=\f{\bH[t-1]\bv}{(N\psi)^{1/2}}
	\bigg\}\,,
	\eeq
where the last identity holds assuming $\bG$ belongs to the event $\COL$ from \eqref{e:col}. Note that 
$\bH[t-1]\bv$ is determined by the
parameter $\delta(J)$ from Definition~\ref{d:J.params}.
We refer to $\ADM$ as the \bemph{admissibility event}, and note $\ADM\subseteq\COL$.
\end{dfn}

In the setting of the perceptron model, the calculation of Lemma~\ref{l:resampling} can be simplified as follows:

\begin{lem}[reduction of column constraints]
\label{l:adm}
If $h:\R^M\to\R$ is any bounded measurable function, then
	\[\E\bigg( h(\bG J)
	\,\bigg|\, \ROW,\COL,(\Gprime)_{\ROW\COL}
	\bigg)
	=\E\bigg( h(\bG J)
	\,\bigg|\, \ROW,\ADM,(\Gprime)_{\ROW\ADM}
	\bigg)\]
where $\bG$ denotes an independent copy of $\Gprime$
and the events $\ROW,\COL,\ADM$ are defined by \eqref{e:row}, \eqref{e:col}, and \eqref{e:adm}.

\begin{proof} Let $V_{\COL\setminus\ROW}$ be the orthogonal complement of $V_\ROW$ inside $V_\ROW+V_\COL$: that is,
	\[
	V_\ROW+V_\COL
	= V_\ROW\obot V_{\COL\setminus\ROW}\,,
	\]
where we use $\obot$ to denote the sum of two orthogonal vector spaces. Note that $V_\ADM$ is a subspace of $V_\COL$ which is orthogonal to $V_\ROW$, so it follows that $V_\ADM$ is also a subspace of $V_{\COL\setminus\ROW}$. 
Let $\proj_{\COL\setminus\ROW}$ denote the orthogonal projection onto $V_{\COL\setminus\ROW}$. 
Note that $V_\ADM$ is a subspace of $V_\PRF$, and $V_\PRF$ is orthogonal to $V_\ROW$. We claim that
	\beq\label{e:adm.orth}
	\proj_{\COL\setminus\ROW}(V_\PRF)
	=V_\ADM\,.\eeq
Since we already noted that $V_\ADM\subseteq V_{\COL\setminus\ROW}$, it suffices to show inclusion in the other direction. The space $V_\PRF$ is spanned by the elements $\eM_a\bv^\st$.
Let
$\bc^{(\ell)}$, $1\le\ell\le t-1$, be any orthonormal basis for the span of the vectors 
$\bn^{(\ell)}$, $1\le \ell\le t-1$. An orthonormal basis for $V_\ADM$ is then given by the matrices
$\bc^{(\ell)}\bv^\st$, $1\le\ell\le t-1$. 
On the other hand, the space $V_\COL$ is spanned by the elements $\bc^{(\ell)}(\eN_i)^\st$. We therefore have
	\begin{align*}
	&\bigg(
	\eM_a\bv^\st
	-\proj_\ADM\Big(\eM_a\bv^\st\Big)
	, \bc^{(\ell)}(\eN_i)^\st
	-\proj_\ROW\Big(\bc^{(\ell)}(\eN_i)^\st\Big)
	\bigg)
	= \bigg(
	\eM_a\bv^\st
	-\proj_\ADM\Big(\eM_a\bv^\st\Big)
	, \bc^{(\ell)}(\eN_i)^\st\bigg) \\
	&\qquad=(\bc^{(\ell)})_a \bv_i
		- \bigg( \sum_{k\le t-1}
		(\eM_a\bv^\st,\bc^{(k)}\bv^\st)
		\bc^{(k)}\bv^\st
		,\bc^{(\ell)}(\eN_i)^\st\bigg)
	=(\bc^{(\ell)})_a \bv_i-(\bc^{(\ell)})_a \bv_i
	= 0\,.
	\end{align*}
It follows that for any $G_\PRF \in V_\PRF$
we have $G_\PRF-\proj_\ADM(G_\PRF)$ orthogonal to $V_{\COL\setminus\ROW}$, which concludes the proof of \eqref{e:adm.orth}. It follows that $V_\PRF = V_\ADM\obot V_{\PRF\setminus\ADM}$ where $V_{\PRF\setminus\ADM}$ is the orthogonal complement of $V_\ADM$ inside $V_\PRF$, and $V_{\PRF\setminus\ADM}$ is orthogonal to $V_{\COL\setminus\ROW}$. As a result, if $\bG$ is an $M\times N$ matrix with i.i.d.\ standard gaussian entries, we can decompose 
$\bG_\PRF=\bG_\ADM+\bG_{\PRF\setminus\ADM}$ where
$\bG_{\PRF\setminus\ADM}=\proj_{\PRF\setminus\ADM}(\bG)$ is independent of $\bG_{\COL\setminus\ROW}$.
It follows that
	\begin{align*}
	&\E\bigg( h(\bG J)
	\,\bigg|\, \ROW,\COL,(\Gprime)_{\ROW\COL}
	\bigg)
	=\E\bigg( h(\bG_\ROW J'+ \bG_\PRF J'')
		\,\bigg|\, \ROW,\ADM,\COL,(\Gprime)_{\ROW\COL}\bigg)\\
	&\qquad=\E\bigg( h(\bG_\ROW J'+ 
		( \bG_\ADM+\bG_{\PRF\setminus\ADM} ) J'')
		\,\bigg|\, \ROW,\ADM,\COL
		,(\Gprime)_{\ROW\COL}
		\bigg) \\
	&\qquad=\E\bigg( h(\bG_\ROW J'+ 
		( \bG_\ADM+\bG_{\PRF\setminus\ADM} ) J'')
		\,\bigg|\, \ROW,\ADM,(\Gprime)_{\ROW\ADM}
		\bigg)
	= \E\bigg( h(\bG J)
	\,\bigg|\, \ROW,\ADM,(\Gprime)_{\ROW\ADM}
	\bigg)\,,
	\end{align*}
as claimed.
\end{proof}
\end{lem}

Further towards the proof of Proposition~\ref{p:first.mmt}, we record the following calculations:

\begin{lem}\label{l:G.R.J}For $J\in\set{-1,+1}^N$, recall the decomposition $J=J'+J''$, and define $\tbX_J\equiv \bG J'/N^{1/2}$. On the event $\ROW$ from \eqref{e:row}, we have
	\begin{align*}
	\tbX_J
	&=\f1{q^{1/2}}\bigg\{ 
		\bh[t]^\st\hpi(J)
	+\acute{\beta} \bn[t-1]^\st\cpi(J)
	\bigg\}\\
	&=\bx[t]^\st\dpi(J)
	+N^{1/2}\bc[t-1]^\st
	\bigg(\dvpi(J) - 
	\Big(1-\|\pi(J)\|^2\Big)^{1/2}
	(\bGam_N)^\st\delta
	\bigg)\,.\end{align*}
In the above,
$\pi(J)$ is given by Definition~\ref{d:pi.vpi};
 $\hpi(J)$
and $\delta(J)$
are given by Definition~\ref{d:J.params};
and $\cpi(J)$,
$\dpi(J)$, and
$\dvpi(J)$ are defined by
Lemma~\ref{l:change.of.basis.APPROX}.

\begin{proof} Fix $J$ and abbreviate $\pi\equiv\pi(J)$, etc. Conditional on the event $\ROW$ from \eqref{e:row}, we have
	\[
	\tbX_J
	\equiv \f{\bG J'}{N^{1/2}}
	\stackrel{\eqref{e:def.hpi}}{=}
	\sum_{s\le t}\hpi_s
		\f{\bG\bmag^{(s)}}{(Nq)^{1/2}}
	\stackrel{\eqref{e:row}}{=}
	\sum_{s\le t}\f{\hpi_s}{q^{1/2}}
	\Big(
	\bh^{(s+1)}+\acute{\beta}\bn^{(s-1)}\Big)\,.
	\]
Recall the notation \eqref{e:gs.h.x} and \eqref{e:gs.n.c.EXACT}, and also that $\bn^{(0)}\equiv\zroM\in\R^M$. Therefore the above can be rewritten as
	\[\tbX_J = \f{\bh[t]^\st\hpi}{q^{1/2}}
	+\f{\acute{\beta} \bn[t-1]^\st \cpi}{q^{1/2}}\,.
	\]
Combining with \eqref{e:gs.h.x} and \eqref{e:gs.n.c.EXACT} gives
	\[
	\tbX_J
	\stackrel{\eqref{e:gs.h.x}}{=}
	\bx[t]^\st\bLam^\st\hpi
	+\f{(1-q)\bn[t-1]^\st \cpi}{q^{1/2}}
	\stackrel{\eqref{e:gs.n.c.EXACT}}{=}
	\bx[t]^\st\bLam^\st\hpi
	+\f{N^{1/2}\psi^{1/2}(1-q)}{q^{1/2}}
		\bc[t-1](\bGam_N)^\st\cpi\,.
	\]
Recalling the notation of Lemma~\ref{l:change.of.basis.APPROX} gives, with $\dpi\equiv\bLam^\st\hpi$ and $\dvpi$ as in \eqref{e:vpi.exact},
	\[
	\tbX_J
	= \bx[t]^\st\dpi
	+N^{1/2}\bc[t-1]^\st
	\bigg(\dvpi - 
	\Big(1-\|\pi\|^2\Big)^{1/2}
	(\bGam_N)^\st\delta
	\bigg)\,.
	\]
This concludes the proof.
\end{proof}
\end{lem}

\begin{lem}\label{l:adm.cgf} Given $J\in\set{-1,+1}^N$, define the cumulant-generating function
	\[
	\tcK_J(\tau)
	\equiv \f1N\log\E\bigg[\exp\bigg\{ N^{1/2}\sum_{\ell\le t-1}
		\tau_\ell
		(\bc^{(\ell)},\bG\bv)
		\bigg\}
		\SAT_J(\bG) \,\bigg|\,\ROW\bigg]
	\]
for $\tau\in\R^{t-1}$. Then, with $L$ as in \eqref{e:L}, the function $\tcK_J$ satisfies
	\[\tcK_J(\tau) - \f{\|\tau\|^2}{2}
	= \f1N \bigg(\oneM,
		L_{\|\pi(J)\|^2}
		\bigg( \tbX_J
		+N^{1/2}\Big(1-\|\pi(J)\|^2\Big)^{1/2} \bc[t-1]^\st \tau
		\bigg)\bigg)
	\equiv \tcL_J(\tau)\,,\]
with
$\pi(J)$ as in Definition~\ref{d:pi.vpi} and
 $\tbX_J$ is as in Lemma~\ref{l:G.R.J}.

\begin{proof} 
Conditional on the event $\ROW$, it follows from Lemma~\ref{l:G.R.J} that $\bG J'/N^{1/2} = \tbX_J \equiv \tbX$. We also have
	\[ \f{\bG J''}{N^{1/2}}
	= \f{\|J''\|}{N^{1/2}}
	\bG \bv
	\equiv
	\Big(1-\|\pi\|^2\Big)^{1/2}
	\bxi
	\,,
	\]
where $\pi\equiv\pi(J)$, and
 $\bxi\equiv\bG\bv$ is distributed as an independent gaussian vector in $\R^N$. It follows that
	\[\tcK_J(\tau) =\f1N
	\sum_{a\le M}\log\E_\xi\bigg[ \exp\bigg\{
	N^{1/2}\sum_{\ell\le t-1}
	\tau_\ell(\bc^{(\ell)})_a\xi
	\bigg\}
	U\bigg( \tbX_a 
	+ \Big(1-\|\pi\|^2\Big)^{1/2}
		\xi\bigg)
	\bigg]\,,\]
where $\xi$ denotes a standard gaussian random variable. Making a change of variable gives
	\[\tcK_J(\tau)
	=
	\f{\|\tau\|^2}2
	+ \f1N\sum_{a\le M}
	\log\E_\xi U\bigg( 
	\tbX_a + \Big(1-\|\pi\|^2\Big)^{1/2}
		\bigg[\xi  +N^{1/2} 
		\sum_{\ell\le t-1} 
			\tau_\ell \bc^{(\ell)}_a
		\bigg]
	 \bigg)\,,\]
from which the result follows.
\end{proof}
\end{lem}

Having collected most of the necessary ingredients, we now prove the main result of this subsection. The proof requires one more slightly technical estimate which we defer to  Proposition~\ref{p:density.bound} in Section~\ref{s:local.clt}.

\begin{proof}[{\hypertarget{proof:p.first.mmt}{Proof of Proposition~\ref{p:first.mmt}}}]
With $\Fprime(t)$ as in \eqref{e:tap.condition.prime} and
 $\SAT_J(\Gprime)$ as in \eqref{e:weight}, let us abbreviate the quantity of interest as
 	\[
	E_J
	\equiv \E\Big(\SAT_J(\Gprime)\,\Big|\,\Fprime(t)\Big)\,.
	\]
By the resampling principle from Lemma~\ref{l:resampling}, we can express
	\[E_J 
	=\E\Big(\SAT_J(\bG)\,\Big|\,\ROW,\COL,
		(\Gprime)_{\ROW\COL}
		\Big)\,,\]
where $\bG$ is an independent copy of $\Gprime$, and $\ROW$ and $\COL$ are the row and column events of Definition~\ref{d:R.C.events}.
Applying Lemma~\ref{l:adm} then gives the further simplification
	\beq\label{e:E.J.given.A.R}
	E_J
	=\E\Big(\SAT_J(\bG)\,\Big|\,\ROW,\ADM,
		(\Gprime)_{\ROW\ADM}\Big)
	\,,\eeq
where $\ADM$ is the admissibility event defined by \eqref{e:adm}.

Let $V_\ROW$ be as in 
Definition~\ref{d:row.col}, and note that an orthonormal basis for $V_\ROW$ is given by the elements $\eM_a(\br^{(s)})^\st$ for $1\le a\le M$, $1\le s\le t$. Denote
	\[
	\bg_\ROW
	\equiv
	\bigg(
	(\bG,\eM_a(\br^{(s)})^\st)
	: 1\le a\le M, 1\le s\le t\bigg)
	\in\R^{Mt}\,.
	\]
Likewise let $V_\PRF$ and $V_\ADM$ be as in Definition~\ref{d:prf.adm.subspaces}: recall that $V_\PRF$ is orthogonal to $V_\ROW$, and $V_\ADM$ is a subpsace of $V_\PRF$. An orthonormal basis for $V_\PRF$ is given by the elements $\eM_a\bv^\st$ for $1\le a\le M$. Denote
	\beq\label{e:g.prf}
	\bg_\PRF
	\equiv
	\bigg(
	(\bG,\eM_a\bv^\st)
	: 1\le a\le M\bigg)
	= \bG\bv \in\R^M\,.
	\eeq
An orthonormal basis for $V_\ADM$ is given by the elements $\bc^{(\ell)}\bv^\st$ for $1\le\ell\le t-1$, and we shall denote
	\beq\label{e:g.adm}
	\bg_\ADM
	\equiv
	\bigg(
	(\bG,\bc^{(\ell)}\bv^\st)
	: 1\le \ell\le t-1\bigg)
	= \bc[t-1] \bG\bv
	\in\R^{t-1}\,.
	\eeq
Lastly, as in the proof of Lemma~\ref{l:adm},
let $V_{\PRF\setminus\ADM}$ be the orthogonal complement of $V_\ADM$ inside $V_\PRF$. Choose an orthonormal basis for $V_{\PRF\setminus\ADM}$, and denote it $\bB_j$ for $1\le j\le M-(t-1)$. We then let
	\beq\label{e:g.prf.minus.adm}
	\bg_\BBB
	\equiv
	\bigg(
	(\bG,\bB_j)
	: 1\le j\le M-(t-1)\bigg)
	\in\R^{M-t+1}\,.
	\eeq
Note that there is an orthogonal transformation of $\R^M$ which maps $\bg_\PRF$ to the pair $(\bg_\ADM,\bg_\BBB)$. In what follows we let $p_\ROW$ denote the probability density function for $\bg_\ROW$, so
	\beq\label{e:gaus.density.r}
	p_\ROW(g_\ROW)
	=\f1{(2\pi)^{Mt/2}}
	\exp\bigg\{-\f{\|g_\ROW\|^2}{2}
	\bigg\}\,.\eeq
Likewise let $p_\ADM$ and $p_\BBB$ denote the densities for $\bg_\ADM$ and $\bg_\BBB$ respectively. Since the three subspaces  $V_\ROW$, $V_\ADM$, and $V_\BBB$ are mutually orthogonal, the joint density of $(\bg_\ROW,\bg_\ADM,\bg_\BBB)$ is simply the product $p_\ROW(g_\ROW) p_\ADM(g_\ADM) p_\BBB(g_\BBB)$.

The weight $\SAT_J(\bG)$, as defined by \eqref{e:weight}, is a function of $\bG J$, which we decomposed in the proof of Lemma~\ref{l:adm.cgf} as a sum of $\bG J'$ and $\bG J''$.
 Note that $\bG J'$ is a function of $\bg_\ROW$,
while $\bG J''$ is a function of $\bg_\PRF$ which in turn is a function of $(\bg_\ADM,\bg_\BBB)$.
Thus \eqref{e:weight} can be rewritten as a function
$\bm{S}_J$ of
$( \bg_\ROW,\bg_\ADM,\bg_\BBB)$: explicitly,
	\[
	\SAT_J(\bG)
	=\prod_{a\le M}
	U\bigg( 
	\sum_{s\le t} 
	\f{(J,\br^{(s)}) }{N^{1/2}} (\bg_\ROW)_{a,s}
	+ \f{\|J''\|}{N^{1/2}} (\bg_\PRF)_a
	\bigg)
	\equiv\bm{S}_J(\bg_\ROW,\bg_\ADM,\bg_\BBB)\,.
	\]
On the event $\ROW$, the value of $\bg_\ROW$ is fixed to a value $\bar{g}_\ROW$:
	\[(\bar{g}_\ROW)_{a,s}
	= (\bG\br[t]^\st)_{a,s}
	\stackrel{\eqref{e:gs.m.r.EXACT}}{=}
	\bigg(\f{\bG\bmag[t]^\st ((\bLam_N)^\st)^{-1}}
		{(Nq)^{1/2}}\bigg)_{a,s}\,,
	\]
where the right-hand side can be computed from \eqref{e:row}. Likewise, on the event $\ADM$, the value of $\bg_\ADM$ is fixed to a value $\bar{g}_\ADM$. We then introduce a parameter $\tau\in\R^{t-1}$, and define
	\beq\label{e:S.tau.fn}
	\SAT_{J,\tau}(\bG)
	\equiv \bm{S}_{J,\tau}
	( \bg_\ROW,\bg_\ADM,\bg_\BBB)
	\equiv
	\bm{S}_J(\bg_\ROW,\bg_\ADM,\bg_\BBB)
	\exp\bigg\{
	N^{1/2} (\tau,\bg_\ADM)\bigg\}
	\,.
	\eeq
Then, for any $\tau\in\R^{t-1}$, we can rewrite \eqref{e:E.J.given.A.R} as
	\begin{align}\nonumber
	E_J 
	&= \E\Big(\SAT_J(\bG)
		\,\Big|\,\ROW,\ADM,(\Gprime)_{\ROW\ADM}
		\Big)
	= \E\bigg(
		\f{\bm{S}_{J,\tau}
		(\bg_\ROW,\bg_\ADM,\bg_\BBB)
		}{\exp(N^{1/2}(\tau,\bar{g}_\ADM))}
		\,\bigg|\,
		(\bg_\ROW,\bg_\ADM)=(\bar{g}_\ROW,\bar{g}_\ADM)
	\bigg)\\
	&= 
	\f1{\exp(N^{1/2}(\tau,\bar{g}_\ADM))}
	\int \bm{S}_{J,\tau}(\bar{g}_\ROW,\bar{g}_\ADM,g_\BBB)
	p_\BBB(g_\BBB)\,dg_\BBB\,.
	\label{e:E.J.integral}
	\end{align}
By contrast, the expected value of $\SAT_{J,\tau}$ given only the row constraints is
	\begin{align}\nonumber
	\bE_J(\tau\,|\,\bar{g}_\ROW)
	&\equiv
	 \E\Big(\SAT_{J,\tau}(\bG)\,\Big|\,\ROW,
	 	(\Gprime)_{\ROW}
		\Big)
	= \E\bigg(\bm{S}_{J,\tau}
		(\bg_\ROW,\bg_\ADM,\bg_\BBB)\,\bigg|\,
		\bg_\ROW=\bar{g}_\ROW\bigg)\\
	&=\int
	p_\ADM(g_\ADM) \int 
	\bm{S}_{J,\tau}(\bar{g}_\ROW,g_\ADM,g_\BBB)
	p_\BBB(g_\BBB)\,dg_\BBB\,g_\ADM
	= \exp(N\tcK_J(\tau))\,,
	\label{e:def.E.tau}
	\end{align}
which was computed in Lemma~\ref{l:adm.cgf} above.
We then let $\ap_{J,\tau}(\cdot\,|\,\bar{g}_\ROW)$ be the probability density function of $g_\ADM$ under the measure that is biased by $\SAT_{J,\tau}(\bG)$, conditional on the event $\ROW$, that is to say,
	\beq\label{e:ap.tau}
	\ap_{J,\tau}(g_\ADM\,|\,\bar{g}_\ROW)\,dg_\ADM
	\equiv \f{\E(\SAT_{J,\tau}(\bG)
		\Ind{\bg_\ADM \in dg_\ADM}\,|\,\ROW)}
		{\E(\SAT_{J,\tau}(\bG)\,|\,\ROW)}
	\equiv \f{p_\ADM(g_\ADM) }{\bE_J(\tau\,|\,\bar{g}_\ROW)}
	\int 
	\bm{S}_{J,\tau}(\bar{g}_\ROW,g_\ADM,g_\BBB)
	p_\BBB(g_\BBB)\,dg_\BBB
	\,dg_\ADM\,.\eeq
Then, for any $\tau\in\R^{t-1}$, we can rewrite \eqref{e:E.J.integral} as
	\beq\label{e:tilt.gaussian}
	E_J
	=\f{ \bE_J(\tau\,|\,\bar{g}_\ROW)
	\cdot \ap_{J,\tau}(\bar{g}_\ADM\,|\,\bar{g}_\ROW)
	}{\exp\{N^{1/2}(\tau,\bar{g}_\ADM)\}
	\cdot
	p_\ADM(\bar{g}_\ADM)}\,.
	\eeq	
We will show in Proposition~\ref{p:density.bound} 
(deferred to Section~\ref{s:local.clt}) that there is a finite constant $\wp_{t,0}$ such that
for any finite constant $\tau_{\max}$, we have the uniform bound
	\beq\label{e:density.bound}
	\max\bigg\{
	\Big\|
		\ap_{J,\tau}(\cdot\,|\,\bar{g}_\ROW)
		\Big\|_\infty
	: J\in\set{-1,+1}^N, \|\pi(J)\| \le \f45,
	\|\tau\| \le \tau_{\max}
	\bigg\}
	\le \wp_{t,0}
	\eeq
with high probability. 
It therefore remains to estimate the other two terms on the right-hand side of \eqref{e:tilt.gaussian}. We then note 
that Definition~\ref{d:ADM.event} implies that, on the event $\ADM$, we have
	\begin{align}\nonumber
	\f{\bar{g}_\ADM}{N^{1/2}}
	=\f{\bg_\ADM}{N^{1/2}}&
	\stackrel{\eqref{e:g.adm}}{=}
		\f{\bc[t-1]\bG\bv}{N^{1/2}}
	\stackrel{\eqref{e:gs.n.c.EXACT}}{=}
	\f{(\bGam_N)^{-1}\bn[t-1]\bG\bv}{N \psi^{1/2}}\\
	&\stackrel{\eqref{e:adm}}{=}
	\f{(\bGam_N)^{-1}\bH[t-1]\bv}{(N\psi)^{1/2}}
	\stackrel{\eqref{e:defn.delta}}{=}
	(\bGam_N)^\st\delta
	\label{e:adm.orth.SECOND}\,.\end{align}
Substituting \eqref{e:adm.orth.SECOND} into the formula for $p_\ADM$ (similar to \eqref{e:gaus.density.r}) gives
	\beq\label{e:p.adm.uncond}
	p_\ADM(\bar{g}_\ADM)
	=\f1{(2\pi)^{(t-1)/2}}
	\exp\bigg\{ -
		\f{N\|(\bGam_N)^\st\delta\|^2}{2} 
	\bigg\}\,.
	\eeq
Meanwhile, it follows by combining
\eqref{e:def.E.tau} and \eqref{e:adm.orth.SECOND} that
	\beq\label{e:remaining.E.factor}
	\f{ \bE_J(\tau\,|\,\bar{g}_\ROW)}
		{\exp\{N^{1/2}(\tau,\bar{g}_\ADM)\}}
	=\exp\bigg\{
	N \Big[
	\tcK_J(\tau) - (\tau,(\bGam_N)^\st\delta)
	\Big]
	\bigg\}\,.\eeq
Substituting \eqref{e:density.bound}, \eqref{e:p.adm.uncond}, and \eqref{e:remaining.E.factor}
into \eqref{e:tilt.gaussian} gives
	\[ \f{E_J}{(2\pi)^{t/2}\cdot
	\wp_{t,0}}
	\le
	\exp\bigg\{ N\bigg[
	\tcK_J(\tau) - (\tau,(\bGam_N)^\st\delta)
		+ \f{\|(\bGam_N)^\st\delta\|^2}{2}
		\bigg]\bigg\}\,.\]
Recalling the calculation of 
$\tcK_J(\tau)$ from Lemma~\ref{l:adm.cgf} gives
	\beq\label{e:tcA}
	\f{E_J}{(2\pi)^{t/2}\cdot
	\wp_{t,0}}
	\le
	\exp\bigg\{ N\bigg[
	\f{\|\tau-(\bGam_N)^\st\delta\|^2}{2}
	 +\tcL_J(\tau)
		\bigg]\bigg\}
	\equiv
	\exp\Big\{N\tcA_J(\tau)\Big\}\,,
	\eeq
where $\tcA_J$ is defined by the last identity.
To simplify the above expression, we will recenter $\tau$ around
	\beq\label{e:bar.tau}
	\bar{\tau}
	\equiv \bar{\tau}(\cpi)
	\equiv - \f{\psi^{1/2}(1-q)
		(\bGam_N)^\st\cpi}
		{q^{1/2}(1-\|\pi\|^2)^{1/2}}
	\stackrel{\eqref{e:vpi.exact}}{=}
	-\f{\dvpi}{(1-\|\pi\|^2)^{1/2}}
		+(\bGam_N)^\st\delta
	\,.\eeq
We then make a change of variables from $\tau$ to $\theta$, via the definition
	\beq\label{e:tau.theta}
	\tau\equiv \bar{\tau} + \f{\theta}{(1-\|\pi\|^2)^{1/2}}\,.
	\eeq
This change of variables results in the simplification
	\[
	\tau-(\bGam_N)^\st\delta
	\stackrel{\eqref{e:tau.theta}}{=}
	\bar{\tau} + \f{\theta}{(1-\|\pi\|^2)^{1/2}}
	-(\bGam_N)^\st\delta
	\stackrel{\eqref{e:bar.tau}}{=}
	\f{\theta-\dvpi}{(1-\|\pi\|^2)^{1/2}}\,.
	\]
The computation of $\tbX_J$ from Lemma~\ref{l:G.R.J} can also be rewritten as
	\beq\label{e:tbX.rewrite}
	\tbX_J
	\stackrel{\eqref{e:bar.tau}}{=}
	\bx[t]^\st\dpi
	-N^{1/2}
	\Big(1-\|\pi\|^2\Big)^{1/2}\bc[t-1]^\st\bar{\tau}\,.
	\eeq
As a result the function $\tcL_J$ from \eqref{l:adm.cgf} can be reparametrized as
	\[
	\tcL_J\bigg(\bar{\tau} + \f{\theta}{(1-\|\pi\|^2)^{1/2}}\bigg) 
	\stackrel{\eqref{e:tbX.rewrite}}{=}
	\f1N\bigg( \oneM,
		L_{\|\pi\|^2}\Big(
			\bx[t]^\st\dpi
			+ N^{1/2}
			\bc[t-1]^\st \theta
			\Big)\bigg)\,.\]
It follows by substituting the above calculations into \eqref{e:tcA}  that
	\[	\tcA_J\bigg(\bar{\tau} + \f{\theta}{(1-\|\pi\|^2)^{1/2}}\bigg) 
	= \f{\|\dvpi-\theta\|^2}{2(1-\|\pi\|^2)}
	+\f1N\sum_{a\le M}L_{\|\pi\|^2}
		\bigg(
		\bx[t]^\st\dpi
		+ N^{1/2}\bc[t-1]^\st\theta
		\bigg)\,.
	\]
The claim follows by taking
$\wp_{t,1}\equiv (\log\wp_{t,0} +t\log(2\pi))/2$.\end{proof}

The above completes the proof of Proposition~\ref{p:first.mmt}, modulo Proposition~\ref{p:density.bound} which is deferred to Section~\ref{s:local.clt}.

\subsection{First moment for partition function}
\label{ss:entropy}

We now collect some of the preceding results to complete the proof of the main result of this section:

\begin{proof}[\hypertarget{proof:t.ubd.slice}{Proof of Theorem~\ref{t:ubd.slice}}]
For any $J\in\set{-1,+1}^N$ we can calculate (abbreviating $\vpi\equiv\vpi(J)$)
	\[
	\f{(\bH^{(t)},J)}{N}
	= \f{(\acute{e}_{t-1})^\st\bH[t-1]J}{N}
	\stackrel{\eqref{e:gs.H.y}}{=}
	\f{\psi^{1/2}(\acute{e}_{t-1})^\st\bGam\by[t-1]J}{N}
	\stackrel{\eqref{e:def.vpi}}{=}
	(\psi^{1/2}\bGam^\st\acute{e}_{t-1},\vpi)
	\stackrel{\eqref{e:def.vpi.star}}{=}
	\f{(\starvpi,
		\vpi)}{1-q}\,.
	\]
It follows by combining with 
Definition~\ref{d:cube.P.Q} that
	\begin{align}\nonumber
	\f{\E(\bZ_\circ(\Gprime)
	 \,|\,\Fprime(t))}{
		\exp\{(\oneN,\log(2\Ch(\bH^{(t)})))\}}
	&=\sum_{J\in\mbH_\circ}
	\bQ(J)
		\bigg(
		\f{\bP(J)}{\bQ(J) 
		\exp\{(\oneN,\log\Ch(\bH^{(t)}))\}}
		\bigg)
		 \E\Big(\SAT_J(\Gprime)\,\Big|\,\Fprime(t)\Big)\\
	&=\sum_{J\in\mbH_\circ}
	\bQ(J)
	\exp\bigg\{-\f{N(\starvpi,\vpi)}{1-q}\bigg\}
	\E\Big(\SAT_J(\Gprime)\,\Big|\,\Fprime(t)\Big)\,.
	\label{e:first.mmt.rewrite.Q}
	\end{align}
Combining Proposition~\ref{p:first.mmt} with 
Lemma~\ref{l:change.of.basis.APPROX} gives, with high probability,
	\[\E\Big(\SAT_J(\Gprime)\,\Big|\,\Fprime(t)\Big)
	\le\f{\|\vpi-\theta\|^2}{2(1-\|\pi\|^2)}
	+\f1N\sum_{a\le M}L_{\|\pi\|^2}
		\bigg(
		\bx[t]^\st\pi
		+ N^{1/2}\bc[t-1]^\st\theta
		\bigg)
	+\ERR_{t,3}\,,
	\]
uniformly over 
$\|\pi(J)\|\le4/5$ and
$\|\theta\|\le\theta_{\max}$.
The claim follows by setting $\theta=\bareps(\vpi-\starvpi)$.
\end{proof}

Recall from \eqref{e:Z.outside} that $\bZ=\bZ_\circ+\bZ_\bullet$ where $\bZ_\circ$ is bounded by Theorem~\ref{t:ubd.slice}. In the remainder of this section we show that the other quantity 
$\bZ_\bullet$ can be bounded by \textit{a priori} estimates. For this purpose we prove a rough estimate on $\pi(J)$ (Lemma~\ref{l:AH.pi}), followed by a more precise estimate on $\vpi(J)$ (Lemma~\ref{l:AH.vpi}). In fact Lemma~\ref{l:AH.vpi} is more precise than what is needed to analyze $\bZ_\bullet$, but 
it will be needed later (in Section~\ref{s:analysis.first.mmt}) in the analysis of $\bZ_\circ$. We first state and prove the estimate for $\pi(J)$:

\begin{lem}\label{l:AH.pi}
Recall $\pi(J)$ from Definition~\ref{d:pi.vpi}
and $\stardpi$ from \eqref{e:def.pi.star}.
For $\bQ$ as in Definition~\ref{d:cube.P.Q}, we have
	\[
	\bQ\bigg(
	\bigg\{ J\in\set{-1,+1}^N
	: \Big\|\pi(J)-\stardpi\Big\|\ge d
	\bigg\} 
	\bigg)
	\le
	\bigg( \f{66t}{q}\bigg)^{t/2}
	\exp\bigg\{ -\f{Nd^2(1-3q^{1/2})}{8} \bigg\}
	\]
for all $|d|\ge 1/N^{1/2}$.
(The bound is vacuous unless $Nd^2$ is large compared to $t\log t$.)

\begin{proof}
Under the measure $\bQ$, the random vector $J-\bmag^{(t)}$ has independent entries of mean zero.
We note also that 
	\beq\label{e:AH.max.i}
	(\mathfrak{m}_i)^2
	\equiv
	\max\bigg\{\Big|J_i-(\bmag^{(t)})_i\Big|^2
	: J_i\in\set{-1,+1}
	\bigg\}
	\le \Big(1+ |(\bmag^{(t)})_i|\Big)^2
	\le 1 + 3|(\bmag^{(t)})_i|
	\le 4\,.\eeq
Thus for any $a\in\R^t$ we can bound
	\[
	V_{\max}(a)
	\equiv \sum_{i\le N}
		\bigg(
		\sum_{s\le t} a_s (\br^{(s)})_i
		\bigg)^2 
		(\mathfrak{m}_i)^2
	\le 4\bigg\|\sum_{s\le t}a_s\br^{(s)}\bigg\|^2
	=4\|a\|^2\,.
	\]
It follows by the Azuma--Hoeffding bound that
	\beq\label{e:AH.pi}
	\bQ\bigg(
	\f1{N^{1/2}}\bigg(
	\sum_{s\le t} a_s \br^{(s)}, J-\bmag^{(t)}
	\bigg) \ge x\bigg)
	\le \exp\bigg\{ - \f{Nx^2}{2V_{\max}(a)}\bigg\}
	\le\exp\bigg\{-\f{Nx^2}{8\|a\|^2}\bigg\}\,.
	\eeq
On the other hand, it follows from Definition~\ref{d:pi.vpi} and
\eqref{e:def.pi.star} that
	\beq\label{e:AH.pi.eval}
	\f1{N^{1/2}}\bigg(
	\sum_{s\le t} a_s \br^{(s)}, J-\bmag^{(t)}
	\bigg)
	=\Big(a, \pi(J)-\stardpi\Big)\,.\eeq
Given $d>0$ and $\epsilon\in(0,1/4]$, note there exists a $(d\epsilon)$-net of $[-4d,4d]^t$  of cardinality at most
	\beq\label{e:d.eps.net}
	\bigg\lceil
		\f{8 t^{1/2}}{\epsilon}\bigg\rceil^t
	\le \bigg( \f{8t^{1/2}}{\epsilon}+1\bigg)^t\,.
	\eeq
If $J$ is any element of $\set{-1,+1}^N$ with
 $d\le \|\pi(J)-\stardpi\|\le 2d$, and $\pi_\textup{net}$ is an element of the $(d\epsilon)$-net at minimal distance from 
  $\pi(J)$, then
 $\|\pi_\textup{net}-\stardpi\|\ge d(1-\epsilon)$, and
	\[
	\Big(\pi_\textup{net}-\stardpi,\pi(J)-\stardpi\Big)
	\ge\Big\|\pi(J)-\stardpi\Big\|^2
		-\bigg|\Big(\pi_\textup{net}-\pi(J),\pi(J)-\stardpi\Big)\bigg|
	\ge d^2(1-2\epsilon)\,.
	\]
Thus, by taking $a=\pi_\textup{net}-\stardpi$
and $\epsilon=q^{1/2}$ in \eqref{e:AH.pi} and \eqref{e:AH.pi.eval}, we obtain
	\begin{align*}
	\bQ \bigg(
		d\le \Big\|\pi(J)-\stardpi\Big\|\le 2d
	\bigg)
	&\le
	\bigg(\f{65t}{q}\bigg)^{t/2}
	\exp\bigg\{
	-\f{Nd^4(1-2q^{1/2})^2}
		{8d^2(1-q^{1/2})^2}
	\bigg\}\\
	&\le\bigg(\f{65t}{q}\bigg)^{t/2}
	\exp\bigg\{ -
		\f{Nd^2(1-3q^{1/2})}{8} \bigg\}\,.
	\end{align*}
Since $4^k\ge 3k$ for all $k\ge0$, 
as long as $Nd^2\ge1$ we can bound
	\[\bQ \bigg(
	\Big\|\pi(J)-\stardpi\Big\|\ge d
	\bigg)
	\le 
	\sum_{k\ge0}\f{(65t/q)^{t/2}}
	{\exp\{ N (2^kd)^2(1-3q^{1/2})/8 \}} 
	\le \f{2(65t/q)^{t/2}}{\exp\{ Nd^2(1-3q^{1/2})/8 \}}
	\,.\]
This proves the claim.
\end{proof}
\end{lem}

The result for $\vpi(J)$ is very similar, although slightly more involved since we require a more precise estimate:

\begin{lem}\label{l:AH.vpi}
Recall $\vpi(J)$ from Definition~\ref{d:pi.vpi},
and $\stardvpi$ from \eqref{e:def.vpi.star}.
For $\bQ$ as in Definition~\ref{d:cube.P.Q}, we have
	\[
	\bQ\bigg(
	\bigg\{ J\in\set{-1,+1}^N
	: \Big\|\vpi(J)-\stardvpi\Big\|\ge d
	\bigg\} 
	\bigg)
	\le \bigg(\f{66t}{q}\bigg)^{t/2} 
	\exp\bigg\{ -\f{Nd^2(1-8q^{1/2})}{2}\bigg\}
	\]
for all $|d| \ge1/N^{1/2}$.
(The bound is vacuous unless $Nd^2$ is large compared to $t\log t$.)

\begin{proof}
For $b\in\R^{t-1}$, denote
	\[
	W_{\max}(b)
	\equiv \f1N
	\sum_{i\le N}
		\bigg(
		\sum_{\ell\le t-1}
			b_\ell (\by^{(\ell)})_i
		\bigg)^2 
		(\mathfrak{m}_i)^2
	\le W_0(b) + 3 W_1(b)\,.\]
Recall the bound \eqref{e:AH.max.i} from the proof of Lemma~\ref{l:AH.pi}; it implies
$W_{\max}\le W_0+3W_1$ where
	\begin{align*}
	W_0(b)&\equiv
	\f1N\bigg\| \sum_{\ell\le t-1}
			b_\ell \by^{(\ell)}\bigg\|^2\,,\\
	W_1(b) &\equiv
	\f1N
	\sum_{i\le N}
	\bigg(
		\sum_{\ell\le t-1}
			b_\ell (\by^{(\ell)})_i
		\bigg)^2 
		|(\bmag^{(t)})_i|\,.
	\end{align*}
It follows from Lemma~\ref{l:AMP.PL} that $W_0(b)\to\|b\|^2$ in probability as $N\to\infty$. Lemma~\ref{l:AMP.PL} also implies
	\[\f{W_1(b)}{\|b\|^2}
	\stackrel{N\to\infty}{\longrightarrow}
	\E\bigg[ 
	\bigg(\rho Z + (1-\rho^2)^{1/2} Z' \bigg)^2
		\Big| \Th(\psi^{1/2}Z)\Big|
	\bigg]
	\equiv w_1(\rho)\,.
	\]
in probability, where $\rho\in[-1,1]$ is a value that can depend on $b$. However we can crudely bound
	\[
	w_1(\rho)
	\le
	\bigg(\E(Z^4) \E[ \Th(\psi^{1/2}Z)^2]\bigg)^{1/2}
	\stackrel{\eqref{e:fp}}{=} 
	(3q)^{1/2}\,.
	\]
It follows by the Azuma--Hoeffding inequality that
	\beq\label{e:AH.vpi}
	\bQ\bigg(
	\f1{N}\bigg(
	\sum_{\ell\le t-1} b_\ell \by^{(\ell)}, J-\bmag^{(t)}
	\bigg) \ge x\bigg)
	\le 
	\exp\bigg\{ - \f{Nx^2}{2W_{\max}(b)}\bigg\}
	\le \exp\bigg\{-\f{Nx^2}{2\|b\|^2
		(1 + 6q^{1/2})}\bigg\}\,.
	\eeq
On the other hand, it follows from Definition~\ref{d:pi.vpi} and \eqref{e:def.vpi.star} that
	\beq\label{e:AH.vpi.eval}
	\f1N\bigg(
	\sum_{\ell\le t-1} b_\ell \by^{(\ell)}, J-\bmag^{(t)}
	\bigg)
	= \Big(b,\vpi(J)-\stardvpi\Big)\,.\eeq
Given $d>0$ and  $\epsilon\in(0,1/4]$, note there exists a $(d\epsilon)$-net of $[-4d,4d]^{t-1}$ 
with cardinality upper bounded by \eqref{e:d.eps.net}.
If $J$ is any element of $\set{-1,+1}^N$ with
 $d\le \|\vpi(J)-\stardvpi\|\le 2d$, and $\vpi_\textup{net}$ is an element of the $(d\epsilon)$-net at minimal distance from 
  $\vpi(J)$, then
 $\|\vpi_\textup{net}-\stardvpi\| \ge d(1-\epsilon)$, and
  	\[
	\Big(\vpi_\textup{net}-\stardvpi,\vpi(J)-\stardvpi\Big)
	\ge \Big\|\vpi(J)-\stardvpi\Big\|^2
	-\bigg|
	 (\vpi_\textup{net}-\vpi(J),\vpi(J)-\stardvpi)
	 \bigg|
	\ge d^2(1-2\epsilon)\,.
	\]
Thus, by taking 
$\epsilon =q^{1/2}$ and
$b=\vpi_\textup{net}-\stardvpi$ in \eqref{e:AH.vpi} and \eqref{e:AH.vpi.eval}, we obtain
	\begin{align*}
	\bQ
	\bigg(
	d\le\Big\|\vpi(J)-\stardvpi\Big\|\le 2d
	\bigg)
	&\le \bigg(\f{65t}{q}\bigg)^{t/2}
	\exp\bigg\{
	-\f{N d^4(1-2q^{1/2})^2}{2 d^2(1-q^{1/2})^2(1+6q^{1/2}) }
	\bigg\}
	\\
	&\le \bigg(\f{65t}{q}\bigg)^{t/2}
	\exp\bigg\{ -\f{Nd^2(1-8q^{1/2})}{2}\bigg\}\,.
	\end{align*}
Since $4^k\ge 3k$ for all $k\ge0$, as long as $Nd^2\ge1$ we can bound
	\[\bQ
	\bigg(
	\Big\|\vpi(J)-\stardvpi\Big\|\ge d
	\bigg)
	\le 
	\sum_{k\ge0}
	\f{(65t/q)^{t/2}}
		{\exp\{N(2^kd)^2(1-8q^{1/2})/2\}} 
	\le
	\f{2(65t/q)^{t/2}}{\exp\{ Nd^2(1-8q^{1/2})/2\}}\,.\]
The claim follows.
\end{proof}
\end{lem}

\section{Technical estimates}
\label{s:technical}

We now collect some technical results which will be used later in the proof. This section is organized as follows:
\begin{itemize}
\item In \S\ref{ss:prelim.bounds} we prove some basic consequences of 
 Assumptions~\ref{a:bdd} and \ref{a:Lip}. 
\item In \S\ref{ss:estimate.rs} we give the 
\hyperlink{proof:p.fp}{proof of Proposition~\ref{p:fp}}, which characterizes the replica symmetric fixed-point solution. As a consequence of this analysis we obtain a rough estimate (Corollary~\ref{c:rs.estimate}) of the replica symmetric formula \eqref{e:rs}, which will be used in later sections. We also \hyperlink{proof:p.RS.U.eps}{prove Proposition~\ref{p:RS.U.eps}}, showing that the replica symmetric formula for $U_\ETA$ converges to the one for $U$ as $\ETA\downarrow0$.
\item In \S\ref{ss:AT} we prove Lemma~\ref{l:AT.conclusion}, which gives the Almeida--Thouless (AT) condition in our setting.

\item In \S\ref{ss:logc} 
we give the 
\hyperlink{proof:p.logc}{proof of Proposition~\ref{p:logc}},
showing that Assumption~\ref{a:Lip} holds if $u\equiv\log U$ is either bounded or concave.
We also give the
\hyperlink{proof:p.AT}{proof of Proposition~\ref{p:AT}} (convergence of the state evolution recursions),
which amounts to checking that AT condition 
derived in Lemma~\ref{l:AT.conclusion} holds for $0<\alpha\le\alpha(U)$. We conclude the section with some further consequences 
(Lemmas \ref{l:F.Lipschitz} and \ref{l:TAP.supnorm}) of Assumption~\ref{a:Lip}.
\end{itemize}
The following notation will be used throughout the paper:
\begin{dfn}\label{d:mu.x.c}
For $c>0$ and $x\in\R$, let $\mu_{x,c}$ denote the probability measure on the real line whose density (with respect to the Lebesgue measure) is given by
	\[
	\f{d\mu_{x,c}}{dz} 
	=\chi_{x,c}(z)\equiv \f{U(x+cz)\varphi(z)}
		{\E_\xi[U(x+c\xi)]}\,.\]
We use $\E_{x,c}$, $\Var_{x,c}$, and $\Cov_{x,c}$ to denote 
expectation,
variance, and covariance under $\mu_{x,c}$. 
\end{dfn}

\subsection{Preliminary bounds}
\label{ss:prelim.bounds}
In this subsection we prove some basic consequences of  Assumptions~\ref{a:bdd} and \ref{a:Lip}. As before, $\xi$ denotes a standard gaussian random variable, and $\E_\xi$ denotes expectation over $\xi$.

\begin{lem}\label{l:uniform.conv.ell.one}
Suppose $U$ satisfies Assumption~\ref{a:bdd}, and let 
$q_{x,c}(z)\equiv U(x+cz)\varphi(z)$ as above. Then, given any 
$\epsilon>0$ and any $L<\infty$, it is possible to choose $\eta'$ small enough such that we have the bound
	\[
	\int \Big| U(x+cz)- U(x'+c'z)\Big|
	\varphi(z)\,dz
	\le\epsilon
	\] 
as long as $c,c'\in[1/3,3]$, $x,x'\in[-L,L]$, and
$\max\set{|x-x'|,|c-c'|}\le \eta'$.

\begin{proof} Given $\epsilon>0$, we can clearly choose $L(\epsilon)$ large enough (depending only on $\epsilon$) such that $L(\epsilon)\ge L$, and
	\beq\label{e:large.z.bound}
	\int_{|z|\ge L(\epsilon)}
		\Big| U(x+cz)- U(x'+c'z)\Big|
	\varphi(z)\,dz
	\le \int_{|z|\ge L(\epsilon)}
		\varphi(z)\,dz
	\le \f{\epsilon}{4}\,.\eeq
If $|z|\le L(\epsilon)$, then the assumptions imply $|x+cz|\le 4L(\epsilon)$ and $|x'+c'z|\le 4L(\epsilon)$, so
	\[\int_{|z|\le L(\epsilon)}
		\Big| U(x+cz)- U(x'+c'z)\Big|
	\varphi(z)\,dz
	\le \int
		\Big| \bar{u}(x+cz)-\bar{u}(x'+c'z)\Big|
	\varphi(z)\,dz
	\]
where $\bar{u}(x)\equiv U(x) \Ind{|x| \le 4 L(\epsilon)}$. Then, since $\bar{u}\in L^1$,
 it is well known that we can choose a function $\tilde{u}$ which is compactly supported and smooth, such that
$\|\bar{u}-\tilde{u}\|_1\le\epsilon/4$
(see e.g.\ \cite[Lem.~2.19]{MR1817225}). Therefore
	\[
	\int\Big| \bar{u}(x+cz)-\tilde{u}(x+cz)\Big|
	\varphi(z)\,dz
	\le  \varphi(0) \int
		\Big| \bar{u}(x+cz)-\tilde{u}(x+cz)\Big|
		\,dz
	= \f{\varphi(0)\|\bar{u}-\tilde{u}\|_1}{c}
	\le \f{\epsilon}{4}\,,
	\]
where this estimate holds for all $x\in\R$ and all $c\ge1/2$. We also have
	\[
	\int
		\Big| \tilde{u}(x+cz)-\tilde{u}(x'+c'z)\Big|
	\varphi(z)\,dz
	\le \|\tilde{u}'\|_\infty
	\int\Big( |x-x'| + |c-c'||z|\Big) \varphi(z)\,dz
	\le
	2 \|\tilde{u}'\|_\infty \eta'\,,
	\]
which can be made at most $\epsilon/4$ by taking
$\eta'=\epsilon/(8\|\tilde{u}'\|_\infty)$. Combining the above estimates gives
	\[
	\int_{|z|\le L(\epsilon)}
		\Big| U(x+cz)- U(x'+c'z)\Big| \le \f{3\epsilon}{4}\,,
	\]
and combining with the estimate
\eqref{e:large.z.bound} for $|z|\ge L(\epsilon)$ gives the conclusion.
\end{proof}
\end{lem}

\begin{lem}\label{l:poly}
Suppose $U$ satisfies Assumption~\ref{a:bdd} . There exists a finite constant $C_1(U)$, depending on $U$ only, such that
	\[\E_{x,c}(|Z|^p)
	= \f{\E_\xi( |\xi|^p U(x+c \xi))}
		{\E_\xi U(x+c\xi)}
	\le C_1(U)
		+ \bigg(\f{1.82 \cdot |x|}{c}\bigg)^p 
	\]
for all $0\le p\le 200$, $1/2\le c\le2$, and $x\in\R$. (We can assume, without loss, $C_1(U)\ge 10$.)

\begin{proof}
It follows from Assumption~\ref{a:bdd} that
$\E_\xi U(c\xi)>0$ for any $c>0$. 
Lemma~\ref{l:uniform.conv.ell.one} gives that $\E_\xi U(c\xi)$ is a continuous function of $1/2\le c\le 2$, so by compactness considerations we must have
	\beq\label{e:bar.c.one.U}
	\bar{c}_1(U)
	\equiv\max\bigg\{2,
	\sup\bigg\{ \f1{\E_\xi U(c\xi)}
		: \f12 \le c\le 2
	\bigg\}
	\bigg\}<\infty
	\eeq
(where we chose $\bar{c}_1(U)\ge2$ for convenience). Next, for any $M>0$, it holds for all $1/2\le c\le 2$ that
	\[
	\E_\xi\Big( U(c\xi) ; |c\xi| \ge M\Big)
	\le \P\bigg(|\xi| \ge \f{M}{c}\bigg)
	\le \f{\varphi(M/c)}{M/c}\,.
	\]
If we take $K \ge K_0(U)=(8\log c_1(U))^{1/2}\ge2$, then for all $1/2\le c\le 2$ we have
	\[
	\E_\xi\Big( U(c\xi) ; |c\xi| \le K\Big)
	\ge \E_\xi U(c\xi) - \f{\varphi(K/2)}{K/2}
	\ge \f1{2\bar{c}_1(U)}\,.
	\]
In what follows let $K(x)\equiv \max\set{K_0(U),|x|}$. Then we can lower bound
	\begin{align}\nonumber
	\E_\xi U(x+c\xi)
	&=\int U(cz) \varphi\bigg(z-\f{x}{c}\bigg)
		\,dz
	= \int U(cz)
	\exp\bigg\{-\f{x^2}{2c^2}
		+\f{xz}{c} \bigg\} \varphi(z) \,dz\\
	&\ge
	\f{\E(  U(c\xi) ; |c\xi| \le K(x) )}
		{\exp\{(3/2) K(x)^2 /c^2\}}
	\ge \f{1/(2\bar{c}_1(U))}
		{\exp\{(3/2) K(x)^2 /c^2\}}\,.
	\label{e:EU.lbd}
	\end{align}
Next we note that for any $M\ge0$
and $\eta'=1/10$ we have
	\begin{align}\nonumber
	\E_\xi( |\xi|^p ; |\xi|\ge M)
	&= \int_{|z|\ge M} \f{|z|^p}{(2\pi)^{1/2}}
		\exp\bigg\{-\f{z^2}{2}\bigg\}
	\,dz\\ \nonumber
	&=\int_{|z|\ge M/(1+\eta')^{1/2}}
	\f{(1+\eta')^{(p+1)/2}|z|^p}{(2\pi)^{1/2}}
	\exp\bigg\{ -\f{(1+\eta')z^2}{2}
	\bigg\}\,dz\\
	&\le 
	\sup\bigg\{ \f{1.05^{p+1}|z|^p}{\exp(z^2/20)}
	: z\in\R 
	\bigg\}
	\P\bigg(|\xi|\ge \f{M}{(1+\eta')^{1/2}}\bigg)
	\le 
	\CSTZRO\f{\varphi(M/1.05)}{M/1.05}\,,
	\label{e:E.xip.ubd}
	\end{align}
where $\CSTZRO\ge5$ is an absolute constant since we restricted $0\le p\le 200$. Combining \eqref{e:EU.lbd} with \eqref{e:E.xip.ubd} gives
	\begin{align*}
	\f{\E_\xi (|\xi|^p U(x+c\xi))}
		{\E_\xi U(x+c\xi)}
	&\le \bigg(\f{1.82 \cdot K(x)}{c}\bigg)^p
	+ \f{\E_\xi( |\xi|^p ; |\xi| \ge 1.82 \cdot K(x)/c)}
		{\E_\xi U(x+c\xi)} \\
	&\le
	 \bigg(\f{1.82\cdot K(x)}{c}\bigg)^p
	 + \CSTZRO\cdot
	 \f{ \varphi( 1.82 \cdot K(x)/(1.05\cdot c) )}
	 	{1.82\cdot K(x)/(1.05\cdot c)}
	\cdot \f{\exp\{(3/2) K(x)^2/c^2\}}
	{ 1/(2\bar{c}_1(U)) }\,.
	\end{align*}
(The first inequality above also uses that $U\le 1$, from Assumption~\ref{a:bdd}.) Recalling again the restrictions $1/2 \le c\le2$
and $0\le p\le 200$, we can simplify the above to obtain
	\begin{align}\nonumber
	\f{\E_\xi (|\xi|^p U(x+c\xi))}
		{\E_\xi U(x+c\xi)}
	&\le \bigg(\f{1.82\cdot K(x)}{c}\bigg)^p
		+ 2\CSTZRO \cdot \bar{c}_1(U) 
	\cdot \varphi\bigg( \f{1.82\cdot K(x)}{1.05\cdot c}\bigg)
	\exp\bigg\{\f{3K(x)^2}{2c^2}\bigg\} \\ \nonumber
	&\le
	\bigg(\f{1.82\cdot|x|}{c}\bigg)^p 
	+\bigg\{ \bigg(\f{1.82 \cdot K_0(U)}{1/2}\bigg)^{200}
		+ \f{2 \cdot \CSTZRO \cdot \bar{c}_1(U)}{(2\pi)^{1/2}}
		 \bigg\}\\ \nonumber
	&\le\bigg(\f{1.82\cdot|x|}{c}\bigg)^p 
	+\bigg\{ \Big(4 K_0(U)\Big)^{200}
		+ \CSTZRO \cdot \bar{c}_1(U)
		 \bigg\}\\ \label{e:def.C.one}
	&\equiv \bigg(\f{1.82\cdot|x|}{c}\bigg)^p 
	 + C_1(U)\,,
	\end{align}
where the last equality defines $C_1(U)$. The above choices guarantee $C_1(U)\ge\CSTZRO\cdot\bar{c}_1(U)\ge 10$.
\end{proof}
\end{lem}

\begin{rmk*}
The bound from Lemma~\ref{l:poly} is reasonably tight.
To see this, consider the function
	\[U(x)=\I\bigg\{|x-a| \le \f{a}{2}\bigg\}\]
for $a>0$. If $U(x+\xi)=1$, then $x+\xi \ge a/2$, so $\xi\ge a/2-x$. In the case that $x\le0$, it implies $|\xi| \ge a/2 + |x|$. It follows that for any $x\le0$ we have
	\[\E_{x,c}(|Z|^p)
	= \f{\E_\xi(|\xi|^p U(x+\xi))}{\E_\xi U(x+\xi)}
	\ge \bigg(\f{a}{2}+|x|\bigg)^p
	\ge \bigg(\f{a}{2}\bigg)^p + |x|^p\,,
	\]
where $a>0$ can be chosen to be arbitrarily large.
\end{rmk*}

Next we combine Assumption~\ref{a:Lip} (which bounds $\Var_{x,c}(Z)$) with the calculations of Lemma~\ref{l:poly}
to obtain bounds on $\Var_{x,c}(Z^2)$ and $\Cov_{x,c}(Z,Z^2)$:

\begin{lem}\label{l:Lip.implications}
Suppose $U$ satisfies Assumptions~\ref{a:bdd} and \ref{a:Lip}, 
and let $C_1(U)$ be as in Lemma~\ref{l:poly}. Then we have
	\begin{align}\label{e:Lip.var.sq}
	\Var_{x,c}(Z^2)
	&\le K_2(U)\cdot\bigg\{ \bigg(\f{1.82\cdot |x|}{c}\bigg)^2 
	+ C_1(U) \bigg\}\,,\\
	\Cov_{x,c}(Z,Z^2)
	&\le
	\f{K_2(U) }{2^{1/2}}
	\cdot \bigg(
	\f{1.82\cdot |x|}{c}
	+  C_1(U)^{1/2}
	\bigg)\,,
	\label{e:Lip.cov}
	\end{align}
for all $1/2\le c\le 2$ and all $x\in\R$.

\begin{proof}
Let $K(x)$ be as in the proof of Lemma~\ref{l:poly}. From the definition of $K_2(U)$ (see Assumption~\ref{a:Lip}),
	\begin{align*}
	(\textup{I})
	&\equiv\f{\E_{\xi,\xi'} [(\xi - \xi')^2(\xi + \xi')^2 U(x+c\xi)U(x+c\xi')
		; |\xi+\xi'| \le  2^{1/2} \cdot 1.82 \cdot K(x)/c]}
		{\E_{\xi,\xi'}[ U(x+c\xi)U(x+c\xi')]} \\
	&\le
	2 K_2(U)\cdot \bigg( \f{1.82 \cdot K(x)}{c} \bigg)^2\,.
	\end{align*}
If $\xi$ and $\xi'$ are independent standard gaussian random variables, then $\xi-\xi'$ and $\xi+\xi'$ are independent gaussian random variables with mean zero and variance $2$. It follows that
	\beq\label{e:E.xip.fourth.ubd}
	\E_{\xi,\xi'}\Big[
		(\xi - \xi')^2 (\xi + \xi')^2 ;|\xi + \xi'| \geq \sqrt{2}M
		\Big]
	= 4\cdot\E_\xi\Big[|\xi|^2 ; |\xi| \ge M\Big]
	\stackrel{\eqref{e:E.xip.ubd}}{\le}
	4\cdot \CSTZRO\f{\varphi(M/1.05)}{M/1.05}\,.
	\eeq
Combining \eqref{e:E.xip.fourth.ubd} with our earlier bound \eqref{e:EU.lbd} gives 
	\begin{align*}
	(\textup{II})
	&\equiv\f{\E_{\xi,\xi'} [(\xi - \xi')^2(\xi + \xi')^2 U(x+c\xi)U(x+c\xi')
		; |\xi+\xi'| \ge  2^{1/2} \cdot 1.82 \cdot K(x)/c]}
		{\E_{\xi,\xi'}[ U(x+c\xi)U(x+c\xi')]}\\
	&\le 4\cdot\CSTZRO
	 \f{ \varphi( 1.82\cdot K(x)/(1.05\cdot c) )}
	 	{1.82\cdot K(x)/(1.05\cdot c)}
	\cdot \f{\exp\{(3/2) K(x)^2/c^2\}}
	{ 1/(2\bar{c}_1(U)) }
	\end{align*}
(The first inequality above also uses that $U\le 1$, from Assumption~\ref{a:bdd}.)
By combining the above bounds for the quantities (\textup{I}) and (\textup{II}), and recalling again that $1/2\le c\le 2$, we obtain
	\begin{align*}
	\Var_{x,c}(Z^2)&=
		\f{\E_{\xi,\xi'} [(\xi - \xi')^2(\xi + \xi')^2 
		U(x+c\xi) U(x+c\xi')]}
		{2\cdot \E_{\xi,\xi'} U(x+c\xi)U(x+c\xi')}\\
	&\le K_2(U)\cdot \bigg( \f{1.82 \cdot K(x)}{c} \bigg)^2
		+ 
		\f{4\cdot \CSTZRO \cdot \bar{c}_1(U)}
			{1.82/1.05}
    	\cdot \varphi\bigg( \f{1.82\cdot K(x)}{1.05\cdot c}\bigg)
    	\exp\bigg\{\f{3K(x)^2}{2c^2}\bigg\} \\
	&\le K_2(U) \cdot\bigg\{\bigg(\f{1.82\cdot |x|}{c}\bigg)^2 
	+ \bigg( \f{1.82 \cdot K_0(U)}{1/2}
	\bigg)^2\bigg\}
	+ \f{4 \cdot \CSTZRO \cdot \bar{c}_1(U) }
		{(1.82/1.05) \cdot (2\pi)^{1/2}}\\
	&\le K_2(U)\cdot\bigg\{ \bigg(\f{1.82\cdot |x|}{c}\bigg)^2 
	+14 \cdot K_0(U)^2
		+  \CSTZRO \cdot \bar{c}_1(U) 
		\bigg\}\\
	&\le K_2(U)\cdot \bigg\{
	\bigg(\f{1.82\cdot |x|}{c}\bigg)^2 
	 + C_1(U)
	 \bigg\}\,,
	\end{align*}
where the second-to-last inequality uses that we took $K_2(U)\ge1$
(see Assumption~\ref{a:Lip}), and the last inequality uses the definition \eqref{e:def.C.one} of $C_1(U)$ from the proof of Lemma~\ref{l:poly}. This proves \eqref{e:Lip.var.sq}. Combining with Assumption~\ref{a:Lip} and the Cauchy--Schwarz inequality gives 
	\begin{align*}
	\Cov_{x,c}(Z,Z^2)
	&\le\bigg\{ \Var_{x,c}(Z)
		\Var_{x,c}(Z^2)\bigg\}^{1/2}
	\le \f{K_2(U)}{2^{1/2}}
	\bigg\{\bigg(\f{1.82\cdot |x|}{c}\bigg)^2 
	+ C_1(U)\bigg\}^{1/2}\\
	&\le
	\f{K_2(U)}{2^{1/2}}
	 \bigg(
	\f{1.82\cdot |x|}{c}
	+ C_1(U)^{1/2}
	\bigg)\,,
	\end{align*}
where the last inequality again uses that $K_2(U)\ge1$. This proves \eqref{e:Lip.cov}.
\end{proof}
\end{lem}

\begin{rmk}\label{r:counterex}
We include here an example of a function $U$ that satisfies Assumption~\ref{a:bdd} but does not satisfy the bound \eqref{e:Lip.var.sq}
(and hence, by Lemma~\ref{l:Lip.implications}, must violate Assumption~\ref{a:Lip}). For $k\ge1$ let $b_k\equiv \exp(-100 \cdot 4^k)$, and let
	\[
	A_k(x)
	\equiv 
	\bigg(\I\Big\{x\in[0,1]\Big\}
		 + \I\Big\{ x \in [2^k-1,2^k]\Big\}
	\bigg) \f{b_k}{\varphi(x)}
	\equiv \f{ b_k f_k(x) }{\varphi(x)}\,.
	\]
Then clearly $A_k$ is a nonnegative measurable function
supported on $[0,1]\cup[2^k-1,2^k]$, with
	\[
	\|A_k\|_\infty
	\le \f{b_k}{\varphi(2^k)}
	= b_k (2\pi)^{1/2} \exp\bigg( \f{4^k}{2}\bigg)
	\le \f{(2\pi)^{1/2}}{\exp( 99 \cdot 4^k)}\,.
	\]
Let $C$ be a large absolute constant, and define $x_k\equiv C 2^k$ and
	\[
	U(x)
	\equiv
	\sum_{k\ge1} A_k(x_k+x)\,.
	\]
From the above bound on $\|A_k\|_\infty$ it is clear that $U$ satisfies Assumption~\ref{a:bdd}. Next we note that
	\begin{align*}
	\f{\E_\xi [\xi^2 A_k(\xi)]}{b_k}
	&= \int z^2 f_k(z) \,dz 
	= \f13 \bigg( (2^k)^3-(2^k-1)^3+1 \bigg)
	= 2^{2k}\bigg(1+ \f{O(1)}{2^k}\bigg)\,,\\
	\f{\E_\xi [\xi^4 A_k(\xi)]}{b_k}
	&= \int z^4 f_k(z) \,dz 
	= \f15 \bigg( (2^k)^5-(2^k-1)^5+1 \bigg)
	= 2^{4k}\bigg(1+ \f{O(1)}{2^k}\bigg)\,.
	\end{align*}
For any $k\ge1$, we have $\E_\xi U(-x_k+\xi) \ge \E_\xi A_k(\xi) = 2b_k$.
For $\ell\ge k+1$ and $0\le p\le 4$, we have
	\begin{align*}
	&\f{\E_\xi[|\xi|^p A_\ell(-x_k+x_\ell+x)]}{b_k}
	\le \f{\|A_\ell\|_\infty \E_\xi(\xi^4)}{b_k} \le
	 \f{3 (2\pi)^{1/2} \exp(100\cdot 4^k)}{\exp(99 \cdot 4^\ell)} \\
	&\qquad
	\le \f{3 (2\pi)^{1/2}}{\exp(4^k[74 \cdot 4^{\ell-k}
		+(25 \cdot4^{\ell-k} -100) ])}
	\le \f{3 (2\pi)^{1/2}}{\exp(74 \cdot 4^\ell)}
	\le \f1{\exp(70 \cdot 4^\ell)}\,.
	\end{align*}
On the other hand, for $1\le\ell\le k-1$ and $0\le p\le 4$, we have
(again taking $C$ large enough)
	\begin{align*}
	&\f{\E_\xi[|\xi|^p A_\ell(-x_k+x_\ell+x)]}{b_k}
	\le \f{\|A_\ell\|_\infty
	\E_\xi[\xi^4 ; |\xi| \ge C 2^k/4]}{b_k}
	\le \f{(2\pi)^{1/2}\exp(100 \cdot 4^k)}{\exp(99 \cdot 4^\ell + C^2 4^k/33)}\\
	&\qquad\le \f{(2\pi)^{1/2}\exp(100 \cdot 4^k)}{\exp(99+C^2 4^k/33)}
	\le  \f{(2\pi)^{1/2}}{\exp(99+C  4^k)} \le \f1{\exp(70 \cdot 4^k)}\,.
	\end{align*}
(In the first inequality above, we used that the support of $A_\ell$ is contained in $[0,2^\ell]$.) Altogether we conclude
	\begin{align*}
	\f{\E_\xi[\xi^2 U(-x_k+\xi)]}{\E_\xi U(-x_k+\xi)}
	&=\bigg(1 + \f{O(1)}{\exp(4^k)}\bigg)
	\f{\E_\xi[\xi^2 A_k(\xi)]}{\E_\xi A_k(\xi)}
	=\bigg(1 + \f{O(1)}{\exp(4^k)}\bigg) \f{2^{2k}}{2}\,,\\
	\f{\E_\xi[\xi^4 U(-x_k+\xi)]}{\E_\xi U(-x_k+\xi)}
	&=\bigg(1 + \f{O(1)}{\exp(4^k)}\bigg)
	\f{\E_\xi[\xi^4 A_k(\xi)]}{\E_\xi A_k(\xi)}
	=\bigg(1 + \f{O(1)}{\exp(4^k)}\bigg) \f{2^{4k}}{2}\,.
	\end{align*}
Recalling the notation of Definition~\ref{d:mu.x.c}, we obtain 
	\[
	\Var_{-x_k,1}(Z^2)
	=\bigg(1 + \f{O(1)}{\exp(4^k)}\bigg) 
	\bigg\{ \f{2^{4k}}{2}-\bigg(\f{2^{2k}}{2}\bigg)^2\bigg\}
	= \Theta( 2^{4k})\,.
	\]
Thus shows that $U$ does not satisfy the bound \eqref{e:Lip.var.sq}, as claimed.\end{rmk}

\subsection{Estimates of the replica symmetric solution}
\label{ss:estimate.rs}
In this subsection we give the \hyperlink{proof:p.fp}{proof of Proposition~\ref{p:fp}}. As a consequence we obtain a rough estimate (Corollary~\ref{c:rs.estimate}) of the replica symmetric formula which will be used later in our analysis.

\begin{lem}\label{l:dq.dpsi}
Suppose $U$ satisfies Assumption~\ref{a:bdd}.
As in Proposition~\ref{p:fp}, let
$\bar{q}(\psi)\equiv\E[\Th(\psi^{1/2}Z)^2]$. Then
	\[
	\max\Big\{0,1-4\psi\Big\} \le
	\f{d\bar{q}}{d\psi} \le 1
	\]
for all $\psi\ge0$.

\begin{proof}
It is clear that $\bar{q}$ is increasing with respect to $\psi\ge0$: indeed,
	\[\f{d\bar{q}}{d\psi}
	= \E\bigg[ \Th(\psi^{1/2} Z)
	\Th'(\psi^{1/2} Z)  \f{Z}{\psi^{1/2}}\bigg]
	>0\,,
	\]
since $\Th'(x)>0$ for all $x\in\R$, and $x\Th(x)\ge 0$ for all $x\in\R$. Integrating by parts gives
	\begin{align*}
	\f{d\bar{q}}{d\psi}
	&= \E \bigg[\Big(\Th'(\psi^{1/2} Z)\Big)^2
		+ \Th(\psi^{1/2} Z) \Th''(\psi^{1/2} Z)
		\bigg] \\
	&= \E\bigg( 1-4\Th(\psi^{1/2} Z)^2 
	+ 3\Th(\psi^{1/2} Z)^4
		\bigg)\,.
	\end{align*}
Note that $x=\Th(\psi^{1/2} Z)^2 \in[0,1]$ almost surely, and $1-4x\le 1-4x+3x^2 \le 1$ for all $x\in[0,1]$, so
	\[
	1\ge \f{d\bar{q}}{d\psi}
	\ge \E\Big( 1-4\Th(\psi^{1/2} Z)^2 \Big)
	\ge 1-4\psi \cdot \E(Z^2) = 1-4\psi\,,
	\]
for all $\psi\ge0$.
\end{proof}
\end{lem}

\begin{lem}\label{l:dpsi.dq}
Suppose $U$ satisfies Assumption~\ref{a:bdd}.
As in Proposition~\ref{p:fp}, let
$\bar{r}(q)\equiv \E[F_q(q^{1/2}Z)^2]$. Then
	\[
	\sup\bigg\{
	\bigg| \f{d\bar{r}}{dq}\bigg|
	: 0\le q\le \f12\bigg\} 
	\le
	\CSTONE \cdot C_1(U)^6\,,
	\]
where $\CSTONE\ge1$ is an absolute constant while $C_1(U)$
is the constant from Lemma~\ref{l:poly}.

\begin{proof}
For convenience we shall rewrite \eqref{e:F} as
	\beq\label{e:F.rewrite}
	F_q(x)
	= \f{\E_\xi U'(x+(1-q)^{1/2}\xi)}{\E_\xi U(x+(1-q)^{1/2}\xi)}\,.
	\eeq
Note that the above makes sense for any $U$ satisfying Assumption~\ref{a:bdd}, without any smoothness assumption, since $U'$ can be interpreted as a distributional derivative (as in e.g.\ \cite[Ch.~6]{MR1817225}). Similarly one can make sense of the distributional derivative $U^{(k)}$ for any integer $k\ge1$.  We can then calculate
	\[\f{d\bar{r}}{dq}
	= \E\bigg[ 2F_q(q^{1/2}Z)
	 \f{d[F_q(q^{1/2}Z)]}{dq}\bigg]
	= \textup{(I)}-\textup{(II)}
	\]
where, abbreviating $U^{(k)} \equiv U^{(k)}(q^{1/2}Z+(1-q)^{1/2}\xi)$, we have
	\begin{align*}
	\textup{(I)}
	&= \E\bigg[ \f{Z}{q^{1/2}} 
		F_q(q^{1/2}Z)\cdot (F_q)'(q^{1/2}Z)\bigg]
	= \E\bigg[ F_q(q^{1/2}Z)\cdot (F_q)''(q^{1/2}Z) + \Big((F_q)'(q^{1/2}Z)\Big)^2\bigg]\,,\\
	\textup{(II)}
	&=
		\E\bigg[
		\f{F_q(q^{1/2}Z)}{(1-q)^{1/2}}
		 \bigg(\f{\E_\xi(\xi U'')}{\E_\xi U}
		- \f{(\E_\xi U') \E_\xi(\xi U')}{(\E_\xi U)^2} \bigg)\bigg]\,.
	\end{align*}
It follows by repeated applications of the inequality $2ab \le a^2+b^2$ that
	\[\bigg|\f{d\bar{r}}{dq}\bigg|
	\le C\sum_{0\le k,p\le3}	\E
	\bigg(\f{\E_\xi[|\xi|^k U(q^{1/2}Z+(1-q)^{1/2}\xi)]}
		{\E_\xi U(q^{1/2}Z+(1-q)^{1/2}\xi)}\bigg)^{2p}
	\]
for all $0\le q\le 1/2$, where $C$ is an absolute constant. It then follows from Lemma~\ref{l:poly} that
	\[\bigg|\f{d\bar{r}}{dq}\bigg|
	\le C\sum_{0\le k,p\le3}\E\bigg[
	\Big( (4q^{1/2}|Z|)^k + C_1(U)\Big)^{2p}\bigg]
	\le \CSTONE \cdot C_1(U)^6
	\]
for all $0\le q\le 1/2$, where $\CSTONE\ge1$ is an absolute constant.
\end{proof}
\end{lem}

\begin{proof}[\hypertarget{proof:p.fp}{Proof of Proposition~\ref{p:fp}}]
We seek a value $q\in[0,1/25]$ that satisfies the fixed-point equation \eqref{e:fp}, i.e., $q=\bar{q}(\alpha\bar{r}(q))$. This is the same as a root $q\in[0,1/25]$ of the function
	\beq\label{e:root.fn}
	\bar{g}(q) = \f{\bar{q}^{-1}(q)}{\alpha} - \bar{r}(q)\,.
	\eeq
Note that $\bar{q}(0)=0$, and it follows from Lemma~\ref{l:dq.dpsi} that 
$\bar{q}'(\psi)\in[4/5,1]$ for all $\psi\le1/20$, so
	\[
	\f45\psi \le \bar{q}(\psi) \le \psi
	\]
for all $\psi\le1/20$. Consequently, if $\bar{q}(\psi)\le1/25$ then we must have $\psi\le1/20$, that is to say,
	\[
	\sup\bigg\{
	(\bar{q})^{-1}(q)
	: q\le \f1{25}
	\bigg\} \le \f1{20}\,.
	\]
It follows from
Lemma~\ref{l:dq.dpsi} that 
$(\bar{q}^{-1})'(q)\in[1,5/4]$ for all $q\le1/25$. Combining with Lemma~\ref{l:dpsi.dq} gives
	\[
	\f{1}{\alpha} - \CSTONE \cdot C_1(U)^6 \le 
	\f{d\bar{g}}{dq} 
	\le \f{5}{4\alpha} + \CSTONE \cdot C_1(U)^6\,,
	\]
where $\CSTONE$ is the absolute constant 
from Lemma~\ref{l:dpsi.dq}. It follows that as long as $\alpha \le\alpha(U)$ as defined by \eqref{e:alpha.U}, 
then for all $0\le q\le 1/ 25$ we will have
	\[\f1{2\alpha} \le \f{d\bar{g}}{dq}
	\le \f{2}{\alpha}\,.\]
At $q=0$ we have
$\bar{g}(0)=-\bar{r}(0)$, and
 it follows by Assumption~\ref{a:bdd} 
combined with Lemma~\ref{l:poly} that
	\[	
	\Big(\E_\xi[\xi U(\xi)]\Big)^2
	\le \bar{r}(0)
	=\bigg( \f{\E_\xi U'(\xi)}{\E_\xi U(\xi)}\bigg)^2
	= \bigg(\f{\E_\xi[\xi U(\xi)]}
		{\E_\xi U(\xi)}\bigg)^2
	\le C_1(U)^2
	\,.\]
It follows that on the interval $0\le q\le 1/25$, the function $\bar{g}$ has a unique root $q$, which must satisfy
	\[
	\f{(\E_\xi[\xi U(\xi)])^2}{2}
	\le \f{q}{\alpha} \le 2 C_1(U)^2
	\]
It follows from the earlier bound on $\psi$ that 
	\[\f{(\E_\xi[\xi U(\xi)])^2}{2}
	\le \f{q}{\alpha}
	\le \f{\psi}{\alpha}
	\le \f{5q}{4\alpha} \le \f{5 C_1(U)^2}{2}
	\,,\]
so this concludes the proof.
\end{proof}

\begin{cor}\label{c:rs.estimate}
If the function $U$ satisfies Assumptions~\ref{a:bdd} and \ref{a:Lip},  then
for all $0<\alpha\le\alpha(U)$ we have
	\[
	\f{\RS(\alpha;U)}{\alpha}
	\ge \f{\annealed(\alpha;U)}{\alpha} 
		- 1.51\cdot C_1(U)^2
	\ge \f{\log 2}{\alpha} -
		1.53 \cdot C_1(U)^2\,,
	\]
where $C_1(U)$ is the constant from Lemma~\ref{l:poly},
and $\alpha(U)$ is given by \eqref{e:alpha.U}.

\begin{proof}
Let $(q,\psi)$ be the solution from Proposition~\ref{p:fp}, and recall from \eqref{e:rs} that
	\[ \RS(\alpha;U) - \log2
	= 
	- \f{\psi(1-q)}{2}
	+\E\bigg\{ \log\Ch(\psi^{1/2}Z)
	+ \alpha L_q(q^{1/2}Z)\bigg\}\,.
	\]
We hereafter abbreviate
	\[
	\bar{\ell}(q)
	\equiv
	\E L_q(q^{1/2}Z)
	= \E\bigg[
	\log \E_\xi U\Big(
	q^{1/2}Z + (1-q)^{1/2}\xi\Big)
	\bigg]\,.
	\]
Since $\Ch(x)\ge1$ for all $x\in\R$, we can lower bound
	\[ \RS(\alpha;U) - \log2
	\ge -\f{\psi}{2}  
	+ \alpha \E L_q(q^{1/2}Z)
	\ge
	 -\f{\psi}{2}  
	+ \alpha \bigg\{ \bar{\ell}(0)- 
	q \sup_{0\le q\le1/2}
	\bigg|\f{d\bar{\ell}}{dq} \bigg|
	\bigg\}
	\,.\]
Similarly as in the proof of Lemma~\ref{l:dpsi.dq}, we can bound
	\beq\label{e:d.ell.d.q}
	\bigg|\f{d\bar{\ell}}{dq} \bigg|
	\le C \sum_{0\le k,p\le2}
	\E\bigg[
	\bigg(\f{\E_\xi[ |\xi|^k U(q^{1/2}Z+(1-q)^{1/2}\xi)]}
	{\E_\xi U(q^{1/2}Z+(1-q)^{1/2}\xi)}
	\bigg)^p\bigg]
	\le \CSTONE \cdot C_1(U)^2
	\eeq
for all $0\le q\le 1/2$, where $\CSTONE$ is an absolute constant (and can be arranged to be the same as the $\CSTONE$ from Lemma~\ref{l:dpsi.dq}). By combining the above bounds we conclude
	\begin{align*}
	&\f{\RS(\alpha;U)
	- (\log2 + \alpha \bar{\ell}(0))}{\alpha}
	\ge -\f{\psi}{2\alpha} - q  \sup_{0\le q\le1/2}
	\bigg|\f{d\bar{\ell}}{dq} \bigg|
	\stackrel{\eqref{e:d.ell.d.q}}{\ge}
	-\f{\psi}{2\alpha} - q  \cdot \CSTONE \cdot C_1(U)^2
	\\
	&\qquad\stackrel{\eqref{e:fp.bounds}}{\ge}
	-3C_1(U)^2
	\bigg(\f12+
	\CSTONE \cdot C_1(U)^2 \alpha
	\bigg)
	\stackrel{\eqref{e:alpha.U}}{\ge}
	-3C_1(U)^2
	\bigg(\f12
	+ \f{1}{e^{10} C_1(U)^4 K_2(U)^4} \bigg)
	\ge -1.51 \cdot C_1(U)^2\,,
	\end{align*}
where the last bound uses that we chose $C_1(U)\ge10$ in the proof of 
Lemma~\ref{l:poly}. Then, recalling \eqref{e:annealed}, we have
	\[
	\annealed(\alpha;U)-\log2
	= \alpha \bar{\ell}(0)
	= -\alpha\log\f1{\E U(Z)}
	\stackrel{\eqref{e:bar.c.one.U}}{\ge}
	- \alpha \log \bar{c}_1(U)
	\ge -\alpha \bar{c}_1(U)
	\ge- \f{\alpha C_1(U)^2}{50}\,,
	\]
using that we also chose $C_1(U)\ge 5 \cdot\bar{c}_1(U)\ge10$ in the proof of 
Lemma~\ref{l:poly}. The claim follows.
\end{proof}
\end{cor}

\begin{lem}\label{l:smoothed.Ktwo}
Suppose $U$ satisfies Assumptions~\ref{a:bdd} and \ref{a:Lip},
and let $U_\ETA=U*\varphi_\ETA$ as in \eqref{e:U.ETA}. Then,
using the notation of Assumption~\ref{a:Lip}, we will have
$K_2(U_\ETA) \le 4 (K_2)'(U)$
for all $\ETA\le1$.

\begin{proof}
Let $\xi,\xi'$ be i.i.d.\ standard gaussian random variables. 
We need to bound the quantity
	\beq\label{e:Ktwo.N.D}
	\f{\E_{\xi,\xi'} [
	(\xi-\xi')^2 U_\ETA(x+c\xi)U_\ETA(x+c\xi')]}
	{\E_{\xi,\xi'} [ U_\ETA(x+c\xi)U_\ETA(x+c\xi')]}
	\equiv \f{N_\ETA(x,c)}{D_\ETA(x,c)}\,.
	\eeq
Let $\zeta,\zeta'$ be independent copies of $\xi,\xi'$, and note that
	\[N_\ETA(x,c)
	= \E_{\xi,\xi',\zeta,\zeta'}\bigg[
	 (\xi-\xi')^2 
	 U(x+c\xi+\ETA\zeta)
	 U(x+c\xi'+\ETA\zeta')\bigg]\,.
	\]
Taking an orthogonal transformation of $(\xi,\zeta)$ gives another pair of i.i.d.\ standard gaussians,
	\[
	\begin{pmatrix}
	X \\ Y
	\end{pmatrix}
	= \f1{(c^2+\ETA^2)^{1/2}} \begin{pmatrix}
	c & \ETA \\ -\ETA & c
	\end{pmatrix} \begin{pmatrix} \xi \\ \zeta\end{pmatrix}\,.
	\]
Likewise we like $(X',Y')$ be the pair obtained by the same transformation applied to $(\xi',\zeta')$. Then note that
	\[
	(\xi-\xi')^2
	=\bigg( \f{c(X-X')-\ETA(Y-Y')}{(c^2+\ETA^2)^{1/2}}\bigg)^2
	\le 2\cdot
	\f{c^2(X-X')^2+\ETA^2(Y-Y')^2 }{c^2+\ETA^2} 
	\,.
	\]
Rewriting $N_\ETA(x,c)$ in terms of the random variables $X,X',Y,Y'$ gives
	\begin{align*}N_\ETA(x,c)
	&\le 2\cdot
	\E_{X,X',Y,Y'}\bigg[
	\bigg(	\f{c^2(X-X')^2+\ETA^2(Y-Y')^2 }{c^2+\ETA^2} 
	\bigg)
		U(x+ (c^2+\ETA^2)^{1/2}X)
		U(x+ (c^2+\ETA^2)^{1/2}X')\bigg] \\
	&= \f{2c^2}{c^2+\ETA^2} N_0(x, (c^2+\ETA^2)^{1/2})
	+\f{4\ETA^2}{c^2+\ETA^2} D_0(x, (c^2+\ETA^2)^{1/2})\,,
	\end{align*}
where $N_0$ and $D_0$ are as in \eqref{e:Ktwo.N.D} but with $U$ in place of $U_\ETA$. If $1/2\le c\le2$ and $\ETA\le1$, then 
$1/2\le (c^2+\ETA^2)^{1/2}\le7/3$, so Assumption~\ref{a:Lip} will give
	\[N_\ETA(x,c)
	\le 
	\bigg(\f{2c^2}{c^2+\ETA^2} (K_2)'(U)
	+\f{4\ETA^2}{c^2+\ETA^2}\bigg) D_0(x, (c^2+\ETA^2)^{1/2})
	\le
	4 (K_2)'(U) D_\ETA(x,c)\,.
	\]
The claim follows.
\end{proof}
\end{lem}

\begin{proof}[\hypertarget{proof:p.RS.U.eps}{Proof of Proposition~\ref{p:RS.U.eps}}]
Recall from Lemma~\ref{l:poly} the constant $C_1(U)$: it depends on the absolute constant $\CSTZRO$, as well as the constant $\bar{c}_1(U)$ defined by \eqref{e:bar.c.one.U}. Let $\xi,\zeta$ be i.i.d.\ standard gaussians, and note
	\[\E_\xi U_\ETA(c\xi)
	= \E_{\xi,\zeta} U\Big(c\xi + \ETA\zeta\Big)
	= \E_\xi U\Big( (c^2+\ETA^2)^{1/2} \xi\Big)\,.
	\]
From this it is clear that $\bar{c}_1(U_\ETA)$ converges to $\bar{c}_1(U)$ as $\ETA\downarrow0$. Next recall from Lemma~\ref{l:smoothed.Ktwo} that if $\ETA\le1$ then we have $K_2(U_\ETA)\le4(K_2)'(U)$. Consequently, 
recalling \eqref{e:alpha.U} and \eqref{e:alpha.p.U}, we have
	\beq\label{e:compare.alphas.U}
	\alpha(U_\ETA)
	\stackrel{\eqref{e:alpha.U}}{\equiv}
	\f1{e^{10}
	\cdot \CSTONE 
		\cdot C_1(U_\ETA)^6
		\cdot K_2(U_\ETA)^4}
	\ge
	\f1{e^{16}
	\cdot \CSTONE 
		\cdot C_1(U)^6
		\cdot (K_2)'(U)^4}
	\stackrel{\eqref{e:alpha.p.U}}{=} \alpha'(U)\,.
	\eeq
This shows that for all $0<\alpha\le\alpha'(U)$, we also have
$\alpha\le\alpha(U_\ETA)$ for all $\ETA$ small enough, which means that the results of Proposition~\ref{p:fp} apply for $U_\ETA$ as well as for $U$. We see from the 
\hyperlink{proof:p.fp}{proof of Proposition~\ref{p:fp}}
that the replica symmetric fixed point $q_\ETA$ for $U_\ETA$ is a root $q_\ETA\in[0,1/25]$ of the function (cf.\ \eqref{e:root.fn})
	\[
	\bar{g}_\ETA(q) = \f{\bar{q}^{-1}(q)}{\alpha}
		- \bar{r}_\ETA(q)\,,
	\]
where $\bar{r}_\ETA$ is defined as in \eqref{e:fp} but with $U_\ETA$ in place of $U$:
	\[
	\bar{r}_\ETA(q)
	= \f1{(1-q)}
	\E\bigg[ \bigg( \f{
	\E_\xi[\xi U_\ETA(Z+(1-q)^{1/}\xi)]
	}{\E_\xi U_\ETA(Z+(1-q)^{1/}\xi)} \bigg)^2\bigg]\,.
	\]
It is clear that $\bar{g}_\ETA$ converges uniformly to $\bar{g}$ over $0\le q \le 1/25$, so $q_\ETA$ converges to $q$, and consequently $\psi_\ETA$ converges to $\psi$. It is then straightforward to deduce from the formula \eqref{e:rs} that $\RS(\alpha;U_\ETA)$ converges to $\RS(\alpha;U)$ as $\ETA\downarrow0$.
\end{proof}

\subsection{Almeida--Thouless condition}
\label{ss:AT}

Recall from Definition~\ref{d:state} the state evolution recursions. 

\begin{lem}\label{l:AMP.well.def}
Suppose $U$ satisfies Assumption~\ref{a:bdd}.
The recursions of Definition~\ref{d:state}
are well-defined: the recursions \eqref{e:def.rho.mu}
lead to $|\rho_s|\le1$ and $|\mu_s|\le1$ for all $s\ge1$,
and the recursions \eqref{e:intro.lm.gm.rec} leads to
$\Lambda_s\in[0,1)$ and $\Gamma_s\in[0,1)$ for all $s\ge0$.

\begin{proof} We abbreviate $F\equiv F_q$ throughout this proof.
We have $\rho_1\equiv\lm_1$ and $\mu_1\equiv\gm_1$ as in \eqref{e:intro.lm.gm.one}, and it follows that
	\[
	0\le (\rho_1)^2 = (\lm_1)^2
	\le \f1q \E \Big[\Th(\psi^{1/2} Z)^2\Big]
	\stackrel{\eqref{e:fp}}{=} 1\,,
	\]
and likewise $0\le (\mu_1)^2= (\gm_1)^2\le1$. 
Then for $s\ge1$ we have $\rho_{s+1}$ and $\mu_{s+1}$ defined by \eqref{e:def.rho.mu}, and it follows by the Cauchy--Schwarz inequality that
	\[
	|\rho_{s+1}|
	\le \f{\E[\Th(\psi^{1/2}Z)^2]}{q}
	\stackrel{\eqref{e:fp}}{=}1\,,\quad
	|\mu_{s+1}|
	\le \f{\alpha\E[F(q^{1/2}Z)^2]}{\psi}
	\stackrel{\eqref{e:fp}}{=}1\,.
	\]
Thus $|\rho_s|\le1$ and $|\mu_s|\le1$ for all $s\ge1$, which confirms that the recursions \eqref{e:def.rho.mu} are well-defined. 

It remains to verify that the quantities $\Lambda_{s-1}$ and $\Gamma_{s-1}$ from \eqref{e:sum.of.squares} are \emph{strictly} smaller than $1$ for all $s\ge1$. The claim holds trivially in the base case $s=0$, since clearly $\Lambda_0=\Gamma_0=0$. We therefore suppose inductively that 
we have $\Lambda_{s-1}<1$ and $\Gamma_{s-1}<1$. This means that the quantities $\lm_s$ and $\gm_s$ are well-defined by the recursions \eqref{e:intro.lm.gm.rec}.
Denote $M_1\equiv q^{1/2}$ and $N_1\equiv (\psi/\alpha)^{1/2}$. Next let $Y_i,X_j$ be a collection of i.i.d.\ standard gaussian random variables, and let
	\begin{align*}
	M_{i+1} &\equiv \Th\bigg( \psi^{1/2}\bigg\{
		\gm_1 Y_1 + \ldots + \gm_{i-1}Y_{i-1}
			+(1-\Gamma_{i-1})^{1/2}Y_i
		\bigg\}\bigg) \\
	N_{j+1} &\equiv F\bigg(q^{1/2}
		\bigg\{
		\lm_1X_1+\ldots+\lm_{j-1}X_{j-1}
			+(1-\Lambda_{i-1})^{1/2}X_j
		\bigg\}
		\bigg)
	\end{align*}
(cf.\ \eqref{e:amp.decomp.H.y} and \eqref{e:amp.decomp.h.x}). This gives well-defined random variables $M_k,N_k$ for all $1\le k\le s+1$, with $\E[(M_k)^2]=q$ and $\E[(N_k)^2]=\psi/\alpha$. If $2\le k<\ell\le s+1$, then
	\begin{align}
	\label{e:M.SCALAR.IDEAL}
	\f{\E(M_k M_\ell)}{q}
	&= \rho\bigg( 
		(\gm_1)^2 + \ldots + (\gm_{k-2})^2 + \gm_{k-1}
		(1-\Gamma_{k-2})^{1/2}\bigg)
	\stackrel{\eqref{e:intro.lm.gm.rec}}{=}
	\rho(\mu_{k-1})
	\stackrel{\eqref{e:def.rho.mu}}{=}
	\rho_k\,,\\
	\f{\E(N_k N_\ell)}{\psi/\alpha}
	&= \mu\bigg(
	(\lm_1)^2 + \ldots + (\lm_{k-2})^2 + \lm_{k-1}
		(1-\Lambda_{k-2})^{1/2}\bigg)
	\stackrel{\eqref{e:intro.lm.gm.rec}}{=}
	\mu(\rho_{k-1})
	\stackrel{\eqref{e:def.rho.mu}}{=}
	\mu_k\,.
	\label{e:N.SCALAR.IDEAL}
	\end{align}
(cf.\ \eqref{e:m.scalar.products}
and \eqref{e:n.scalar.products}).
Now let $R_i,C_i$ ($i\ge1$) be the Gram--Schmidt orthogonalization of the random variables $M_i,N_i$:
	\begin{align}\label{e:gs.M.R.IDEAL}
	R_{i+1}
	&= \f1{r_{i+1}}
	\bigg\{ M_{i+1}-\sum_{j\le i} \E(M_{i+1}R_j)R_j
	\bigg\}\,,\\
	C_{i+1}&=\f1{c_{i+1}}
	\bigg\{ N_{i+1}-\sum_{j\le i} \E(N_{i+1}C_j)C_j
	\bigg\}
	\label{e:gs.N.C.IDEAL}
	\end{align}
where $r_{i+1}$ and $c_{i+1}$ are the normalizing constants such that $\E[(R_{i+1})^2]=1$ and  $\E[(C_{i+1})^2]=1$. To see that $r_{s+1}$ is a well-defined positive number, we apply the inductive hypothesis 
$\Gamma_{s-1}<1$: then follows from the above definition (together with the fact that $\Th$ is a non-constant function) that $M_{s+1}$ depends non-trivially on $Y_s$. On the other hand, the random variables $R_j$ for $j\le s$ can depend only on $Y_1,\ldots,Y_{s-1}$. It follows that the random variable
	\[M_{s+1}-\sum_{j\le s} \E(M_{s+1}R_j)R_j\]
has strictly positive variance, so $r_{s+1}$ is well-defined and positive.
Likewise, using the inductive hypothesis $\Lambda_{s-1}<1$ together with the fact that $F$ is non-constant, we deduce that $c_{s+1}$ is also well-defined and positive. Next, since we see from above that the quantities $\E(M_kM_\ell)$ and $\E(N_kN_\ell)$
depend only on $\min\set{k,\ell}$,
it follows that there is a value $l_j$ such that
$\E(M_{i+1}R_j)=q^{1/2}l_j$ for all $i\ge j$,
and likewise there is a value $y_j$ such that
$\E(N_{i+1}C_j)=(\psi/\alpha)^{1/2}y_j$ for all $i\ge j$. As in \eqref{e:sum.of.squares}, let us abbreviate
	\[
	L_i \equiv \sum_{j\le i} (l_j)^2\,,\quad
	Y_i \equiv \sum_{j\le i} (y_j)^2\,.
	\]
It follows by the above calculations that
	\begin{align*}
	l_{i+1}
	&=
	\f{\E(M_{i+2}R_{i+1})}{q^{1/2}}
	\stackrel{\eqref{e:gs.M.R.IDEAL}}{=} 
	\f{q^{1/2}}{r_{i+1}}
	\bigg\{
	\f{\E(M_{i+2}M_{i+1})}{q}
	-\sum_{j\le i} (l_j)^2
	\bigg\}
	\stackrel{\eqref{e:M.SCALAR.IDEAL}}{=}
	\f{\rho_{i+1}-L_i}
		{(1-L_i)^{1/2}}\,,\\
	y_{i+1}
	&= \f{\E(N_{i+2}C_{i+1})}{(\psi/\alpha)^{1/2}}
	\stackrel{\eqref{e:gs.N.C.IDEAL}}{=} 
	\f{(\psi/\alpha)^{1/2}}{r_{i+1}}
	\bigg\{
	\f{\E(N_{i+2}N_{i+1})}{\psi/\alpha}
	-\sum_{j\le i} (y_j)^2
	\bigg\}
	\stackrel{\eqref{e:N.SCALAR.IDEAL}}{=}
	\f{\mu_{i+1}-Y_i}
		{(1-Y_i)^{1/2}}\,.
	\end{align*}
Recalling that $r_{i+1}$ and $c_{i+1}$ are positive for all $i\le s$, we deduce
	\[
	0< \f{r_{i+1}}{q^{1/2}}
	= \f1{q^{1/2}}\E\bigg[\bigg(
		M_{i+1}-\sum_{j\le i} \E(M_{i+1}R_j)R_j\bigg)^2\bigg]^{1/2}
	= (1-L_i)^{1/2}\,,
	\]
and similarly $0<(1-Y_i)^{1/2}$, which implies $L_i,Y_i\in[0,1)$. We see moreover that the sequences $l_i,m_i$ satisfy the same recursions
\eqref{e:intro.lm.gm.rec} as the sequences $\lm_i,\mu_i$, which implies $l_i=\lm_i$ and $m_i=\mu_i$ for all $i\ge1$. This proves that $\Lambda_i=L_i$ and $\Gamma_i=Y_i$ both lie in $[0,1)$ for all $i\ge1$.
Therefore the recursions \eqref{e:intro.lm.gm.rec} give well-defined quantities  $\lm_s$ and $\gm_s$ for all $s\ge1$, as desired.
\end{proof}
\end{lem}

\begin{lem}[Almeida--Thouless condition]
\label{l:AT.conclusion}
Suppose $U$ satisfies Assumption~\ref{a:bdd}, and moreover that
	\beq\label{e:AT.condition}
	\AT(\alpha;U)\equiv
	\alpha \cdot \bigg\{
		\E\Big( (F_q)'(q^{1/2} Z)^2\Big)
		\bigg\}\bigg\{
		\E\Big(\Th'(\psi^{1/2}Z)^2\Big)
		\bigg\} \le1\,.\eeq
In this case, the recursions of Definition~\ref{d:state}
lead to $\Gamma_s\to1$ and $\Lambda_s\to1$ as $s\to\infty$.

\begin{proof}
We begin with a general observation.
Let $Z,\xi,\xi'$ be i.i.d.\ standard gaussians. Suppose $f:\R\to\R$ is any function with at most polynomial growth, and consider the function
	\[
	r_f(t)\equiv
	\E\bigg[f\Big(t^{1/2} Z + (1-t)^{1/2}\xi\Big)
	f\Big(t^{1/2} Z + (1-t)^{1/2}\xi'\Big)\bigg]\,,
	\]
which is defined for $0\le t\le 1$. Write 
$Z(t)\equiv t^{1/2}Z +(1-t)^{1/2}\xi$, and note that
	\[
	r_f(t)
	=\E\bigg[ \bigg(
	 \E_\xi f\Big(t^{1/2}Z +(1-t)^{1/2}\xi\Big)\bigg)^2\bigg]
	=\E\bigg[ \Big(
	 \E_\xi f(Z(t))\Big)^2\bigg]\ge0\,.
	\]
Next we differentiate with respect to $t$ and apply gaussian integration by parts to obtain
	\begin{align*}
	(r_f)'(t)
	&=  \E \bigg\{
	\Big(\E_\xi f(Z(t))\Big)
	\E_\xi\bigg[ f'(Z(t)) \bigg( \f{Z}{t^{1/2}}
		-\f{\xi}{(1-t)^{1/2}}\bigg)\bigg]
		\bigg\}\\
	&=\E\bigg[ \Big(\E_\xi f'(Z(t))\Big)^2\bigg]
	= r_{f'}(t)\ge0\,.
	\end{align*}
It follows moreover that $(r_f)''(t) =  (r_{f'})'(t)=r_{f''}(t)\ge0$ for all $0\le t\le1$, so $r_f$ is convex.

Now, returning to the state evolution recursions 
from Definition~\ref{d:state}, we will consider $r_S$ and $r_T$ for
	\[
	S(x)\equiv\bigg(\f{\alpha}{\psi}\bigg)^{1/2}
		F(q^{1/2}x)\,,\quad
	T(x)\equiv \bigg(\f1q\bigg)^{1/2}
		\Th(\psi^{1/2}x)\,.
	\]
Denote $r_{ST}\equiv r_S\circ r_T$.
Note that the fixed point equation \eqref{e:fp} implies $r_S(1)=1$ and $r_T(1)=1$, so $r_{ST}(1)=1$. We also have
from \eqref{e:intro.lm.gm.one}
that $r_S(0)=\mu_1$,
while $r_T(0)=\rho_1=0$; so if
 $\mu_1=0$ then $r_{ST}(0)=0$. However, 
the condition \eqref{e:AT.condition} is equivalent to
 $(r_{ST})'(1)<1$, which implies
 $(r_{ST})'(t)<1$ for all $t\in[0,1]$, and consequently
	\[
	1-r_{ST}(0)=r_{ST}(1)-r_{ST}(0)
	= \int_0^1 (r_{ST})'(t)\,dt
	\le (r_{ST})'(1)<1\,.
	\]
This shows that if $(r_{ST})'(1)<1$ then we must have $r_{ST}(0)=r_S(0)=(\mu_1)^2>0$. 

Next we argue that $\rho_2\ne0$. To this end, for the function $T$ we can directly calculate that for all $0\le t\le 1$,
	\[
	(r_T)'(t)=r_{T'}(t) \ge r_{T'}(0) =\Big(\E T'(Z)\Big)^2
	= \f{\psi}{q} \Big(\E\Th'(\psi^{1/2}Z)\Big)^2
	\stackrel{\eqref{e:fp}}{=} \f{\psi(1-q)^2}{q}>0\,,
	\]
so $r_T(t)$ is strictly increasing. Thus, in the case $\mu_1>0$ we obtain
	\[
	\rho_2
	\stackrel{\eqref{e:def.rho.mu}}{=}
	\rho(\mu_1)
	=r_T\Big((\mu_1)^{1/2}\Big)
	>r_T(0)=0\,.
	\]
Since $T$ is an odd function, in the case $\mu_1<0$ we obtain
	\[
	\rho_2\stackrel{\eqref{e:def.rho.mu}}{=}
	\rho(\mu_1)
	=-r_T\Big((-\mu_1)^{1/2}\Big)
	<-r_T(0)=0\,.
	\]
In both cases we obtain $\rho_2\ne0$ as claimed. 

To conclude, note that
 Lemma~\ref{l:AMP.well.def} implies that $\Lambda_s\uparrow\Lambda_\infty\le1$ and $\Gamma_s\uparrow\Gamma_\infty\le1$ as $s\to\infty$. 
If $(r_{ST})'(1)<1$, the above considerations give
$(\Gamma_\infty)^2\ge(\gm_1)^2=(\mu_1)^2>0$, as well as
	\[
	(\Lambda_\infty)^2
	\ge (\lm_2)^2
	\stackrel{\eqref{e:intro.lm.gm.rec}}{=}
	\bigg( \f{\rho_2-\Lambda_1}{(1-\Lambda_1)^{1/2}}
	\bigg)^2
	\stackrel{\eqref{e:intro.lm.gm.one}}{=}
	(\rho_2)^2>0\,.
	\]
Clearly we must also have $\lm_s\to0$ and $\gm_s\to0$ as $s\to\infty$, so
	\[
	\begin{pmatrix}\rho_s \\ \mu_s\end{pmatrix}
	\stackrel{\eqref{e:intro.lm.gm.rec}}{=}
	\begin{pmatrix}
	\Lambda_{s-1}+\lm_s(1-\Lambda_{s-1})^{1/2}\\
	\Gamma_{s-1}+\gm_s(1-\Gamma_{s-1})^{1/2}
	\end{pmatrix}
	\stackrel{s\to\infty}{\longrightarrow}
	\begin{pmatrix}
		\Lambda_\infty \\ \Gamma_\infty\end{pmatrix}
	>\begin{pmatrix}0\\0\end{pmatrix}\,.\]
Thus, for $s$ large enough, we can express
	\[
	\mu_{s+1}
	\stackrel{\eqref{e:def.rho.mu}}{=}
		\mu(\rho_s)
	\stackrel{\eqref{e:def.rho.mu}}{=}
		\mu(\rho(\mu_{s-1})
	=r_{ST}\Big((\mu_{s-1})^{1/2}\Big)\,\,
	\]
which shows that $\Gamma_\infty$ must be a fixed point of $r_{ST}$, and $\Lambda_\infty=r_T(\Gamma_\infty)$. Since we saw above that $r_{ST}(t)$ is convex on the interval $0\le t\le 1$, if $(r_{ST})'(1) \le1$ then the only fixed point of $r_{ST}(t)$ on the interval $0\le t\le1$ occurs at $t=1$, and thus we obtain $\Lambda_\infty=\Gamma_\infty=1$.
\end{proof}
\end{lem}

\subsection{Logconcavity}
\label{ss:logc}

In this subsection we review the \hyperlink{proof:p.logc}{proof of Proposition~\ref{p:logc}} which follows from well-known results on logconcave measures. We then state and prove Lemmas~\ref{l:F.Lipschitz} and \ref{l:TAP.supnorm}, which give some further consequences of Assumption~\ref{a:Lip}.  We also present the \hyperlink{proof:p.AT}{proof of Proposition~\ref{p:AT}}. 

\begin{thm}[{\cite{MR1097258}}]
\label{t:maurey}
Suppose $U$ satisfies Assumption~\ref{a:bdd}
and is logconcave.
Recall that $\varphi$ denotes the standard gaussian density on $\R$, and let $\mu$ be the probability measure on $\R$ whose density
(with respect to Lebesgue measure) is
	\[
	\f{d\mu}{dz} = \f{U(z)\varphi(z)}{\E_\xi U(\xi)}\,.
	\]
Then for any measurable subset $B\subseteq\R$ we have the concentration bound
	\[
	\int \exp\bigg(\f{d(z,B)^2}{4}
	\bigg)\,d\mu(z)
	\le \f1{\mu(B)}\,,
	\]
where $d(z,B)$ denotes the minimum distance from $z$ to $B$.
\end{thm}

Theorem~\ref{t:maurey} is obtained as a consequence of the Pr\'ekopa--Leindler inequality
(or functional Brunn--Minkowski inequality)
\cite{MR315079,MR404557,MR0430188}
from convex geometry; see also \cite{MR1800062} and \cite[Thm.~3.1.4]{MR2731561}. In this paper we use Theorem~\ref{t:maurey} only in the \hyperlink{proof:p.logc}{proof of Proposition~\ref{p:logc}}, which is not needed for the main result Theorem~\ref{t:main}. See \S\ref{sss:lit.spherical} for a discussion of results on the positive spherical perceptron which use convex geometry in more essential ways. By well-known arguments, Theorem~\ref{t:maurey} can be used to deduce the following:

\begin{thm}[{see e.g. \cite[Thm.~3.1.4]{MR2731561}}]
\label{t:talagrand} In the same setting as Theorem~\ref{t:maurey},
if $f:\R\to\R$ is Lipschitz, then
	\[
	\int \f{( f(y)-f(z))^{2k}}{{(16k)^k}}
		\,d\mu(y)\,d\mu(z)
	\le
	\int \exp
	\bigg\{ \f{(f(y)- f(z))^2}{16} 
	\bigg\}\,d\mu(y)\,d\mu(z)\le4
	\]
for any integer $k\ge1$.
\end{thm}

Note that the concentration bounds from Theorems~\ref{t:maurey}~and~\ref{t:talagrand} rely on the strong logconcavity of the gaussian density $\varphi(x)$, and the bounds hold uniformly over all logconcave functions $U$. As a consequence we obtain:

\begin{proof}[\hypertarget{proof:p.logc}{Proof of Proposition~\ref{p:logc}}] Suppose $U$ satisfies Assumption~\ref{a:bdd}. If $U$ is bounded away from zero or compactly supported, then Assumption~\ref{a:Lip} holds by trivial calculations. In the case that $U$ is logconcave, Assumption~\ref{a:Lip} 
follows from the above result Theorem~\ref{t:talagrand}.
\end{proof}

\begin{lem}\label{l:F.Lipschitz}
If $U$ satisfies Assumption~\ref{a:bdd} and \ref{a:Lip}, then the function $F_q$ of \eqref{e:F} satisfies
	\[\|(F_q)'\|_\infty
	\le \f1{1-q}\bigg(
	\f{K_2(U)}{2} + 1
	\bigg)\,.\]
Therefore $F_q$ is Lipschitz for any $q\in[0,1)$.

\begin{proof}
From \eqref{e:F} and \eqref{e:F.rewrite} we calculate
	\beq\label{e:F.prime.prelim}
	(F_q)'(x)
	=\f{\E_\xi U''(x+(1-q)^{1/2}\xi)}
		{\E_\xi U(x+(1-q)^{1/2}\xi)}
	-\bigg(\f{\E_\xi U'(x+(1-q)^{1/2}\xi)}
		{\E_\xi U(x+(1-q)^{1/2}\xi)}\bigg)^2
	\,.\eeq
Applying gaussian integration by parts gives
	\begin{align*}
	(F_q)'(x)
	&=\f1{1-q}\bigg\{ 
		\f{\E_\xi[(\xi^2-1) U(x+(1-q)^{1/2}\xi)]}
		{\E_\xi U(x+(1-q)^{1/2}\xi)}
		-\bigg(\f{\E_\xi[\xi U(x+(1-q)^{1/2}\xi)]}
		{\E_\xi U(x+(1-q)^{1/2}\xi)}\bigg)
		\bigg\} \\
	&=\f1{1-q}
	\bigg\{
	\f12
		\f{\E_{\xi,\xi'}[ (\xi-\xi')^2
		U(x+c\xi)U(x+c\xi')]}
		{\E_{\xi,\xi'}[U(x+c\xi)U(x+c\xi')]}
	  -1\bigg\}\,.
	\end{align*}
The result follows from Assumption~\ref{a:Lip}.
\end{proof}
\end{lem}

\begin{proof}[\hypertarget{proof:p.AT}{Proof of Proposition~\ref{p:AT}}]
In view of Lemma~\ref{l:AT.conclusion},
it suffices to check that the condition
\eqref{e:AT.condition}
holds for $0<\alpha\le\alpha(U)$. 
By Lemma~\ref{l:F.Lipschitz}
and the fact that
 $\Th'(x)\in(0,1)$ for all $x\in\R$, we can bound 
	\begin{align*}
	\AT(\alpha;U)
	&\le
	\f{\alpha}{(1-q)^2}
	\bigg(
	\f{K_2(U)}{2} + 1
	\bigg)^2
	\le \f{(3/2)^2 \cdot  K_2(U)^2 \cdot \alpha}{(1-q)^2}
	\stackrel{\eqref{e:fp.bounds}}{\le}
	 3 \cdot K_2(U)^2 \cdot \alpha \\
	&\stackrel{\eqref{e:alpha.U}}{\le}
	\f{3}{e^{10}  C_1(U)^6 K_2(U)^2} <1 \,,
	\end{align*} 
having used that
$C_1(U)\ge10$
and $K_2(U)\ge1$.
\end{proof}

\newcommand{\MAXEVENT}{\bm{\Omega}}

\begin{lem}\label{l:TAP.supnorm}
Suppose $U$ satisfies Assumption~\ref{a:bdd} and \ref{a:Lip}. Let $\Fprime(t)$ be as in \eqref{e:tap.condition.prime}. Then
	\[
	\max\bigg\{
	\|\bh^{(\ell)}\|_\infty,
	\|\bn^{(\ell)}\|_\infty,
	\|\bH^{(s)}\|_\infty,
	\|\bmag^{(s)}\|_\infty
	:s\le t,\ell\le t-1
	\bigg\}
	\le N^{0.01}
	\]
with probability $1-o_N(1)$.

\begin{proof} Note that Lemma~\ref{l:AMP.PL} implies, for all $1\le \ell\le t-1$ and all $1\le s\le t$,
	\[\lim_{N\to\infty} \f1N
		\sum_{i\le N}
		\bigg( \f{(\bH^{(s)})_i}{\psi^{1/2}}\bigg)^{101}
	=\lim_{N\to\infty} \f1M \sum_{a\le M}
	\bigg( \f{(\bh^{(\ell)})_a}{q^{1/2}}\bigg)^{101}
	= \E(Z^{101})\,,
	\]
where the convergence holds in probability. It follows that the event
	\[
	\MAXEVENT
	\equiv\bigg\{
	\max\bigg\{
	\f1N\sum_{i\le N}
		\bigg( \f{(\bH^{(s)})_i}{\psi^{1/2}}\bigg)^{101},
	\f1M \sum_{a\le M}
	\bigg( \f{(\bh^{(\ell)})_a}{q^{1/2}}\bigg)^{101}
	:s\le t,\ell\le t-1
	\bigg\}
	\le 2\E(Z^{101})
	\bigg\}
	\]
occurs with probability $1-o_N(1)$.
We claim that 
$\MAXEVENT$ implies the desired bounds.
Indeed, $\MAXEVENT$ clearly implies
	\[
	\max\bigg\{
	\|\bmag^{(s)}\|_\infty
	 : s\le t\bigg\}
	\le
	\max\bigg\{
	\|\bH^{(s)}\|_\infty : s\le t\bigg\}
	\le \psi^{1/2} \Big(2N \E(Z^{101})\Big)^{1/101}
	\le N^{1/100}\,.
	\]
In the above, the first inequality uses
that $\bmag^{(s)}=\Th(\bH^{(s)})$ and
$|\Th(x)|\le|x|$; and the last bound
holds for $N$ large enough (depending on $\psi$). Similarly, $\MAXEVENT$ implies
	\[
	\max\bigg\{
	\|\bh^{(\ell)}\|_\infty
	:\ell\le t-1\bigg\}
	\le q^{1/2}
	\Big(N\alpha\E(Z^{101})\Big)^{1/101}
	\le N^{1/100}\,,
	\]
where the last bound holds for $N$ large enough (depending on $\alpha$, $q$). Finally, it follows using Lemma~\ref{l:poly} that
	\[
	\Big|(\bn^{(\ell)})_a\Big|
	=\Big| F((\bh^{(\ell)})_a)\Big|
	\stackrel{\eqref{e:F}}{=}
		\bigg| \f1{(1-q)^{1/2}}
		\f{\E_\xi [\xi U(
		(\bh^{(\ell)})_a + (1-q)^{1/2}\xi
		)]}{\E_\xi U(
		(\bh^{(\ell)})_a + (1-q)^{1/2}\xi
		)}\bigg|
	\le \f{C_1(U) + 4	|(\bh^{(\ell)})_a|
	}{(1-q)^{1/2}}\,,
	\]
and combining with the previous bound 
on $\|\bh^{(\ell)}\|_\infty$ gives
	\begin{align*}
	\max\bigg\{
	\|\bn^{(\ell)}\|_\infty
		: \ell\le t-1\bigg\}
	&=\max\bigg\{ \|F(\bh^{(\ell)})\|_\infty
		: \ell\le t-1\bigg\}\\
	&\le \f{C_1(U) 
		+ 4	q^{1/2}(N\alpha\E(Z^{101}))^{1/101}
	}{(1-q)^{1/2}}
	\le N^{1/100}\,,
	\end{align*}
where the last bound holds for $N$ large enough. This proves the claim.
\end{proof}
\end{lem}

\section{Analysis of first moment}
\label{s:analysis.first.mmt}

In this section we finish analyzing the conditional first moment bound (Theorem~\ref{t:ubd.slice}) obtained in Section~\ref{s:first.mmt}. This leads to the
\hyperlink{proof:t.ubd}{proof of Theorem~\ref{t:ubd}},
our main result on the conditional first moment. From this we can deduce the \hyperlink{proof:t.main.ubd}{upper bound in Theorem~\ref{t:main}}, as presented at the end of this section. For the reader's convenience, we begin by reviewing some important notations. Recall 
from \eqref{e:def.pi.star} that
	\beq\label{e:pi.star}
	\starpi
	\equiv q^{1/2}\bLam^\st \hat{e}_t
	= q^{1/2}
	\begin{pmatrix}\lambda_1 \\ \vdots \\
		\lambda_{t-1} \\
		(1-\Lambda_{t-1})^{1/2}) \end{pmatrix}
		\in\R^t\,.
	\eeq
Recall also from
\eqref{e:def.vpi.star} that we defined
	\beq\label{e:vpi.star}
	\starvpi
	\equiv (1-q)\psi^{1/2} \bGam^\st \acute{e}_{t-1}
	\equiv (1-q)\psi^{1/2} 
		\begin{pmatrix}\gm_1 \\
		\vdots \\ \gm_{t-2} \\
			(1-\Gamma_{t-2})^{1/2})
			\end{pmatrix}\in\R^{t-1}\,.
	\eeq
Given $\pi\in\R^t$ with $\|\pi\|^2\le1$, we denote $c(\pi)\equiv(1-\|\pi\|^2)^{1/2}$.
Next, as in the statement of Theorem~\ref{t:ubd.slice}, given a parameter $\bareps\in\R$
(see \eqref{e:eps.U} below), we let
	\beq\label{e:bX.rewrite}
	\bX(\pi,\vpi)
	\equiv
	\bx[t]^\st\starpi
	+\bigg\{ \bx[t]^\st(\pi-\starpi)
	+ N^{1/2}\bareps\bc[t-1]^\st
		(\vpi-\starvpi)
		\bigg\}
	\in\R^M\,.
	\eeq
We then recall the function $L$ from \eqref{e:L}, and use it to define
	\beq\label{e:Lsum}
	\cL(\pi,\vpi)
	\equiv
	\f1N \sum_{a\le M} L_{\|\pi\|^2}
	(\bX_a(\pi,\vpi))
	= \f1N \sum_{a\le M}
	\log \E_\xi U(\bY_a(\pi,\vpi))\,,
	\eeq
The bound in Theorem~\ref{t:ubd.slice} is expressed in terms of the function
	\beq\label{e:Psi.rewrite}
	\Psi(\pi,\vpi) 
	\equiv
	\f{\|\vpi
	-\bareps(\vpi-\starvpi)\|^2}
		{2c(\pi)^2}
	- \f{(\starvpi,\vpi)}{1-q}
	+ \cL(\pi,\vpi)\,.\eeq
Recall \eqref{e:Z.outside} that we decomposed $\bZ(\Gprime)=\bZ_\circ(\Gprime)+\bZ_\bullet(\Gprime)$.
The rest of this section is organized as follows:
\begin{itemize}
\item In \S\ref{ss:AH.conclusion}
we use Lemmas~\ref{l:AH.pi} and \ref{l:AH.vpi}
to prove Corollary~\ref{ss:AH.conclusion}, which gives a bound on $\bZ_\bullet(\Gprime)$.
This takes care of the case 
$(\pi,\vpi)\notin\bm{N}_\circ$
(see \eqref{e:N.circ}), so in the rest of the section we restrict to
$(\pi,\vpi)\in\bm{N}_\circ$.

\item In \S\ref{ss:stationary}
we prove
Lemmas~\ref{l:dPsi.dpi} and \ref{l:dPsi.dvpi},
which show that the point $(\starpi,\starvpi)$,
as defined by \eqref{e:pi.star} and \eqref{e:vpi.star},
is approximately a stationary point of the function $\Psi$ of \eqref{e:Psi.rewrite}.

\item In \S\ref{ss:hess}
we prove Proposition~\ref{p:hess.Psi}, which bounds
$\Hess\Psi$ for $(\pi,\vpi)\in\bm{N}_\circ$.

\item In \S\ref{ss:rs.ubd} we
combine the results described above
to conclude the
\hyperlink{proof:t.ubd}{proof of Theorem~\ref{t:ubd}}.
We then use this to conclude the \hyperlink{proof:t.main.ubd}{proof of the upper bound in Theorem~\ref{t:main}}.
\end{itemize}
Lastly, we now fix the parameter
\beq\label{e:eps.U}
\bareps = e^5 C_1(U) \alpha^{1/2}
\stackrel{\eqref{e:alpha.U}}{\le}
\f{1}{C_1(U)^2 K_2(U) }\,.
\eeq	
However, this choice of $\bareps$ will not become important until Lemma~\ref{l:hess.P} below.

\subsection{Azuma--Hoeffding bounds}
\label{ss:AH.conclusion}

\begin{cor}\label{c:AH.bound}
If $U$ satisifes Assumptions~\ref{a:bdd} and \ref{a:Lip}, then with high probability we have
	\[
	\E\Big( \bZ_\bullet(\Gprime)
	\,\Big|\,\Fprime(t)\Big)
	\le 
	\exp\bigg\{ N \Big(
		\RS(\alpha;U) - C_1(U)^2 \alpha
		\Big)\bigg\}
	\]
for $\bZ_\bullet(\Gprime)$ as defined by \eqref{e:Z.outside}.

\begin{proof} Recall from Corollary~\ref{c:rs.estimate} that for $\alpha\le\alpha(U)$ we have
	\[\f{\RS(\alpha;U)-\log2}{\alpha}
	\ge-1.53\cdot C_1(U)^2\,,\]
where $C_1(U)\ge10$
is the constant from Lemma~\ref{l:poly}. 
Recalling \eqref{e:Z.outside}, 
we will first bound the case where $\|\vpi(J)-\starvpi\|$ is large. To this end, denote
	\[\bZ_d(\Gprime)
	\equiv\sum_J
	\I\bigg\{ d\le\|\vpi(J)-\starvpi\|\le 2d\bigg\}
	\bQ(J)
	\SAT_J(\Gprime)\,.\]
Thanks to Assumption~\ref{a:bdd}, in Proposition~\ref{p:first.mmt} we also have the trivial bound $\E(\SAT_J(\Gprime)\,|\,\Fprime(t))\le1$.
Substituting this into the calculation \eqref{e:first.mmt.rewrite.Q}
from the \hyperlink{proof:t.ubd.slice}{proof of Theorem~\ref{t:ubd.slice}} gives
	\[
	\f{\E(\bZ_d(\Gprime)\,|\,\Fprime(t))}{
		\exp\{(\oneN,\log(2\Ch(\bH^{(t)})))\}}
	\le
	\sum_{J : d\le\|\vpi(J)-\starvpi\|\le 2d}
	\bQ(J)
	\exp\bigg\{-\f{N(\starvpi,\vpi(J))}{1-q}\bigg\}\,.
	\]
By Lemma~\ref{l:AMP.PL} combined with Jensen's inequality, we have
	\begin{align}\nonumber
	\lim_{N\to\infty}\f{(\oneN,\log\Ch(\bH^{(t)}))}{N}
	&=\E \log\Ch(\psi^{1/2}Z) 
	\le \log \E\Ch(\psi^{1/2}Z) \\
	&=\log \E \exp(\psi^{1/2}Z)
	= \f{\psi}{2}
	\stackrel{\eqref{e:fp.bounds}}{\le}
	\f{3C_1(U)^2 \alpha}{2}
	\,,\label{e:log.cosh.rough.bound}
	\end{align}
where the convergence holds in probability as $N\to\infty$. It follows that, with high probability,
	\[
	\f{\exp\{(\oneN,\log(2\Ch(\bH^{(t)})))\}}
	{\exp\{N\RS(\alpha;U)\}}
	\le
	\exp\bigg\{ N
		\bigg( 1.53 +1.51\bigg)
	C_1(U)^2\alpha
		\bigg\} 
	\le \exp\bigg\{ 3.05\cdot
		NC_1(U)^2\alpha\bigg\}
	\,.
	\]
Next, it follows from \eqref{e:Gamma.matrix} and \eqref{e:def.vpi.star} that
	\[
	\|\starvpi\|
	\stackrel{\eqref{e:def.vpi.star}}{=}
	\psi^{1/2}(1-q)
	\|\bGam^\st\acute{e}_{t-1}\|
	\stackrel{\eqref{e:Gamma.matrix}}{=}
	\psi^{1/2}(1-q)
	\bigg(
	\sum_{\ell\le t-2}(\gm_\ell)^2 + 1-\Gamma_{t-2}
	\bigg)^{1/2}
	=\psi^{1/2}(1-q)\,.
	\]
Note also that if $d\le\|\vpi(J)-\starvpi\|\le 2d$ then
	\[
	-\f{(\starvpi,\vpi(J))}{1-q}
	= -\f{\|\starvpi\|^2+(\starvpi,\vpi(J)-\starvpi)}
		{1-q}
	\le \f{\|\starvpi\|\|\vpi(J)-\starvpi\|}{1-q}
	\le \f{2d\|\starvpi\|}{1-q}
	= 2\psi^{1/2}d\,.
	\]
Combining the above bounds gives, with high probability,
	\begin{align*}
	\f{\E(\bZ_d(\Gprime)\,|\,\Fprime(t))}{\exp\{N\RS(\alpha;U)\}}
	&\le\exp\bigg\{
	N\bigg(
	3.05\cdot C_1(U)^2\alpha+
	2\psi^{1/2} d
	\bigg)
	\bigg\}
	\bQ\bigg(d\le\Big\|\vpi(J)-\starvpi\Big\|\le 2d\bigg) \\
	&\le
	\exp\bigg\{ -N
	\bigg(
	\f{d^2}{2.01}- 2\psi^{1/2} d 
				-3.05\cdot C_1(U)^2\alpha
				+ o_N(1)
\bigg)
	\bigg\}\,,
	\end{align*}
where the last inequality is by
Lemma~\ref{l:AH.vpi}. If we take
$d \ge d_0
\equiv 8 \cdot C_1(U)\alpha^{1/2}$, then we obtain
 	\begin{align*}
	\f{\E(\bZ_d(\Gprime)\,|\,\Fprime(t))}
		{\exp\{N\RS(\alpha;U)\}}
	&\le\exp\bigg\{ -N C_1(U)^2\alpha\bigg(
		\f{8^2}{2.01}
		-2 \cdot 3^{1/2} \cdot 8 - 3.05 - o_N(1)
	\bigg)
	\bigg\} \\
	&\le \exp\bigg\{- 1.1
		\cdot N C_1(U)^2\alpha\bigg\}\,.
	\end{align*}
This concludes our analysis of the case where 
$\|\vpi(J)-\starvpi\|$ is large, so we next turn to the case that
$\|\pi(J)-\starpi\|$ is large. To this end, let us denote
	\[
	\bZ'(\Gprime)
	\equiv
	\sum_J
	\I\bigg\{ \f{\|\vpi(J)-\starvpi\|}
		{C_1(U)\alpha^{1/2}}
		\le 8,
	\f{\|\pi(J)-\starpi(J)\|}
		{C_1(U)\alpha^{1/2}}
		\ge 16
	\bigg\}
	\SAT_J(\Gprime)\,.
	\]
It follows from the previous bounds that
	\begin{align*}
	\f{\E(\bZ'(\Gprime)\,|\,\Fprime(t))}
		{\exp\{N\RS(\alpha;U)\}}
	&\le \exp\bigg\{N\bigg(
		3.05\cdot C_1(U)^2\alpha+
		2\psi^{1/2}
		\cdot 8
		\cdot C_1(U)\alpha^{1/2}
	\bigg)\bigg\}
	\bQ\bigg( \f{\|\pi(J)-\starpi(J)\|}
		{C_1(U)\alpha^{1/2}} \ge 16 \bigg)
		\\
	&\le
	\exp\bigg\{ N C_1(U)^2\alpha\bigg(
	3.05 + 2 \cdot 3^{1/2}\cdot 8
	- \f{16^2 }{8.01} + o_N(1)
	\bigg)\bigg\} \\
	&\le
	\exp\bigg\{ - 1.2 \cdot
		N C_1(U)^2\alpha\bigg\}\,,
	\end{align*}
where the second-to-last inequality is by
Lemma~\ref{l:AH.pi}.
Recalling the 
 definition \eqref{e:Z.outside} of $\bZ_\circ(\Gprime)$,
we have
	\[\bZ_\bullet(\Gprime)
	\le \bZ'(\Gprime) + 
	\sum_{k\ge0} \bZ_{2^k d_0}(\Gprime)\,,
	\]
where $d_0=8 \cdot C_1(U)\alpha^{1/2}$ as above.
It follows by combining the above bounds that
	\[
	\f{\E(\bZ_\bullet(\Gprime)\,|\,\Fprime(t))}
		{\exp\{N\RS(\alpha;U)\}}
	\le \exp\bigg\{ -NC_1(U)^2\alpha\bigg\}
	\]
with high probability, which proves the claim.
\end{proof}
\end{cor}

\subsection{Stationarity at replica symmetric value}
\label{ss:stationary}

In this subsection we show that the function $\Psi(\pi,\vpi)$ from \eqref{e:Psi.rewrite}
is approximately stationary at the point $(\starpi,\starvpi)$.

\begin{lem}\label{l:dPsi.dpi}
Suppose $U$ satisfies Assumption~\ref{a:bdd} and \ref{a:Lip}. Then for all $1\le s\le t$ we have
	\[
	\f{\partial\Psi}{\partial\pi_s}
	(\starpi,\starvpi)
	\simeq0\,,
	\]
where $\simeq$ indicates convergence in probability as $N\to\infty$.

\begin{proof} 
Recalling \eqref{e:F}, \eqref{e:F.rewrite}, and \eqref{e:F.prime.prelim}, we can rewrite
	\beq\label{e:F.prime}
	(F_q)'(x)
	=\f{\E_\xi U''(x+(1-q)^{1/2}\xi)}
		{\E_\xi U(x+(1-q)^{1/2}\xi)}
	- (F_q(x))^2\,.\eeq
Recall from \eqref{e:bX.rewrite} the definition of $\bY\equiv\bY(\pi,\vpi)$. We then calculate
	\begin{align*}
	\f{\partial\cL}{\partial\pi_s}
	&\stackrel{\eqref{e:Lsum}}{=} \f1N\sum_{a\le M}\bigg\{
	\f{\pd\bX_a}{\pd\pi_s}
	\f{\E_\xi U'(\bY_a)}{\E_\xi U(\bY_a)}
	+ \f{\partial c}{\partial \pi_s}
	\f{\E_\xi[ \xi U'(\bY_a)]}{\E_\xi U(\bY_a)}
	\bigg\} \\
	&\stackrel{\eqref{e:F.rewrite}}{=}
	\f1N\sum_{a\le M}\bigg\{
	(\bx^{(s)})_a
	F_{\|\pi\|^2}(\bX_a)
	- \pi_s
	\f{\E_\xi U''(\bY_a)}{\E_\xi U(\bY_a)}
	\bigg\}\\
	&\stackrel{\eqref{e:F.prime}}{=}
	\f1N\bigg\{
	(\bx^{(s)},F_{\|\pi\|^2}(\bX_a))
	-\pi_s (\oneM, (F_{\|\pi\|^2})'(\bX))
	-\pi_s \|F_{\|\pi\|^2}(\bX)\|^2
	\bigg\}
	\end{align*}
It follows from \eqref{e:pi.star} that
$\|\starpi\|^2=q$, and $c_*\equiv c(\starpi)=(1-q)^{1/2}$. We also note that
	\[
	\bh^{(t+1)}
	= \bh[t]^\st\hat{e}_t
	\stackrel{\eqref{e:gs.h.x}}{=}
	q^{1/2}\bx[t]^\st\bLam^\st\hat{e}_t
	\stackrel{\eqref{e:pi.star}}{=}
	\bx[t]^\st\pi_*
	\stackrel{\eqref{e:bX.rewrite}}{=}
	\bX(\pi_*,\starvpi)
	\equiv\bX_*\,.
	\]
It follows using \eqref{e:gs.h.x} and Lemma~\ref{l:AMP.PL} that at $(\starpi,\starvpi)$ we have
	\[
	\f{(\bx^{(s)},F_{\|\pi_*\|^2}(\bX_*))}{N}
	\simeq
	\alpha \Lambda_{t,s}\E Z F_q(q^{1/2}Z)
	\stackrel{\eqref{e:pi.star}}{=}
	\pi_{*,s} \alpha \E (F_q)'(q^{1/2}Z)\,,
	\]
having again used gaussian integration by parts at the last step. As a consequence
	\[
	\f{\partial\cL}{\partial\pi_s}(\starpi,\starvpi)
	\simeq
	 - \pi_{*,s}
	\alpha \E\Big[F_q(q^{1/2}Z)^2\Big]
	\stackrel{\eqref{e:fp}}{=} 
	-\pi_{*,s}\psi\,.
	\]
Substituting this into \eqref{e:Psi.rewrite} gives
	\[\f{\partial\Psi}{\partial\pi_s}(\starpi,\starvpi)
	\simeq
	\f{\|\starvpi\|^2\pi_{*,s}}{(1-\|\starpi\|^2)^2}
	-\pi_{*,s}\psi
	\stackrel{\eqref{e:vpi.star}}{=} 0\,,
	\]
as claimed.
\end{proof}
\end{lem}

\begin{lem}\label{l:dPsi.dvpi}
Suppose $U$ satisfies Assumption~\ref{a:bdd} and \ref{a:Lip}. Then for all $1\le\ell\le t-2$ we have
	\[
	\f{\partial\Psi}{\partial\vpi_\ell}
	(\starpi,\starvpi)
	\simeq0\,,\]
where $\simeq$ indicates convergence in probability as $N\to\infty$. For $\ell=t-1$ we have
	\[
	\f{\partial\Psi}{\partial\vpi_\ell}
	(\starpi,\starvpi)
	\simeq
	\bareps\psi^{1/2}
	\bigg(\gamma_{t-1} - (1-\Gamma_{t-2})^{1/2}\bigg)
	\,,
	\]
where the right-hand side is $o_t(1)$.

\begin{proof}
Similarly to the proof of Lemma~\ref{l:dPsi.dpi}, we calculate
	\[
	\f{\partial\cL}{\partial\vpi_\ell}
	\stackrel{\eqref{e:Lsum}}{=}
	\f1N\sum_{a\le M}
	\f{\partial\bX_a}{\partial\vpi_\ell}
	\f{\E_\xi U'(\bY_a)}{\E_\xi U(\bY_a)}
	= \f{\bareps
		(\bc^{(\ell)},F_{\|\pi\|^2}(\bX))}
		{N^{1/2}}\,.
	\]
It follows by recalling Lemma~\ref{l:AMP.PL} that
	\[
	\f{\partial\cL}{\partial\vpi_\ell}(\starpi,\starvpi)
	=  \f{\bareps
		(\bc^{(\ell)},\bn^{(t+1)})}
		{N^{1/2}}
	\simeq
	\bareps\psi^{1/2} \gamma_\ell\,.
	\]
Substituting this into \eqref{e:Psi.rewrite} gives
	\[\f{\partial\Psi}{\partial\vpi_\ell}
		(\starpi,\starvpi)
	\simeq\bareps\bigg\{-\f{\vpi_{*,\ell}}{1-q} 
		+\psi^{1/2} \gamma_\ell
	\bigg\}\,,
	\]
and combining with \eqref{e:vpi.star} gives the claim.
\end{proof}
\end{lem}

\subsection{Hessian calculation} 
\label{ss:hess}

In this subsection we analyze the Hessian of the function $\Psi$ from \eqref{e:Psi.rewrite} to prove:

\begin{ppn}\label{p:hess.Psi}
If $U$ satisfies Assumptions~\ref{a:bdd} and \ref{a:Lip}, then the function $\cL$ of \eqref{e:Lsum} satisfies
	\[
	\Hess\Psi(\pi,\vpi) 
	= \begin{pmatrix}
		\Psi_{\pi,\pi} & \Psi_{\pi,\vpi} \\
		\Psi_{\pi,\vpi} & \Psi_{\vpi,\vpi}
		\end{pmatrix} \bigg|_{(\pi,\vpi)}
	\preccurlyeq 
	\begin{pmatrix}
	e^7 C_1(U)^2K_2(U) \alpha I & 0\\
	0&1-1.9\bareps\end{pmatrix}
	\]
for all $(\pi,\vpi)\in\bm{N}_\circ$ (as defined by \eqref{e:N.circ}), 
for $0\le \alpha\le \alpha(U)$ as defined by \eqref{e:alpha.U}, and
 $\bareps=\bareps(\alpha;U)$
as in \eqref{e:alpha.U}.
\end{ppn}

The \hyperlink{proof:p.hess.Psi}{proof of Proposition~\ref{p:hess.Psi}} is given at the end of this subsection. We divide the analysis into several steps. Define
	\begin{align}\label{e:diag.A}
	A_c(x)
	&=\f{\E_\xi U''(x+c\xi)}{\E_\xi U(x+c\xi)}
	-\bigg(\f{\E_\xi U'(x+c\xi)}{\E_\xi U(x+c\xi)}
		\bigg)^2
	=(F_{1-c^2})'(x)\\
	B_c(x)
	&\equiv
	\f{\E_\xi[\xi U''(x+c\xi)]}{\E_\xi U(x+c\xi)}
	-\f{\E_\xi[\xi U'(x+c\xi)]}{\E_\xi U(x+c\xi)}
	\f{\E_\xi U'(x+c\xi)}{\E_\xi U(x+c\xi)}
	\label{e:offdiag.b}
	\end{align}
Define the $M$-dimensional vectors $\mathbf{A}\equiv A_{c(\pi)}(\bX)$
and $\mathbf{B}\equiv B_{c(\pi)}(\bX)$. Next let
	\begin{align}\label{e:bar.a}
	a_c(x)
	&\equiv \f{\E_\xi[\xi U'(x+c\xi)]}
		{\E_\xi U(x+c\xi)}\,,\\
	b_c(x)
	&\equiv
	\f{\E_\xi[\xi^2 U''(x+c\xi)]}
		{\E_\xi U(x+c\xi)}
	-\bigg(\f{\E_\xi[\xi U'(x+c\xi)]}
		{\E_\xi U(x+c\xi)}\bigg)^2
	\label{e:bar.b}\,,
	\end{align}
and define the scalars
$\bar{a}\equiv(\oneM,a_{c(\pi)}(\bX))$
and $\bar{b}\equiv(\oneM,b_{c(\pi)}(\bX))$.
	
\begin{lem} \label{l:compute.Hess}
For the function $\cL$ defined by \eqref{e:Lsum} we have
	\begin{align}
	\label{e:L.pi.pi}
	\cL_{\pi,\pi} 
	&=\f1N \bigg\{ \bx[t] (\diag \mathbf{A}) \bx[t]^\st
		+ \bigg(\bx[t] \mathbf{B} (\nabla c)^\st
		+ (\nabla c)(\bx[t] \mathbf{B})^\st
		\bigg)
		 + \bar{a}\cdot \Hess c
			+\bar{b}
			\cdot
			(\nabla c)(\nabla c)^\st
			\bigg\}\,,\\
	\label{e:L.pi.vpi}
	\cL_{\pi,\vpi}
	&= 
	\f{\bareps}{N^{1/2}}
	\bigg\{ \bx[t] (\diag \mathbf{A}) 
	\bc[t-1]^\st
	+ (\nabla c) \bc[t-1]\mathbf{B}
	\bigg\}\,,\\
	\cL_{\vpi,\vpi}
	&= \bareps^2 
	\bigg\{ \bc[t-1]
	(\diag\mathbf{A})
	\bc[t-1]^\st\bigg\}
	\label{e:L.vpi.vpi}
	\end{align}
for $\mathbf{A}$, $\mathbf{B}$, $\bar{a}$, and $\bar{b}$ as defined above.

\begin{proof}
Note that $\bY$ is linear in $\vpi$, with first derivative
	\[
	\f{\partial\bY_a}{\partial\vpi_\ell}
	=\f{\partial\bX_a}{\partial\vpi_\ell}
	= N^{1/2}\bareps (\bc^{(\ell)})_a\,.
	\]
It follows by differentiating \eqref{e:Lsum} twice that
	\begin{align*}
	\f{\partial\cL^2}{\partial\vpi_k\partial\vpi_\ell}
	&= \f1N\sum_{a\le M}
		\bigg\{ \f{\E_\xi[U''(\bY_a)
	\f{\partial\bY_a}{\partial\vpi_k}
	\f{\partial\bY_a}{\partial\vpi_\ell}
	]}{\E_\xi U(\bY_a)}
	- \bigg(\f{\E_\xi[U'(\bY_a)
	\f{\partial\bY_a}{\partial\vpi_k}
	]}{\E_\xi U(\bY_a)}
	\bigg)\bigg(\f{\E_\xi[U'(\bY_a)
	\f{\partial\bY_a}{\partial\vpi_\ell}
	]}{\E_\xi U(\bY_a)}
	\bigg)
	\bigg\} \\
	&=\f1N\sum_{a\le M} \mathbf{A}_a
	\f{\partial\bX_a}{\partial\vpi_k}
	\f{\partial\bX_a}{\partial\vpi_\ell}
	=\bareps^2\bigg\{
	\bc[t-1](\diag \mathbf{A})
	 \bc[t-1]^\st\bigg\}_{k,\ell}\,,\end{align*}
which verifies \eqref{e:L.vpi.vpi}.
On the other hand we note that $\bY$ depends on $\pi$ both through $\bX$ and through $c(\pi)$, and
	\[
	\f{\partial\bY_a}{\partial\pi_s}
	= \f{\partial\bX_a}{\partial\pi_s}
		+ \f{\partial c}{\partial\pi_s} \xi\,.
	\]
We use this to calculate the mixed partial
	\begin{align*}
	\f{\partial\cL^2}{\partial\pi_s\partial\vpi_\ell}
	&=
	\f1N \bigg\{\sum_{a\le M} \mathbf{A}_a
		 \f{\partial\bX_a}{\partial\pi_s}
		  \f{\partial\bX_a}{\partial\vpi_\ell}
		 + 
		 \f{\partial c}{\partial\pi_s}
		\sum_{a\le M} \mathbf{B}_a
		\f{\partial\bX_a}{\partial\vpi_\ell}
		\bigg\}\\
	&= \f{\bareps}{N^{1/2}}
	\bigg( \bx[t] (\diag \mathbf{A}) 
	\bc[t-1]^\st
	+ (\nabla c) \bc[t-1] \mathbf{B}
	\bigg)_{s,\ell}\,,
	\end{align*}
which verifies \eqref{e:L.pi.vpi}.
Finally, a similar calculation gives
	\[\f{\partial\cL^2}{\partial\pi_r\partial\pi_s}
	=\f1N\bigg\{
	\sum_{a\le M} \mathbf{A}_a
		 \f{\partial\bX_a}{\partial\pi_s}
		  \f{\partial\bX_a}{\partial\vpi_\ell}
		 + 
		 \f{\partial c}{\partial\pi_s}
		\sum_{a\le M} \mathbf{B}_a
		\f{\partial\bX_a}{\partial\vpi_\ell}
		+
		\bar{a} \f{\partial^2 c}{\partial\pi_r\partial\pi_s}
		+\bar{b}
		\f{\partial c}{\partial\pi_r}
		\f{\partial c}{\partial\pi_s}
	\bigg\}\,,\]
which implies \eqref{e:L.pi.pi}.
\end{proof}
\end{lem}

We now proceed to bound the quantities
defined above.

\begin{lem}\label{l:bound.B}
Suppose $U$ satisfies Assumptions~\ref{a:bdd} and \ref{a:Lip}. With the notation from  \eqref{e:offdiag.b}, we have
	\[
	\f{\|\mathbf{B}\|}{M^{1/2}}
	=\f{\|B_{c(\pi)}(\bX)\|}{M^{1/2}}
	\le K_2(U)
	\bigg( 2.5	\cdot C_1(U) 
	+5.8 \cdot \f{\|\bX\|}{M^{1/2}}\bigg)\,.
	\]
for all $0.95\le c\le1$.

\begin{proof}
Recalling the notation of Definition~\ref{d:mu.x.c},
we first use gaussian integration by parts to rewrite \eqref{e:offdiag.b} as
	\[
	B_c(x)
	=\f1{c^2}
	\bigg\{	
	\Cov_{x,c}(Z^2,Z)
	-2\cdot \E_{x,c}(Z)
	\bigg\}\,.
	\]
It follows by combining Lemmas~\ref{l:poly}  and \ref{l:Lip.implications} that for all $0.95\le c\le 1$, 
	\[|B_c(x)|
	\le \f1{c^2}\bigg\{
	2\bigg(C_1(U) + \f{1.82 \cdot|x|}{0.95} \bigg)
	+
	\f{K_2(U) }{2^{1/2}}
	\bigg(
	\f{1.82\cdot |x|}{0.95}
	+ C_1(U)^{1/2}
	\bigg)
	\bigg\}\,.\]
Recall also that we assumed (without loss) $C_1(U)\ge10$ and $K_2(U)\ge1$. Therefore
	\begin{align*}|B_c(x)|
	&\le
	\f{K_2(U)}{0.95^2}
	\bigg\{
	\bigg(2 
	+ \f{1}{(2\cdot 10)^{1/2}}
	\bigg) C_1(U)
	+\bigg(2+\f1{2^{1/2}}\bigg)
	\f{1.82\cdot |x|}{0.95}
	\bigg\}\\
	&\le K_2(U)
	\bigg(  2.5 \cdot C_1(U)
	+ 5.8\cdot |x|
	\bigg)\,.
	\end{align*}
where the last bound again uses that $C_1(U)\ge10$.
The claim follows.
\end{proof}
\end{lem}

\begin{lem}\label{l:bound.a}
Suppose $U$ satisfies Assumption~\ref{a:bdd}. With the notation from \eqref{e:bar.a}, we have
	\[\f{|\bar{a}|}{M}
	=\f{| (\oneM,a_{c(\pi)}(\bX))|}{M}
	\le 1.1\cdot C_1(U)
		+ 3.7\cdot \f{\|\bX\|^2}{M}
	\]
for all $0.95\le c\le1$.

\begin{proof}We use gaussian integration by parts to rewrite \eqref{e:bar.a} as
	\[
	a_c(x)
	=\f{\E_\xi[\xi^2U(x+c\xi)]}
		{\E_\xi U(x+c\xi)}-1\,.
	\]
It follows by Lemma~\ref{l:poly}
(which uses only Assumption~\ref{a:bdd}) that
for all $0.95\le c\le1$,
	\[|a_c(x)| \le
	1 + \bigg(C_1(U) + \bigg( \f{1.82\cdot x}{0.95} \bigg)^2\bigg)
	\le 1.1\cdot C_1(U) 
		+ 3.7\cdot x^2
		\,,\]
where the last bound uses that we took $C_1(U)\ge10$.
\end{proof}
\end{lem}

\begin{lem}\label{l:bound.b}
Suppose $U$ satisfies Assumptions~\ref{a:bdd} and \ref{a:Lip}. With the notation from \eqref{e:bar.b}, we have
	\[
	\f{|\bar{b}|}{M}
	=\f{| (\oneM,b_{c(\pi)}(\bX))|}{M}
	\le K_2(U)\bigg( 4.6\cdot C_1(U)
		+ 17 \cdot
		\f{\|\bX\|^2}{M}\bigg)
	\]
for all $0.95\le c\le1$.

\begin{proof}
We use gaussian integration by parts to rewrite \eqref{e:bar.b} as
	\begin{align*}
	b_c(x)
	&=\f1{c^2}
	\bigg\{
	\f{\E_\xi[(\xi^4-5\xi^2+2)U(x+c\xi)]}
		{\E_\xi U(x+c\xi)}
	-\bigg( \f{\E_\xi[\xi^2 U(x+c\xi)]}
		{\E_\xi U(x+c\xi)}-1 \bigg)^2
	\bigg\} \\
	&=\f1{c^2}\bigg\{ \Var_{x,c}(Z^2)
	 -3 \cdot \E_{x,c}(Z^2) +1
	\bigg\}\,.
	\end{align*}
It follows by combining Lemmas~\ref{l:poly} and
 \ref{l:Lip.implications} that for all $0.95\le c\le1$, 
	\begin{align*}
	|b_c(x)|
	&\le \f1{c^2}\bigg\{
	K_2(U) \bigg\{ \bigg(\f{1.82\cdot x}{c}\bigg)^2 
	 + C_1(U)\bigg\}
	+3\bigg(C_1(U) + \bigg(\f{1.82\cdot x}{0.95}\bigg)^2\bigg)+1
	\bigg\} \\
	&\le \f{K_2(U)}{0.95^2}\bigg\{
	\bigg( 4 + \f1{10}\bigg)C_1(U)  
	+4\cdot \bigg( \f{1.82\cdot x}{0.95}\bigg)^2
	\bigg\}\,,
	\end{align*}
where the last bound uses that we took 
$C_1(U)\ge10$ and
$K_2(U)\ge1$. The claim follows.
\end{proof}
\end{lem}

\begin{cor}\label{c:pair.hess} If $U$ satisfies Assumptions~\ref{a:bdd} and \ref{a:Lip}, then the function $\cL$ of \eqref{e:Lsum} satisfies
	\[
	\Hess\cL(\pi,\vpi) 
	= \begin{pmatrix}
		\cL_{\pi,\pi} & \cL_{\pi,\vpi} \\
		\cL_{\pi,\vpi} & \cL_{\vpi,\vpi}
		\end{pmatrix} \bigg|_{(\pi,\vpi)}
	\preccurlyeq K_2(U)
	\begin{pmatrix}
	 C_1(U)^2 \alpha I & 0 \\
	0 & 5\cdot\bareps^2 I
	\end{pmatrix}
	\]
for all $(\pi,\vpi)\in\bm{N}_\circ$ (as defined by \eqref{e:N.circ}), for $0\le \alpha\le \alpha(U)$ as defined by \eqref{e:alpha.U}, and $\bareps=\bareps(\alpha;U)$ as in \eqref{e:alpha.U}.

\begin{proof}
We will bound each of the
terms computed in Lemma~\ref{l:compute.Hess}.
Let $u$ denote any vector in $\R^t$
and let $v$ denote any vector in $\R^{t-1}$. It follows from Lemma~\ref{l:AMP.PL} that, with high probability,
	\beq\label{e:norm.x.bound}
	\sup
	\bigg\{
	\f{\|\bx[t]^\st u\|^2}{M}
	=\f1M\bigg\|\sum_{s\le t} u_s\bx^{(s)}\bigg\|^2
	: u\in\R^{t-1},\|u\|=1
	\bigg\}
	\le 2\,.\eeq
Next, it follows from
\eqref{e:pi.star} that $\|\starpi\|=q^{1/2}$, so the restriction $(\pi,\vpi)\in\bm{N}_\circ$ (see \eqref{e:N.circ}) implies
	\beq\label{e:N.circ.pi}
	\f{\|\nabla c\|}{2}
	\le \|\pi\|
	\le q^{1/2} 
		+ {16\cdot C_1(U)}\alpha^{1/2}
	\stackrel{\eqref{e:fp.bounds}}{\le}
	18\cdot C_1(U)\alpha^{1/2}
	\stackrel{\eqref{e:alpha.U}}{\le}
	\f{18}{e^5 \cdot
		C_1(U)^2 K_2(U)^2 }
	\le \f1{e^6\cdot  K_2(U)^2}\,,
	\eeq
where the last bound uses that we assumed (without loss) $C_1(U)\ge10$ and $K_2(U)\ge1$. Therefore  we certainly have 
$c(\pi)=(1-\|\pi\|^2)^{1/2} \ge  0.95$.
It follows from \eqref{e:diag.A} and Lemma~\ref{l:F.Lipschitz}
(which uses Assumption~\ref{a:Lip}) that
	\beq\label{e:bound.A}
	\|\mathbf{A}\|_\infty
	\le \|A_{c(\pi)}\|_\infty
	= \|(F_{\|\pi\|})'\|_\infty
	\le \f1{c(\pi)^2}
	\bigg\{ \f{K_2(U)}{2}+1\bigg\}
	\le \f{K_2(U)/2 + 1}{0.95^2}
	\le 1.7\cdot K_2(U)\,. 
	\eeq
It follows that, with high probability,
it holds for all unit vectors $u,v$ that
	\begin{align}
	\label{e:bound.L.pi.pi.first}
	u^\st\bigg(
	\f1N \bx[t](\diag \mathbf{A})\bx[t]^\st\bigg)
		u
	&\le 
	\f{\|\mathbf{A}\|_\infty}{N}
	\cdot
	\Big\|\bx[t]^\st u\Big\|^2
	\le 3.4\cdot K_2(U) \alpha\,\\
	\label{e:bound.L.pi.vpi.first}
	u^\st\bigg(
	\f1{N^{1/2}}
	\bx[t]
	(\diag \mathbf{A})\bc[t-1]^\st\bigg)v
	&\le
	\f{\|\mathbf{A}\|_\infty}{N^{1/2}}
	\cdot\Big\|\bx[t]^\st u\Big\|
	\cdot\Big\|\bc[t-1]^\st v\Big\|
	\le 2.5\cdot K_2(U)\alpha^{1/2} \,,\\
	v^\st \bigg(
	 \bc[t-1]
		(\diag\mathbf{A})\bc[t-1]^\st
	\bigg)v
	&\le
	\|\mathbf{A}\|_\infty
	\cdot\Big\|\bc[t-1]^\st v\Big\|^2
	\le 1.7\cdot K_2(U)\,.
	\label{e:bound.L.vpi.vpi}
	\end{align}
Next, recalling \eqref{e:bX.rewrite}, for all $(\pi,\vpi)\in\bm{N}_\circ$ we have 
	\beq\label{e:bX.bound}
	\f{\|\bX(\pi,\vpi)\|}{M^{1/2}}
	\stackrel{\eqref{e:norm.x.bound}}{\le}
	2^{1/2}\|\pi\|
	+\f{\bareps\|\vpi-\starvpi\|}{\alpha^{1/2}}
	\stackrel{\eqref{e:N.circ}}\le
	2^{1/2}\|\pi\|
	+16\cdot \bareps C_1(U) 
	\stackrel{\eqref{e:N.circ.pi}}{\le}
	0.33 \cdot C_1(U)\,,\eeq
having used that $C_1(U)\ge10$ and
$|\bareps|\le  1/50$.
Combining \eqref{e:bX.bound} with 
Lemma~\ref{l:bound.B} gives 
	\[
	\f{\|\mathbf{B}\|}{M^{1/2}}
	\le K_2(U) \bigg( 2.5	\cdot C_1(U) 
	+ 5.8 \cdot
		\f{\|\bX\|}{M^{1/2}}\bigg)
	\le 4.5\cdot C_1(U) K_2(U)\,.
	\]
Combining the above with \eqref{e:norm.x.bound} and \eqref{e:N.circ.pi} gives
that with high probability, for all unit vectors $u,v$ we have
	\begin{align}\nonumber
	u^\st
	\bigg(
	\f1N
	\bx[t] \mathbf{B}(\nabla c)^\st
	\bigg)
	u
	&\le\f{\|\bx[t]^\st u\|
	\|\mathbf{B}\|
	\|\nabla c\|}{N}
	\le
	\f{(2M)^{1/2}}{N}
	\bigg(4.5\cdot C_1(U) K_2(U) M^{1/2}\bigg)
	\bigg(36\cdot C_1(U)\alpha^{1/2}\bigg)
	\\
	&\stackrel{\eqref{e:alpha.U}}{\le}
	\f{2^{1/2}
		\cdot 4.5
		\cdot 36}
		{ e^5 \cdot C_1(U)}\cdot\alpha
	\le 0.16 \cdot  \alpha\,,
	\label{e:bound.L.pi.pi.second} 
	\end{align}
again using that $C_1(U)\ge10$. Similarly, with high probability, it holds for all unit vectors $u,v$ that 
	\begin{align}\nonumber
	u^\st
	\bigg(
	\f1{N^{1/2}}
	(\nabla c) \bc[t-1] \mathbf{B}
	\bigg)v
	&\le \f{\|\mathbf{B}\| \|\nabla c\|}
		{N^{1/2}}
	\le
	\bigg(4.5\cdot C_1(U) K_2(U) \alpha^{1/2}\bigg)
	\bigg(36\cdot C_1(U)\alpha^{1/2}\bigg)\\
	&\stackrel{\eqref{e:alpha.U}}{\le}
	\f{4.5 \cdot 36}
		{e^5\cdot C_1(U)}
		\cdot \alpha^{1/2}
	\le 0.11\cdot \alpha^{1/2} \,.
	\label{e:bound.L.pi.vpi.second}
	\end{align}
Next, combining \eqref{e:bX.bound} with Lemma~\ref{l:bound.b} gives 
	\[\f{|\bar{b}|}{M}
	\le
	K_2(U)\bigg(4.6\cdot C_1(U)
		+ 17 \cdot
	\Big(0.33\cdot C_1(U)\Big)^2\bigg)
	\le 2.4\cdot C_1(U)^2 K_2(U)\,.\]
Combining the above with \eqref{e:N.circ.pi} gives,
for any unit vector $u$,
	\beq\label{e:final.bar.b}
	u^\st\bigg(\f1N\bar{b}(\nabla c)(\nabla c)^\st
	\bigg)u
	\le
	\f{|\bar{b}|
	\|\nabla c\|^2}{N}
	\le\bigg( 2.4\cdot C_1(U)^2 K_2(U)\bigg)
	\cdot
	\bigg(\f{18}{e^5 \cdot
		C_1(U)^2 K_2(U)^2 }\bigg)^2
	\alpha
	\le \f{\alpha}{e^7}\,.
	\eeq
Finally, we note that the Hessian of $c(\pi)=(1-\|\pi\|^2)^{1/2}$ can be calculated as
	\[
	\Hess c(\pi)
	= -\f1{c(\pi)}
	\bigg\{ I + \f{\pi\pi^\st}{c(\pi)^2}\bigg\}\,.
	\]
We can bound the above in operator norm by
	\[\|\Hess c(\pi)\|
	\le \f1{0.95} \bigg( 1 + \f{\|\pi\|^2}{0.95^2}\bigg)
	\stackrel{\eqref{e:N.circ.pi}}{\le}
	\f1{0.95}
	\bigg(1 + \f{(1/e^6)^2}{0.95^2}
		\bigg)\le 1.1\,.\]
Combining \eqref{e:bX.bound} with Lemma~\ref{l:bound.a} gives
	\[
	\f{|\bar{a}|}{M}
	\le  1.1 \cdot C_1(U)
	+3.7\cdot \Big(0.33\cdot C_1(U)\Big)^2
	\le 0.6\cdot C_1(U)^2\,,
	\]
so altogether we obtain, for any unit vector $u$,
	\beq\label{e:final.bar.a}
	u^\st\bigg(\f1N\bar{a}\cdot\Hess c\bigg) u
	\le 1.1\cdot 0.6\cdot C_1(U)^2 \alpha
	\le 0.7\cdot
	C_1(U)^2 \alpha
	\,.\eeq
To conclude, we note that
substituting \eqref{e:bound.L.vpi.vpi}
into \eqref{e:L.vpi.vpi} implies
	\[
	\f{\|\cL_{\vpi,\vpi}\|}{\bareps^2}
	\le 1.7 \cdot K_2(U)\,.
	\]
Substituting \eqref{e:bound.L.pi.vpi.first} and \eqref{e:bound.L.pi.vpi.second}
into \eqref{e:L.pi.vpi} implies
	\[
	\f{\|\cL_{\pi,\vpi}\|}{\alpha^{1/2} \bareps}
	\le 2.5\cdot K_2(U)+ 0.11
	\le 2.7 \cdot  K_2(U)\,.
	\]
Finally, substituting 
\eqref{e:bound.L.pi.pi.first},
\eqref{e:bound.L.pi.pi.second}, 
\eqref{e:final.bar.b}, and
\eqref{e:final.bar.a}
into \eqref{e:L.pi.pi} gives
	\[\f{\|\cL_{\pi,\pi}\|}{\alpha}
	\le 3.4 \cdot  K_2(U)
	+2\cdot 0.16 + \f1{ e^7}  + 0.7\cdot C_1(U)^2
	\le 0.8\cdot C_1(U)^2K_2(U) \,.
	\]
Consequently, for any vector $x\equiv (\dot{x},\ddot{x})$ where $\dot{x}\in\R^t$ and $\ddot{x}\in\R^{t-1}$, we have
	\begin{align*}
	\f{|x^\st (\Hess\cL) x|}{ K_2(U)}
	&\le  0.8\cdot C_1(U)^2 \alpha \|\dot{x}\|^2
	+ 1.7\cdot \bareps^2\|\ddot{x}\|^2
	+2 \cdot 2.7\cdot
		\alpha^{1/2}\bareps
		\|\dot{x}\|\|\ddot{x}\|\\
	&\le
	\bigg(  0.8 \cdot C_1(U)^2 + 2.7 \bigg)
	\alpha \|\dot{x}\|^2
	+
	\Big( 1.7 + 2.7 \Big)
	\bareps^2\|\ddot{x}\|^2\,.
	\end{align*}
The claim follows.
\end{proof}
\end{cor}

Recalling \eqref{e:Psi.rewrite}, let us now denote
	\beq\label{e:cP}
	\cP(\pi,\vpi)
	\equiv\f{\|\vpi
	-\bareps(\vpi-\starvpi)\|^2}
		{2c(\pi)^2}
	- \f{(\starvpi,\vpi)}{1-q}\,,
	\eeq
so that $\Psi=\cP+\cL$.

\begin{lem}\label{l:hess.P}
If $U$ satisfies Assumptions~\ref{a:bdd} and \ref{a:Lip}, then the function $\cP$ of \eqref{e:cP} satisfies
	\[
	\Hess\cP(\pi,\vpi) 
	= \begin{pmatrix}
		\cP_{\pi,\pi} & \cP_{\pi,\vpi} \\
		\cP_{\pi,\vpi} & \cP_{\vpi,\vpi}
		\end{pmatrix} \bigg|_{(\pi,\vpi)}
	\preccurlyeq 
	\begin{pmatrix}
	1080\cdot C_1(U)^2\alpha I & 0\\
	0 & (1-1.95\cdot\bareps) I
	\end{pmatrix}
	\]
for all $(\pi,\vpi)\in\bm{N}_\circ$ (as defined by \eqref{e:N.circ}), for $0\le \alpha\le \alpha(U)$ as defined by \eqref{e:alpha.U}, and
 $\bareps=\bareps(\alpha;U)$
as in \eqref{e:alpha.U}.

\begin{proof} 
We first calculate the mixed partial derivatives
	\begin{align*}
	\cP_{\pi,\pi}
	&=
	\f{ \| \vpi-\bareps(\vpi-\starvpi) \|^2 }{c(\pi)^4}
	\bigg\{ 1+ \f{4\pi\pi^\st}{c(\pi)^2}
		\bigg\}\,,\\
	\cP_{\pi,\vpi}
	&= \f{2(1-\bareps)}{c(\pi)^4}
	\pi\Big( \vpi-\bareps(\vpi-\starvpi)\Big)^\st
	\,,\\
	\cP_{\vpi,\vpi}
	&=\f{(1-\bareps)^2}{c(\pi)^2}I
	= \f{(1-\bareps)^2}{1-\|\pi\|^2}I\,.
	\end{align*}
We have from \eqref{e:vpi.star} that $\|\starvpi\|=(1-q)\psi^{1/2}\le\psi^{1/2}$. Then, for $(\pi,\vpi)\in\bm{N}_\circ$
(as defined by \eqref{e:N.circ}) we must have
	\beq\label{e:N.circ.vpi}
	\|\vpi\|
	\le
	\psi^{1/2}
	+16\cdot C_1(U)\alpha^{1/2}
	\stackrel{\eqref{e:fp.bounds}}{\le} 
	18\cdot C_1(U)\alpha^{1/2}
	\stackrel{\eqref{e:alpha.U}}{\le}
	\f{18}{ e^5\cdot C_1(U)^2}
	\le \f1{e^6}\eeq
(very similarly to \eqref{e:N.circ.pi}). 
It follows using \eqref{e:N.circ.pi} and \eqref{e:N.circ.vpi} that
	\begin{align*}
	\|\cP_{\pi,\pi}\|
	&\le \f{(18 \cdot C_1(U))^2\alpha}{0.95^4}
	\bigg(
	1+ \f{4 \cdot (1/e^6)^2}{0.95^2}
	\bigg)
	\le e^6 \cdot C_1(U)^2 \alpha\,,\\
	\|\cP_{\pi,\vpi}\|
	&\le \f{2( 18\cdot C_1(U))^2\alpha}{0.95^4}
	\le 720 \cdot C_1(U)^2 \alpha
	\stackrel{\eqref{e:eps.U}}{\le}
	\f{\bareps^2}{e^3}
	\,,\\
	\|\cP_{\vpi,\vpi}\|
	&\le \f{(1-\bareps)^2}
		{1 - ( 18 \cdot C_1(U))^2 \alpha}
	\le (1-\bareps)^2 + 
	2 \cdot (18 \cdot C_1(U))^2 \alpha
	\stackrel{\eqref{e:eps.U}}{\le}
	1-2\bareps + 1.03 \cdot \bareps^2\,.
	\end{align*}
Consequently, for any vector $x\equiv (\dot{x},\ddot{x})$ where $\dot{x}\in\R^t$ and $\ddot{x}\in\R^{t-1}$, we have
	\begin{align*}
	|x^\st(\Hess\cP)x|
	&\le 360\cdot  C_1(U)^2 \alpha\|\dot{x}\|^2
	+ 2\cdot 720\cdot  C_1(U)^2 \alpha \|\dot{x}\|
		\|\ddot{x}\|
	+ \bigg(
	1-2\bareps + 1.03 \cdot \bareps^2\bigg)
		\|\ddot{x}\|^2\\
	&\le
	C_1(U)^2
	\Big( 360 + 720  \Big)
	 \alpha \|\dot{x}\|^2
	+\bigg(
	1-2\bareps + 
	1.03 \cdot \bareps^2
	+ \f{\bareps^2}{e^3}
	\bigg)
		\|\ddot{x}\|^2\,.
	\end{align*}
The claim follows.
\end{proof}
\end{lem}

\begin{proof}[\hypertarget{proof:p.hess.Psi}{Proof of Proposition~\ref{p:hess.Psi}}]
It follows by combining
Corollary~\ref{c:pair.hess} and Lemma~\ref{l:hess.P} that
	\[
	\Hess\Psi
	=\Hess\cP+\Hess\cL
	\preccurlyeq
	\begin{pmatrix} C_1(U)^2
	(K_2(U) + 1080 )\alpha I & 0\\
	0& (1-1.95\cdot \bareps 
	+ 5 \cdot K_2(U) \bareps^2)I
	\end{pmatrix}\,.
	\]
We use the choice of $\bareps$ from \eqref{e:eps.U} to bound
	\[
	5 \cdot K_2(U) \bareps^2
	\stackrel{\eqref{e:eps.U}}
	\le 
	\f{5\bareps}{C_1(U)^2} 
	\le 0.05 \cdot \bareps\,,
	\]
and the claim follows.
\end{proof}

\subsection{Replica symmetric upper bound} 
\label{ss:rs.ubd}

In this subsection we give the 
\hyperlink{proof:t.ubd}{proof of Theorem~\ref{t:ubd}}.
We then use this to conclude the \hyperlink{proof:t.main.ubd}{proof of the upper bound in Theorem~\ref{t:main}}.

\begin{proof}[\hypertarget{proof:t.ubd}{Proof of Theorem~\ref{t:ubd}}] Recall from \eqref{e:Z.outside} that we decomposed $\bZ(\Gprime)=\bZ_\circ(\Gprime)+\bZ_\bullet(\Gprime)$. For $\bZ_\circ(\Gprime)$, we will analyze the bound from Theorem~\ref{t:ubd.slice}. Note that Lemma~\ref{l:AMP.PL} implies
	\beq\label{e:log.cosh.limit}
	\f{(\oneN,\log(2 \Ch(\bH^{(t)})))}{N}
	\stackrel{N\to\infty}{\longrightarrow}
	\log + \E\log \Ch(\psi^{1/2}Z)
	\eeq
in probability. Recalling 
\eqref{e:pi.star}, \eqref{e:vpi.star}, and
\eqref{e:Psi.rewrite}, and 
applying Lemma~\ref{l:AMP.PL} again, we have
	\beq\label{e:Psi.star}
	\Psi(\starpi,\starvpi)
	= -\f{\|\starvpi\|^2}{2(1-q)}
		+\cL(\starpi,\starvpi)
	\stackrel{N\to\infty}{\longrightarrow}
	-\f{\psi(1-q)}{2}
		+ \alpha \E L_q(q^{1/2}Z)\eeq
in probability, for $L$ as in \eqref{e:L}. It follows
by comparing
\eqref{e:log.cosh.limit}
and \eqref{e:Psi.star}
 with \eqref{e:rs} that
	\beq\label{e:first.mmt.gives.rs}
	\f{(\oneN,\log(2 \Ch(\bH^{(t)})))}{N}
	+ \Psi(\starpi,\starvpi)
	\stackrel{N\to\infty}{\longrightarrow}
	\RS(\alpha;U)
	\eeq
in probability. Next, it follows by combining Lemmas~\ref{l:AH.pi} and \ref{l:AH.vpi} that
	\begin{align}\nonumber
	&\bQ\bigg(
	\bigg\{ J\in\set{-1,+1}^N
	: \Big\|\pi(J)-\starpi\Big\|\ge d_1
	\textup{ and }
	\Big\|\vpi(J)-\starvpi\Big\|\ge d_2
	\bigg\} 
	\bigg) \\
	&\qquad\le
	\exp\bigg\{ -N
	\bigg[
	\vartheta \f{(1-3q^{1/2})}{8}	(d_1)^2
	+ (1-\vartheta) \f{(1-8q^{1/2})}{2}
		(d_2)^2
	+ o_N(1)
	\bigg]\bigg\}
	\label{e:AH.quadratic}
	\end{align}
for any $\vartheta\in[0,1]$.
On the other hand,
if $(\pi,\vpi)\in\bm{N}_\circ$ (as defined by \eqref{e:N.circ}) with $\|\pi-\starpi\|\le d_1$
and $\|\vpi-\starvpi\|\le d_2$, then it follows by combining Lemmas~\ref{l:dPsi.dpi} and \ref{l:dPsi.dvpi} 
with Proposition~\ref{p:hess.Psi} that
	\begin{align}\nonumber
	\Psi(\pi,\vpi)-\Psi(\starpi,\starvpi)
	&\le \nabla\Psi(\starpi,\starvpi)
		\begin{pmatrix}
		\pi-\starpi\\\vpi-\starvpi
		\end{pmatrix}
	+
	\f{e^7C_1(U)^2 K_2(U) \alpha}{2}(d_1)^2
	+\f{( 1-1.9\bareps)}{2}
		(d_2)^2\\
	&\le
	o_N(1)+ o_t(1)
	+\f{e^7 C_1(U)^2 K_2(U) \alpha }{2}(d_1)^2
	+\f{( 1-1.9\bareps)}{2}
		(d_2)^2\,.
	\label{e:Psi.quadratic}
	\end{align}
Let us take $\vartheta=4\alpha^{1/2}$.
Then,  for $d_1\le \|\pi-\starpi\|\le (1+\alpha)^{1/2}d_1$, combining the 
$\|\pi-\starpi\|^2$ terms in \eqref{e:AH.quadratic} and \eqref{e:Psi.quadratic} results in
	\[
	-\f{4\alpha^{1/2}(1-3q^{1/2})}{8}
	+\f{e^7 C_1(U)^2 K_2(U)\alpha(1+\alpha)}{2}
	\stackrel{\eqref{e:alpha.U}}{\le}
	\bigg(
	-1+3q^{1/2}
	+\f{e^7 (1+\alpha)}
		{e^5 C_1(U) }
	\bigg)
	\f{\alpha^{1/2}}{2}
	\le -\f{\alpha^{1/2}}{10}\,.
	\]
For $d_2\le \|\vpi-\starvpi\|\le (1+\alpha)^{1/2}d_2$,
combining the 
$\|\vpi-\starvpi\|^2$ terms in \eqref{e:AH.quadratic} and \eqref{e:Psi.quadratic} results in
	\begin{align*}
	&-\f{(1-4\alpha^{1/2})(1-8q^{1/2})}{2}
	+ \f{(1+\alpha)(1-1.9\bareps)}{2}
	\stackrel{\eqref{e:fp.bounds}}{\le}
	\bigg(2 + 4 \cdot 3^{1/2} C_1(U) 
	+ \alpha^{1/2}
	\bigg)\alpha^{1/2}
		-\f{1.9\cdot\bareps}{2} \\
	&\qquad\le 8 C_1(U) \alpha^{1/2} - 
		\f{1.9\cdot\bareps}{2}
	\stackrel{\eqref{e:eps.U}}{\le}
	\bigg( 8 - \f{1.9 \cdot e^5 }{2}
	\bigg) C_1(U) \alpha^{1/2}
	\le -1000 \cdot \alpha^{1/2}
	\,.
	\end{align*}
Substituting the above bounds into the result of Theorem~\ref{t:ubd.slice} gives,
with high probability,
	\[
	\f{\E(\bZ_\circ(\Gprime) \,|\,\Fprime(t))}
		{\exp\{N(\RS(\alpha;U) 
		+ o_t(1))
		\}}
	\le 
	\sum_{k_1,k_2\ge0}
	\exp\bigg\{- \f{N\alpha^{1/2}}{10}
	\sum_{i=1}^2 (d_i)^2 (1+\alpha)^{k_i}
	\bigg\}
	\le O(1)\,.
	\]
The result follows by combining with
the bound on $\bZ_\bullet(\Gprime)$ from Corollary~\ref{c:AH.bound}.
\end{proof}

\begin{proof}[\hypertarget{proof:t.main.ubd}{Proof of Theorem~\ref{t:main} upper bound}]
It follows from Theorem~\ref{t:ubd}
and Markov's inequality that for any $\epsilon>0$,
	\[\P\bigg(
	\f1N \log \bZ(\Gprime) \ge \RS(\alpha;U) + \epsilon
	\,\bigg|\,\Fprime(t)\bigg)
	\le
	\f{\exp(N o_t(1))}{\exp(N\epsilon)}\,,
	\]
with high probability over the randomness of $\Fprime(t)$.
It follows that
	\[
	\P\bigg(\f1N \log \bZ(\Gprime) \ge \RS(\alpha;U) + \epsilon\bigg)
	\le o_N(1)
	+ \f{\exp(N o_t(1))}{\exp(N\epsilon)}\,.
	\]
The left-hand side does not depend on $t$, so it follows that
	\[\limsup_{N\to\infty}
	\f1N \log \bZ \le \RS(\alpha;U) \]
in probability, which gives the upper bound in Theorem~\ref{t:main}.
\end{proof}

\section{Second moment conditional on AMP}
\label{s:second.mmt}

In this section we give the \hyperlink{proof:t.ubd}{proof of Theorem~\ref{t:lbd}}, our main result on the conditional second moment. From this we will deduce the
\hyperlink{proof:t.main.lbd.bdd}{lower bound in Theorem~\ref{t:main} in the bounded case}, as explained at the end of this section. The \hyperlink{proof:t.main.conclusion}{lower bound in the general case} will be treated in Section~\ref{s:conc.log.Z}. Recalling \eqref{e:N.circ},
we now restrict further to 
	\beq\label{e:N.star}
	\bm{N}_*
	\equiv
	\bigg\{
	(\pi,\vpi)
	: \max\Big\{
	\|\pi(J)-\starpi\|,
	\|\vpi(J)-\starvpi\|
	\Big\} \le o_N(1)
	\bigg\}\,,\eeq
so $\bm{N}_*\subseteq\bm{N}_\circ$.
Then, analogously to \eqref{e:H.circ}, we let
	\[\mbH_*
	\equiv \bigg\{ J\in\set{-1,+1}^N:
		(\pi(J),\vpi(J))\in\bm{N}_*\bigg\}\,,
	\]
so $\mbH_*\subseteq\mbH_\circ$. Analogously to
\eqref{e:Z.outside}, we let
	\beq\label{e:Z.star}
	\bZ_*(\bG)
	\equiv
	\sum_{J\in\mbH_*}\SAT_J(\bG)
	\le \bZ_\circ(\bG)
	\le \bZ(\bG)\,.
	\eeq
We will prove Theorem~\ref{t:lbd} for the random variable
	\beq\label{e:Z.bar}
	\bar{\bZ}(\bG)
	\equiv
	\sum_{J\in\mbH_*}\SAT_J(\bG)
	\I\bigg\{
	\f{\|\bG\bv_J\|^2}{M} \le 5C_1(U)^2
	\bigg\}
	\le \bZ_*(\bG)\,,
	\eeq
where $\bv_J=J''/\|J''\|$ as in Definition~\ref{d:J.params}, and $C_1(U)$ is the constant from Lemma~\ref{l:poly}.
The remainder of this section is organized as follows:
\begin{itemize}
\item In \S\ref{ss:first.mmt.lbd} we 
\hypertarget{proof:restricted.first.mmt.lbd}{prove the 
first moment lower bound
\eqref{e:lbd.restricted.first.mmt}}, which is the first assertion of Theorem~\ref{t:lbd}.

\item In \S\ref{ss:pair} we introduce a parameter $\lm=\lm(J,K)$ (Definition~\ref{d:lambda}) which captures the correlation of a pair of configurations $J,K\in\set{-1,+1}^N$. We then prove
Theorem~\ref{t:slice.second.mmt} which
gives a preliminary bound on the second moment contribution from pairs with small $\lm$
(see \eqref{e:pair.small.lm}). We also prove Corollary~\ref{c:small.lambda} which bounds 
 the second moment contribution from pairs with larger $\lm$.

\item In \S\ref{ss:second.mmt.analysis}
we further analyze the bound obtained in Theorem~\ref{t:slice.second.mmt}. We show in Proposition~\ref{p:stationary.lambda} that the bound is approximately stationary at $\lm=0$, and then in
Corollary~\ref{c:pair.hess} we control the second derivative of the bound with respect to $\lm$.

\item In \S\ref{ss:second.mmt.conclusion}
we combine the results of the preceding sections to 
\hyperlink{proof:t.lbd}{conclude the proof of Theorem~\ref{t:lbd}}. From this we \hyperlink{proof:t.main.lbd.bdd}{deduce the lower bound of Theorem~\ref{t:main}
in the case $\|u\|<\infty$}.
\end{itemize}
The calculation of this section follows a similar outline as that of Sections~\ref{s:first.mmt} and \ref{s:analysis.first.mmt}, so we will point out the parallels throughout. As before, we let $\bG$ be an independent copy of $\Gprime$.

\subsection{First moment lower bound}
\label{ss:first.mmt.lbd}

In this subsection we 
\hyperlink{proof:restricted.first.mmt.lbd}{prove \eqref{e:lbd.restricted.first.mmt}}, the first assertion of Theorem~\ref{t:lbd}. To this end, we begin with the following result which essentially says that the upper bound of
Theorem~\ref{t:ubd.slice} is tight in the case 
$(\pi,\vpi)=(\starpi,\starvpi)$.

\begin{ppn}\label{p:restricted.first.mmt}
Suppose $U$ satisfies Assumptions~\ref{a:bdd} and \ref{a:Lip}.
Let $\Fprime(t)$ be as in \eqref{e:tap.condition.prime}.
For $\bZ_*$ as in \eqref{e:Z.star} we have
	\[\E\Big(\bZ_*(\Gprime)\,\Big|\,\Fprime(t)\Big)
	\ge
	\exp\bigg\{N\Big(\RS(\alpha;U) 
	- o_t(1)
	\Big)\bigg\}\]
with high probability.

\begin{proof} Recall from the 
\hyperlink{proof:p.first.mmt}{proof of Proposition~\ref{p:first.mmt}} that
	\beq\label{e:tilt.gaussian.rewrite}
	E_J
	\equiv \E\Big(\SAT_J(\Gprime)\,\Big|\,\Fprime(t)\Big)
	\stackrel{\eqref{e:tilt.gaussian}}{=}
	\f{ \bE_J(\tau\,|\,\bar{g}_\ROW)
	\cdot 
	\ap_{J,\tau}(\bar{g}_\ADM\,|\,\bar{g}_\ROW)
	}{\exp\{N^{1/2}(\tau,\bar{g}_\ADM)\}
	\cdot
	p_\ADM(\bar{g}_\ADM)}\,.
	\eeq
If $J\in\mbH_*$, then
it follows from Lemma~\ref{l:change.of.basis.APPROX} that
$\cpi(J)\simeq\cpi_*\equiv q^{1/2}\acute{e}_{t-1}$, and
	\beq\label{e:delta.star.is.zero}
	\f{\bar{g}_\ADM}{N^{1/2}}
	\stackrel{\eqref{e:adm.orth.SECOND}}{=}
	(\bGam_N)^\st \delta(J)
	\simeq 
	\f1{(1-q)^{1/2}} \bigg\{
	\starvpi
	- \f{\psi^{1/2}}{q^{1/2}}(1-q)
		\bGam^\st \cpi_*\bigg\}
	\stackrel{\eqref{e:def.vpi.star}}{=} \mathbf{0}\in\R^{t-1}\,.
	\eeq
Substituting this into the result of 
Proposition~\ref{p:density.bound} gives
	\beq\label{e:local.clt.ideal}
	\f{\ap_{J,\bar{\tau}}(\bar{g}_\ADM\,|\,\bar{g}_\ROW)}{\psi^{1/2} |\det \bGam_N|}
	\simeq
	\f{\bg_{J,\bar{\tau}}(\bar{g}_\ADM)}
	{\psi^{1/2} |\det \bGam_N|}
	\stackrel{\eqref{e:CLT.transformed}}{=}
	g_{J,\bar{\tau}}
		\bigg(- (N\psi)^{1/2}
		\bigg[ \bGam_N \bar{\tau}
			+ c(\pi)
			\f{\bn[t-1]F_{\|\pi\|^2}
				(\bX_{J,\bar{\tau}})}{N\psi^{1/2}}
			+ o_N(1) \bigg]
			\bigg)\,,
	\eeq
for $\bX_{J,\bar{\tau}}$ as defined by \eqref{e:bX.J.tau.rewrite}.
To evaluate the right-hand side above, note that $J\in\mbH_*$ implies
	\[
	\f{\bGam_N\bar{\tau}}{\psi^{1/2}(1-q)^{1/2}}
	\equiv
	\f{\bGam_N\bar{\tau}(J)}{\psi^{1/2}(1-q)^{1/2}}
	\stackrel{\eqref{e:bar.tau}}{\simeq}
	- \f{\bGam\starvpi}{\psi^{1/2}(1-q)}
	\stackrel{\eqref{e:def.vpi.star}}{=}
	-\bGam\bGam^\st\acute{e}_{t-1}
	\stackrel{\eqref{e:Gamma.matrix}}{=}
	-\begin{pmatrix}
	\mu_1 \\ \vdots \\ \mu_{t-2} \\ 1
	\end{pmatrix}
	\,.
	\]
It also implies $\bX_{J,\bar{\tau}} \simeq\bh^{(t+1)}$, and consequently
	\[c(\pi)
			\f{\bn[t-1]F_{\|\pi\|^2}(  \bX_{J,\bar{\tau}} )}{N \psi(1-q)^{1/2}}
	\simeq
	\f{\bn[t-1] F_q(\bh^{(t+1)})}{N \psi}
	\stackrel{\eqref{e:n.scalar.products}}{\simeq}
	\begin{pmatrix}
	\mu_1 \\ \vdots \\ \mu_{t-2} \\ \mu_{t-1}
	\end{pmatrix}\,.
	\]
Note moreover that
Proposition~\ref{p:AT} and
Lemma~\ref{l:AT.conclusion} together imply
$\mu_{t-1}=1-o_t(1)$. Substituting these calculations into \eqref{e:local.clt.ideal} gives (cf.\ \eqref{e:density.bound})
	\beq\label{e:local.CLT.LBD}
	\ap_{J,\bar{\tau}}(\bar{g}_\ADM\,|\,\bar{g}_\ROW)
	\simeq
	\psi^{1/2}|\det\bGam_N|
	g_{J,\tau}\left( (N\psi)^{1/2}
	\begin{pmatrix}
	o_N(1) \\ \vdots \\ o_N(1) \\ o_t(1)
	\end{pmatrix}
	\right)
	= \exp\{No_t(1)\}\,.\eeq
Substituting \eqref{e:p.adm.uncond}, \eqref{e:remaining.E.factor}, and \eqref{e:local.CLT.LBD} into \eqref{e:tilt.gaussian} gives (cf.\ \eqref{e:tcA})
	\[
	E_J
	= \exp\bigg\{
	N\Big[ \tcA_J(\bar{\tau}) + o_t(1)\Big]
	\bigg\}\,.
	\]
It then follows from the 
\hyperlink{proof:t.ubd.slice}{proof of Theorem~\ref{t:ubd.slice}} that
(cf.\ \eqref{e:first.mmt.rewrite.Q})
	\beq\label{e:pair.lbd.rewrite}
	\f{\E(\bZ_*(\Gprime)
	 \,|\,\Fprime(t))}{
		\exp\{(\oneN,\log(2\Ch(\bH^{(t)})))\}}
	=
	\bQ(\mbH_*)
	\exp\bigg\{
	N\Big[\Psi(\starpi,\starvpi) + o_t(1)\Big]
	\bigg\}\,.
	\eeq
We have $\bQ(\mbH_*)\simeq1$ by the law of large numbers, so the claim follows by recalling
\eqref{e:first.mmt.gives.rs}.
\end{proof}
\end{ppn}

To finish the  proof of \eqref{e:lbd.restricted.first.mmt}, it remains only to account for the restriction on $\|\bG\bv\|$ in \eqref{e:Z.bar}:

\begin{proof}[\hypertarget{proof:restricted.first.mmt.lbd}{Proof of 
first moment lower bound
\eqref{e:lbd.restricted.first.mmt}}]
We begin with an easy large deviations calculation.
If $\zeta$ is a standard gaussian random variable, then
it is well known that
$\zeta^2/2$ is a gamma random variable with shape parameter $1/2$, and moment-generating function
	\[
	\E\exp\bigg( \f{\theta\zeta^2}{2}\bigg)
	= \int_0^\infty \f{e^{-(1-\theta) x}}
		{2 \pi^{1/2} x^{1/2}}\,dx
	= \f1{(1-\theta)^{1/2}}\,,
	\]
for any $\theta<1$. If $\bze$ is a standard gaussian random vector in $\R^M$, then for any $L>1$ we have
	\begin{align}\nonumber
	\P\bigg( \f{\|\bze\|^2}{M} \ge L\bigg)
	&\le 
	\exp\bigg\{
	-\f{M}{2}
	\inf
	\bigg\{ 
	\log(1-\theta) + L \theta
	:\theta\in[0,1) \bigg\}
	\bigg\}\\
	&\le \exp\bigg\{
	-\f{M}{2}
	\Big[ L -\log L - 1\Big]
	\bigg\}\,.
	\label{e:chi.sq.large.devs}
	\end{align}
Now, recalling \eqref{e:Z.star} and \eqref{e:Z.bar}, let us take $L \equiv L_1(U) \equiv 5C_1(U)^2 \ge 500$ and define
	\beq\label{e:large.devs.Z}
	\bZ''(\Gprime)
	\equiv
	\bZ_*(\Gprime)-\bar{\bZ}(\Gprime)
	\equiv \sum_J
	\SAT_J(\Gprime)
	\I\bigg\{ \f{\|\Gprime\bv_J\|^2}{M} >  L\bigg\}\,.
	\eeq
It follows from Lemmas~\ref{l:resampling} and \ref{l:adm} that
	\[
	\E\Big(\bZ''(\Gprime)\,\Big|\,\Fprime(t)\Big)
	=
	\E\Big(\bZ''(\bG)\,\Big|\,\ROW,\ADM,
		(\Gprime)_{\ROW\ADM} \Big)
	\le
	\sum_J \P\bigg(
		\f{\|\bG\bv_J\|^2}{M} > L
	\,\bigg|\,\ROW,\ADM,
		(\Gprime)_{\ROW\ADM}
		\bigg)\,.
	\]
Recall from Definition~\ref{d:row.col}
that $V_\COL$ denotes the span of the vectors $\bc^{(\ell)}$ for $\ell\le t-1$. Let us decompose $\bG\bv_J\in\R^M$ as $(\bG\bv_J)^\parallel+(\bG\bv_J)^\perp$ where $(\bG\bv_J)^\parallel$ is the orthogonal projection of $\bG\bv_J$ onto $V_\COL$. Conditional on the events $\ROW$ and $\ADM$, 
$(\bG\bv_J)^\parallel$ is fixed by the admissibility condition (see \eqref{e:adm.orth}), while 
$(\bG\bv_J)^\perp$ behaves as an independent standard gaussian random vector in the orthogonal complement of $V_\COL$. It follows that, conditional  on $\ROW$ and $\ADM$, $\|\bG\bv_J\|^2/M$ is equidistributed as
	\[
	\f{N\|(\bGam_N)^\st\delta\|}{M}
	+ \f{\|\bze'\|^2}{M}
	= o_N(1) + \f{\|\bze'\|^2}{M}\,,
	\]
where $\bze'$ is a standard gaussian random vector in $\R^{M-\ell-1}$. 
It follows by applying 
\eqref{e:chi.sq.large.devs} that
	\begin{align*}
	\E\Big(\bZ''(\bG) \,\Big|\,\FF(t)\Big)
	&\le 2^N
	\P\bigg( \f{\|\bze\|^2}{M} \ge L\bigg)
	\le \exp\bigg\{
	N\bigg[ \log 2 - \f{5\alpha C_1(U)^2}{3}
		\bigg]
	\bigg\}\\
	&\le
	\exp\bigg\{
	N\bigg[ \RS(\alpha;U)
	- \f{\alpha C_1(U)^2}{10}
	\bigg]
	\bigg\}\,,
	\end{align*}
where the last bound uses the result of Corollary~\ref{c:rs.estimate}. 
Combining with the result of Proposition~\ref{p:restricted.first.mmt} gives
	\[
	\E\Big( \bar{\bZ}(\Gprime) \,\Big|\,\Fprime(t)\Big)
	\ge
	\E\Big( \bZ_*(\Gprime) \,\Big|\,\Fprime(t)\Big)
	-\E\Big( \bZ''(\Gprime) \,\Big|\,\Fprime(t)\Big)
	\ge \exp \bigg\{ N\Big(
		 \RS(\alpha;U) - o_t(1)
		 \Big)\bigg\}\,,
	\]
with high probability.
\end{proof}

\subsection{Expected weight of a correlated pair}
\label{ss:pair}

\begin{dfn}\label{d:measure.Q.J}
Recall the function $\bm{S}_J(\bg_\ROW,\bg_\ADM,\bg_\BBB)$ from \eqref{e:S.tau.fn}. Moreover recall that by 
\eqref{e:g.adm} and
\eqref{e:g.prf.minus.adm} combined, the pair
$(\bg_\ADM,\bg_\BBB)$ is equivalent to $\bg_\PRF\equiv\bG\bv$. We let $\Q_J(\cdot)$ denote the measure on $\R^M$ such that
	\[
	\Q_J(B)
	= \f{\E (\SAT_J(\bG) \Ind{\bG\bv\in B} \,|\,\ROW,\
	\ADM)}{\E (\SAT_J(\bG) \,|\,\ROW,\
	\ADM)}
	= \f{\E(\bm{S}_J(\bar{g}_\ROW,\bar{g}_\ADM,\bg_\BBB)
		\Ind{(\bar{g}_\ADM,\bg_\BBB) \in B}
		)}{\E(\bm{S}_J\bar{g}_\ROW,\bar{g}_\ADM,\bg_\BBB)
		)}\,.
	\]
Note that $\Q_J$ depends on $\bar{g}_\ROW$ and $\bar{g}_\ADM$, where $\bar{g}_\ROW$ does not depend on $J$, but $\bar{g}_\ADM$ does. 
\end{dfn}

\begin{dfn}[analogous to Definition~\ref{d:pi.vpi}]
\label{d:lambda}
Let $J,K\in\set{-1,+1}^N$. Recall from 
Definition~\ref{d:pi.vpi}
that we decompose $J=J'+J''$ where $J'$ is the orthogonal projection of $J$ onto the span of the vectors $\bmag^{(s)}$, $1\le s\le t$. Analogously decompose $K=K'+K''$.  Recall that $\bv\equiv J''/\|J''\|$, and define analogously $\bv_K\equiv K''/\|K''\|$. Then let
	\[
	\lm(J,K)
	\equiv 
	\bigg(
	\f{J''}{\|J''\|},
	\f{K''}{\|K''\|}
	\bigg)
	= (\bv,\bv_K)\,,
	\]
so clearly we have $-1\le\lm(J,K)\le1$.
We further denote
	\beq\label{e:K.w}
	\bw \equiv 
	\f{K'' - (K'',\bv)\bv}{\|K'' - (K'',\bv)\bv\|}
	\simeq \f{\bv_K - \lm \bv}{(1-\lm^2)^{1/2}}
	\,,
	\eeq
so $\bw$ is a unit vector in $\R^N$ orthogonal to $\bv$.
\end{dfn}

\begin{dfn}[analogous to Definition~\ref{d:prf.adm.subspaces}]
\label{d:pair.prf.adm.subspaces}
Given $\Fprime(t)$ as in \eqref{e:tap.condition.prime},  and $J,K\in\set{-1,+1}^N$,
recall from Definition~\ref{d:lambda} that we decompose
$J=J'+J''$ and $K=K'+K''$, and define corresponding unit vectors $\bv$ and $\bw$. Let
	\begin{align*}
	V_{\PRF(K)}
	&\equiv \spn\bigg\{
	\eM_a\bw^\st : 1\le a\le M
	\bigg\}\,,\\
	V_{\ADM(K)}
	&\equiv \spn\bigg\{
	\bn^{(\ell)}\bw^\st : 1\le \ell\le t-1
	\bigg\}\,.\end{align*}
Note $V_{\ADM(K)}$ is a subspace of $V_{\PRF(K)}$, and is also a subspace of the space $V_\COL$ from Definition~\ref{d:row.col}.
Let $\proj_{\ADM(K)}$ denote the orthogonal projection onto $V_{\ADM(K)}$, and note that 
$(\Gprime)_{\ADM(K)}s\equiv\proj_{\ADM(K)}(\Gprime)$ is measurable with respect to $\Fprime(t)$.
\end{dfn}

\begin{dfn}[analogous to Definition~\ref{d:ADM.event}]
\label{d:pair.ADM.event} As before, let $\bG$ be an independent copy of $\Gprime$. Let
	\beq\label{e:pair.adm}
	\ADM(K)\equiv
	\Big\{\proj_{\ADM(K)}(\bG)=(\Gprime)_{\ADM(K)}\Big\}
	\stackrel{\eqref{e:col}}{=} 
	\bigg\{
	\f{\bn[t-1]	\bG\bw}{N\psi^{1/2}}
	=\f{\bH[t-1]\bw}{(N\psi)^{1/2}}
	\bigg\}\,,
	\eeq
where the last identity holds assuming $\bG$ belongs to the event $\COL$ from \eqref{e:col}.
\end{dfn}

\begin{dfn}[extension of Definition~\ref{d:cube.P.Q}]
\label{d:pair.P.Q}
We now let $\bP$ denote the uniform probability measure over pairs $(J,K)\in(\set{-1,+1}^N)^2$, and let $\bQ$ be the probability measure on the same space which is given by
	\[
	\f{d\bQ}{d\bP}
	=\f{\exp\{(\bH^{(t)},J+K)\}}
		{\exp\{2\cdot(\oneN,\log\Ch\bH^{(t)})\}}\,.
	\]
Note that $J$ and $K$ are independent under $\bQ$, and each has mean $\bmag^{(t)}$.
\end{dfn}

\begin{ppn}[analogous to Proposition~\ref{p:first.mmt}]
\label{p:ubd.pair.lambda} For $\bze\in\R^M$ define
	\[
	\cA_2(\lm\,|\,\bze)
	\equiv
	\f{\psi(1-q)}{2(1-\lm^2)}
	 + \f1N\bigg(\oneN,
	L_{q+\lm^2(1-q)}
	\bigg(
	\bh^{(t+1)}
	+ (1-q)^{1/2}\lm\bze
	\bigg)\bigg)\,.
	\]
Then, for $J,K\in\mbH_*$, we have
	\[
	\f{\E(\SAT_J(\Gprime)\SAT_K(\Gprime)
	\Ind{\|\Gprime\bv\|^2/M \le L}
	\,|\,\Fprime(t))}
	{\E(\SAT_J(\Gprime)
	\,|\,\Fprime(t))}
	\le 
	\int 
	\I\bigg\{\f{\|\bze\|^2}{M} \le L\bigg\}
	\exp\{N\cA_2(\lm\,|\,\bze)\}
	 \,\Q_J(d\bze)
	\]
with $\Q_J$ as in Definition~\ref{d:measure.Q.J},
and $\lm=\lm(J,K)$ as in Definition~\ref{d:lambda}.
\end{ppn}

In preparation for the 
\hyperlink{proof:p.ubd.pair.lambda}{proof of Proposition~\ref{p:ubd.pair.lambda}}, we record the following calculation:

\begin{lem}[analogous to Lemma~\ref{l:adm.cgf}]
\label{l:pair.adm.cgf}
For $J,K\in\set{-1,+1}^N$
and $\bze\in\R^M$, define the cumulant-generating function
	\[\tcK_{K|J}(\tau\,|\,\bze)
	\equiv
	\f1N\log
	\E\bigg(\SAT_K(\bG)
	\exp\bigg\{N^{1/2}
	\tau^\st\bc[t-1]\bG\bw
	\bigg\}
	\,\bigg|\,\ROW, (\Gprime)_\ROW,
	\bG\bv=\bze\bigg)
	\]
for $\tau\in\R^{t-1}$. Next, with $L$ as in \eqref{e:L}
and with $\tbX_K$ as defined by Lemma~\ref{l:G.R.J},
 define
	\beq\label{e:pair.tcL}\tcL_{K|J}(\tau\,|\,\bze)
	\equiv
	\f1N\bigg(\oneN,
	L_{\|\pi(K)\|^2(1-\lm^2)+\lm^2}
	\bigg( \tbX_K
	+ c(\pi(K)) 
	\Big[\lm\bze + (1-\lm^2)^{1/2}
		N^{1/2} \bc[t-1]^\st\tau
	 \Big]
	\bigg)\bigg)\,,
	\eeq
where $c(\pi(K))\equiv (1-\|\pi(K)\|^2)^{1/2}$.
Then the function $\tcK_{K|J}$ satisfies
	\[\tcK_{K|J}(\tau\,|\,\bze)
	= \f{\|\tau\|^2}{2}
	+ \tcL_{K|J}(\tau\,|\,\bze)\,.
	\]

\begin{proof}
Conditional on the event $\ROW$, it follows from Lemma~\ref{l:G.R.J} that $\bG K'/N^{1/2} = \tbX_K\equiv\tbX$.
We also have
	\beq\label{e:decomp.G.Kpp}
	\f{\bG K''}{N^{1/2}}
	= \f{\|K''\|}{N^{1/2}}
	\bG \bv_K
	= c(\pi(K))
	\Big(\lambda\bze +(1-\lm^2)^{1/2}\bxi
		\Big)\,,\eeq
where $\bxi=\bG\bw$ is distributed as an independent gaussian vector in $\R^N$. Thus
	\[\tcK_{K|J}(\tau\,|\,\bze)
	= \f1N\sum_{a\le M}
	\log \E_\xi\bigg[
	\exp\bigg\{
	N^{1/2}\sum_{\ell\le t-1}\tau_\ell
		(\bc^{(\ell)})_a \xi
	\bigg\}
	U\bigg(\tbX_a + c(\pi(K))
	\Big\{
		\lambda\bze_a
		+(1-\lm^2)^{1/2} \xi\Big\}
	\bigg)
	\bigg]\,,
	\]
where $\xi$ denotes a standard gaussian random variable. Making a change of variable gives
	\[\tcK_{K|J}(\tau\,|\,\bze)
	= \f{\|\tau\|^2}{2}
	+\f1N\sum_{a\le M}
	\log \E_\xi U
	\bigg(
	\tbX_a
	+ c(\pi(K))
	\bigg\{\lambda\bze_a
	+(1-\lm^2)^{1/2} 
	\bigg[\xi+ N^{1/2}\sum_{\ell\le t-1}
		\tau_\ell (\bc^{(\ell)})_a
	\bigg]
	\bigg\}
	\bigg)\,,
	\]
from which the result follows.
\end{proof}
\end{lem}

\begin{proof}[\hypertarget{proof:p.ubd.pair.lambda}{Proof of Proposition~\ref{p:ubd.pair.lambda}}]
We follow a very similar outline as in the 
\hyperlink{proof:p.first.mmt}{proof of Proposition~\ref{p:first.mmt}}. 
As in \eqref{e:large.devs.Z} above, let us write $L=5C_1(U)^2$.
Given $\Fprime(t)$ as in \eqref{e:tap.condition.prime} and  $J,K\in\set{-1,+1}^N$,
we abbreviate the quantity of interest as	
	\beq\label{e:def.E.J.K}
	E_{J,K}\equiv
	\E\bigg(\SAT_J(\Gprime)\SAT_K(\Gprime)
	\I\bigg\{
	\f{\|\Gprime\bv\|^2}{M}\le L
	\bigg\}
	\,\bigg|\,\Fprime(t)\bigg)\,.
	\eeq
It follows by the obvious generalization of Lemma~\ref{l:resampling} that
	\[
	E_{J,K}
	=\E\bigg(\SAT_J(\bG)\SAT_K(\bG)
	\I\bigg\{ \f{\|\bG\bv\|^2}{M}
	\le L\bigg\}
	\,\bigg|\, \ROW,\COL,(\Gprime)_{\ROW\COL}
	\bigg)\,,\]
where $\bG$ is an independent copy of $\Gprime$. Next, the obvious generalization of Lemma~\ref{l:adm} gives the simplification
	\[
	E_{J,K}
	=\E\bigg(\SAT_J(\bG)\SAT_K(\bG)
	\I\bigg\{ \f{\|\bG\bv\|^2}{M}\le L
		\bigg\}
	\,\bigg|\, \ROW,\ADM,\ADM(K)\bigg)\,,
	\]
where $\ADM$  and $\ADM(K)$ are as in Definition~\ref{d:ADM.event} and 
Definition~\ref{d:pair.ADM.event} respectively. By the law of iterated expectations, 
	\beq\label{e:pair.iterated.expectations}
	E_{J,K}
	= \E\bigg(\SAT_J(\bG)
	\I\bigg\{ \f{\|\bG\bv\|^2}{M}
	\le L\bigg\}
	\E\Big[\SAT_K(\bG)
	\,\Big|\,\ROW,\bG\bv,\ADM(K)\Big]
	\,\bigg|\, \ROW,\ADM\bigg)\,.
	\eeq
We therefore first consider the calculation of
	\beq\label{e:defn.E.K.zeta}
	E_{K|J}(\bze)
	\equiv 
	\E\bigg(\SAT_K(\bG)
	\,\bigg|\,\ROW,\bG\bv=\bze,\ADM(K)\bigg)
	\eeq
(where we assume that $\bze$ satisfies the constraints imposed by $\ADM$).

Towards the calculation of \eqref{e:defn.E.K.zeta},
recall the notation of Definition~\ref{d:pair.prf.adm.subspaces}, and let $V_{\PRF(K)\setminus\ADM(K)}$ be the orthogonal complement of $V_{\ADM(K)}$ inside $V_{\PRF(K)}$. 
Analogously to \eqref{e:g.prf}
and \eqref{e:g.adm}, define
$\bg_{\PRF(K)}$ and $\bg_{\ADM(K)}$, for instance
	\beq\label{e:g.adm.K}\bg_{\ADM(K)}
	\equiv 
	\bigg(
	(\bG,\bc^{(\ell)}\bw^\st)
	: 1\le \ell\le t-1\bigg)
	=\bc[t-1]\bG\bw \in\R^{t-1}\,.
	\eeq
Choose an orthonormal basis for $V_{\PRF(K)\setminus\ADM(K)}$, and denote it $\bm{B}_j(K)$ for $1\le j\le M-(t-1)$. Analogously to \eqref{e:g.prf.minus.adm}, let
	\[
	\bg_{\BBB(K)}
	\equiv
	\bigg(
	(\bG,\bB_j(K))
	: 1\le j\le M-(t-1)\bigg)
	\in\R^{M-t+1}\,.
	\]
Note that there is an orthogonal transformation of $\R^M$ which maps $\bg_{\PRF(K)}$ to the pair $(\bg_{\ADM(K)},\bg_{\BBB(K)})$.

The weight $\SAT_K(\bG)$, as defined by \eqref{e:weight}, is a function of $\bG K$, 
which we decomposed in the proof of Lemma~\ref{l:pair.adm.cgf} as a sum of $\bG K'$ and $\bG K''$. Recall that $\bG K'$ is a function of $\bg_\ROW$. Meanwhile
  (see e.g.\ \eqref{e:decomp.G.Kpp}) $\bG K''$ is a linear combination of $\bg_\PRF=\bG\bv=\bze$ and $\bg_{\PRF(K)}=\bG\bw$, where $\bg_{\PRF(K)}$ is equivalent to the pair $(\bg_{\ADM(K)},\bg_{\BBB(K)})$ as noted above. Thus $\SAT_K(\bG)$ can be rewritten as a function $\bm{S}_{K|J}$ of $(\bg_\ROW,\bg_\PRF,\bg_{\ADM(K)},\bg_{\BBB(K)} )$: explicitly,
	\[
	\SAT_K(\bG)
	=\prod_{a\le M}
	U\bigg( 
	\sum_{s\le t} 
	\f{(K,\br^{(s)}) }{N^{1/2}} (\bg_\ROW)_{a,s}
	+ \f{\|K''\|}{N^{1/2}}
	\Big( \lambda (\bg_\PRF)_a
	+(1-\lm^2)^{1/2}
	(\bg_{\PRF(K)})_a \Big)
	\bigg)
	\equiv\bm{S}_{K|J}
		(\bg_\ROW,\bg_\PRF,
		\bg_{\ADM(K)},\bg_{\BBB(K)} )\,,
	\]
with $\lm\equiv\lm(J,K)$ as given by  Definition~\ref{d:lambda}. On the event $\ADM(K)$, the value of $\bg_{\ADM(K)}$ is fixed to a value $\bar{g}_{\ADM(K)}$. We then introduce a parameter $\tau\in\R^{t-1}$, and define (analogously to \eqref{e:S.tau.fn})
	\[
	\SAT_{K|J,\tau}(\bG)
	\equiv
	\bm{S}_{K|J,\tau}(\bg_\ROW,\bg_\PRF,
		\bg_{\ADM(K)},\bg_{\BBB(K)} )
	\equiv 
	\bm{S}_{K|J}(\bg_\ROW,\bg_\PRF,
		\bg_{\ADM(K)},\bg_{\BBB(K)} )
	\exp\bigg\{
	N^{1/2}(\tau,\bg_{\ADM(K)})
	\bigg\}\,.
	\]
Then, analogously to \eqref{e:E.J.integral}, for any $\tau\in\R^{t-1}$ we can rewrite \eqref{e:defn.E.K.zeta} as
	\begin{align}\nonumber
	E_{K|J}(\bze)
	&= \E\bigg( \f{ \bm{S}_{K|J,\tau}(\bg_\ROW, \bg_\PRF,
		\bg_{\ADM(K)},\bg_{\BBB(K)} )}
			{\exp\{N^{1/2}
				(\tau,\bar{g}_{\ADM(K)})\}}
		\,\bigg|\,
			(\bg_\ROW,\bg_\PRF,\bg_{\ADM(K)})=
				(\bar{g}_\ROW,\bze,\bar{g}_{\ADM(K)})
			\bigg) \\
	&=\f1{\exp\{N^{1/2}(\tau,\bar{g}_{\ADM(K)})\}}
	\int 
	\bm{S}_{K|J,\tau}(\bar{g}_\ROW, \bze,
		\bar{g}_{\ADM(K)}, g_{\BBB(K)} )
			p_{\BBB(K)}(g_{\BBB(K)})
			\,dg_{\BBB(K)}\,.
	\label{e:E.K.zeta.integral}
	\end{align}
By contrast, the expected value of $\SAT_{K|J,\tau}$ given only the row constraints is (cf.\ \eqref{e:def.E.tau})
	\begin{align}\nonumber
	\bE_{K|J}(\tau\,|\,\bar{g}_\ROW,\bze)
	&\equiv
	 \E\bigg( \SAT_{K|J,\tau}(\bG)\,\bigg|\,
	\ROW,\bG\bv=\bze\bigg)
	=\E\bigg( \bm{S}_{K|J,\tau}(\bg_\ROW, \bg_\PRF,
		\bg_{\ADM(K)},\bg_{\BBB(K)} )
		\,\bigg|\,
			(\bg_\ROW,\bg_\PRF)=
				(\bar{g}_\ROW,\bze)
			\bigg)  \\ \nonumber
	&= \int p_{\ADM(K)}(g_{\ADM(K)})
		\int 
	\bm{S}_{K|J,\tau}(\bar{g}_\ROW, \bze,
		g_{\ADM(K)}, g_{\BBB(K)} )
			p_{\BBB(K)}(g_{\BBB(K)})
			\,dg_{\BBB(K)}
			\,dg_{\ADM(K)}\\
	&=\exp\bigg\{ N\tcK_{K|J}(\tau\,|\,\bze)\bigg\}\,.
	\label{e:def.E.tau.pair}
	\end{align}
Then, analogously to \eqref{e:ap.tau}, we define the probability density function
	\begin{align}\nonumber
	\ap_{K|J,\tau}(g_{\ADM(K)}\,|\,\bar{g}_\ROW,\bze)
		\,dg_{\ADM(K)}
	&\equiv
	\f{\E(\SAT_{K|J,\tau}(\bG)
		\Ind{
		\bg_{\ADM(K)} \in dg_{\ADM(K)}
		} \,|\,\ROW,\bG\bv=\bze)}
		{\E(\SAT_{K|J,\tau}(\bG)\,|\,\ROW,\bG\bv=\bze)}
	\\
	&= \f{p_{\ADM(K)}(g_{\ADM(K)})}
		{\bE_{K|J}(\tau\,|\,\bar{g}_\ROW,\bze) }
	\int \bm{S}_{K|J,\tau}(
		\bar{g}_\ROW,\bze,g_{\ADM(K)},
			g_{\BBB(K)})
		p_{\BBB(K)}(g_{\BBB(K)})
	\,dg_{\BBB(K)}
	\label{e:ap.tau.two}
	\end{align}
Then it follows similarly to \eqref{e:tilt.gaussian} that we can rewrite \eqref{e:E.K.zeta.integral} as
	\beq\label{e:pair.tilt.gaussian}
	E_{K|J}(\bze)
	= \f{\bE_{K|J}(\tau\,|\,\bar{g}_\ROW,\bze) 
		 \cdot \ap_{K|J,\tau}(\bar{g}_{\ADM(K)}
		\,|\,\bar{g}_\ROW,\bze)}
		{
	{\exp\{N^{1/2}(\tau,\bar{g}_{\ADM(K)})\}}
		\cdot
		p_{\ADM(K)}(\bar{g}_{\ADM(K)})}\,.
	\eeq
We will show in 
Proposition~\ref{p:pair.density.bound}
(deferred to Section~\ref{s:local.clt})
 that (cf.\ \eqref{e:density.bound})
	\beq\label{e:pair.density.bound}
	\max \bigg\{ 
		\Big\| \ap_{K|J,\tau}
		(\cdot\,|\,\bar{g}_\ROW)
		\Big\|_\infty
	: J\in\set{-1,+1}^N, \|\pi(J)\| \le \f45,
	|\lm|\le\f45,
	\|\tau\| \le \tau_{\max}
	\bigg\}
	 \le\wp_{t,2}
	\,.\eeq
It therefore remains to estimate the other two terms on the right-hand side of \eqref{e:pair.tilt.gaussian}.
We then note 
that Definition~\ref{d:pair.ADM.event} implies that,
on the event $\ADM(K)$, we have
(cf.\ \eqref{e:adm.orth.SECOND})
	\begin{align}\nonumber
	\f{\bar{g}_{\ADM(K)}}{N^{1/2}}
	&\stackrel{\eqref{e:g.adm.K}}{=} \f{\bc[t-1]\bG\bw}{N^{1/2}}
	\stackrel{\eqref{e:gs.n.c.EXACT}}{=}
	\f{(\bGam_N)^{-1}\bn[t-1]\bG\bw}{N \psi^{1/2}}
	\stackrel{\eqref{e:pair.adm}}{=}
	\f{(\bGam_N)^{-1}\bH[t-1]\bw}{(N \psi)^{1/2}}\\
	&\stackrel{\eqref{e:K.w}}{=} 
	\f{(\bGam_N)^{-1}\bH[t-1]}{(N \psi)^{1/2}}
	\bigg( \f{\bv_K-\lm\bv}{(1-\lm^2)^{1/2}} \bigg)
	\stackrel{\eqref{e:defn.delta}}{=}
	\f{(\bGam_N)^\st [\delta(K)-\lm\delta(J)]}
		{(1-\lm^2)^{1/2}}
	= o_N(1)\,,
	\label{e:pair.adm.orth}
	\end{align}
where the last estimate holds thanks to the restriction $J,K\in\mbH_*$ (see \eqref{e:delta.star.is.zero}).
Substituting \eqref{e:pair.adm.orth} into the formula for $p_{\ADM(K)}$ (similar to \eqref{e:p.adm.uncond}) gives
	\beq\label{e:pair.p.adm.uncond}
	p_{\ADM(K)}(\bar{g}_{\ADM(K)})
	= \f1{(2\pi)^{(t-1)/2}} \exp\bigg\{
	-\f{N}{2}\bigg\|
	\f{(\bGam_N)^\st[\delta(K)-\lm\delta(J)]}
		{(1-\lm^2)^{1/2}}
		\bigg\|^2
		\bigg\}
	= \exp\{ N \cdot o_N(1)\}\,.
	\eeq
Meanwhile, it follows by combining
\eqref{e:def.E.tau.pair} and \eqref{e:pair.adm.orth} that (cf.\ \eqref{e:remaining.E.factor})
	\beq\label{e:pair.remaining.E.factor}
	E_{K|J}(\bze)
	\le \exp\bigg\{
	N \bigg[\tcK_{K|J}(\tau\,|\,\bze) - 
	\bigg(\tau,\f{(\bGam_N)^\st [\delta(K)-\lm\delta(J)]}
		{(1-\lm^2)^{1/2}}
\bigg) + o_N(1) \bigg]
	\bigg\}
	\eeq
Substituting \eqref{e:pair.density.bound},
\eqref{e:pair.p.adm.uncond}, and
\eqref{e:pair.remaining.E.factor}
into \eqref{e:pair.tilt.gaussian},
and combining with Lemma~\ref{l:pair.adm.cgf}, gives
(cf.\ \eqref{e:tcA})
	\beq\label{e:pair.tcA}
	\f{E_{K|J}(\bze)}{ \wp_{t,2} (2\pi)^{t/2}}
	\le
	\exp\bigg\{
	N\bigg[
	 \f12\bigg\|
	\tau-\f{(\bGam_N)^\st[\delta(K)-\lm\delta(J)]}
		{(1-\lm^2)^{1/2}}
	\bigg\|^2 + \tcL_{K|J}(\tau\,|\,\bze)
	\bigg]
	\bigg\}
	\equiv \exp\Big\{N
		\tcA_{K|J}(\tau\,|\,\bze)\Big\}\,,
	\eeq
where $\tcA_{K|J}$ is defined by the last identity.
To simplify the above expression, we set 
$\tau = \bar{\tau}(\lm)$ where 
	\[
	\bar{\tau}(\lm)
	\equiv
	\f{\starvpi}{(1-q)^{1/2}(1-\lm^2)^{1/2}}
	\stackrel{\eqref{e:vpi.star}}{=} - \f{\psi^{1/2}(1-q)^{1/2}}
		{(1-\lm^2)^{1/2}}
		\bGam^\st\acute{e}_{t-1}\,.
	\]
Substituting this into \eqref{e:pair.tcL},
and recalling the definition of $\tbX_K$ from Lemma~\ref{l:G.R.J}, we obtain
	\beq\label{e:pair.cL.def}
	\tcL_{K|J}(\bar{\tau}(\lm)\,|\,\bze)
	\simeq
	\f1N\bigg(\oneN, L_{
	q(1-\lm^2)+\lm^2}\Big( 
		\bh^{(t+1)}
		+ (1-q)^{1/2} \lm \bze \Big)\bigg)
	 + o_N(1)
	\equiv \cL_2(\lm\,|\,\bze)\,,
	\eeq
where $\cL_2$ is defined by the last identity.
 By substituting the above into
\eqref{e:pair.tcA}, we see that the quantity from 
\eqref{e:defn.E.K.zeta} can be upper bounded by
	\[
	E_{K|J}(\bze) \le
	 \exp\bigg\{ N\bigg[ \f{\psi(1-q)}{2(1-\lm^2)}
		+\cL_2(\lm\,|\,\bze) + o_N(1)
		\bigg]\bigg\}
	=
	\exp\bigg\{N\Big[
	\cA_2(\lm\,|,\bze) + o_N(1)\Big]
	\bigg\}
	\,,\]
for $\cA_2(\lm\,|,\bze)$ as in the statement of the proposition. By comparing 
\eqref{e:pair.iterated.expectations} with 
\eqref{e:defn.E.K.zeta}, we see that
	\[
	E_{J,K}
	=\E(\SAT_J\,|\,\ROW,\ADM)
	\int
	\I\bigg\{\f{\|\bze\|^2}{M} \le L\bigg\}
	 E_{K|J}(\bze)\,\Q_J(d\bze)\,,
	\]
so the claim follows.
\end{proof}

Analogously to \eqref{e:H.circ}, we now define
	\[
	(\mbH_*)^{2,\circ}
	\equiv \bigg\{
	(J,K)\in(\mbH_*)^2
	: \f{|\lm(J,K)|}{\alpha^{1/2}}
	\le 10\cdot C_1(U)
	\bigg\}\,,
	\]
so $(\mbH_*)^{2,\circ}$ is a subset of $(\mbH_*)^2$.
Then decompose
$\bar{\bZ}^2(\Gprime)\equiv
\bar{\bZ}^{2,\circ}(\Gprime)+\bar{\bZ}^{2,\bullet}(\Gprime)$
where (cf.\ \eqref{e:Z.outside})
	\beq\label{e:pair.small.lm}
	\bar{\bZ}^{2,\circ}(\Gprime)
	\equiv
	\sum_{(J,K)\in(\mbH_*)^{2,\circ}}
	\SAT_J(\Gprime)
	\SAT_K(\Gprime)
	\I\bigg\{
	\f{\|\Gprime\bv_J\|^2}{M} \le L,
	\f{\|\Gprime\bv_K\|^2}{M} \le L
	\bigg\}\,.
	\eeq
We bound $\bar{\bZ}^{2,\circ}(\Gprime)$ as follows:

\begin{thm}[analogous to Theorem~\ref{t:ubd.slice}]
\label{t:slice.second.mmt}
Suppose $U$ satisfies Assumptions~\ref{a:bdd} and \ref{a:Lip}, and let $\Fprime(t)$ be as in \eqref{e:tap.condition.prime}.
Recalling Proposition~\ref{p:ubd.pair.lambda}
and \eqref{e:pair.cL.def}, let
$\Psi_2(\lm\,|\,\bze)$ be defined by
	\[\Psi_2(\lm\,|\,\bze)
	-\Psi(\starpi,\starvpi) 
	\equiv
	-\psi(1-q)
		+ \cA_2(\lm\,|\,\bze)
	=-\psi(1-q)
	+ \f{\psi(1-q)}{2(1-\lm^2)}
		+\cL_2(\lm\,|\,\bze)
	\,.
	\]
For $\bar{\bZ}^{2,\circ}(\Gprime)$ as defined by \eqref{e:pair.small.lm}, we have
	\[
	\f{\E(\bar{\bZ}^{2,\circ}(\Gprime)
	\,|\,\Fprime(t))}
		{\exp\{2 \cdot(\oneN,\log(2 \Ch(\bH^{(t)})))\}}
	\le
	\sum_{(J,K)\in (\mbH_*)^{2,\circ}}
	\bQ(J,K)
	\int
	\exp\{N\Psi_2(\lm\,|\,\bze)\}
	\,\Q_J(d\bze)
	\]
for $\Q_J$ as in Definition~\ref{d:measure.Q.J}
and $\bQ$ as in Definition~\ref{d:pair.P.Q}.

\begin{proof}
We follow the
\hyperlink{proof:t.ubd.slice}{proof of 
Theorem~\ref{t:ubd.slice}}.
Suppose $J,K\in\mbH_*$ with $\lm=\lm(J,K)$
as given by Definition~\ref{d:lambda}.
Recalling Definition~\ref{d:pi.vpi},
the restriction $J\in\mbH_*$ implies
	\[
	\f{\bH[t-1]J}{N\psi^{1/2}}
	\stackrel{\eqref{e:gs.H.y}}{=}
	\f{\bGam\by[t-1]J}{N}
	\stackrel{\eqref{e:def.vpi}}{=}
	\bGam\vpi(J)
	\simeq
	\bGam\starvpi
	\stackrel{\eqref{e:def.vpi.star}}{=}
	\psi^{1/2}(1-q)\bGam\bGam^\st\acute{e}_{t-1}\,.
	\]
It follows that, for all $J\in\mbH_*$,
	\[
	\f{(\bH^{(t)},J)}{N}
	\simeq
	\psi(1-q)(\bGam\bGam^\st)_{t-1,t-1}
	\stackrel{\eqref{e:Gamma.matrix}}{=}
	\psi(1-q)\,.
	\]
Since $(\bH^{(t)},\bH^{(t+1)})/(N\psi) = 1-o_t(1)$, we conclude that, for all $J\in\mbH_*$.
	\[
	\f{(\bH^{(t+1)},J)}{N}
	\simeq \psi(1-q) - o_t(1)\,.
	\]
Let $E_{J,K}$ be as in \eqref{e:def.E.J.K}. Combining with
Definition~\ref{d:pair.P.Q} gives
	\begin{align}\nonumber
	&\f{\E(\bar{\bZ}^{2,\circ}(\Gprime)
	\,|\,\Fprime(t))}
		{\exp\{2 \cdot(\oneN,\log(2 \Ch(\bH^{(t)})))\}}
	\le\sum_{(J,K)\in (\mbH_*)^{2,\circ}}
		\bQ(J,K)\bigg(
	\f{\bP(J,K)/\bQ(J,K)}{
		\exp\{2 \cdot(\oneN,\log\Ch(\bH^{(t)}))\}}
		\bigg) E_{J,K} \\ \nonumber
	&\qquad\le\sum_{(J,K)\in (\mbH_*)^{2,\circ}}
	\bQ(J,K)
	\f{ E_J \cdot E_{J,K}/E_J}
	{\exp\{2 N[\psi(1-q) - o_t(1)]\}}
\\
	&\qquad\le
	\exp\Big\{N\Psi(\pi_*,\vpi_*)\Big\}
	\sum_{(J,K)\in (\mbH_*)^{2,\circ}}
	\bQ(J,K)
	\f{E_{J,K}/E_J}{\exp\{N[\psi(1-q) - o_t(1)]\}}
	\label{e:barZ.sq.circ.ubd}\,,
	\end{align}
where the last bound is by 
the calculation \eqref{e:pair.lbd.rewrite}
from the proof of Proposition~\ref{p:restricted.first.mmt}.
Combining with Proposition~\ref{p:ubd.pair.lambda} gives
the claim.
\end{proof}
\end{thm}

\begin{cor}[analogous to Corollary~\ref{c:AH.bound}]
\label{c:small.lambda}
Suppose $U$ satisfies Assumptions~\ref{a:bdd} and \ref{a:Lip}, and let $\Fprime(t)$ be as in \eqref{e:tap.condition.prime}. We then have the bound
	\[
	\E\Big(\bar{\bZ}^{2,\bullet}(\Gprime)
		\,\Big|\,\Fprime(t)\Big)
	\le \exp\bigg\{
	2N\bigg(\RS(\alpha;U) - 0.1\cdot
		C_1(U)^2\alpha\bigg)
	\bigg\}
	\]
for $\bar{\bZ}^{2,\bullet}(\Gprime)
=\bar{\bZ}^2(\Gprime)-\bar{\bZ}^{2,\circ}(\Gprime)$ as defined by \eqref{e:pair.small.lm}.

\begin{proof}
For $E_{J,K}$ as in \eqref{e:def.E.J.K}
we also have trivially
$E_{J,K}\le1$,
and combining this with the calculation 
\eqref{e:barZ.sq.circ.ubd}
gives
	\[
	\f{\E(\bar{\bZ}^{2,\bullet}(\Gprime)
	\,|\,\Fprime(t))}
		{\exp\{2 \cdot(\oneN,\log(2 \Ch(\bH^{(t)})))\}}
	\le
	\sum_{(J,K)\in(\mbH_*)^{2,\bullet}}
	\f{\bQ(J,K) }
	{\exp\{ 2N
	[\psi(1-q) - o_t(1)]\}}\,.
	\]
Combining Proposition~\ref{p:fp} with
Corollary~\ref{c:rs.estimate} and
 \eqref{e:log.cosh.rough.bound} gives, with high probability,
	\begin{align*}
	\f{\E(\bar{\bZ}^{2,\bullet}(\Gprime)
	\,|\,\Fprime(t))}
		{\exp\{2 N\RS(\alpha;U)\}}
	&\le\bQ\Big((\mbH_*)^{2,\bullet}\Big) 
	\cdot \exp\bigg\{
	N\cdot 2\bigg[
	3
	+ 1.53
	+ \f{3}{2} + o_t(1) \bigg]
	 C_1(U)^2\cdot \alpha
	\bigg\} \\
	&\le
	\bQ\Big((\mbH_*)^{2,\bullet}\Big)
	\cdot
	\exp\bigg\{ 12.1 \cdot N
		\cdot C_1(U)^2 \alpha \bigg\}\,.
	\end{align*}
For any $J\in\set{-1,+1}^N$, it follows by the Azuma--Hoeffding inequality that
	\[
	\bQ\bigg(
	K\in\set{-1,+1}^N :
	\bigg|\f{(J-\bmag^{(t)},K-\bmag^{(t)})}{N}\bigg| \ge x
	\bigg)
	\le
	2\exp\bigg\{
	-\f{Nx^2}{8}
	\bigg\}
	\]
for any $x\ge0$.
Recalling Definition~\ref{d:lambda}, it follows that
for any $J\in\mbH_*$,
	\beq\label{e:AH.lambda}
	\bQ\bigg(
	K\in\mbH_*
	: |\lm(J,K)| \ge l
	\bigg)
	\le 2
	\exp\bigg\{
	-\f{N(1-q + o_N(1)) l^2}{8}
	\bigg\}
	\eeq
for any $l\ge0$. Taking $l=10\cdot C_1(U)\alpha^{1/2}$
and summing over $J$ gives
	\[\bQ\Big((\mbH_*)^{2,\bullet}\Big)
	\le
	\sum_{J\in\mbH_*}
	\bQ\bigg(
	K\in\mbH_*
	: |\lm(J,K)| \ge 10\cdot C_1(U) \alpha^{1/2}
	\bigg)
	\le \exp\bigg\{
	-12.4\cdot N \cdot  C_1(U)^2 \alpha
	\bigg\}\,.
	\] 
It follows by combining the above bounds that
	\[
	\f{\E(\bar{\bZ}^{2,\bullet}(\Gprime)
	\,|\,\Fprime(t))}
		{\exp\{2 N\RS(\alpha;U)\}}
	\le 
	\exp\bigg\{
	-0.3\cdot N \cdot C_1(U)^2\alpha
	\bigg\}\,,
	\]
which concludes the proof.
\end{proof}
\end{cor}

\subsection{Analysis of second moment}
\label{ss:second.mmt.analysis}

In this subsection we analyze the bound from Theorem~\ref{t:slice.second.mmt}.

\begin{ppn}[analogous to Lemmas~\ref{l:dPsi.dpi} and \ref{l:dPsi.dvpi}]\label{p:stationary.lambda}
For $\cA_2(\lm\,|\,\bze)$ as defined by Proposition~\ref{p:ubd.pair.lambda}, we have
	\[
	\f{d\cA_2(\lm\,|\,\bze)}{d\lm}
	\bigg|_{\lm=0}
	=o_N(1) + o_t(1) \alpha L^{1/2}
	\]
provided that $\bze=\bG\bv$ is compatible with 
the admissibility condition \eqref{e:adm},
and satisfies the bound $\|\bze\|^2 \le ML$ where $L=L_1(U)=5C_1(U)^2$ as defined above.

\begin{proof}
For $\cL_2(\lm\,|\,\bze)$ as defined by \eqref{e:pair.cL.def}, we have
	\begin{align*}
	\f{d\cA_2(\lm\,|\,\bze)}{d\lm}
	\bigg|_{\lm=0}
	&=
	\f{d\cL_2(\lm\,|\,\bze)}{d\lm}
	\bigg|_{\lm=0}
	= \f{(1-q)^{1/2}}{N}
	\sum_{a\le M}
	\f{\E_\xi U'( (\bh^{(t+1)})_a + (1-q)^{1/2} \xi)}
		{\E U((\bh^{(t+1)})_a
			 + (1-q)^{1/2} \xi)}
	\bze_a \\
	&\stackrel{\eqref{e:F.rewrite}}{=}
	\f{(1-q)^{1/2}
	(F_q(\bh^{(t+1)}),\bze)}{N}
	= \f{(1-q)^{1/2} (\bn^{(t+1)},\bze)}{N}\,.
	\end{align*}
On the other hand, it follows from the admissibility condition that $\bze=\bG\bv$ must satisfy
	\beq\label{e:adm.rewrite.bze}
	\f{\bc[t-1]\bze}{N^{1/2}}
	\stackrel{\eqref{e:adm.orth}}{=}
	(\bGam_N)^\st \delta(J)
	= o_N(1)\,,
	\eeq
where the last step is by the restriction $J\in\mbH_*$.
The span of the vectors $\bc^{(\ell)}$ for $\ell\le t-1$
--- which is the same as the span of the vectors
$\bn^{(\ell)}$ for $\ell\le t-1$
--- does not contain $\bn^{(t+1)}$, but recall from \eqref{e:n.scalar.products} that
	\[
	\f{(\bn^{(t+1)},\bn^{(t-1)})}{N\psi}
	\stackrel{\eqref{e:n.scalar.products}}{=}
	\mu_{t-1}\,.
	\]
It follows from Proposition~\ref{p:AT}
and Lemma~\ref{l:AT.conclusion}
that $\mu_{t-1}=1-o_t(1)$, so we can decompose
	\[
	\f{\bn^{(t+1)}}{(N\psi)^{1/2}}
	= \bc^\parallel+\bc^\perp
	\]
where $\bc^\parallel$ lies in the span of the vectors $\bc^{(\ell)}$ and has norm $1-o_t(1)$, while $\bc^\perp$ is orthogonal to the vectors $\bc^{(\ell)}$ and has norm $o_t(1)$. It follows that
	\[
	\bigg|
	\bigg(\f{\bn^{(t+1)}}{N \psi^{1/2}},\bze\bigg)
	\bigg|
	= \f{|(\bc^\parallel + \bc^\perp,\bze)|}{N^{1/2}}
	\stackrel{\eqref{e:adm.rewrite.bze}}{\le}
	o_N(1) 
	+\f{\|\bc^\perp\| \|\bze\|}{N^{1/2}}
	\le
	o_N(1) + o_t(1) (\alpha L)^{1/2}\,,
	\]
having used Cauchy--Schwarz together with 
the assumption $\|\bze\|^2\le ML$. In conclusion we find
	\[\f{d\cA_2(\lm\,|\,\bze)}{d\lm}
	\bigg|_{\lm=0} 
	=
	o_N(1) + o_t(1) \alpha L^{1/2}\,,
	\]
as claimed.
\end{proof}
\end{ppn}

\begin{lem}[analogous to  Corollary~\ref{c:pair.hess}]
\label{l:hess.second.mmt}
Suppose $U$ satisfies Assumptions \ref{a:bdd} and \ref{a:Lip}. Recall $K_2(U)$ from Assumption~\ref{a:Lip}, and $C_1(U)$ from Lemma~\ref{l:poly}.
For $\cL_2(\lm\,|\,\bze)$ as defined by \eqref{e:pair.cL.def}, we have  the bound 
		\[\bigg|
	\f{d^2\cL_2(\lm\,|\,\bze)}{d\lm^2}\bigg|
	\le 420 \cdot 
	C_1(U)^2 K_2(U) \cdot \alpha\,
	\]
as long as $|\lm|\le 4/5$
and $\|\bze\|^2/M\le L= L_1(U) = 5C_1(U)^2$.

\begin{proof}
Let $e(\lm)\equiv(1-q)^{1/2}(1-\lm^2)^{1/2}$. Denote
	\[
	\bY\equiv\bY(\lm;\bze)
	=\bigg\{ \bh^{(t+1)} + (1-q)^{1/2}\lm\bze\bigg\}
	+(1-q)^{1/2}(1-\lm^2)^{1/2} \xi\oneM
	\equiv\bX(\lm;\bze)
		+e(\lm) \xi\oneM\,.
	\]
Then the function $\cL_2(\lm\,|\,\bze)$ from Proposition~\ref{p:ubd.pair.lambda} can be rewritten as
	\[
	\cL_2(\lm\,|\,\bze)
	\equiv \f1N\sum_{a\le M}
	\log\E_\xi U(\bY_a(\lm,\bze))\,.
	\]
Recalling the notation of \eqref{e:diag.A}, \eqref{e:offdiag.b}, \eqref{e:bar.a}, and \eqref{e:bar.b}, let us now define
$\mathbf{A}\equiv A_c(\bX)$,
$\mathbf{B}\equiv B_c(\bX)$,
$\bar{a}\equiv(\oneM,a_c(\bX))$,
and $\bar{b}\equiv(\oneM,b_c(\bX))$,
for $c=e(\lm)$ and $\bX=\bX(\lm;\bze)$. With this notation, we have
	\[
	\f{d^2 \cL_2(\lm\,|\,\bze)}{d\lm^2}
	= \f1N\sum_{a\le M}
	\bigg\{
	\f{ \E_\xi U'(\bY_a) \f{d^2\bY_a}{d\lm^2}
	+\E_\xi U''(\bY_a)
		(\f{d\bY_a}{d\lm})^2}
	{\E_\xi U(\bY_a)}
	- \bigg(\f{\E_\xi U'(\bY_a)
		\f{d\bY_a}{d\lm}}
	{\E_\xi U(\bY_a)}\bigg)^2
	\bigg\}\,.\]
We decompose the above as $(\textup{I})+ (\textup{II})+ (\textup{III})+(\textup{IV})$ where (cf.\ Lemma~\ref{l:compute.Hess})
	\begin{align*}
	(\textup{I})
	&\equiv
	\f1N (\bX'(\lm))^\st (\diag\mathbf{A}) \bX'(\lm)
	= \f{1-q}{N} \bze^\st (\diag\mathbf{A}) \bze\,,\\
	(\textup{II})
	&\equiv
	\f{2}{N}e'(\lm)
	\mathbf{B}^\st \bX'(\lm)
	=-\f{2(1-q)^{1/2}\lm}{N(1-\lm^2)^{1/2}}
	\mathbf{B}^\st \bze\,,\\
	(\textup{III})
	&\equiv \f1N e''(\lm)\bar{a}
	= -\f{(1-q)^{1/2}}{N(1-\lm^2)^{3/2}} \bar{a}\,,\\
	(\textup{IV})
	&\equiv \f1N e'(\lm^2) \bar{b}
	= \f{(1-q)\lm^2}{N(1-\lm^2)}
	\bar{b}\,.
	\end{align*}
We bound each of the above terms, assuming $|\lm|\le 4/5$. Applying \eqref{e:bound.A} gives
	\[
	|(\textup{I})|
	\le
	\f1N \|\mathbf{A}\|_\infty \|\bze\|^2
	\le
	\f1N 1.7 \cdot K_2(U) \|\bze\|^2
	\le 
	1.7 \cdot K_2(U) \alpha L
	= 8.5 \cdot C_1(U)^2 K_2(U)
	\cdot \alpha
	\,.
	\]
It follows from the above definition of $\bX\equiv\bX(\lm;\bze)$ that
	\beq\label{e:pair.bound.bX}
	\f{\|\bX\|}{M^{1/2}}
	\le \f{\|\bh^{(t+1)}\| + \|\bze\|}{M^{1/2}}
	\le 2q^{1/2}+ L^{1/2}
	\le 2.5 \cdot C_1(U)
	\,,\eeq
with high probability. Combining
\eqref{e:pair.bound.bX} with Lemma~\ref{l:bound.B} gives 
	\begin{align*}
	|(\textup{II})|
	&\le
	\f{8}{3N}
	\|\mathbf{B}\|\|\bze\|
	\le
	\f{8M^{1/2}\|\bze\|}{3N}
	K_2(U) \bigg(2.5 \cdot C_1(U) 
	+ 5.8\cdot \f{\|\bX\|}{M^{1/2}}
	\bigg) \\
	&\le 
	\f{8\alpha L^{1/2}}{3} K_2(U)
	\bigg( 2.5 \cdot C_1(U) 
	+ 5.8\cdot 2.5 \cdot C_1(U)
	\bigg)
	\le 105 \cdot C_1(U)^2 K_2(U) \cdot\alpha
	\end{align*}
Next, combining
\eqref{e:pair.bound.bX} with Lemma~\ref{l:bound.a} gives 
	\begin{align*}
	|(\textup{III})|
	&\le \f{4.7}{N}|\bar{a}|
	\le \f{4.7\cdot M}{N}\bigg\{
	1.1\cdot C_1(U) + 3.7 \cdot\f{\|\bX\|^2}{M}\bigg\}\\
	&\le
	\f{4.7\cdot M}{N}\bigg\{
	1.1\cdot C_1(U) + 3.7 \cdot 2.5^2 \cdot C_1(U)^2\bigg\}
	\le 110 \cdot C_1(U)^2 \cdot \alpha\,.
	\end{align*}
Finally, 
combining
\eqref{e:pair.bound.bX} with Lemma~\ref{l:bound.b} gives 
	\begin{align*}
	|(\textup{IV})|
	&\le \f{1.8}{N}|\bar{b}|
	\le  1.8 \cdot K_2(U) \bigg( 4.6\cdot C_1(U)
		+ 17 \cdot \f{\|\bX\|^2}{M}\bigg) \cdot \alpha
	 \\
	&\le 1.8\cdot K_2(U)
	\bigg( 4.6\cdot C_1(U)
	+ 17\cdot \Big( 2.5 \cdot C_1(U)\Big)^2
	\bigg) \cdot \alpha
	\le 193 \cdot C_1(U)^2 K_2(U) \cdot \alpha\,.
	\end{align*}
Combining the above bounds gives the claim.
\end{proof}
\end{lem}

\begin{cor}[analogous to Proposition~\ref{p:hess.Psi}]
\label{c:pair.hessian}
Suppose $U$ satisfies Assumptions \ref{a:bdd} and \ref{a:Lip}. For $\Psi_2(\lm\,|\,\bze)$ as in the statement of 
Theorem~\ref{t:slice.second.mmt}, we have the bound
	\[\bigg|
	\f{d^2\Psi_2(\lm\,|\,\bze)}{d\lm^2}\bigg|
	\le 610 \cdot C_1(U)^2 K_2(U) \cdot \alpha\,,
	\]
as long as $|\lm|\le 4/5$
and $\|\bze\|^2/M\le L_1(U) = 5C_1(U)^2$.

\begin{proof}
It follows from the definition that
	\[
	\bigg|\f{d^2\Psi_2(\lm\,|\,\bze)}{d\lm^2}\bigg|
	\le \f{\psi(1-q)(1+3\lm^2)}{(1-\lm^2)^3}
	+\bigg|\f{d^2\cL_2(\lm\,|\,\bze)}{d\lm^2}\bigg|
	\le 62.6\cdot\psi
	+\bigg|\f{d^2\cL_2(\lm\,|\,\bze)}{d\lm^2}\bigg|\,.
	\]
Applying Proposition~\ref{p:fp} and Lemma~\ref{l:hess.second.mmt} gives
	\[\bigg|\f{d^2\cA_2(\lm\,|\,\bze)}{d\lm^2}\bigg|
	\le\bigg\{
	62.6 \cdot 3 \cdot C_1(U)^2
	+ 420 \cdot C_1(U)^2 K_2(U)  \bigg\} \cdot \alpha\,.
	\]
The claim follows.
\end{proof}
\end{cor}

\subsection{Conclusion of second moment}
\label{ss:second.mmt.conclusion}

In this concluding subsection we 
finish the \hyperlink{proof:t.lbd}{proof of Theorem~\ref{t:lbd}},
and use it to \hyperlink{proof:t.main.lbd.bdd}{deduce the lower bound of Theorem~\ref{t:main} in the case $\|u\|_\infty<\infty$}.

\begin{proof}[\hypertarget{proof:t.lbd}{Proof of Theorem~\ref{t:lbd} (conclusion)}]
Recall that the \hyperlink{proof:restricted.first.mmt.lbd}{proof of the first moment lower bound \eqref{e:lbd.restricted.first.mmt}}
 was already given at the end of \S\ref{ss:first.mmt.lbd}.
It therefore remains to show the second moment upper bound \eqref{e:ubd.restricted.second.mmt}, and for this 
we follow the \hyperlink{proof:t.ubd}{proof of Theorem~\ref{t:ubd}}. Recall from \eqref{e:pair.small.lm} that we decomposed
$\bar{\bZ}^2(\Gprime)$ as the sum of
$\bar{\bZ}^{2,\circ}(\Gprime)$
and $\bar{\bZ}^{2,\bullet}(\Gprime)$.
For $\bar{\bZ}^{2,\circ}(\Gprime)$, we will analyze the bound from Theorem~\ref{t:slice.second.mmt}.
We note that at $\lm=0$ we have
	\[
	\Psi_2(0\,|\,\bze)-\Psi(\starpi,\starvpi)
	= 
	- \f{\psi(1-q)}{2} + \cL_2(0\,|\,\bze)
	\stackrel{\eqref{e:pair.cL.def}}{\simeq}
	- \f{\psi(1-q)}{2} + 
	\alpha\E L_q(q^{1/2}Z)
	\stackrel{\eqref{e:Psi.star}}{\simeq}
	\Psi(\starpi,\starvpi)\,.
	\]
It follows by combining with
\eqref{e:log.cosh.limit} that
	\[
	(\oneN,\log(2 \Ch(\bH^{(t)})))
	+ \f{\Psi(\lm\,|\,\bze)}{2}
	\stackrel{N\to\infty}{\longrightarrow}
	 \RS(\alpha;U)
	\]
Next, for $|\lm|\le4/5$, it follows by combining
Proposition~\ref{p:stationary.lambda} 
and Corollary~\ref{c:pair.hessian} that
	\begin{align*}
	\Psi_2(\lm\,|\,\bze)
	-\Psi_2(0\,|\,\bze)
	&\le \f{d\Psi_2(\lm\,|\,\bze)}{d\lm}
		\bigg|_{\lm=0} \cdot \lm
		+  
		\max
	\bigg\{
	\bigg|\f{d^2\Psi_2(\lm\,|\,\bze)}{d\lm^2}\bigg|
	: |\lm|\le\f45
	\bigg\}
	\cdot\f{	\lm^2}{2} \\
	&\qquad
	\le o_t(1) \lm + 610\cdot
	C_1(U)^2 K_2(U) \cdot \alpha \cdot \f{\lm^2}{2}\,.
	\end{align*}
Recalling that $\alpha\le \alpha(U)$ as defined by 
\eqref{e:alpha.U}, the above can be simplified as
	\[\Psi_2(\lm\,|\,\bze)
	-\Psi_2(0\,|\,\bze)
	\stackrel{\eqref{e:alpha.U}}{\le}
	o_t(1)
	+ \f{610\cdot\lm^2}{2
	e^{10}  C_1(U)^4 K_2(U)^3}
	\le o_t(1)+
	\f{\lm^2}{e^{13}}\,.
	\]
Substituting this into the bound from
Theorem~\ref{t:slice.second.mmt}
 gives
 	\[
	\f{\E(\bar{\bZ}^{2,\circ}(\Gprime)
	\,|\,\Fprime(t))}
		{\exp\{2N
		[ \RS(\alpha;U) + o_t(1)
		]\}}
	\le
	\sum_{(J,K)\in (\mbH_*)^{2,\circ}}
	\bQ(J,K)
	\exp\bigg\{\f{N\lm^2}{e^{13}}\bigg\}\,.
	\]
It follows by combining with \eqref{e:AH.lambda} 
that the right-hand side is bounded by a constant.
Finally, we recall that $\E(\bar{\bZ}^{2,\bullet}(\Gprime)$ was bounded by Corollary~\ref{c:small.lambda},
so the claim follows.
\end{proof}

\begin{proof}[\hypertarget{proof:t.main.lbd.bdd}{Proof of Theorem~\ref{t:main} lower bound assuming $\|u\|_\infty<\infty$}] It follows from the first bound from Theorem~\ref{t:lbd} that
	\beq\label{e:restricted.first.mmt.lbd.half}
	\E\bigg(\bar{\bZ} \I\bigg\{
	\bar{\bZ} \ge \f{\E(\bar{\bZ}\,|\,\FF(t))}{2}
	\bigg\} \,\bigg|\,\FF(t)\bigg)
	\ge
	\f{\E(\bar{\bZ}\,|\,\FF(t))}{2}
	\stackrel{\eqref{e:lbd.restricted.first.mmt}}
		{\ge}
	\f{\exp\{ N(
		 \RS(\alpha;U) -o_t(1))\}}{2}\,,
	\eeq
with high probability over the randomness of $\FF(t)$.
On the other hand, the Cauchy--Schwarz inequality gives
	\[
	\E\bigg(\bar{\bZ} \I\bigg\{
	\bar{\bZ} \ge \f{\E(\bar{\bZ}\,|\,\FF(t))}{2}
	\bigg\} \,\bigg|\,\FF(t)\bigg)^2
	\le
	\E(\bar{\bZ}^2\,|\,\FF(t)) \cdot
	\P\bigg(
	\bar{\bZ} \ge \f{\E(\bar{\bZ}\,|\,\FF(t))}{2}
	\,\bigg|\,\FF(t)\bigg)\,.
	\]
Combining the above with the second bound from Theorem~\ref{t:lbd} gives, again with high probability,
	\begin{align}\nonumber
	\P\bigg(
	\bar{\bZ} \ge \f{\exp\{ N(
		 \RS(\alpha;U) - o_t(1))\}}{2}
	\,\bigg|\,\FF(t)\bigg)
	&\stackrel{\eqref{e:restricted.first.mmt.lbd.half}}{\ge}
	\f{\exp\{ 2N(\RS(\alpha;U)-o_t(1))\}/4}
		 {\E(\bar{\bZ}^2\,|\,\FF(t)) }\\
	&\stackrel{\eqref{e:ubd.restricted.second.mmt}}{\ge}
	\f{1/4}{\exp(2N o_t(1))}\,.
	\label{e:second.mmt.lbd.cond.on.t}
	\end{align}
Next let $\P^j$ denote probability conditional on 
the first $j$ rows of $\bG$, and let $\E^j$ denote expectation with respect to $\P^j$. Then, as in the proof of \cite[Propn.~9.2.6]{MR3024566}, we take the martingale decomposition
	\[\f1N\bigg\{\log\bZ-\E\log\bZ\bigg\}
	= \sum_{j\le M} \f1N\bigg\{
		\E^j\log\bZ
		-\E^{j-1}\log\bZ\bigg\}
	\equiv \sum_{j\le M} X_j\,.\]
To bound $X_j$, let $\bZ_j$ denote the normalized partition function without the $j$-th factor,
	\beq\label{e:Z.j}
	\bZ_j
	\equiv \sum_J 
	\prod_{\substack{a\le M, \\
		 a \ne j}}
		U \bigg(\f{(\bg^a,J)}{N^{1/2}}\bigg)\,.\eeq
Since $\bZ_j$ does not depend on the $j$-th row of $\bG$, we can rewrite
	\[N X_j
	=\E^j \log \f{\bZ}{\bZ_j}
	-\E^{j-1} \log \f{\bZ}{\bZ_j}\,.\]
By Assumption~\ref{a:bdd} and the uniform bound on $u\equiv\log U$, we have
	\[\f1{\exp(\|u\|_\infty)}
	\le \f{\bZ}{\bZ_j}\le1\,,
	\]
which implies $|NX_j|\le\|u\|_\infty$ almost surely. It follows from the Azuma--Hoeffding bound that
	\beq\label{e:azuma.log.Z}
	\P\bigg(
	\Big|\log \bZ -\E\log\bZ\Big|
	\ge N\epsilon\bigg)
	\le 2\exp\bigg\{
	\f{-N\epsilon^2}{2\alpha(\|u\|_\infty)^2}\bigg\}
	\equiv \f{2}{\exp(N s(\epsilon))}\,.
	\eeq
On the other hand, if we fix any $\epsilon>0$, then \eqref{e:second.mmt.lbd.cond.on.t} implies
	\beq\label{e:second.mmt.lbd.cond.err.t}
	\P\bigg(\f1N\log\bZ
	\ge \RS(\alpha;U) - \epsilon - \f{\log2}{N}
	\bigg) \ge o_N(1)+
	\f{1/4}{\exp(2N o_t(1) )}
	\ge \f{1/4}{\exp(Ns(\epsilon)/2)}\,.
	\eeq
Note that \eqref{e:azuma.log.Z} and \eqref{e:second.mmt.lbd.cond.err.t} contradict one another unless
	\beq\label{e:E.log.Z.lbd}
	\f1N\E\log \bZ\ge
	\RS(\alpha;U) - 2\epsilon - \f{\log2}{N}\,.
	\eeq
It follows using \eqref{e:azuma.log.Z}
again that, for $N$ large enough,
	\begin{align*}
	&\P\bigg(\f1N\log\bZ
		\le \RS(\alpha;U) - 4\epsilon\bigg)
	\stackrel{\eqref{e:E.log.Z.lbd}}{\le}  
	\P\bigg(\log\bZ-\E\log\bZ
		\le  - N\epsilon\bigg)
	\stackrel{\eqref{e:azuma.log.Z}}{\le} 
	 o_N(1)\,.
	\end{align*}
In the above, the left-hand side does not depend on $t$, so it follows that
	\[\liminf_{N\to\infty}
	\f1N \log \bZ \ge \RS(\alpha;U) \]
in probability. This gives the lower bound in Theorem~\ref{t:main} in the case $\|u\|_\infty<\infty$.
\end{proof}

\section{Local central limit theorem}
\label{s:local.clt}

In this section we state and prove Proposition~\ref{p:density.bound}
(used in the proofs 
\hyperlink{proof:p.first.mmt}{of Proposition~\ref{p:first.mmt}} and
Proposition~\ref{p:restricted.first.mmt})
and Proposition~\ref{p:pair.density.bound}
(used in the 
\hyperlink{proof:p.ubd.pair.lambda}{proof of Proposition~\ref{p:ubd.pair.lambda}}).
Recall the calculation of $\tbX_J$ from  Lemma~\ref{l:G.R.J}. Given $J\in\set{-1,+1}^N$ and $\tau\in\R^{t-1}$, we define
	\begin{align}\nonumber
	\bX\equiv\bX_{J,\tau}
	&= \tbX_J + N^{1/2} c(\pi(J))
	\bc[t-1]^\st\tau
	= \f{\bh[t]^\st\hpi}{q^{1/2}}
	+\f{(1-q)\bn[t-1]^\st\cpi}{q^{1/2}}
	+N^{1/2} c(\pi(J))
	\bc[t-1]^\st\tau \\
	&\stackrel{\eqref{e:gs.n.c.EXACT}}{=} 
	\f{\bh[t]^\st\hpi}{q^{1/2}}
	+\bn[t-1]^\st\bigg\{
	\f{(1-q)\cpi}{q^{1/2}}
	+ \f{c(\pi(J))}{\psi^{1/2}}
	((\bGam_N)^\st)^{-1}\tau
	\bigg\}\,.\label{e:bX.J.tau.rewrite}
	\end{align}
Let $\bze_a$ ($a\le M$) be independent scalar random variables, such that $\bze_a$ has density
given by (cf.\ Definition~\ref{d:mu.x.c})
	\beq\label{e:bxi.a.density}
	\chi_{\bX_a,c}(z)\equiv \f{U(\bX_a+cz)\varphi(z)}
		{\E_\xi U(\bX_a+c\xi)}\,,
	\eeq
where $\bX\equiv\bX_{J,\tau}$ as above, and $c\equiv c(\pi(J))\equiv (1-\|\pi(J)\|^2)^{1/2}$. 
Note that
	\beq\label{e:adm.zeta} \E\bze_a
	= \f{\E_\xi[\xi U(\bX_a+c\xi)]}
		{\E_\xi U(\bX_a+c\xi)}
	\stackrel{\eqref{e:F}}{=}
	c(\pi(J)) F_{\|\pi(J)\|^2}(\bX_a)
	\,.\eeq
Let $\bn_a\in\R^{t-1}$ denote the $a$-th column of the matrix $\bn[t-1]$, and consider the random variable
	\beq\label{e:local.clt.W}
	\bW
	\equiv \f1{N^{1/2}}
	\sum_{a\le M} 
	(\bze_a-\E\bze_a)
	 \bn_a
	= \f{\bn[t-1](\bze-\E\bze)}{N^{1/2}}
	\in \R^{t-1}\,.\eeq
Let $P_{J,\tau}$ denote the law of $\bW$.
We will compare $P_{J,\tau}$ with the 
gaussian distribution on $\R^{t-1}$ that has mean zero and covariance
	\beq\label{e:cov.Sigma}
	\Sigma
	\equiv
	\Sigma_{J,\tau}
	\equiv \f1N\sum_{a\le M}
		(\Var \bze_a) \bn_a(\bn_a)^\st
	\in\R^{(t-1)\times(t-1)}\,.\eeq
(We bound the singular values of $\Sigma_{J,\tau}$
in Lemma~\ref{l:nondegenerate.Sigma} below.) The majority of this section is occupied with proving the following result:

\begin{ppn}[local central limit theorem]\label{p:local.CLT}
Suppose $U$ satisfies Assumptions~\ref{a:bdd} and \ref{a:Lip}.
Recall that $P_{J,\tau}$ is the law of the random variable $\bW$ from \eqref{e:local.clt.W}. For any finite constant $\tau_{\max}$, it holds with  high probability
that for all $J\in\set{-1,+1}^N$
and all $\|\tau\|\le \tau_{\max}$, the measure $P_{J,\tau}$ has a bounded continuous density $p_{J,\tau}$. Moreover, again with high probability,
	\[\sup\bigg\{
	\|p_{J,\tau}-g_{J,\tau}\|_\infty:
	J\in\set{-1,+1}^N, \|\pi(J)\|\le \f45,
	\|\tau\|\le \tau_{\max}
	\bigg\}
	\le \f1{(2\pi)^{t-1}N^{0.35}}
	\le \f1{N^{0.3}}\,,
	\]
where $g_{J,\tau}$ denotes the density of the centered gaussian distribution on $\R^{t-1}$ with covariance $\Sigma\equiv\Sigma_{J,\tau}$. 
\end{ppn}

At the end of this section we will show that Proposition~\ref{p:local.CLT} readily implies  the required results Propositions~\ref{p:density.bound} and \ref{p:pair.density.bound}.
Towards the \hyperlink{proof:p.local.CLT}{proof of  Proposition~\ref{p:local.CLT}}, we introduce some notation. Write $p_a$ for the density function of the random variable $\bze_a-\E\bze_a$, so
in the notation of \eqref{e:bxi.a.density} we have
	\[p_a(z)
	= \chi_{\bX_a,c}(z+\E\bze_a)
	= \chi_{\bX_a,c}\bigg(
		z + \f{\E_\xi[\xi U(\bX_a+c\xi)]}
			{\E_\xi U(\bX_a+c\xi)}
	\bigg)\,.\]
The characteristic function of the random variable $\bW$ from \eqref{e:local.clt.W} (i.e., the Fourier transform of the measure $P_{J,\tau}$) is given by the function
	\beq\label{e:hat.p}
	\hat{p}(\sss)
	\equiv\hat{p}_{J,\tau}(\sss)
	\equiv \E \exp(\ii(\sss,\bW))
	= \prod_{a\le M}
	\E\exp\bigg\{
	\f{\ii(\sss,\bn_a)
		(\bze_a-\E\bze_a)}{N^{1/2}}\bigg\}
	= \prod_{a\le M}
	\hat{p}_a\bigg(
	\f{(\sss,\bn_a)}{N^{1/2}}\bigg)\,,
	\eeq
where $\hat{p}_a$ denotes the Fourier transform of $p_a$. The Fourier transform of the gaussian density $g\equiv g_{J,\tau}$ 
 is given by
	\beq\label{e:hat.g}
	\hat{g}(\sss)
	\equiv\hat{g}_{J,\tau}(\sss)
	\equiv \exp\bigg\{-\f{(\sss,\Sigma\sss)}{2}\bigg\}
	= \prod_{a\le M}
		\exp\bigg\{-\f{(\sss,\bn_a)^2
			\Var\bze_a}{2N}\bigg\}
	\equiv\prod_{a\le M} \hat{g}_a(\sss)\,.\eeq
With $\hat{p}\equiv\hat{p}_{J,\tau}$ as in \eqref{e:hat.p}
and $\hat{g}\equiv\hat{g}_{J,\tau}$ as in \eqref{e:hat.g}, we define
	\begin{align}
	\label{e:fourier.ONE}
	I_1(J,\tau) & \equiv\int
		\Big|\hat{p}_{J,\tau}(\sss) 
	-\hat{g}_{J,\tau}(\sss)\Big|
		\I\Big\{\|\sss\|\le N^{0.01}\Big\}
	\,d\sss\,,\\
	\label{e:fourier.TWO}
	I_2(J,\tau,\epsilon_2)
		&\equiv
	\int\Big|\hat{p}_{J,\tau}(\sss)
		-\hat{g}_{J,\tau}(\sss)\Big| 
	\I\Big\{
		N^{0.01} \le \|\sss\| \le \epsilon_2 N^{1/2}
	\Big\} \,d\sss\,, \\
	\label{e:fourier.THREE}
	I_3(J,\tau,\epsilon_2)
		&\equiv
	\int \Big| \hat{p}_{J,\tau}(\sss)
			-\hat{g}_{J,\tau}(\sss)\Big| 
		\I\Big\{ \|\sss\| \ge 
		\epsilon_2 N^{1/2}
	\Big\} \,d\sss
		\,.
	\end{align}
In the analysis below we show that the integrals $I_j(J,\tau)$ can be bounded uniformly over  $J\in\set{-1,+1}^N$ such that $\|\pi(J)\|\le4/5$,
and any bounded range of vectors $\tau$. The remainder of this section is organized as follows:
\begin{itemize}
\item In \S\ref{ss:fourier.low.med}
	we bound the quantities $I_1$ and $I_2$ from
	\eqref{e:fourier.ONE} and \eqref{e:fourier.TWO}.
\item In \S\ref{ss:tap.nondeg}, in preparation for bounding $I_3$ from  \eqref{e:fourier.THREE},
we prove rough estimates concerning the nondegeneracy of the vectors arising from the AMP iteration.

\item In \S\ref{ss:fourier.high.freq}
	we bound $I_3$ from \eqref{e:fourier.THREE}.
\item In \S\ref{ss:local.clt.conclusion}
	we combine the bounds from the preceding sections to finish the proof of Proposition~\ref{p:local.CLT}.
	We then state and prove
Proposition~\ref{p:density.bound}
and \ref{p:pair.density.bound}.
\end{itemize}
The analysis of this section is based on standard methods; see e.g.\ 
\cite{MR0388499,borovkov2017generalization}.

\subsection{Fourier estimates at low and intermediate frequency}
\label{ss:fourier.low.med}

In this subsection we prove Lemmas \ref{l:fourier.ONE} and \ref{l:fourier.TWO},
bounding the quantities $I_1$ and $I_2$ from \eqref{e:fourier.ONE} and \eqref{e:fourier.TWO}.

\begin{lem}\label{l:nondegenerate.Sigma}
Suppose $U$ satisfies Assumption~\ref{a:bdd} and \ref{a:Lip},
and let $\Sigma$ be as in \eqref{e:cov.Sigma}. Given any $\tau_{\max}<\infty$, there is a positive constant $\iota_1$, depending on $t$ and on $\tau_{\max}$, such that we have the bounds
	\begin{align}\label{e:Sigma.N.ubd}
	\inf\bigg\{
	(u,\Sigma_{J,\tau}u)
	: \|u\|=1, J\in\set{-1,+1}^2,
	\|\pi(J)\|\le\f45,
	\|\tau\|\le \tau_{\max}
	\bigg\} \ge \iota_1\,, \\
	\sup\bigg\{
	(u,\Sigma_{J,\tau}u)
	: \|u\|=1, J\in\set{-1,+1}^2,
	\|\pi(J)\|\le\f45,
	\|\tau\|\le \tau_{\max}
	\bigg\} \le \f1{\iota_1}\,,
	\label{e:Sigma.N.lbd}
	\end{align}
with probability $1-o_N(1)$.
	
\begin{proof}
Abbreviate $v_a\equiv\Var\bze_a$, and note that Assumption~\ref{a:Lip} gives, with $c=c(\pi(J))
\equiv(1-\|\pi(J)\|^2)^2$,
	\beq\label{e:var.W.a.ubd}
	v_a
	= v(\bX_a,c)
	\equiv\f12 \f{\E_\xi[(\xi-\xi')^2
		U(\bX_a+c\xi)U(\bX_a+c\xi')]}
		{\E_{\xi,\xi'}[U(\bX_a+c\xi)U(\bX_a+c\xi')]}
	\le \f{K_2(U)}{2}\,.
	\eeq
It follows that, for any unit vector $u\in\R^{t-1}$, and with $\varsigma_t$ as defined by Remark~\ref{r:singval.bound}, we have
	\begin{align*}
	(u,\Sigma u)
	&= \f1N\sum_{a\le M} v_a (\bn_a,u)^2
	\le \f{K_2(U)}{2N}
	\sum_{a\le M} (\bn_a,u)^2\\
	&= \f{K_2(U)}{2N}
	\Big\|\bn[t-1]^\st u\Big\|^2
	\stackrel{\eqref{e:gs.n.c.EXACT}}{=}
	\f{K_2(U) \psi}{2}
	\Big\|(\bGam_N)^\st u\Big\|^2
	\le\f{K_2(U) \psi \varsigma_t}{2}\,,
	\end{align*}
which proves \eqref{e:Sigma.N.ubd}. Next, for any $L$, let $M(L)\subseteq[M]$ denote the subset of indices $a\le M$ satisfying the condition
	\beq\label{e:mathfrak.a}
	\mathfrak{m}_a
	\equiv
	\max\bigg\{
		|(\bh^{(s)})_a|,|(\bn^{(\ell)})_a|
		: s\le t,\ell\le t-1
		\bigg\} \le L\,.\eeq
It follows from Lemma~\ref{l:AMP.PL} that with high probability we can bound
	\[
	\max\bigg\{
		\f1M\sum_{a\le M}((\bh^{(s)})_a)^4,
		\f1M\sum_{a\le M}((\bn^{(\ell)})_a)^4
		:s\le t,\ell\le t-1
	\bigg\} \le \wp_4
	\]
for a constant $\wp_4$. As a result, for any finite $L$, we can bound
	\[
	\f1M \sum_{a\le M}
	\I\Big\{|(\bh^{(s)})_a|\ge L\Big\}
	\le 
	\f1M \sum_{a\le M} 
	\f{((\bh^{(s)})_a)^4}{L^4}
	\le \f{\wp_4}{L^4}\,,
	\]
and similarly with $\bn^{(\ell)}$ in place of $\bh^{(s)}$. It follows using the Cauchy--Schwarz that
	\begin{align*}
	\f1M
	\sum_{a\le M}
	(\bn^{(\ell)}_a)^2
	\I\Big\{|(\bh^{(s)})_a|\ge L\Big\}
	&\le
	\bigg(\f1M\sum_{a\le M}
	(\bn^{(\ell)}_a)^4\bigg)^{1/2}
	\bigg(\f1M\sum_{a\le M}
	\I\Big\{|(\bh^{(s)})_a|\ge L\Big\}
	\bigg)^{1/2}
	\le \f{\wp_4}{L^2}\,,\\
	\f1M
	\sum_{a\le M}
	(\bn^{(\ell)}_a)^2
	\I\Big\{|(\bn^{(j)})_a|\ge L\Big\}
	&\le
	\bigg(\f1M\sum_{a\le M}
	(\bn^{(\ell)}_a)^4\bigg)^{1/2}
	\bigg(\f1M\sum_{a\le M}
	\I\Big\{|(\bn^{(j)})_a|\ge L\Big\}
	\bigg)^{1/2}
	\le \f{\wp_4}{L^2}\,,
	\end{align*}
where the bounds hold for all $s\le t$ and all $j,\ell\le t-1$.
Combining these bounds gives, with $\mathfrak{m}_a$ as defined in \eqref{e:mathfrak.a},
	\beq\label{e:bound.outside.M.L}
	\f1M
	\sum_{a\notin M(L)}
	(\bn^{(\ell)}_a)^2
	= 
	\f1M
	\sum_{a\le M}
	(\bn^{(\ell)}_a)^2
	\I\Big\{|\mathfrak{m}_a|\ge L\Big\}
	\le \f{2t\wp_4}{L^2}\,.
	\eeq
Next, it follows from the definition \eqref{e:bX.J.tau.rewrite} 
of $\bX\equiv\bX_{J,t}$ that for all $a\le M$,
	\begin{align}\nonumber
	|\bX_a|
	&\le \f{\|\hpi\|_\infty}{q^{1/2}}
	\sum_{s\le t} |(\bh^{(s)})_a|
	+\bigg( \f{\|\hpi\|_\infty}{q^{1/2}}
	+\f{\varsigma_t\|\tau\|}{\psi^{1/2}}
	\bigg)\sum_{\ell\le t-1} |(\bn^{(\ell)})_a|\\
	&\le \f{\varsigma_t
		(1+\|\tau\|)}{(q\psi)^{1/2}}
	\sum_{s\le t}\bigg( |(\bh^{(s)})_a|
		+|(\bn^{(s)})_a|\bigg)\,.
	\label{e:bX.J.tau.bound}
	\end{align}
If $\|\tau\|\le\tau_{\max}$ where (without loss) $\tau_{\max}\ge1$, then we obtain
	\beq\label{e:bound.X.a.on.M.L}
	\max\bigg\{
	|\bX_a| : a\in M(L)\bigg\}
	\le 
	\f{4\varsigma_t\tau_{\max} tL}{(q\psi)^{1/2}}
	\equiv L'\,.\eeq
It follows using Assumption~\ref{a:bdd} that
for any finite $L'$ we must  have
	\[
	\inf\bigg\{
	v(x,c) : \f12\le c\le1, |x|\le L'
	\bigg\} \ge \epsilon(L') >0\,.
	\]
It follows that, for any unit vector $u\in\R^{t-1}$, we have the lower bound
	\begin{align*}
	(u,\Sigma u)
	&\ge \f{\epsilon(L')}{N}
		\sum_{a\in M(L)}
		 (\bn_a,u)^2
	\ge \f{\epsilon(L')}{N}
	\bigg\{
	\Big\|\bn[t-1]^\st u\Big\|^2
	-\sum_{a\notin M(L)} (\bn_a,u)^2
	\bigg\} \\
	&\stackrel{\eqref{e:gs.n.c.EXACT}}{\ge}
	\epsilon(L')
	\bigg\{ \psi
	\Big\| (\bGam_N)^\st u\Big\|^2
	-\f1N \sum_{a\notin M(L)} \|\bn_a\|^2
	\bigg\}
	\stackrel{\eqref{e:bound.outside.M.L}}{\ge}
	\epsilon(L')
	\bigg\{
	\f{\psi}{\varsigma_t}
	- \f{2t^2\wp_4}{L^2}
	\bigg\}
	\ge \f{\epsilon(L')\psi}{2\varsigma_t}\,,
	\end{align*}
where the last inequality can be arranged by taking $L$ large enough (note that $L$ depends on $t$, and $L'$ depends on $L$). This proves the second assertion \eqref{e:Sigma.N.lbd}.
\end{proof}
\end{lem}

\begin{lem}[Taylor expansion of characteristic function]
\label{l:cgf.taylor.expand}
Suppose $U$ satisfies Assumptions~\ref{a:bdd} and \ref{a:Lip}.
Let $\hat{p}_a$ be as in \eqref{e:hat.p}, and recall that it depends on both $J$ and $\tau$. It holds with high probability that
	\[\max\bigg\{
	\bigg|
	\hat{p}_a(\sss)
	-\bigg( 1 
		- \f{(\sss,\bn_a)^2\Var\bze_a}{2N}
		\bigg)\bigg|
	: J\in\set{-1,+1}^N,
	\|\pi(J)\|\le\f45,
	\|\tau\|\le N^{0.01}
	\bigg\}
	\le\f{\|\sss\|^3}{N^{1.4}}
	\]
for all $\sss\in\R^{t-1}$ and all $a\le M$.

\begin{proof} 
It is well-known that for all $x\in\R$ we have	\[
	\bigg|
	e^{ix}-\bigg(1 + ix - \f{x^2}{2}\bigg)\bigg|
	\le \f{|x|^3}{6}\,.
	\]
We also note that Lemma~\ref{l:poly} implies the third moment bound
	\[\E\bigg( \Big|\bze_a-\E\bze_a\Big|^3\bigg)
	\le 8\E(|\bze_a|^3)
	=\f{8\E_\xi[|\xi|^3U(\bX_a+c\xi)]}
		{\E_\xi U(\bX_a+c\xi)}
	\le 8\bigg( C_1(U) + (8|\bX_a|)^3
	 \bigg)\,.
	\]
As a consequence, for all $\sss\in\R^{t-1}$ we have
	\begin{align}\nonumber
	\bigg|
	\hat{p}_a(\sss)
	-\bigg( 1 
		- \f{(\sss,\bn_a)^2\Var\bze_a}{2N}
		\bigg)\bigg|
	&\le \f{|(\sss,\bn_a)|^3}{6N^{3/2}}
	\E\bigg(\Big|\bze_a-\E\bze_a\Big|^3\bigg)\\
	&\le
	\f{4|(\sss,\bn_a)|^3}{3N^{3/2}}
	\bigg(
	C_1(U)+(8|\bX_a|)^3\bigg)
	\label{e:cgf.taylor}
	\,.\end{align}
By combining Lemma~\ref{l:TAP.supnorm} with the bound \eqref{e:bX.J.tau.bound} and the restriction $\|\tau\|\le N^{0.01}$, we must have $\|\bX\|_\infty\le N^{0.021}$ with high probability. Therefore, with high probability,
	\[
	\f{4|(\sss,\bn_a)|^3}{3N^{3/2}}
	\bigg(C_1(U)+ (8|\bX_a|)^3\bigg)
	\le \f{4\|\sss\|^3 t^{3/2} N^{0.03}}
		{3N^{3/2}}
		\bigg( C_1(U)
	 + (8|\bX_a|)^3
	 \bigg)
	\le \f{\|\sss\|^3}{N^{1.4}}\,.
	\]
Combining with \eqref{e:cgf.taylor} concludes the proof.
\end{proof}
\end{lem}

\begin{lem}[low-frequency estimate]
\label{l:fourier.ONE}
Suppose $U$ satisfies Assumption~\ref{a:bdd} and \ref{a:Lip}. In the notation of \eqref{e:fourier.ONE}, we have
	\[\max\bigg\{
	I_1(J,\tau)
	: J\in\set{-1,+1}^N,
	\|\pi(J)\|\le\f45,
	\|\tau\| \le N^{0.01}
	\bigg\} \le \f1{N^{0.38}}
	\]
with probability $1-o_N(1)$.

\begin{proof} 
Recall from 
\eqref{e:var.W.a.ubd} that $\Var\bze_a\le K_2(U)/2$.
Combining with Lemma~\ref{l:TAP.supnorm} gives, with high probability,
	\[	\f{(\sss,\bn_a)^2\Var\bze_a}{2N}
	\le\f{ \|\sss\|^2 t(\|\bn_a\|_\infty)^2 K_2(U)}{4N}
	\le\f{ \|\sss\|^2 t K_2(U)}{4N^{0.98}}
	\le\f{\|\sss\|^2}{N^{0.97}}\,,
	\]
We have $|\log(1-x)+x|\le x^2$ for all $x$ small enough, so if $\|\sss\|\le N^{0.01}$, then combining with Lemma~\ref{l:cgf.taylor.expand} gives
	\begin{align*}
	\bigg|
	\log \hat{p}_a(\sss)
	+ \f{(\sss,\bn_a)^2\Var\bze_a}{2N}\bigg|
	&\le
	\f{\|\sss\|^3}{N^{1.4}}
	+ \bigg(
	\f{\|\sss\|^2}{N^{0.97}}
	+\f{\|\sss\|^3}{N^{1.4}}
	\bigg)^2\\
	&\le
	\f{\|\sss\|^3}{N^{1.4}}
	+\bigg(\f{2\|\sss\|^2}{N^{0.97}}\bigg)^2
	\le\f{2\|\sss\|^3}{N^{1.4}}
	\le \f1{N^{1.39}}\,.
	\end{align*}
Summing the above over $a\le M$ gives
that the multiplicative error between $\hat{p}(\sss)$ and $\hat{g}(\sss)$ is small for all $\|\sss\|\le N^{0.01}$. Therefore, with high probability, we have the bound
	\[
	I_1(J,\tau)
	\le
	\int
	\hat{g}(\sss)
	\bigg\{ \exp\bigg( \f{M}{N^{1.39}}\bigg)-1\bigg\}
	\,d\sss
	\le \f{(2\pi)^{(t-1)/2}}
		{N^{0.385} (\det \Sigma)^{1/2}}
	\le \f1{N^{0.38}}\,,
	\]
uniformly over all $J\in\set{-1,+1}^N$ and all $\|\tau\|\le N^{0.01}$.
\end{proof}
\end{lem}

\begin{lem}[moderate-frequency estimate]
\label{l:fourier.TWO} Suppose $U$ satisfies Assumption~\ref{a:bdd} and \ref{a:Lip}. With the notation of \eqref{e:fourier.TWO}, 
for any finite constant $\tau_{\max}$,
we can choose $\epsilon_2$ depending on $\tau_{\max}$ such that
	\[\sup
	\bigg\{
	I_2(J,\tau,\epsilon_2)
	: J\in\set{-1,+1}^N,
	\|\pi(J)\|\le\f45,
	\|\tau\|\le \tau_{\max}
	\bigg\}
	\le \f1{\exp(N^{0.01})}
	\]
with probability $1-o_N(1)$.

\begin{proof}
It follows from
the bound \eqref{e:cgf.taylor} in
the proof of
 Lemma~\ref{l:cgf.taylor.expand} that,
with high probability, we have
	\[|\hat{p}(\sss)|
	\le \hat{g}(\sss)
	\exp\bigg\{
	\f{4\|\sss\|^3}{3N^{3/2}}
	\sum_{a\le M}
	\|\bn_a\|^3
	\Big( C_1(U)+(8|\bX_a|)^3 \Big)	\bigg\}
	\]
for all $\sss\in\R^{t-1}$. Recall the bound \eqref{e:bX.J.tau.bound} on $\bX=\bX_{J,\tau}$. If we assume without loss of generality that 
$\tau_{\max}\ge1$, then combining \eqref{e:bX.J.tau.bound} with
 Lemma~\ref{l:AMP.PL} gives, with high probability,
	\[
	\sup\bigg\{\f{4}{3N} \sum_{a\le M}\|\bn_a\|^3
	\Big( C_1(U)+(8|\bX_a|)^3 \Big)
	: J\in\set{-1,+1}^n,\|\tau\|\le \tau_{\max}\bigg\}
	\le 
	(\tau_{\max})^3 \wp_1\,,
	\]
where $\wp_1$ is a finite constant. On the other hand, by Lemma~\ref{l:nondegenerate.Sigma}, with high probability
	\[\hat{g}(\sss)
	\le \exp\bigg\{ -\f{\iota_1\|\sss\|^2}{2}
	\bigg\}\,.
	\]
It follows that, with high probability,
	\[
	|\hat{p}(\sss)|
	\le \exp\bigg\{
	-\f{\iota_1\|\sss\|^2}{2}
	+\f{
	(\tau_{\max})^3\wp_1
	\|\sss\|^3}{N^{1/2}} 
	\bigg\}\,.
	\]
To ensure that the quadratic term dominates the cubic term, we
restrict to $\|\sss\|\le \epsilon_2 N^{1/2}$ where
	\[\epsilon_2
	\equiv
	\f{\iota_1}{4(\tau_{\max})^3\wp_1}\,.
	\]
For this choice of $\epsilon_2$ we find that, with high probability, we have the bound
	\[
	I_2(J,\tau,\epsilon)
	\le \int_{\|\sss\|\ge N^{0.01}}
	\exp\bigg\{ -\f{\iota_1\|\sss\|^2}{4} \bigg\}
		\,d\sss
	\le \f1{\exp(N^{0.01})}
	\]
uniformly over all $J\in\set{-1,+1}^N$ and $\|\tau\|\le\tau_{\max}$.
\end{proof}
\end{lem}

\subsection{Non-degeneracy of TAP iterates}
\label{ss:tap.nondeg}

In this subsection we prove some preliminary results 
which will be used in \S\ref{ss:fourier.high.freq} to estimate the quantity $I_3$ from \eqref{e:fourier.THREE}.

\newcommand{\bb}{\mathbf{b}}
\begin{lem}\label{l:submatrix}
If $B$ is any $k\times M$ matrix such that $BB^\st=I_k$, then $B$ has a $k\times k$ submatrix $U$ such that
	\[|\det U|
	\ge \bigg(\f{k!}{M^k}\bigg)^{1/2}
	\,.\]

\begin{proof}
We argue by induction on $k$. If $k=1$ then $B$ consists of a single row which is a unit vector in $\R^M$, so clearly $B$ must have an entry with absolute value at least $1/M^{1/2}$. Now suppose $k\ge2$ and that the claim has been proved up to $k-1$. Denote the columns of $B$ as $\bb_1,\ldots,\bb_M$ where each $\bb_a\in\R^k$. Since
	\[
	\sum_{a\le M}\|\bb_a\|^2
	=\tr(BB^\st)=k\,,
	\]
there must exist at least one index $a\le M$ with
	\[\|\bb_a\|^2 \ge \f{k}{M}\,.\]
We assume without loss that $a=1$. Let $O$ be a $k\times k$ orthonormal matrix such that
$O\bb_1 = \|\bb_1\|e_1$, where $e_1$ denotes the first standard basis vector in $\R^k$. Let $\bar{B}\equiv OB$, and note that $\bar{B}\bar{B}^\st=OBB^\st O^\st=I_k$, so $\bar{B}$ also has orthonormal rows. We can further decompose
	\[\bar{B}= OB
	= \begin{pmatrix}
	\|\bb_1\| & * \\ 
	\mathbf{0} & \tilde{B}
	\end{pmatrix}
	\]
where $\mathbf{0}$ denotes the zero vector in $\R^{k-1}$, and $\tilde{B}$ is a $(k-1)\times (M-1)$ matrix with orthonormal rows. It follows from the inductive hypothesis that $\tilde{B}$ has a $(k-1)\times(k-1)$ submatrix $\tilde{U}$ with
	\[|\det\tilde{U}|\ge
	\bigg( \f{(k-1)!}{(M-1)^{k-1}} \bigg)^{1/2}
	\,.\]
As a result, 
$\bar{B}$ has a $k\times k$ submatrix $\bar{U}$
with
	\[
	|\det \bar{U}|
	=\bigg| \det
	\begin{pmatrix}
	\|\bb_1\| & * \\ 
	\mathbf{0} & \tilde{U}
	\end{pmatrix}\bigg|
	\ge \|\bb_1\| \cdot |\det\tilde{U}|
	\ge \f{k^{1/2}}{M^{1/2}}
	\bigg( \f{(k-1)!}{(M-1)^{k-1}} \bigg)^{1/2}
	\ge\bigg( \f{k!}{M^k} \bigg)^{1/2}\,.
	\]
The claim follows by noting that $U=O^\st\bar{U}$ is a submatrix of the original matrix $B$.
\end{proof}
\end{lem}

\begin{cor}\label{c:submatrix.APPROX}
If $B$ is any $k\times M$ matrix such that 
$\|BB^\st-I_k\|_\infty \le 1/(3k)$, then $B$ has a $k\times k$ submatrix $U$ with
	\beq\label{e:det.U.one.lbd}
	|\det U|
	\ge \f13 \bigg(\f{k!}{M^k}\bigg)^{1/2}
	\,.\eeq
(In the above, as elsewhere, $\|\cdot\|_\infty$ denotes the entrywise maximum absolute value of the matrix.)

\begin{proof}
Denote the rows of $B$ as $\bu^1,\ldots,\bu^k$ where each $\bu^\ell\in\R^M$. Consider the 
Gram--Schmidt orthogonalization of these vectors: for each $\ell\le k$, we decompose
	\[
	\bu^\ell
	\equiv\bu^{\ell,\parallel}+\bu^{\ell,\perp}
	\equiv\sum_{j\le\ell-1} c_{\ell,j}\bu^j
		+\bu^{\ell,\perp}
	\]
where $\bu^{\ell,\parallel}$ is the orthogonal projection of $\bu^\ell$ onto the span of $\bu^1,\ldots,\bu^{\ell-1}$. Then for all $j\le\ell-1$ we must have
	\[
	0 = (\bu^{\ell,\perp},\bu^j)
	= (\bu^\ell,\bu^j)
		-\sum_{i\le\ell-1}
		\Ind{i\ne j} c_{\ell,i}(\bu^i,\bu^j)
		-c_{\ell,j}\|\bu^j\|^2\,.
	\]
Abbreviate $\epsilon\equiv \epsilon(k) \equiv 1/(3k)$, so the assumptions imply that
 $\|\bu^j\|^2\ge1-\epsilon$ while 
$|(\bu^i,\bu^j)|\le\epsilon$ for all $i\ne j$.
Rearranging the above gives an upper bound for $|c_{\ell,j}|$ in terms of the other coefficients $c_{\ell,i}$ ($i\ne j$). If we further denote $c_{\max}\equiv\max\set{|c_{\ell,j}| : \ell\le k, j\le \ell-1}$, then we have
	\[
	c_{\max}
	\le \f{\epsilon}{1-\epsilon}
	\bigg\{1 + (k-2) c_{\max}\bigg\}\,.
	\]
Rearranging the inequality gives the bound
	\[
	c_{\max}
	\le \f{\epsilon}{1-\epsilon}\bigg/
	\bigg(1 - \f{\epsilon(k-2)}{1-\epsilon}\bigg)
	= \f{\epsilon}{1-\epsilon(k-1)}
	\le \f{3\epsilon}{2}\,.
	\]
From this bound we can deduce that for all $\ell\le k$ we have
	\[
	\|\bu^{\ell,\parallel}\|^2
	=\bigg\|\sum_{j\le\ell-1}c_{\ell,j}\bu^j\bigg\|^2
	\le 
	\bigg(\f{3\epsilon}{2}\bigg)^2
	\bigg\{(\ell-1)(1+\epsilon)
	+(\ell-1)(\ell-2)\epsilon\bigg\}
	\le 
	\bigg(\f{3\epsilon}{2}\bigg)^2 \f{4k}{3}
	= \epsilon\,.
	\]
It follows that
$\|\bu^{\ell,\perp}\|^2= \|\bu^\ell\|^2-\|\bu^{\ell,\parallel}\|^2\ge 1-2\epsilon$.
(We also have trivially
$\|\bu^{\ell,\perp}\|^2\le\|\bu^\ell\|^2\le1+\epsilon$.)
Let $R$ denote the Gram--Schmidt matrix, so $R$ is $k\times k$ lower triangular with entries
	\[
	R_{\ell,j} =
	\f1{\|\bu^{\ell,\perp}\|}
	\bigg\{
	\Ind{\ell=j} - \Ind{\ell<j}c_{\ell,j}
	\bigg\}\,.
	\]
Since $R$ is lower triangular, its determinant is simply the product of its diagonal entries, so
	\[
	\f1{(1+\epsilon)^k} \le \det R
	= \prod_{\ell\le k} \f1{\|\bu^{\ell,\perp}\|}
	\le \f1{(1-2\epsilon)^k}\,.
	\]
By construction, $\acute{B}=RB$ is a $k\times M$ matrix with orthonormal rows, so Lemma~\ref{l:submatrix} implies that $\acute{B}$ has a $k\times k$ submatrix $\acute{U}$ with
	\[|\det \acute{U}| \ge
		\bigg(\f{k!}{M^k}\bigg)^{1/2}\,.\]
Therefore $U=R^{-1}\acute{U}$ is a $k\times k$ submatrix of the original matrix $B$, with
	\[
	|\det U|
	= \f{|\det \acute{U}|}{\det R}
	\ge (1-2\epsilon)^k
	\bigg(\f{k!}{M^k}\bigg)^{1/2}
	= \bigg(1-\f{2}{3k}\bigg)^k
	\bigg(\f{k!}{M^k}\bigg)^{1/2}
	\ge\f13
	\bigg(\f{k!}{M^k}\bigg)^{1/2}\,,
	\]
where the bound holds for all $k\ge1$.
\end{proof}
\end{cor}

\begin{lem}\label{l:c.submatrix}
Suppose $U$ satisfies Assumption~\ref{a:bdd} and \ref{a:Lip}. Recall from \eqref{e:gs.n.c.EXACT} that the matrix $\bc[t-1]$ is $(t-1)\times M$ with orthonormal rows. Let $M(L)\subseteq[M]$ be as defined in the proof of Lemma~\ref{l:nondegenerate.Sigma} (see \eqref{e:mathfrak.a}). If $B=B(L)$ is the submatrix of $\bc[t-1]$ with column indices in $M(L)$, then 
with high probability it satisfies
	\[\|B\|_\infty\le \f{t\varsigma_t L}{(N\psi)^{1/2}}
	\,.\]
It is possible to choose $L=L(t)$ large enough such that, with high probability, 
$\|BB^\st-I_{t-1}\|_\infty\le 1/(4t)$.

\begin{proof} It follows using \eqref{e:gs.n.c.EXACT} that for each $\ell\le t-1$,
	\beq\label{e:bound.c.by.n}
	\Big((\bc^{(\ell)})_a\Big)^2
	=\bigg(
	\sum_{j\le t-1} \f{((\bGam_N)^{-1})_{\ell,j}
		(\bn^{(j)})_a}{(N\psi)^{1/2}}
		\bigg)^2
	\le \f{t(\varsigma_t)^2}{N\psi}
	\sum_{j\le t-1} \Big((\bn^{(j)})_a\Big)^2\,.
	\eeq
Applying \eqref{e:bound.c.by.n} for $a\in M(L)$ gives the claimed bound on $\|B(L)\|_\infty$. On the other hand, by applying \eqref{e:bound.c.by.n} for $a\notin M(L)$ and combining with the bound
\eqref{e:bound.outside.M.L} from
the proof of Lemma~\ref{l:nondegenerate.Sigma}, we find, with high probability,
	\[
	\sum_{a\notin M(L)}
	\Big((\bc^{(\ell)})_a\Big)^2
	\le\f{t(\varsigma_t)^2}{N\psi}
	\sum_{j\le t-1} \sum_{a\notin M(L)}
	\Big((\bn^{(j)})_a\Big)^2
	\stackrel{\eqref{e:bound.outside.M.L}}{\le}
	\f{M t^2(\varsigma_t)^2}{N\psi}
	\cdot \f{2t\wp_4}{L^2}\,,
	\]
which can be made $\le 1/(4t)$ by choosing $L$ large enough. Then, for any $\ell,j\le t-1$, we have
	\[\bigg|
	\sum_{a\in M(L)}
	(\bc^{(\ell)})_a(\bc^{(j)})_a
	-\Ind{\ell=j}\bigg|
	= \bigg|\sum_{a\notin M(L)}(\bc^{(\ell)})_a(\bc^{(j)})_a
	\bigg|
	\le \f1{4t}\,,
	\]
which shows that the matrix $B=B(L)$ satisfies $\|BB^\st-I_{t-1}\|_\infty\le 1/(4t)$ as desired.
\end{proof}
\end{lem}

\begin{cor}\label{c:n.submatrices}
Let $M(L)\subseteq[M]$ be as in Lemma~\ref{l:c.submatrix}, where $L=L(t)$.
With high probability, the matrix $\bn[t-1]$ has disjoint $(t-1)\times(t-1)$ submatrices $A_1,\ldots,A_{\lfloor N^{0.9}\rfloor}$, all involving only columns indexed by $M(L)$, such that
each $A_i$ has minimal singular value lower bounded by a positive constant $\iota_2$
(depending on $t$).

\begin{proof}Let $B=B(L)$ be the submatrix of $\bc[t-1]$ guaranteed by  Lemma~\ref{l:c.submatrix}, so
	\[\|B\|_\infty
	\le \f{t\varsigma_t L}{(N\psi)^{1/2}}\,,\quad
	\|BB^\st-I\|\le\f1{4t}
	\,.\]
Then $B$ satisfies the conditions of Corollary~\ref{c:submatrix.APPROX}, so it has a $(t-1)\times(t-1)$ submatrix $U_1$ 
satisfying the determinant lower bound
\eqref{e:det.U.one.lbd}. Let $B_1$ be the matrix obtained by deleting $U_1$ from $B$. Then for all $N$ large enough we have
	\[
	\Big\|(B_1)(B_1)^\st-I_{t-1}\Big\|_\infty
	\le
	\f1{4t} + t\bigg(\f{t\varsigma_t L}{(N\psi)^{1/2}}\bigg)^2
	\le \f1{3t}\,.
	\]
Thus $B_1$ also satisfies the conditions of Corollary~\ref{c:submatrix.APPROX}, so it has a $(t-1)\times (t-1)$ submatrix $U_2$ which also satisfies the determinant lower bound \eqref{e:det.U.one.lbd}. Repeating the same argument, we see that with high probability the original matrix $B$ 
has disjoint $(t-1)\times(t-1)$ submatrices $U_1,\ldots,U_{\lfloor N^{0.9}\rfloor}$, all satisfying 
\eqref{e:det.U.one.lbd}. 
Recalling \eqref{e:gs.n.c.EXACT}, the corresponding 
submatrices of $\bn[t-1]$ are given by
$A_i\equiv(N\psi)^{1/2} \bGam_N U_i$, and
	\[
	|\det A_i|
	\ge \f{(N\psi)^{(t-1)/2} 
	|\det U_i|}{\varsigma_t}
	\ge 
	\bigg(\f{N\psi}{M}\bigg)^{(t-1)/2} 
	\f{((t-1)!)^{1/2}}{3\varsigma_t}
	\equiv\hat{\iota}_2\,.
	\]
Take any $A=A_i$, and denote its singular values
$\sigma_1\ge \ldots\ge \sigma_{t-1}\ge0$. 
Note that
$\sigma_1\le t\|A\|_\infty\le tL$, where the last bound holds since $A$ only involves columns of $\bn[t-1]$ indexed by $a\in M(L)$ (as in Lemma~\ref{l:c.submatrix}). Then
	\[
	\sigma_{t-1}
	\ge \f{|\det A|}{(\sigma_1)^{t-2}}
	\ge 
	\f{\hat{\iota}_2}{(tL)^{t-2}}
	\equiv \iota_2\,.
	\]
This concludes the proof.
\end{proof}
\end{cor}

\subsection{Fourier estimates at high frequency}
\label{ss:fourier.high.freq}

The main result of this subsection is the following lemma:

\begin{lem}[high-frequency estimate]
\label{l:fourier.THREE} Suppose $U$ satisfies Assumption~\ref{a:bdd} and \ref{a:Lip}. With the notation of \eqref{e:fourier.THREE},
it holds for any $\tau_{\max}<\infty$ and any
 $\epsilon_2>0$ that
	\[\max\bigg\{
	I_3(J,\tau,\epsilon_2)
	: J\in\set{-1,+1}^N,
	\|\pi(J)\|\le\f45,
	\|\tau\| \le \tau_{\max}
	\bigg\}
	\le \f1{\exp(N^{0.8})}
	\]
with probability $1-o_N(1)$.
\end{lem}

Towards the \hyperlink{proof:l.fourier.THREE}{proof of Lemma~\ref{l:fourier.THREE}}, 
recall that the random variable $\bze_a$ has density given by \eqref{e:bxi.a.density}. Thus
	\beq\label{e:relate.p.and.chi}
	\hat{p}_a(s)
	= \E\exp\bigg\{
		\ii s \Big(\bze_a-\E\bze_a\Big)\bigg\}
	=
	\f{\hat{\chi}_{\bX_a,c}(s)}{\exp(\ii s\E\bze_a)}\,.
	\eeq
We also denote $q_{x,c}(z)\equiv U(x+cz)\varphi(z)$, and note that
	\beq\label{e:relate.chi.and.q}
	\hat{\chi}_{x,c}(s)
	=\f{\hat{q}_{x,c}(s)}{\E_\xi U(x+c\xi)}
	=\f{\hat{q}_{x,c}(s)}{\hat{q}_{x,c}(0)}
	\,.\eeq
Note that Jensen's inequality implies
	\begin{align}\nonumber
	\Big| \hat{q}_{x,c}(s)
	-\hat{q}_{x',c'}(s)\Big|
	&=\bigg| \int e^{isz}
		\Big( U(x+cz)-U(x'+c'z)\Big)
		\varphi(z)\,dz\bigg|\\
	&\le
	\int \Big| U(x+cz)- U(x'+c'z)\Big|
	\varphi(z)\,dz\,,
	\label{e:cgf.jensen}\end{align}
and the last expression is bounded by Lemma~\ref{l:uniform.conv.ell.one}.

\begin{cor}\label{c:riemann.lebesgue}
Suppose $U$ satisfies Assumption~\ref{a:bdd}. 
Given any $\epsilon>0$ and any $L<\infty$, it is possible to choose $K$ large enough
(depending on $\epsilon$ and $L$) such that
	\[
	\sup\bigg\{ |\hat{\chi}_{x,c}(s)|
	: \f12\le c\le 2, |x|\le L, |s|\ge K\bigg\}
	 \le \epsilon\,,
	\]
where $\chi_{x,c}$ is as defined by \eqref{e:bxi.a.density}.

\begin{proof}
Recall from \eqref{e:relate.chi.and.q} the relation
	\[\hat{\chi}_{x,c}(s)
	=\f{\hat{q}_{x,c}(s)}{\E_\xi U(x+c\xi)}
	=\f{\hat{q}_{x,c}(s)}{\hat{q}_{x,c}(0)}\,.
	\]
By Assumption~\ref{a:bdd}, the denominator
$\hat{q}_{x,c}(0)=\E_\xi U(x+c\xi)$ is strictly positive for any given $x\in\R$, $c>0$. On the other hand, it follows from Lemma~\ref{l:uniform.conv.ell.one}
and \eqref{e:cgf.jensen} that 
$\hat{q}_{x,c}(0)$ is continuous in $(x,c)$. It follows that
	\[
	\inf\bigg\{
	\E_\xi U(x+c\xi)
	:\f12\le c\le 2, |x|\le L
	\bigg\}<\infty
	\]
for any finite $L$.  Therefore it suffices to show the claim with $q_{x,c}$ in place of $\chi_{x,c}$. By Lemma~\ref{l:uniform.conv.ell.one} again, given any $\epsilon>0$, we can choose $\eta'$ small enough such that
	\[
	\Big\|\hat{q}_{x,c}-\hat{q}_{x',c'}\Big\|_\infty
	\le
	\int \Big| U(x+cz)- U(x'+c'z)\Big|
	\varphi(z)\,dz
	\le \f{\epsilon}{2}
	\] 
as long as $c,c'\in[1/2,2]$, $x,x'\in[-L,L]$, and
$\max\set{|x-x'|,|c-c'|} \le \eta'$.
Let $\set{x_i}$ be a finite $\eta'$-net of $[-L,L]$, and let $\set{c_j}$ be a finite $\eta'$-net of $[1/2,2]$. It follows by the Riemann--Lebesgue lemma that there exists $K$ finite such that
	\[
	\sup\bigg\{ \max_{i,j} |\hat{q}_{x_i,c_j}(s)|
	: |s| \ge K\bigg\}
	\le \f{\epsilon}{2}\,.
	\]
For any $|x|\le L$ and $1/2\le c\le 2$,
we can find $x_i,c_j$ with $\max\set{|x-x_i|,|c-c_j|}\le\eta'$, so
	\[
	\sup\bigg\{
	\hat{q}_{x,c}(s)
	: |s|\ge K\bigg\} \le \epsilon
	\]
by combining the previous bounds. This concludes the proof.
\end{proof}
\end{cor}

\begin{cor}\label{c:uniform.cgf.bound}
Suppose $U$ satisfies Assumption~\ref{a:bdd}. Let
$\chi_{x,c}$ be as defined by \eqref{e:bxi.a.density}. Then
	\[\sup\bigg\{ |\hat{\chi}_{x,c}(s)|
	: \f12\le c\le 2, |x|\le L, |s|\ge 
		\epsilon\bigg\}
	\le 1-\epsilon' <1
	\]
for any finite $L$ and any $\epsilon>0$, where $\epsilon'$ is a small positive constant depending on $U$, $L$, and $\epsilon$.

\begin{proof}
By Lemma~\ref{c:riemann.lebesgue}, we can choose $K$ large enough such that
	\[\sup\bigg\{ |\hat{\chi}_{x,c}(s)|
	: \f12\le c\le 2, |x|\le L, |s|\ge K\bigg\}
	 \le\f12\,.\]
For any given $x,c$, let $\bze$ be a random variable with density $\chi_{x,c}$. For any $s\ne0$,
	\[
	|\chi_{x,c}(s)|
	=\bigg\{\Big( \E\cos(s\bze) \Big)^2
	+\Big( \E\sin(s\bze) \Big)^2\bigg\}^{1/2}
	<1
	\]
by Jensen's inequality. 
It follows from Lemma~\ref{l:uniform.conv.ell.one}
and \eqref{e:cgf.jensen}
that $\chi_{x,c}(s)$ is continuous in $(x,c,s)$, so
	\[\sup\bigg\{ |\hat{\chi}_{x,c}(s)|
	: \f12\le c\le 2, |x|\le L, \epsilon\le |s|\le K\bigg\}<1\]	
by compactness considerations.
The claim follows.
\end{proof}
\end{cor}

\begin{proof}[\hypertarget{proof:l.fourier.THREE}{Proof of Lemma~\ref{l:fourier.THREE}}]
For any subset of indices $T=\set{i(1),\ldots,i(t-1)}\subseteq[M]$ denote
	\[
	\phi_T(\sss)
	\equiv
	\prod_{\ell\le t-1} \hat{p}_{i(\ell)}(\sss_\ell)
	\]
for $\sss\in\R^{t-1}$. It follows from \eqref{e:relate.p.and.chi} and Plancherel's identity that
the $L^2$ norm of the function $\hat{p}_a(s)$
is the same as the $L^2$ norm of the function
$\chi_{\bX_a,c}(s)$ defined by \eqref{e:bxi.a.density}. 
We also note that Assumption~\ref{a:bdd} implies
	\[
	\|\chi_{x,c}\|^2
	=  \int \f{U(x+cz)^2 \varphi(z)^2}{(\E_\xi U(x+c\xi))^2}\,dz
	\le \f1{(1\pi)^{1/2}}
	\int \f{U(x+cz) \varphi(z)}{(\E_\xi U(x+c\xi))^2}\,dz
	= \f1{\E_\xi U(x+c\xi)}\,.
	\]
By compactness considerations (similarly as for \eqref{e:bar.c.one.U}), we must have
	\[
	\inf
	\bigg\{
	\E_\xi U(x+c\xi)
	: \f12 \le c\le 2, |x| \le L'
	\bigg\} \ge \bar{c}_1(U,L')\,.
	\]
If $T\subseteq M(L)$ (as defined by Lemma~\ref{l:c.submatrix}),
then it follows by combining the above with \eqref{e:bound.X.a.on.M.L} that
	\[
	\|\phi_T\|_2
	= \prod_{\ell\le t-1}
		\|p_{i(\ell)}\|_2
	\le
	\bigg(
	\sup\bigg\{
	\|\chi_{x,c}\|_2
	: \f12 \le c\le 2, |x|\le L
	\bigg\} \bigg)^{t-1}
	\le \bigg(\f1{\bar{c}_1(U,L')}\bigg)^t
	\equiv \wp_5 \,.
	\]
Now let $A_1,\ldots,A_{\lfloor N^{0.9}\rfloor}$ be the submatrices of $\bn[t-1]$ guaranteed
(with high probability) by Corollary~\ref{c:n.submatrices}. Let $T_i$ denote the subset of column indices involved in $A_i$, and note
	\[
	|\hat{p}(\sss)|
	\le
	\prod_{i\le \lfloor N^{0.9}\rfloor}
	\bigg|
	\phi_{T_i}\bigg(\f{(A_i)^\st\sss}{N^{1/2}}\bigg)
	\bigg|\,.
	\]
Moreover, each individual factor $\phi_{T_i}$ has modulus at most one. Combining with the preceding $L^2$ bound gives
	\[\int \bigg|\phi_{T_i}
	\bigg(\f{(A_i)^\st\sss}{N^{1/2}}\bigg)
		\bigg|^2\,
		d\sss
	= \f{N^{(t-1)/2}(\|\phi_{T_i}\|_2)^2}{|\det A_i|}
	\le \f{N^{(t-1)/2}(\wp_5)^2}{(\iota_2)^{t-1}}
	\,.
	\]
It follows using the Cauchy--Schwarz inequality that
	\[
	\int \bigg| 
		\phi_{T_1}\bigg(\f{(A_1)^\st\sss}{N^{1/2}}\bigg)
		\phi_{T_2}\bigg(\f{(A_2)^\st\sss}{N^{1/2}}\bigg)
		\bigg|\,d\sss
	\le\f{N^{(t-1)/2}(\wp_5)^2}{(\iota_2)^{t-1}}\,.
	\]
On the other hand, if $\|\sss\|\ge\epsilon_2N^{1/2}$, then
the least singular value bound from Corollary~\ref{c:n.submatrices} implies
	\[
	\max\bigg\{\f{|(\sss,\bn_a)|}{N^{1/2}}
	: a\in T_i\bigg\}
	= 
	\f{\|(A_i)^\st\sss\|_\infty}{N^{1/2}}
	\ge
	\f{\|(A_i)^\st\sss\|}{(Nt)^{1/2}}
	\ge 
	\f{\iota_2\epsilon_2}{t^{1/2}}\,.
	\]
Recall again that for $a\in M(L)$, $|\bX_a|$ is bounded by \eqref{e:bound.X.a.on.M.L}. Combining with the result of
Corollary~\ref{c:uniform.cgf.bound} gives
	\[
	\bigg| \phi_{T_i}\bigg(\f{(A_i)^\st\sss}{N^{1/2}}\bigg)\bigg|
	\le\sup\bigg\{
	|\hat{\chi}_{x,c}(s)| :
	\f12 \le c\le 2,
	|x|\le L', 
	|s|\ge \f{\iota_2\epsilon_2}{t^{1/2}}
	\bigg\} \le 1-\epsilon'<1\,.
	\]
To conclude we note that the quantity 
$I_3(J,\tau,\epsilon_2)$ from \eqref{e:fourier.THREE}
can be bounded by $I_{3,g}+I_{3,p}$ 
where $I_{3,g}$ is the integral of $\hat{g}_{J,\tau}$,
while $I_{3,p}$ is the integral of $\hat{p}_{J,\tau}$.
By \eqref{e:hat.g} and Lemma~\ref{l:nondegenerate.Sigma}, we have
with high probability
	\[I_{3,g}
	\equiv \int \Big| \hat{g}_{J,\tau}(\sss)\Big| 
		\I\Big\{ \|\sss\| \ge 
		\epsilon_2 N^{1/2}
	\Big\} \,d\sss
	\le  \f1{\exp(N^{0.9})}\,.
	\]
By the previous calculations, we also have with high probability
	\begin{align*}
	I_{3,p}
	&\le
	\bigg\{ \int\bigg| \prod_{i=1,2}
		 \phi_{T_i}\bigg(\f{(A_i)^\st\sss}{N^{1/2}}\bigg)
		 \bigg|
		 \,d\sss\bigg\}
	\cdot
	\sup\bigg\{
	\prod_{i=3}^{N^{0.9}}
	\bigg|
	\phi_{T_i}\bigg(\f{(A_i)^\st\sss}{N^{1/2}}\bigg)
	\bigg|
	: |\sss| \ge \epsilon_2N^{1/2}
	\bigg\}\\
	&\le \f{N^{(t-1)/2}(\wp_5)^2}{(\iota_2)^{t-1}}
	(1-\epsilon')^{N^{0.85}}
	\le \f1{\exp(N^{0.8})}\,.
	\end{align*}
This concludes the proof.
\end{proof}

\subsection{Conclusion of local CLT}
\label{ss:local.clt.conclusion}

In this concluding subsection we \hyperlink{proof:p.local.CLT}{prove the local CLT  
Proposition~\ref{p:local.CLT}}, and apply it to deduce
 Propositions~\ref{p:density.bound} and \ref{p:pair.density.bound}.

\begin{proof}[\hypertarget{proof:p.local.CLT}{Proof of Proposition~\ref{p:local.CLT}}]
Recall that $\hat{p}_{J,\tau}$ and $\hat{g}_{J,\tau}$ are defined by \eqref{e:hat.p} and \eqref{e:hat.g}. It follows by combining Lemmas 
\ref{l:fourier.ONE},
\ref{l:fourier.TWO}, and 
\ref{l:fourier.THREE} that
for any finite constant $\tau_{\max}$, we have
	\[
	\sup\bigg\{
	\int\Big| \hat{p}_{J,\tau}(\sss)
		-\hat{g}_{J,\tau}(\sss)\Big|\,d\sss
	: J\in\set{-1,+1}^N,
	\|\pi(J)\|\le\f45,
	\|\tau\|\le \tau_{\max}
	\bigg\}
	\le \f1{N^{0.35}} 
	\]
with high probability. Inverting the Fourier transform shows that, with high probability, the random variable $\bW$ from \eqref{e:local.clt.W} has a bounded continuous density function $p_{J,\tau}$, which satisfies
	\[\sup\bigg\{
	\|p_{J,\tau}-g_{J,\tau}\|_\infty:
	J\in\set{-1,+1}^N,
	\|\pi(J)\|\le\f45,
	\|\tau\|\le \tau_{\max}
	\bigg\}
	\le \f1{(2\pi)^{t-1}N^{0.35}}
	\le \f1{N^{0.3}}\,,
	\]
as claimed.
\end{proof}

We now define the transformed gaussian density
	\beq\label{e:CLT.transformed}
	\bg_{J,\tau}(z)
	\equiv
	\psi^{1/2} |\det \bGam_N|
	g_{J,\tau}
		\bigg(\psi^{1/2}\bGam_N (z- N^{1/2}\tau)
			-\f{\bn[t-1]\E\bze}{N^{1/2}}\bigg)\,,
	\eeq
where $\E\bze$ is as in \eqref{e:adm.zeta}.

\begin{ppn}[density bound for first moment]
\label{p:density.bound}
Suppose $U$ satisfies Assumptions~\ref{a:bdd} and \ref{a:Lip}. Then we have
	\[
	\sup\bigg\{
	\Big\|\ap_{J,\tau}(\cdot\,|\,\bar{g}_\ROW)
	-	\bg_{J,\tau}(\cdot)\Big\|_\infty
	: J\in\set{-1,+1}^N,
	\|\tau\|\le\tau_{\max}
	\bigg\}
	\le \f1{N^{0.25}}
	\]
with high probability,
where $\ap_{J,\tau}(\cdot\,|\,\bar{g}_\ROW)$ is as in \eqref{e:ap.tau}, while $\bg_{J,\tau}$ is as in \eqref{e:CLT.transformed}.

\begin{proof} Recall that 
Proposition~\ref{p:local.CLT} above estimates the density $p_{J,\tau}$ of the random variable $\bW$ from \eqref{e:local.clt.W},
	\beq\label{e:local.clt.W.repeat}
	\bW
	=\f{\bn[t-1](\bze-\E\bze)}{N^{1/2}}
	\in \R^{t-1}\,,
	\eeq
where each $\bze_a$ has density given by \eqref{e:bxi.a.density}.
On the other hand, let $\bxi\in\R^M$ be a random vector with independent coordinates, such that $\bxi_a$ has density
	\[
	\tilde{p}_a(z)\cong
	U\Big( (\tbX_J)_a + cz\Big)
	\exp\bigg\{ N^{1/2} 
		\tau^\st\bc[t-1] \eM_a
		z\bigg\}
	\varphi(z)\,,\]
where $c=c(\pi(J))$,
$\tbX_J$ is as in Lemma~\ref{l:G.R.J}, and $\cong$ denotes equality up to a normalizing constant. We see from
\eqref{e:ap.tau} that $\ap_{J,\tau}(\cdot\,|\,\bar{g}_\ROW)$ is the density of the random variable $\bc[t-1]\bxi$, for $\bxi$ as we have just described. Note that
	\begin{align*}
	&\tilde{p}_a\bigg(
	z + N^{1/2}\tau^\st\bc[t-1] \eM_a
	\bigg)
	\cong
	U\bigg( (\tbX_J)_a + c\Big\{
	z + N^{1/2}\tau^\st\bc[t-1] \eM_a
	\Big\}\bigg)
	\varphi(z) \\
	&\qquad\stackrel{\eqref{e:bX.J.tau.rewrite}}{=}
	U\Big((\bX_{J,\tau})_a
	+cz\Big)\varphi(z)
	\stackrel{\eqref{e:bxi.a.density}}{\cong}
	\chi_{\bX_a,c}(z)\,,
	\end{align*}
so it follows that $\bxi-N^{1/2}\bc[t-1]^\st\tau$ is equidistributed as $\bze$ for $\bze$ as in \eqref{e:local.clt.W.repeat}. Thus
$\ap_{J,\tau}(\cdot\,|\,\bar{g}_\ROW)$ is the same as the density of
	\begin{align*}
	\bc[t-1]
	\Big(\bze+ N^{1/2}\bc[t-1]^\st\tau\Big)
	&\stackrel{\eqref{e:gs.n.c.EXACT}}{=}
	\f{(\bGam_N)^{-1}\bn[t-1]\bze}{(N\psi)^{1/2}}
		+ N^{1/2}\tau\\
	&\stackrel{\eqref{e:ap.tau}}{=}
	\f{(\bGam_N)^{-1}\bW}{\psi^{1/2}}
	+\f{(\bGam_N)^{-1}\bn[t-1]\E\bze}{(N\psi)^{1/2}}
	+N^{1/2}\tau\,.
	\end{align*}
It follows by making a change of variables that
	\[
	\ap_{J,\tau}(z\,|\,\bar{g}_\ROW)
	=\psi^{1/2} |\det \bGam_N|
	p_{J,\tau}
		\bigg(\psi^{1/2}\bGam_N (z- N^{1/2}\tau)
			-\f{\bn[t-1]\E\bze}{N^{1/2}}
			\bigg)\,.
	\]
Comparing with \eqref{e:CLT.transformed}, we have
	\[
	\Big\|\ap_{J,\tau}(\cdot\,|\,\bar{g}_\ROW)
	-	\bg_{J,\tau}(\cdot)\Big\|_\infty
	= \psi^{1/2} |\det \bGam_N|
	\Big\| p_{J,\tau}-g_{J,\tau} \Big\|_\infty\,,
	\]
so the result follows from Proposition~\ref{p:local.CLT}.
\end{proof}
\end{ppn}

\begin{ppn}[density bound for second moment]
\label{p:pair.density.bound}
Suppose $U$ satisfies Assumptions~\ref{a:bdd} and \ref{a:Lip}. Then
the bound \eqref{e:pair.density.bound} holds
with high probability,
where $\ap_{K|J,\tau}(\cdot\,|\,\bar{g}_\ROW,\bze)$ is as in
\eqref{e:ap.tau.two}.

\begin{proof} Through we abbreviate $c=c(\pi(K))$.
First we slightly modify 
the definition from
\eqref{e:bxi.a.density}: let $\bsi\in\R^M$ be a random vector with independent coordinates, such that each $\bsi_a$
has density given by
$\chi_{\bX_a,e(\lm)}$ for
	\[\bX=\bX_{K|J,\tau}(\bze)
	\equiv\tbX_K + c\cdot \bigg( \lm\bze
	+(1-\lm^2)^{1/2} N^{1/2} \bc[t-1]^\st\tau
	\bigg)
	\]
and $e(\lm)= c\cdot(1-\lm^2)^{1/2}$.
In this definition, $\tbX_K$ is as in Lemma~\ref{l:pair.adm.cgf}, and $\lm=\lm(J,K)$. We define also (cf.\ \eqref{e:local.clt.W})
	\beq\label{e:local.clt.W.modified}
	\bW'
	\equiv
	 \f{\bn[t-1](\bsi-\E\bsi)}{N^{1/2}}
	 \in\R^{t-1}\,.\eeq
On the other hand, let $\bxi\in\R^M$ be a random vector with independent coordinates, such that $\xi_a$ has density
	\[
	\tilde{p}_a(z)
	\cong
	U\bigg(
	(\tbX_K)_a + c\Big(\lm
		\bze_a + (1-\lm^2)^{1/2}z\Big)
	\bigg) \exp\bigg\{
	N^{1/2}\tau^\st \bc[t-1] \eM_a z
	\bigg\}\varphi(z)\,.
	\]
We see from \eqref{e:ap.tau.two} that
$\ap_{K|J,\tau}(\cdot\,|\,\bar{g}_\ROW,\bze)$
is the density of the random variable $\bc[t-1]\bxi$. Note that
	\begin{align*}
	\tilde{p}_a\bigg(z+N^{1/2} \tau^\st \bc[t-1]\eM_a\bigg)
	&\cong
	U\bigg( \bX_{K|J,\tau}(\bze)
	+c \cdot (1-\lm^2)^{1/2}z
	\bigg) \varphi(z)
	\cong \chi_{\bX_a,e(\lm)}\,,
	\end{align*}
which implies that
$\bxi-N^{1/2}\bc[t-1]^\st\tau$
is equidistributed as $\bsi$.
Thus $\ap_{K|J,\tau}(\cdot\,|\,\bar{g}_\ROW,\bze)$ is the same as the density of
	\[\bc[t-1]\Big(
	\bsi+N^{1/2}\bc[t-1]^\st\tau\Big)
	= \f{(\bGam_N)^{-1}\bW'}{\psi^{1/2}}
	+\f{(\bGam_N)^{-1}\bn[t-1]\E\bsi}{(N\psi)^{1/2}}
	+N^{1/2}\tau
	\]
for $\bsi$ and $\bW'$ as defined above. It follows by a minor modification of Proposition~\ref{p:local.CLT} (replacing $\bW$ 
from \eqref{e:local.clt.W}
with $\bW'$ from \eqref{e:local.clt.W.modified}) 
that $\ap_{K|J,\tau}(\cdot\,|\,\bar{g}_\ROW,\bze)$ can be uniformly approximated by a gaussian density. The claim follows.
\end{proof}
\end{ppn}

\section{Concentration of partition function}\label{s:conc.log.Z}

In this section we prove Propositions~\ref{p:poly.conc}, \ref{p:Z.eps.approx.Z}, and \ref{p:conc.smoothed}; and use these to \hyperlink{proof:t.main.conclusion}{conclude the proof of Theorem~\ref{t:main}}. The section is organized as follows:

\begin{itemize}
\item As commented earlier, both Propositions~\ref{p:poly.conc} and \ref{p:Z.eps.approx.Z} rely on a bound for near-isotropic gaussian processes, Proposition~\ref{p:gaus.poly}, which is proved in \S\ref{ss:near.isotropic}.
See Remark~\ref{r:gaus.poly} for further discussion of this result. 
\item In \S\ref{ss:poly.conc} we give the
	 \hyperlink{proof:p.poly.conc}{proof of Proposition~\ref{p:poly.conc}}.
\item In \S\ref{ss:approx.arg} we give the 
	\hyperlink{proof:p.conc.smoothed}{proof of Proposition~\ref{p:conc.smoothed}},
and use this to deduce that the free energy of the smoothed model \eqref{e:tilde.Z} is given by the replica symmetric formula
(Corollary~\ref{c:smoothed.free.energy}).

\item In \S\ref{ss:conclusion.main.thm} we 
give the 
	\hyperlink{proof:p.Z.L.approx.Z}{proof of Proposition~\ref{p:Z.eps.approx.Z}},
	and
	\hyperlink{proof:t.main.conclusion}{conclude the proof of Theorem~\ref{t:main}}.
\end{itemize}
Recall from \S\ref{ss:intro.unbdd} that Assumption~\ref{a:bdd} implies
\eqref{e:positive.measure.set}, where we can assume without loss that $E(U)\subseteq[-E_{\max}(U),E_{\max}(U)]$ for some finite $E_{\max}(U)$.

\subsection{Bounds for near-isotropic gaussian processes}\label{ss:near.isotropic}

The following is a variant of
\cite[Cor.~8.2.5]{MR3024566}:

\begin{ppn}\label{p:gaus.poly} Let $c\in(0,1/12]$. Let $\bv^1,\ldots,\bv^n$ be unit vectors in $\R^n$ such that $(\bv^i,\bv^j)\le c$ for all $i \ne j$.  Then
	\[
	\P\bigg(
	\f1n
	\Big|\Big\{i \le n : (\bg,\bv^i) \in E(U) \Big\}
	\Big| \le \gamma
	\bigg) 
	\le 
	\gamma^{1/(25c)}
	\]
for all $\log(5/c)/(\log n)\le\gamma\le \gamma_0 = \gamma_0(|E(U)|,E_{\max}(U))$
and $n$ large enough.
\end{ppn}

\begin{rmk}\label{r:gaus.poly}
We point out that there are two main differences between \cite[Cor.~8.2.5]{MR3024566} and Proposition~\ref{p:gaus.poly}. First, \cite[Cor.~8.2.5]{MR3024566} considers the event
$\set{(\bg,\bv^i)\ge a}$, and
the proof relies crucially on Gordon's inequality. By contrast,  Proposition~\ref{p:gaus.poly} considers the event
$\set{(\bg,\bv^i)\in E(U)}$, where it does not seem possible to apply standard gaussian comparison inequalities. As a result we rely on more ad hoc
arguments which yield a weaker bound, in the sense that \cite[Cor.~8.2.5]{MR3024566} holds for $\gamma$ polynomially small in $n$ while 
Proposition~\ref{p:gaus.poly}  holds only for $\gamma$ decaying logarithmically in $n$.
\end{rmk}

The \hyperlink{proof:p.gaus.poly}{proof of Proposition~\ref{p:gaus.poly}} is given at the end of this subsection. We begin with some preparatory lemmas:

\begin{lem}[used in proof of Lemma~\ref{l:gs.reindex}]\label{l:K.inds.close}
Let $c\in(0,1)$ and denote $\eta'(c)=1/\log(4/c)$.
For any $K\in\mathbb{N}$ there exists
$n_0(c,K)<\infty$ such that the following holds for all $n\ge n_0(c,K)$: if $\bv^1,\ldots,\bv^n$ are unit vectors in $\R^n$ and $m\le \eta'(c) \log n$, then there must exist $K$ distinct indices $m<i_1<\ldots<i_K\le n$ such that
	\[
	\max \bigg\{
	\Big\|P_m\Big(\bv^{i_a}-\bv^{i_b}\Big)\Big\|: a,b\le K
	\bigg\}
	\le c\,,
	\]
where $P_m$ denotes the orthogonal projection onto the span of $\set{\bv^1,\ldots,\bv^m}$.

\begin{proof}
Suppose for contradiction that for all $m<i_1<\ldots<i_K\le n$ we have
	\beq\label{e:K.inds.contradict}
	\max \bigg\{
	\Big\|P_m\Big(\bv^{i_a}-\bv^{i_b}\Big)\Big\|: a,b\le K
	\bigg\}>c\,.\eeq
Let $U$ denote the \emph{disjoint} union of $U_{m+1},\ldots,U_n$, where $U_i$ is a \emph{copy} of
	\[B_m\bigg(P_m\bv^i, \f{c}{2}\bigg)
	\equiv\bigg(
	\spn\Big\{\bv^1,\ldots,\bv^m\Big\}\bigg)
	\cap B\bigg(P_m\bv^i, \f{c}{2}\bigg)\,.
	\]
Note that if $x\in B(P_m\bv^i,c/2)$ then
$\|x\| \le \|P_m\bv^i\| + c/2 \le 3/2$, so we have a natural mapping $i : U \to B_m(\textbf{0},3/2)$.
By the assumption \eqref{e:K.inds.contradict},
each point in $B_m(\textbf{0},3/2)$ has at most $K-1$ distinct preimages under the mapping $i$, so
	\[
	(n-m)\vol  B\bigg(P_m\bv^i,\f{c}{2}\bigg) = \vol U
	\le (K-1) \vol B_m\bigg(\textbf{0},\f32\bigg)\,.
	\]
If $m'=\dim\spn\set{\bv^1,\ldots,\bv^m}\le m$, then it follows that
	\[
	n-m \le (K-1) \bigg( \f{3/2}{c/2}\bigg)^{m'}
	\le K \bigg(\f{3}{c}\bigg)^{\eta'(c)\log n}
	= K \exp\bigg\{ \f{\log(3/c)}{\log(4/c)} \log n\bigg\}\,,
	\]
which yields a contradiction for $n$ large enough
(depending on $c$ and $K$).
\end{proof}
\end{lem}

\begin{lem}\label{l:gs.reindex}
Let $c\in(0,1)$ and denote $\eta'(c)=1/\log(4/c)$.
There exists $n_0(c)<\infty$ such that the following holds for all $n\ge n_0(c)$:
if $\bv^1,\ldots,\bv^n$ are unit vectors in $\R^n$ with $(\bv^i,\bv^j)\le c$ for all $i\ne j$, then the vectors can be re-indexed in such a way that
	\[
	\max\bigg\{ 
	\Big\|P_m\bv^{m+1}\Big\|
	: 1\le m\le\eta'(c)\log n\bigg\}
	\le(3c)^{1/2}\,,
	\]
where $P_m$ denotes the orthogonal projection onto the span of $\set{\bv^1,\ldots,\bv^m}$.
(The claim is non-trivial only if $c<1/3$.)

\begin{proof} We shall assume the vectors are indexed such that for all $1\le\ell\le n$ we have
	\beq\label{e:reindex}
	\Big\|P_{\ell-1}\bv^\ell\Big\|
	= \min \bigg\{ \Big\|P_{\ell-1}\bv^k\Big\|
		: \ell\le k \le n\bigg\}\,.
	\eeq
Now suppose for the sake of contradiction that for some $m\le \eta'(c)\log n$ we have 
	\beq\label{e:m.contradiction}
	\Big\|P_m\bv^{m+1}\Big\|
	\stackrel{\eqref{e:reindex}}{=}
	\min \bigg\{ \Big\|P_m\bv^k\Big\|
		: m+1 \le k \le n\bigg\}>(3c)^{1/2}\,.
	\eeq
Take $K=2+\lceil 1/c\rceil$. By Lemma~\ref{l:K.inds.close},
for all $n$ large enough we can
 find indices
$m<i_1<\ldots<i_K\le n$ such that
	\beq\label{e:K.inds.close.app}\max \bigg\{
	\Big\|P_m\Big(\bv^{i_a}-\bv^{i_b}\Big)\Big\|: a,b\le K
	\bigg\} \le c\,.\eeq
As a consequence, for any $a\ne b$ where $a,b\le K$, we have
	\begin{align*}
	\bigg(
	(I-P_m)\bv^{i_a},(I-P_m)\bv^{i_b}
	\bigg)
	&= (\bv^{i_a},\bv^{i_b})-(P_m\bv^{i_a},P_m\bv^{i_b})
	\\
	&= -\|P_m\bv^{i_a}\|^2
	+\bigg\{ (\bv^{i_a},\bv^{i_b})
	- \Big(P_m\bv^{i_a},P_m(\bv^{i_b}-\bv^{i_a})\Big) 
	\bigg\} \le -c\,,
	\end{align*}
where the last bound uses
\eqref{e:m.contradiction}, \eqref{e:K.inds.close.app},
and the assumption that $(\bv^i,\bv^j)\le c$ for all $i\ne j$. If we let 
	\[
	\bx^a\equiv \f{(I-P_m)\bv^{i_a}}
		{\|(I-P_m)\bv^{i_a}\|}\,,\]
then the above implies that
$(\bx^a,\bx^b) \le -c$ for all $a\ne b$. It follows that
	\[
	0\le \bigg\|\sum_{a\le K} \bx^a\bigg\|^2
	=\sum_{a,b\le K}(\bx^a,\bx^b)
	\le K\Big(1 - c(K-1)\Big)\,,
	\]
which gives a contradiction since we chose $K\ge 2+1/c$.
\end{proof}
\end{lem}

\begin{lem}\label{l:gs.poly.bd}
Let $\bv^1,\ldots,\bv^m$ be unit vectors in $\R^n$
(for any $m,n$) such that
	\[\max\bigg\{
	\Big\|P_{\ell-1}\bv^\ell\Big\|
	 : 1\le\ell\le m\bigg\}
	\le c' \le \f12\,,\]
where $P_{\ell-1}$ denotes the orthogonal projection onto the span of $\set{\bv^1,\ldots,\bv^{\ell-1}}$. Let $\bg$ be a standard gaussian random vector in $\R^n$.
There exists $\gamma_0=\gamma_0(|E(U)|,E_{\max}(U))>0$ such that
	\[
	\P\bigg(
	\f1m
	\Big|\Big\{i \le m : (\bg,\bv^i) \in E(U)\Big\}
	\Big| \le \gamma
	\bigg) 
	\le 
	\gamma^{1/(8(c')^2)}
	\]
for all $1/m\le\gamma\le\gamma_0$.

\begin{proof}
We shall assume without loss that $m\gamma$ is integer-valued. Let $u_i \equiv (\bg,\bv^i)$, so that $(u_i)$ defines a (centered) gaussian random vector indexed by $i\le n$. For each $i$ we can decompose $u_i\equiv \zeta_i + \xi_i$ where $\zeta_i\equiv(\bg,P_{i-1}\bv^i)$; at the first step $\zeta_1=0$. Define a parameter
	\beq\label{e:s.a.b}
	s \equiv s(U)  \le \max\bigg\{10, E_{\max}(U),
	\bigg( \Big|\log |E(U)| \Big|\bigg)^{1/2}
	\bigg\}\,,\eeq
and define the random subset of indices
	\[B \equiv \bigg\{i\le m 
		: |\zeta_i| \le s\bigg\}\,.
	\]
Let $\Omega_\gamma$ denote the event of interest,
	\[\Omega_\gamma
	\equiv\bigg\{ 
	\f1m
	\Big|\Big\{i \le m : u_i\in E(U)\Big\}
	\Big| \le \gamma\bigg\}\,.
	\]
On the event $\Omega_\gamma$ there must be a subset $A\subseteq[m]$ of size $m\gamma$ such that $u_i\notin E(U)$ for all $i\notin A$. Therefore
	\beq\label{e:fix.coords.in.a.b}
	\P(\Omega_\gamma)
	\le \P\bigg(|B|\le \f{m}{2}\bigg)
	+ \sum_{|A|=m\gamma}
	\P\bigg( u_i\notin E(U)
		\ \forall i\notin A; |B|>\f{m}{2}\bigg)\,. 
	\eeq
To bound the above we will consider a fixed subset $A$, without loss $A=\set{m-m\gamma+1,\ldots,m}$. Define
	\[
	\mathscr{G}_\ell
	\equiv \sigma\bigg(
	(\zeta_i, \xi_i) : 1\le i\le \ell\bigg)\,.
	\]
Let $\tau_0\equiv0$ and define the increasing sequence
	\[
	\tau_\ell
	\equiv \inf\bigg\{ i >\tau_{\ell-1}:
		i\le m,|\zeta_i| \le s \bigg\}\,.
	\]
Note that since $\zeta_\ell\in\mathscr{G}_{\ell-1}$,
the $\tau_\ell$ are stopping times with respect to the filtration $\mathscr{G}_\ell$. We take the usual convention that $\inf\varnothing\equiv\infty$, so the set of finite stopping times corresponds exactly to the set $B$. Let $f(i)\equiv \Ind{u_i\notin E(U)}$. It follows from the assumption that $\xi_i$ has the law of a gaussian random variable which is independent of $\mathscr{G}_i$, and has variance between $1-(c')^2\ge3/4$ and $1$. Therefore we have
	\begin{align*}
	p_\ell&\equiv\E\bigg( \Ind{\tau_\ell<\infty} f(\tau_\ell)
		\,\bigg|\, \mathscr{G}_{\tau_\ell-1}\bigg)
	=\Ind{\tau_\ell<\infty}
	\P\bigg( u_{\tau_\ell}=\zeta_{\tau_\ell}+\xi_{\tau_\ell}
		\notin E(U)
		\,\bigg|\, \mathscr{G}_{\tau_\ell-1}\bigg)  \\
	&
	\le \max
	\bigg\{ \P\bigg(  Z \notin \f{E(U)-x}{\lambda} \bigg)
	: \bigg(\f34\bigg)^{1/2} \le \lambda \le1,
	|x|\le s
	\bigg\}\,.
	\end{align*}
To bound the above, note that the 
set $\lm^{-1}(E(U)-x)$ has Lebesgue measure at least $|E(U)|$ (since $\lm\le1$), and is contained in the interval $[-5s/2,5s/2]$ (by the assumption $s\ge E_{\max}(U)$ from
\eqref{e:s.a.b}, together with the restriction $\lambda\ge(3/4)^{1/2}$). It follows that
	\[
	p_\ell
	\le 1- |E(U)| \varphi\bigg(\f{5s}{2}\bigg)
	\le 1- \f1{(2\pi)^{1/2}}
		\exp\bigg\{-\f{7s^2}{2}\bigg\}\,,
	\]
where the last bound uses the assumption $s^2\ge|\log|E(U)||$ from \eqref{e:s.a.b}. It then follows by iterated expectations that
	\begin{align*}
	&\P\bigg( u_i\notin E(U)
	\ \forall i \notin A; 
	|B|>\f{m}{2}
	\bigg)
	\le
	\E\bigg[ \prod_{j\le m/2} 
		\Ind{\tau_j<\infty} f(\tau_j)\bigg]
	\\
	&\qquad\le
	\E\bigg[
	\bigg( \prod_{j\le m/2-1}
	\Ind{\tau_j<\infty} f(\tau_j)\bigg)
	\E\bigg(
	\Ind{\tau_{\lceil m/2\rceil}<\infty}
	f(\tau_{\lceil m/2\rceil})
	\,\bigg|\,
	\mathscr{G}_{\tau_{\lceil m/2\rceil}-1}\bigg)\bigg]\\
	&\qquad\le
	 \bigg( 
	 1- \f1{(2\pi)^{1/2}}
		\exp\bigg\{-\f{7s^2}{2}\bigg\}
	 \bigg)^{m/2}
	\le
	\exp\bigg\{
	-\f{m \exp(-7s^2/2)}{2(2\pi)^{1/2}}
	\bigg\}
	\,.
	\end{align*}
Substituting this bound into 
\eqref{e:fix.coords.in.a.b} and accounting for the number of choices of $A$ gives
	\[
	\P(\Omega_\gamma)
	\le \P\bigg(|B|\le \f{m}{2}\bigg)
	+
	\exp\bigg\{ m
	\bigg[\mathcal{H}(\gamma) 
	- \f{\exp(-7s^2/2)}{2(2\pi)^{1/2}}
		\bigg]\bigg\}\,,
	\]
where $\mathcal{H}$ denotes the binary entropy function, and satisfies $\mathcal{H}(\gamma)\le\gamma\log(e/\gamma)$. If we take $\gamma=\exp(-4s^2)$, then
	\[	
	\mathcal{H}(\gamma) 
	- \f{\exp(-7s^2/2)}{2(2\pi)^{1/2}}
	\le \f1{\exp(7s^2/2)}
	\bigg(\f{1+4s^2}{\exp(s^2/2)}
	-\f1{2(2\pi)^{1/2}} \bigg)
	\le \f{-1}{6\exp(7s^2/2)}\,,
	\]
where the last bound uses the assumption $s\ge10$ from \eqref{e:s.a.b}. It follows that
	\beq\label{e:p.Omega.gamma.simp}
	\P(\Omega_\gamma)
	\le \P\bigg(|B|\le \f{m}{2}\bigg)
	+ \exp\bigg\{-\f{m}{6\exp(7s^2/2)}\bigg\}\,,
	\eeq
and it remains to bound the probability that $|B|\le m/2$. To this end, note each $\zeta_i$ is a gaussian random variable with variance at most $(c')^2$, so
	\[
	\P(|\zeta_i| \ge s)
	\le \P( c'|Z| \ge s)
	\le \f{c'}{s}
	\exp\bigg\{ -\f{s^2}{2(c')^2} \bigg\}\,.
	\]
It follows by Markov's inequality and the preceding bound that
	\begin{align*}
	\P\bigg(|B| \le \f{m}{2}\bigg)
	&=\P\bigg( |B^c|
	\ge \f{m}{2} 
	\bigg)
	\le 2\max\bigg\{ \P(|\zeta_i| \ge s)
		: i\le m\bigg\}\\
	&\le  \f{2c'}{s}
	\exp\bigg\{ -\f{s^2}{2(c')^2} \bigg\}
	\le \f12 \exp\bigg\{ -\f{s^2}{2(c')^2} \bigg\}
	\,,
	\end{align*}
where the last bound follows trivially from the bounds $c'\le1$ and $s\ge10$ (from \eqref{e:s.a.b}). If $m\ge1/\gamma$, then
	\[
	\f{s^2}{2} \cdot 6\exp(7s^2/2)
	= \f{3s^2}{\exp(s^2/2)} \cdot \f1\gamma
	\le \f1\gamma \le m\,,
	\]
so that \eqref{e:p.Omega.gamma.simp} is dominated by the first term. It follows that
	\[
	\P(\Omega_\gamma)
	\le
	\exp\bigg\{ -\f{s^2}{2(c')^2} \bigg\}
	=\gamma^{1/(8(c')^2)}\,,
	\]
provided $\gamma=\exp(-4s^2)$ for $s$ satisfying \eqref{e:s.a.b}, and $m\ge1/\gamma$. This concludes the proof.
\end{proof}
\end{lem}

\begin{proof}[\hypertarget{proof:p.gaus.poly}{Proof of Proposition~\ref{p:gaus.poly}}]
As in Lemma~\ref{l:gs.reindex}, let $\eta'(c)=1/\log(4/c)$. Let 
	\[
	m=\bigg\lfloor \f12 \eta'(c)\log n \bigg\rfloor\,,\quad
	L =\bigg\lfloor \f{n-n^{1/2}}{m} \bigg\rfloor\,.
	\]
By repeatedly applying Lemma~\ref{l:gs.reindex}, we see that there exists a re-indexing of $\bv^1,\ldots,\bv^n$ such that
	\[
	\max\bigg\{
	\Big\|P_{\ell m,i-1}\bv_{\ell m +i}\Big\|^2
	: 0\le\ell \le L-1,1\le i\le m
	\bigg\}
	\le (3c)^{1/2}
	\equiv c' \le\f12\,,
	\]
where $P_{\ell m,i-1}$ denotes the orthogonal projection onto the span of $\set{\bv^{\ell m+1},\ldots,\bv^{\ell m+i-1}}$. Let
	\[
	N_\ell
	\equiv \bigg|\Big\{ 1\le i\le m
	:(\bg,\bv^{\ell m+i})\in E(U)
	\Big\}\bigg|\,.
	\]
Note that if $N_\ell\ge 2m\gamma$ for at least $n/(2m)$ indices $0\le\ell\le L-1$, then we will have $(\bg,\bv^i)\in E(U)$ for at least $n\gamma$ indices $1\le i\le n$. It follows by combining with Markov's inequality that
	\[
	\P\bigg(
	\f1n
	\Big|\Big\{i \le n : (\bg,\bv^i) \in E(U)\Big\}
	\Big| \le \gamma
	\bigg) 
	\le\P\bigg(
	\sum_{\ell\le L} \Ind{N_\ell \le 2m\gamma}
	\ge \f{n}{3m}
	\bigg)
	\le
	\f{3m}{n} 
	\sum_{0\le\ell\le L-1} \P(N_\ell \le 2m\gamma)\,.
	\]
Applying Lemma~\ref{l:gs.poly.bd} gives, for 
$1/m\le 2\gamma\le \gamma_0=\gamma_0(|E(U)|,E_{\max}(U))$,
	\[
	\P\bigg(
	\f1n
	\Big|\Big\{i \le n : (\bg,\bv^i) \in E(U)\Big\}
	\Big| \le \gamma
	\bigg) 
	\le 4 (2\gamma)^{1/(24c)}\,.
	\]
The claim follows. 
\end{proof}

\subsection{Polynomial concentration of free energy} 
\label{ss:poly.conc}

In this subsection we give the \hyperlink{proof:p.poly.conc}{proof of Proposition~\ref{p:poly.conc}}. 
Towards this end, we first state and prove Lemma~\ref{l:gibbs.poly.bound} below.
This is an adaptation of \cite[Propn.~8.2.6]{MR3024566} (see also \cite[Lem.~9.2.2]{MR3024566}), using 
Proposition~\ref{p:gaus.poly} in place of \cite[Cor.~8.2.5]{MR3024566}.

\begin{lem}\label{l:gibbs.poly.bound}
Let $\mu$ be any probability measure on $\set{-1,+1}^N$ with weights proportional to $w(J)$ such that
 $0\le w(J)\le 1/2^N$ for all $J\in\set{-1,+1}^N$, and
	\[
	W = \sum_J w(J) \ge e^{-N\tau}
	\]
for $\tau=\exp(-12)$. If $\P$ denotes the law of a standard gaussian vector $\bg$ in $\R^N$, then
	\[
	\P\bigg(\mu\bigg(\bigg\{
		J\in\set{-1,+1}^N: \f{(\bg,J)}{N^{1/2}} \in E(U)
		\bigg\}\bigg)
		\le\f{\gamma}{4}\bigg)
	\stackrel{\eqref{e:Upsilon.bound}}{\le}
	\gamma^{11/2}
	\,,
	\]
for $\exp(14)/N\le\gamma\le\gamma_0=\gamma_0(|E(U)|,E_{\max}(U))$ and $N$ large enough.

\begin{proof} First, it follows by a direct application of \cite[Lem.~9.2.1]{MR3024566} that
since $W\ge \exp(-N\tau)$, we have
	\beq\label{e:small.overlap}
	\mu^{\otimes 2}
	\bigg( \bigg\{ (J^1,J^2) \in \set{-1,+1}^{2N} :
	\f{(J^1,J^2)}{N}\ge (8\tau)^{1/2}
	\bigg\}
	\bigg) \le \f1{\exp(2N\tau)}\,.\eeq
We then proceed to adapt the proof of
\cite[Propn.~8.2.6]{MR3024566}. Let
	\[
	Q_n
	\equiv
	\bigg\{J^{1:n}\equiv(J^1,\ldots,J^n)
	\in\set{-1,+1}^{nN}
	: \f{(J^k,J^\ell)}{N} \le (8\tau)^{1/2}
	\ \forall 1\le k<\ell \le n
	\bigg\}\,.
	\]
It follows from \eqref{e:small.overlap}
(and taking a union bound over all $1\le k<\ell \le n$) that
	\beq\label{e:n.overlaps.small}
	\mu^{\otimes n}(Q_n)
	\stackrel{\eqref{e:small.overlap}}{\ge}
	1- \f{n^2}{2\exp(2N\tau)}
	\ge \f12\,,
	\eeq
where the last inequality holds provided $n \le \exp(N\tau)$. Next define
	\[
	\Omega_\gamma(J^{1:n})
	\equiv \bigg\{\bg
		: \f1n\bigg|
		\bigg\{ \ell\le n:
		\f{(\bg,J^\ell)}{N^{1/2}} \in E(U)
		\bigg\} \le \gamma
		\bigg\}\,.
	\]
If we take $c=(8\tau)^{1/2}$, then $c\le1/12$ by the assumption $\tau=\exp(-12)$, and so
Proposition~\ref{p:gaus.poly} implies that 
for every $J^{1:n}\in Q_n$ we have the bound
	\beq\label{e:apply.gaus.tail.bound}
	\P\Big( \Omega_\gamma(J^{1:n})\Big)
	\le \gamma^{1/(25c)}\,,
	\eeq
for $\log(5/c)/(\log n)\le\gamma\le\gamma_0$ and $n$ large enough. Define the random variable
	\[\Upsilon_\gamma
	\equiv \sum_{J^{1:n}\in Q_n}
	\mu^{\otimes n}(J^{1:n})
	\I\Big\{\bg \in \Omega_\gamma(J^{1:n})\Big\}\,,
	\]
and note that Markov's inequality combined with \eqref{e:apply.gaus.tail.bound} gives
	\beq\label{e:Upsilon.bound}
	\P\bigg(\Upsilon_\gamma\ge \f14\bigg)
	\le \f{\E \Upsilon_\gamma}{1/4}
	= 4 
	\sum_{J^{1:n}\in Q_n}
	\mu^{\otimes n}(J^{1:n})
	\P\Big( \Omega_\gamma(J^{1:n})\Big)
	\stackrel{\eqref{e:apply.gaus.tail.bound}}
		{\le}
		4\gamma^{1/(25c)}\,.
	\eeq
On the other hand, we can lower bound
	\begin{align*}
	\Gm &\equiv  \mu\bigg(\bigg\{
		J\in\set{-1,+1}^N: \f{(\bg,J)}{N^{1/2}} \in E(U)
		\bigg\}\bigg)
	=\sum_{J^{1:n}}
	\mu^{\otimes n}(J^{1:n})
	\f1n
	\bigg|\bigg\{
	\ell\le n:\f{(\bg,J^\ell)}{N^{1/2}} \in E(U)\bigg\}
	\bigg| \\
	&\ge 
	\gamma
	\sum_{J^{1:n}\in Q_n}
	\mu^{\otimes n}(J^{1:n})
	\I\Big\{\bg \notin \Omega_\gamma(J^{1:n})\Big\}
	= \gamma\bigg(
	\mu^{\otimes n}(Q_n) - \Upsilon_\gamma
	\bigg)
	\stackrel{\eqref{e:n.overlaps.small}}{\ge} \gamma
	\bigg(\f12- \Upsilon_\gamma
	\bigg)\,.
	\end{align*}
As a consequence, if $\Gm\le\gamma/4$, we must have $\Upsilon_\gamma\ge1/4$. It follows that
	\[
	\P\bigg(\Gm \le\f{\gamma}{4}\bigg)
	\le
	\P\bigg(\Upsilon_\gamma \ge\f14\bigg) 
	\stackrel{\eqref{e:Upsilon.bound}}{\le}
	4\gamma^{1/(25c)}\,,
	\]
again for $\log(5/c)/(\log n)\le\gamma\le\gamma_0$ and $n$ large enough.
Recall moreover that for \eqref{e:n.overlaps.small} to hold we must have 
 $n\le \exp(N\tau)$, so we must ultimately require
	\[
	\gamma_0\ge\gamma
	\ge\f{\log(5/c)}{N\tau}
	= \f{\log(5/(8\tau)^{1/2})}{N\tau}\,.
	\]
The claim follows by recalling $\tau=\exp(-12)$ and $c=(8\tau)^{1/2}$.
\end{proof}
\end{lem}

We now proceed to prove Proposition~\ref{p:poly.conc}.
This is an adaptation of the proof of
\cite[Propn.~9.2.6]{MR3024566}, using 
the above result Lemma~\ref{l:gibbs.poly.bound} in place of \cite[Propn.~8.2.6]{MR3024566}.

\begin{proof}[\hypertarget{proof:p.poly.conc}{Proof of Proposition~\ref{p:poly.conc}}]
As in the
\hyperlink{proof:t.main.lbd.bdd}{proof of Theorem~\ref{t:main} in the bounded case}, let $\P^j$ denote probability conditional on 
the first $j$ rows of $\bG$,
 and let $\E^j$ denote expectation with respect to $\P^j$. Then, as in the proof of \cite[Propn.~9.2.6]{MR3024566}, we let $\bW\equiv \bZ/2^N$ and decompose
	\[\f1N\bigg\{\log_{N\tau} \bW
	-\E \log_{N\tau} \bW\bigg\}
	= \sum_{j\le M} \f1N\bigg\{ \E^j \log_{N\tau}\bW
	-\E^{j-1} \log_{N\tau}\bW\bigg\}
	\equiv \sum_{j\le M} X_j\,.\]
To bound $X_j$, recall \eqref{e:Z.j} and denote
	\[
	\bW_j
	\equiv \f{\bZ_j}{2^N}
	\equiv \sum_J w_j(J)
	\equiv \sum_J \f1{2^N}
	\prod_{\substack{a\le M, \\
		 a \ne j}}
		U \bigg(\f{(\bg^a,J)}{N^{1/2}}\bigg)
	\,.\]
Note that $0\le \bW\le \bW_j\le1$.  Since $\bW_j$ does not depend on the $j$-th row of $\bG$, we can rewrite
	\[N X_j
	=\E^j\bigg( \log_{N\tau}\bW-\log_{N\tau}\bW_j\bigg)
	-\E^{j-1}\bigg( \log_{N\tau}\bW-\log_{N\tau}\bW_j
		\bigg)\,.\]
Recall that $0\le \bW\le \bW_j$, so
if $\bW_j \le e^{-N\tau}$ then
$\log_{N\tau}\bW_j=-N\tau=\log_{N\tau}\bW$. It follows that
	\[
	L_j\equiv \log_{N\tau}\bW_j-\log_{N\tau}\bW
	= \I\Big\{\bW_j\ge e^{-N\tau}\Big\}
	\bigg(\log_{N\tau}\bW_j-\log_{N\tau}\bW\bigg)
	\in[0,N\tau]\,.
	\]
Recall that $\P_j$ denotes probability conditional on all rows of $\bG$ except the $j$-th one, and note $\E^{j-1}=\E_j\E^j$ where $\E_j$ is expectation with respect to $\P_j$. We can rewrite $X_j=-\dot{x}_j+\ddot{x}_j$ where
	\begin{align}\nonumber
	N\dot{x}_j
	&\equiv\E^j\bigg[
		\Big(\log_{N\tau}\bW_j-\log_{N\tau}\bW\Big)
			; \bW_j\ge e^{-N\tau}\bigg]
		=\E^j L_j\in[0,N\tau] \,,\\
	N\ddot{x}_j
	&\equiv\E^{j-1}\bigg[
		\Big(\log_{N\tau}\bW_j-\log_{N\tau}\bW\Big)
			; \bW_j\ge e^{-N\tau}\bigg]
		=\E^{j-1} L_j
		=\E_j(N\dot{x}_j)\,.
	\label{e:interpolate.X.k}
	\end{align}
For comparison let 
$\bar{X}_j =-\dot{z}_j+\ddot{z}_j$ where
$\dot{z}_j \equiv \dot{x}_j-\dot{e}_j$
and $\ddot{z}_j \equiv \ddot{x}_j-\ddot{e}_j$, for
	\begin{align*}
	N\dot{e}_j
	&\equiv\E^j\bigg[L_j;
		\f{\bW}{\bW_j} < \f{\delta'e^{14}}{4N}
			\bigg]
		=\E^j\bigg[
		L_j ; 
		\f{\bW}{\bW_j} < \f{\delta'e^{14}}{4N}
		\bigg]\in[0,N\tau]\,,\\
	N\ddot{e}_j
	&\equiv\E^{j-1}\bigg[ L_j;
			\f{\bW}{\bW_j} < \f{\delta'e^{14}}{4N}
			\bigg]
		= \E^{j-1}\bigg[
		L_j ; \f{\bW}{\bW_j} < \f{\delta'e^{14}}{4N}
		\bigg]
		=\E_j(N\dot{e}_j)\,.
	\end{align*}
Let $G_j$ be as in Lemma~\ref{l:gibbs.poly.bound}, and note the assumption $U(x)\ge\delta'\Ind{x\in E(U)}$ implies
	\[
	\f{\bW}{\bW_j}
	\ge \delta'
	G_j\bigg(\bigg\{
		J\in\set{-1,+1}^N: \f{(\bg^j,J)}{N^{1/2}} \in E(U)
		\bigg\}\bigg) \equiv \delta'\Gamma_j\,.
	\]
Since $0\le L_j\le N\tau$, we can use Markov's inequality to bound
	\beq\label{e:cond.exp.error}
	0\le 
	\E_j(N\dot{e}_j)
	=N\ddot{e}_j
	\le
	N\tau \,\E^{j-1}\bigg[
	\I\Big\{\bW_j\ge e^{-N\tau}\Big\}
	\P_j\bigg(\f{\bW}{\bW_j} 
	< \f{\delta'e^{14}}{4N}\bigg)
	\bigg]
	\le N\tau
	\bigg(\f{ \delta' e^{14}}{N}\bigg)^{11/2}
	\eeq
where the last inequality is by Lemma~\ref{l:gibbs.poly.bound}. It follows using Markov's inequality again that
	\beq\label{e:X.tilde.X.err}
	\P\bigg( \sum_{j\le M}
		\Big|X_j-\bar{X}_j\Big| \ge \f1{2N^2}
	\bigg)
	\le 2N^2 \sum_{j\le M} \E 
		\Big(\dot{e}_j+\ddot{e}_j\Big)
	\stackrel{\eqref{e:cond.exp.error}}{\le} N^3
	\bigg(\f{\delta' e^{14}}{N}\bigg)^{11/2}\,.
	\eeq
It remains to bound the random variables $\bar{X}_j=-\dot{z}_j+\ddot{z}_j$. Using Jensen's inequality, 	\[
	\exp( N\ddot{z}_j)
	\le \E^{j-1}\exp( N \dot{z}_j)
	\le1+\E^{j-1} \bigg[\f{\bW_j}{\bW}
	; \bW_j\ge e^{-N\tau},
			\f{\bW_j}{\bW} \le \f{4N}{\delta'e^{14}}\bigg]\,.\]
It then follows by using Lemma~\ref{l:gibbs.poly.bound} again that the above can be bounded by
	\begin{align}\nonumber
	&1+\E^{j-1} \bigg[
		\f{\bW_j}{\bW}
	; \bW_j\ge e^{-N\tau}, \f{\bW_j}{\bW} 
		\le \f{4N}{\delta'e^{14}}\bigg]
	\le 1+ \int_0^{4N/(\delta'e^{14})}
		\P^{j-1} \bigg(
		\bW_j\ge e^{-N\tau};\f{\bW_j}{\bW} \ge u\bigg)
			\,du \\
	&\qquad
	\le1+ \f{4}{\delta'\gamma_0}
	+ \int_{4/(\delta'\gamma_0)}^\infty
		\bigg( \f{4}{\delta'u}\bigg)^{11/2}\,du
	\le C_0
	\equiv C_0(|E(U)|,E_{\max}(U),\delta')\,.
	\label{e:integrate.tail.bound}
	\end{align}
It follows that we can choose $\lambda_0$ small enough (depending on $C_0$) such that for all $0\le\lambda\le\lambda_0$,
	\begin{align*}
	\E^{j-1}\bigg[ \exp(N\lambda|\bar{X}_j|)\bigg]
	&\le \exp(N\lambda\ddot{z}_j) \cdot
	\E^{j-1} \bigg[ \exp(N\lambda\dot{z}_j)\bigg] \\
	&\le \exp(N\lambda\ddot{z}_j) \cdot
	\bigg(\E^{j-1} \exp(N\dot{z}_j)\bigg)^\lambda
	\le (C_0)^{2\lambda}\le2\,.
	\end{align*}
It follows by the martingale version of Bernstein's inequality (see e.g.\ \cite[eq.~(A.41)]{MR3024566}) that
	\[
	\P\bigg(\bigg|\sum_{j\le M}
		\bar{X}_j\bigg| \ge t\bigg)
	\le 2\exp\bigg( -\f{Nt\lambda}{2}
	\min \bigg\{ 1,\f{t\lambda}{2} \bigg\} \bigg)
	\]
for all $t\ge0$. In particular, taking $t=(\log N)/N^{1/2}$ gives
	\beq\label{e:bernstein.conclusion}
	\P\bigg(\bigg|\sum_{j\le M}
		\bar{X}_j\bigg| \ge \f{\log N}{N^{1/2}} \bigg)
	\le \exp\bigg(-\f{\lambda^2(\log N)^2}{2}\bigg)\,.
	\eeq
The claimed bound follows by combining \eqref{e:X.tilde.X.err} with \eqref{e:bernstein.conclusion}.
\end{proof}

\subsection{Exponential concentration for smoothed model}
\label{ss:approx.arg} 

In this subsection we give the
\hyperlink{proof:p.conc.smoothed}{proof of Proposition~\ref{p:conc.smoothed}}, showing concentration for the log-partition function of the smoothed model \eqref{e:tilde.Z}.

\begin{thm}[{Pisier \cite{MR864714}}] \label{t:pisier}
If $f:\R^n\to\R$ is $C^1$, and $X$ and $Y$ are independent standard gaussian random variables in $\R^n$, then for any convex function $g:\R\to\R$ it holds that
	\[\E g\Big(f(X)-f(Y)\Big)
	\le \E g\bigg(\f{\pi}{2} (\nabla f(X),Y)\bigg)\,.
	\]
In particular, taking $g(x)=\exp(sx)$ for any real number $s$ gives
	\[
	\E \exp \bigg\{ s\Big(f(X)-f(Y)\Big)\bigg\}
	\le 
	\E \exp\bigg\{ \f{s^2\pi^2}{8} \|\nabla f(X)\|^2\bigg\}\,.
	\]
In the case that $\nabla f$ is bounded, this recovers the standard theorem \cite{MR0458556} (see also \cite{MR399402})
on concentration of Lipschitz functionals of gaussian random variables.
\end{thm}

\newcommand{\smax}{s_{\max}}

We also recall that if $\bG$ is an $M\times N$ matrix with i.i.d.\ standard gaussian entries and $M\le N$, then the maximum singular value $ s_{\max}(\bG)$ satisfies the tail bound
	\beq\label{e:svbound}
	\P\bigg( s_{\max}(\bG) \ge (\BIGCTWO N)^{1/2} + t\bigg)
	\le \f{2}{\exp(\csttwo t^2)}
	\eeq
for all $t\ge0$, where $\csttwo$ and $\BIGCTWO$ are absolute constants.
See for instance \cite[Propn.~2.4]{MR2827856} where the result is 
in fact stated more generally for matrices with independent subgaussian entries (with mean zero and unit variance). 
From this bound it is straightforward to deduce the following:

\begin{lem}\label{l:sv.bound}
If $\bG$ is an $M\times N$ matrix with i.i.d.\ standard gaussian entries and $M\le N$, then
we can take $\csttwo=1/\BIGCTWO \le 1$
in the bound \eqref{e:svbound}. With this choice of constants, we have
	\[
	\E\exp\Big( \vartheta \smax(\bG)^2\Big)
	\le 16 N \exp(2\vartheta \BIGCTWO N)
	\]
for all $0\le\vartheta \le \csttwo/2 = 1/(2\BIGCTWO)$.

\begin{proof} It follows by a change of variables that
	\begin{align*}
	E(\vartheta)
	&\equiv\E\exp\Big( \vartheta \smax(\bG)^2\Big)
	= \int_0^\infty 
	\P\bigg( \exp( \vartheta \smax(\bG)^2) \ge x\bigg)\,dx \\
	&=2\vartheta \int_0^\infty
	u \exp(\vartheta u^2)   \cdot
	\P\bigg( \smax(\bG) \ge u
		\bigg)  \,du
	\le \textup{(I)}+\textup{(II)}\,,\end{align*}
where $\textup{(I)}$ is the contribution to the integral from $u\le (\BIGCTWO N)^{1/2}$, while $\textup{(II)}$ is the contribution from $u\ge (\BIGCTWO N)^{1/2}$. We then have the trivial bound
	\begin{align*}
	\textup{(I)}
	&\le 2\vartheta \int_0^{(\BIGCTWO N)^{1/2}}
		u\exp(\vartheta u^2)   \,du
	\le
	2\vartheta \exp(\vartheta \BIGCTWO N)
	\int_0^{(\BIGCTWO N)^{1/2}} u\,du\\
	&= \vartheta \BIGCTWO N
	\exp(\vartheta \BIGCTWO N)
	\le \f{N}{2}
	\exp(\vartheta \BIGCTWO N)
	\,,
	\end{align*}
where the last inequality uses
the assumption $\vartheta \le \csttwo/2 = 1/(2\BIGCTWO)$.
For the other term, it follows from the singular value tail bound \eqref{e:svbound}
(and again using $\vartheta \le \csttwo/2= 1/(2\BIGCTWO)$) that
	\begin{align*}\textup{(II)}
	&\le 4\vartheta\int_0^\infty
	\bigg((\BIGCTWO N)^{1/2}+u\bigg)
	\exp\bigg\{
	\vartheta \Big((\BIGCTWO N)^{1/2}+u\Big)^2 
	-\csttwo u^2
	\bigg\}\,du \\
	&\le4\vartheta
	\exp(\vartheta \BIGCTWO N)
	\int_0^\infty
	\bigg((\BIGCTWO N)^{1/2}+u\bigg)
	\exp\bigg\{ 2\vartheta (\BIGCTWO N)^{1/2} u
	- \f{\csttwo u^2}{2}\bigg\} \,du\,.\end{align*}
Completing the square and making another change of variables gives
	\begin{align*}
	\textup{(II)}
	&\le 4\vartheta
	\exp\bigg\{\bigg(1+\f{2\vartheta }{\csttwo}\bigg)\vartheta \BIGCTWO N
	\bigg\}
	\int_{-\infty}^\infty
	\bigg| u +\bigg( 1 + 
		\f{2\vartheta}{\csttwo}
		\bigg) (\BIGCTWO N)^{1/2}
		\bigg|
	\exp\bigg\{ - \f{\csttwo u^2}{2}\bigg\}\,du\\
	&\le 
	\f{4\vartheta}{(\csttwo)^{1/2}}
	\exp(2 \vartheta\BIGCTWO N)
	\int_{-\infty}^\infty
	\bigg| \f{u}{(\csttwo)^{1/2}} +2 (\BIGCTWO N)^{1/2}
		\bigg|
	\exp\bigg\{ - \f{u^2}{2}\bigg\}\,du\\
	&\le
	\exp(2 \vartheta\BIGCTWO N)
	\f{4\vartheta(2\pi)^{1/2}}
		{\csttwo}
		\Big( 1 + 2N^{1/2}\Big)
	\le 6(2\pi)^{1/2}\cdot N \exp(2 \vartheta\BIGCTWO N)
	\,.
	\end{align*} 
Combining the bounds for \textup{(I)} and \textup{(II)} gives the claimed bound.
\end{proof}
\end{lem}

\begin{lem}\label{l:exp.gradient.squared}
Suppose $U$ satisfies Assumption~\ref{a:bdd}, and let $\bZ(\ETA)$ be as defined by \eqref{e:tilde.Z}. If $f=\log\bZ(\ETA)$ viewed as a function of the gaussian disorder $\bG$, then
there exists a finite constant $C_1(U;\ETA)$ such that
	\[
	\E\exp\bigg( s^2 \|\nabla f(\bG)\|^2\bigg)
	\le 16N \cdot
	\exp \bigg\{
	N\cdot 6 \BIGCTWO   C_1(U;\ETA)^2 s^2
	\bigg\}
	\] 
for all $|s| \le (\csttwo)^{1/2}/(2C_1(U;\ETA))$,
where $\BIGCTWO$ and $\csttwo$ are the constants from Lemma~\ref{l:sv.bound}.

\begin{proof}
Recall that $U_\ETA\equiv U*\varphi_\ETA$, and denote $u_\ETA\equiv\log U_\ETA$. Denote the probability meausure
	\[
	\mu_\ETA(J)
	\equiv \f{w_\ETA(J)}{\bZ(\ETA)}
	= \f{1}{\bZ(\ETA)}
	\prod_{a\le M} U_\ETA\bigg(\f{(\bg^a,J)}{N^{1/2}}\bigg)
	= \f{1}{\bZ(\ETA)}
	\prod_{a\le M} U_\ETA(\Delta_a)\,,
	\]
where we abbreviate $\Delta_a = (\bg^a,J)/N^{1/2}$. Then
	\[
	\bigg|\f{df}{dg_{a,i}}\bigg|
	=\bigg| \sum_J \mu_\ETA(J) (u_\ETA)'(\Delta_a) \f{J_i}{N^{1/2}}\bigg|
	\le \f1{N^{1/2}} \sum_J w(J) |(u_\ETA)'(\Delta_a)| \,.
	\]
Note that if $(u_\ETA)'$ were uniformly bounded, then $f$ would be $A$-Lipschitz with $A=\|(u_\ETA)'\|_\infty/N^{1/2}$, and the desired exponential concentration for $\log\bZ(\ETA)$ would follow from standard concentration theorems for Lipschitz functionals of gaussians.
Since $(u_\ETA)'$ may be unbounded, we cannot conclude that $f$ is Lipschitz. However, we note that
	\[(u_\ETA)'(x)
	=\f{\E_\xi[\xi U(x+\ETA\xi)]}{\ETA\E_\xi U(x+\ETA\xi)}
	\le C_1(U;\ETA)\Big(1+|x|\Big)\,,\]
where the last bound holds by an obvious extension of Lemma~\ref{l:poly} (using Assumption~\ref{a:bdd}). Therefore
	\[
	\bigg|\f{df}{dg_{a,i}}\bigg|
	\le\f{C_1(U;\ETA)}{N^{1/2}}
		\bigg(1+\sum_J \mu_\ETA(J) |\Delta_a|\bigg)
	= \f{C_1(U;\ETA)}{N^{1/2}}
		\bigg(1 + \<|\Delta_a|\>_\ETA\bigg)\,,
	\]
where $\<\cdot\>_\ETA$ denotes expectation over $\mu_\ETA$. It follows that
	\begin{align*}
	\|\nabla f\|^2
	&\le C_1(U;\ETA)^2
	\sum_{a\le M}  \bigg(1 + \<|\Delta_a|\>_\ETA\bigg)^2
	\le 2C_1(U;\ETA)^2
	\sum_{a\le M}  \bigg(1 + (\<|\Delta_a|\>_\ETA)^2\bigg) \\
	&\le 2C_1(U;\ETA)^2
	\sum_{a\le M} \bigg( 1+ \<(\Delta_a)^2\>_\ETA\bigg)
	= 2C_1(U;\ETA)^2
	\bigg\{  M + \f{\|\bG J\|^2}{N}  \bigg\}\\
	&\le
	2C_1(U;\ETA)^2
	\bigg\{  M + s_{\max}(\bG)^2 \bigg\}
	\,,
	\end{align*}
where $s_{\max}(\bG)$ denotes the maximum singular value of $\bG$, as above. Taking the expectation over $\bG$ and applying Lemma~\ref{l:sv.bound} gives
	\begin{align*}
	\E\exp \bigg(s^2 \|\nabla f(\bG)\|^2\bigg)
	&\le
	16N \cdot
	\exp \bigg\{
	2\Big( M 
	+
	2 \BIGCTWO N\Big)
	C_1(U;\ETA)^2 s^2
	\bigg\}\,,
	\end{align*}
where the bound holds provided $|s| \le (\csttwo)^{1/2}/(2C_1(U;\ETA))$.
The result follows by recalling that we assumed $M\le N$
and $\BIGCTWO\ge1$.
\end{proof}
\end{lem}

\begin{proof}[\hypertarget{proof:p.conc.smoothed}{Proof of Proposition~\ref{p:conc.smoothed}}]
Let $\Gprime$ be an independent copy of $\bG$.
It follows by Theorem~\ref{t:pisier} and Lemma~\ref{l:exp.gradient.squared} that
	\begin{align*}
	&\E \exp \bigg\{ s\Big(f(\bG)-\E f(\bG)\Big)\bigg\}
	\le \E \exp \bigg\{ s\Big(f(\bG)-f(\Gprime)\Big)\bigg\} \\
	&\qquad\le 
	\E \exp\bigg\{ \f{s^2\pi^2}{8} \|\nabla f(\bG)\|^2\bigg\}
	\le 16N\cdot \exp\bigg\{
	N\cdot 8 \BIGCTWO   C_1(U;\ETA)^2 s^2
	\bigg\}\,,
	\end{align*}
for all $|s| \le (\csttwo)^{1/2}/(3C_1(U;\ETA))$.
Thus, for $x\ge0$, it holds for $0\le s\le (\csttwo)^{1/2}/(3C_1(U;\ETA))$ that
	\begin{align*}
	\P\bigg(f(\bG)-\E f(\bG) \ge Nx\bigg)
	&\le\E \exp\bigg\{
	 s\Big( f(\bG)-\E f(\bG)\Big)
	 -Nsx\bigg\} \\
	&\le 16N \cdot
	\exp\bigg\{ N\Big(
	8 \BIGCTWO   C_1(U;\ETA)^2 s^2
	- sx
	\Big)\bigg\}\,.
	\end{align*}
A similar bound holds for $x\le0$. In any case it is clear that we can take $s$ small enough to obtain exponential decay. In particular, for $x\ge0$ small enough we can let
	\[s = \f{x}{16 \BIGCTWO \cdot C_1(U;\ETA)^2}
		\le 
		\f{(\csttwo)^{1/2}}{3C_1(U;\ETA)}
		\,,\]
where the last bound holds for $x \le 5 (\BIGCTWO)^{1/2} C_1(U;\ETA)$. 
This results in the bound
	\[\P\bigg(\Big|f(\bG)-\E f(\bG)\Big| \ge Nx\bigg)
	\le32N \cdot
	\exp\bigg\{ -\f{Nx^2}{ 32 \BIGCTWO C_1(U;\ETA)^2 }	\bigg\}\,,
	\]
which concludes the proof.
\end{proof}

\begin{cor}\label{c:smoothed.free.energy}
Suppose $U$ satisfies Assumptions~\ref{a:bdd} and \ref{a:Lip}, and let $\bZ(\ETA)$ be as in \eqref{e:tilde.Z}. Then
	\[
	\lim_{N\to\infty}\f1N\log\bZ(\ETA)
	=\RS(\alpha;U_\ETA)
	\]
for all $0<\alpha\le\alpha'(U)$.
 
\begin{proof} Recall from \eqref{e:compare.alphas.U} that if $0<\alpha\le\alpha'(U)$, then we will also have $\alpha\le\alpha(U_\ETA)$ for $\ETA$ small enough. The upper bound on $\bZ(\ETA)$ follows from the upper bound in Theorem~\ref{t:main},
\hyperlink{proof:t.main.ubd}{which was already proved at the end of Section~\ref{s:analysis.first.mmt}}. For the lower bound on $\bZ(\ETA)$, we argue similarly as in the 
\hyperlink{proof:t.main.lbd.bdd}{proof of the Theorem~\ref{t:main} lower bound for the case $\|u\|_\infty<\infty$}, but using the concentration result from Proposition~\ref{p:conc.smoothed} in place of the Azuma--Hoeffding bound. To this end, let $\bar{\bZ}(\ETA)$ be defined as $\bar{\bZ}$ from \eqref{e:Z.bar}, but with $U_\ETA$ in place of $U$. It follows from Theorem~\ref{t:lbd}
(by the same calculation leading to \eqref{e:second.mmt.lbd.cond.on.t}) that,
with high probability,
	\[\P\bigg(
	\f1N \log \bar{\bZ}(\ETA) \ge 
		 \RS(\alpha;U_\ETA) - o_t(1)
	\,\bigg|\,\FF(t)\bigg)
	\ge \f1{\exp(No_t(1))}\,.\]
On the other hand, it follows from Proposition~\ref{p:conc.smoothed} that,
again with high probability,
	\[
	\P\bigg(
	\f1N\log \bZ \ge \f1N\E\log \bZ(\ETA) + x\,\bigg|\,\FF(t)
	\bigg) \le 35N \cdot
	\exp\bigg\{ -\f{Nx^2}{ 35 \BIGCTWO C_1(U;\ETA)^2 }	\bigg\}
	\]
for sufficiently small $x>0$. The above two bounds are in contradiction with one another unless
	\[
	\f1N\E\log \bZ(\ETA) 
	\ge
	\RS(\alpha;U_\ETA) - o_N(1)\,.
	\]
It then follows by another application of Proposition~\ref{p:conc.smoothed} that
	\begin{align*}
	\P\bigg(
	\f1N\log \bZ(\ETA) \le \RS(\alpha;U_\ETA) - o_N(1) - x
	\bigg)
	&\le \P\bigg(
	\f1N\log \bZ(\ETA) \le \f1N\E\log \bZ(\ETA) - x
	\bigg) \\
	&\le35N \cdot
	 \exp\bigg\{ -\f{Nx^2}{35 \BIGCTWO C_1(U;\ETA)^2 }\bigg\}
	\end{align*}
for sufficiently small $x>0$. This yields the lower bound for $\bZ(\ETA)$ and concludes the proof.
\end{proof}
\end{cor}

\subsection{Comparison with smoothed model and conclusion}
\label{ss:conclusion.main.thm}

In this subsection we prove Proposition~\ref{p:Z.eps.approx.Z} which gives the comparison between the quantities $\bZ$ and $\bZ(\ETA)$ from \eqref{e:Z} and \eqref{e:tilde.Z}. We then \hyperlink{proof:t.main.conclusion}{conclude the proof of the main theorem}.

\begin{proof}[\hypertarget{proof:p.Z.L.approx.Z}{Proof of Proposition~\ref{p:Z.eps.approx.Z}}] 
Some of the steps below are similar to the steps in the \hyperlink{proof:p.poly.conc}{proof of Proposition~\ref{p:poly.conc}}.  Let
	\[\bV_k \equiv
	\f1{2^N} \sum_J 
		\bigg\{\prod_{a\le k}
		U_\ETA\bigg( \f{(\bg^a,J)}{N^{1/2}}\bigg)\bigg\}
		\bigg\{ \prod_{k<a\le M}
		U\bigg( \f{(\bg^a,J)}{N^{1/2}}\bigg)
		\bigg\}\,.\]
Recall $\bW\equiv\bZ/2^N$, and write $\tilde{\bW}\equiv\bZ(\ETA)/2^N$. Note $\bV_0=\bW$, and $\bV_M=\tilde{\bW}$. Let us also define
	\beq
	\label{e:cavity.k}
	\bV_{k,\circ} \equiv
	\f1{2^N} \sum_J 
		\bigg\{\prod_{a< k}
		U_\ETA\bigg( \f{(\bg^a,J)}{N^{1/2}}\bigg)\bigg\}
		\bigg\{ \prod_{k<a\le M}
		U\bigg( \f{(\bg^a,J)}{N^{1/2}}\bigg)
		\bigg\}
		\equiv \sum_J w_{k,\circ}(J)\,.\eeq
Note that $\bV_{k,\circ} \ge \max\set{\bV_{k-1},\bV_k}$. We can then decompose
	\[\f1N\E\bigg(
	\log_{N\tau} \tilde{\bW}-\log_{N\tau} \bW
	\bigg)
	=\f1N\sum_{k\le M}
	\E\bigg(
	\log_{N\tau} \bV_k-\log_{N\tau} \bV_{k-1}
	\bigg)
	= \sum_{k\le M} y_k
	\]
(compare with \eqref{e:interpolate.X.k}). Let $G_{k,\circ}$ be the probability measure on $\set{-1,+1}^N$ with weights proportional to $w_{k,\circ}(J)$ as defined by \eqref{e:cavity.k}.
Write $\<\cdot\>_{k,\circ}$ for expectation with respect to $G_{k,\circ}$. Abbreviate
	\[
	U_k \equiv 
	U\bigg( \f{(\bg^k,J)}{N^{1/2}}\bigg)\,,\quad
	\tilde{U}_k
	\equiv U_\ETA\bigg( \f{(\bg^k,J)}{N^{1/2}}\bigg)
	\,.
	\]
Recalling that $U(x) > \delta'\Ind{x\in E(U)}$ (from \eqref{e:positive.measure.set}), we have
	\[\f{\bV_k}{\bV_{k,\circ}} =
	\<U_k\>_{k,\circ}
	\ge\delta'
	G_{k,\circ}
	\bigg(\bigg\{
		J\in\set{-1,+1}^N: \f{(\bg^k,J)}{N^{1/2}} \in E(U)
		\bigg\}\bigg)
	\equiv \delta'\Gamma_{k,\circ}\,.
	\]
For $\ETA$ small enough we will also have $U_\ETA(x) > \delta'\Ind{x\in E(U)}$, so we can also bound
	\[
	\f{\bV_{k-1}}{\bV_{k,\circ}} =
	\<\tilde{U}_k\>_{k,\circ}
	\ge\delta'
	G_{k,\circ}
	\bigg(\bigg\{
		J\in\set{-1,+1}^N: \f{(\bg^k,J)}{N^{1/2}} \in E(U)
		\bigg\}\bigg)
	= \delta'\Gamma_{k,\circ}\,.
	\]
We also have from \cite[Lem.~8.3.10]{MR3024566} that
if $x,y,z\le1$, then
	\[
	\bigg|\log_A(xz)-\log_A(yz)\bigg|
	\le \bigg|\log_A x-\log_A y\bigg|
		\cdot \mathbf{1}\Big\{z \ge e^{-A}\Big\}\,.
	\]
Combining the above bounds gives
	\begin{align*}
	|Ny_k|
	&\le \E\bigg|
		\log_{N\tau}\Big(\bV_{k,\circ} \<U_k\>_{k,\circ}\Big)
		-\log_{N\tau}\Big(\bV_{k,\circ} \<\tilde{U}_k\>_{k,\circ}\Big)
		\bigg| \\
	&\le\E\bigg[\bigg| \log_{N\tau}\<U_k\>_{k,\circ}
		-\log_{N\tau}\<\tilde{U}_k\>_{k,\circ}\bigg|
		; \bV_{k,\circ}\ge e^{-N\tau}\bigg] \le \textup{(I)}+ \textup{(II)}\,,
	\end{align*}
for \textup{(I)} and \textup{(II)} defined by
	\begin{align*}
	\textup{(I)}
	&\equiv N\tau
	\P\bigg(
	\bV_{k,\circ}\ge e^{-N\tau},
	\delta'\Gamma_{k,\circ} < \f{e^{14}}{4N}
	\bigg)\,,\\
	\textup{(II)} &\equiv
	 \E\bigg[
	\bigg|
	\log\bigg(1+
	 \f{\<\tilde{U}_k\>_{k,\circ}
	 -\<U_k\>_{k,\circ}}{\delta'\Gamma_{k,\circ}}
	\bigg)\bigg|
	; \bV_{k,\circ}\ge e^{-N\tau},
	\Gamma_{k,\circ} \ge \f{e^{14}}{4N}
	\bigg]\,.
	\end{align*}
Combining with Lemma~\ref{l:gibbs.poly.bound} gives
(similarly to \eqref{e:cond.exp.error})
	\[
	\textup{(I)}
	\le N\tau \bigg(\f{e^{14}}{N}\bigg)^{11/2}\,.
	\]
Meanwhile, using
the bound $\log(1+x)\le x$
together with the Cauchy--Schwarz inequality gives
	\begin{align*}
	\textup{(II)}
	&\le
	\E \bigg[\bigg|\f{\<\tilde{U}_k\>_{k,\circ}
	 -\<U_k\>_{k,\circ}}{\delta'\Gamma_{k,\circ}}
	 \bigg| ; \bV_{k,\circ}\ge e^{-N\tau},
	\Gamma_{k,\circ} \ge \f{e^{14}}{4N} \bigg] \\
	&\le \f1{\delta'}
	\bigg\{ \E\bigg[ \Big(\<\tilde{U}_k-U_k\>_{k,\circ}\Big)^2\bigg]
	 \cdot
	\E\bigg[ \f1{(\Gamma_{k,\circ})^2}
		; \bV_{k,\circ}\ge e^{-N\tau},
		\Gamma_{k,\circ} \ge \f{e^{14}}{4N} \bigg]
	 \bigg\}^{1/2}\,.
	\end{align*}
For the first factor we note that
	\[
	\E\bigg[ \Big(\<\tilde{U}_k-U_k\>_{k,\circ}\Big)^2\bigg]
	\le
	\bigg\<\E\Big[ (\tilde{U}_k-U_k)^2\Big]\bigg\>_{k,\circ}
	=\E\bigg[ \Big(U_\ETA(\xi)-U(\xi)\Big)^2\bigg]
	\le o_\ETA(1)\,.
	\]
For the second factor, applying Lemma~\ref{l:gibbs.poly.bound} again gives
(similarly to \eqref{e:integrate.tail.bound})
	\[
	\E\bigg[ \f1{(\Gamma_{k,\circ})^2}
		; \bV_{k,\circ}\ge e^{-N\tau},
		\Gamma_{k,\circ} \ge \f{e^{14}}{4N} \bigg]
	\le
	\int_0^{(4N/e^{14})^2}
	\P\bigg( \f{e^{14}}{4N} \le \Gamma_{k,\circ} 
		\le \f1{y^{1/2}}\bigg)\,dy
	\]
	\[
	\le\bigg( \f{\gamma_0}{4}\bigg)^2
	+ \int_{(\gamma_0/4)^2}^{(4N/e^{14})^2}
	\bigg(\f{4}{y^{1/2}}\bigg)^{11/2}\,dy
	\le \f{ e^{12}}{(\gamma_0)^{7/2}}\,,
	\]
where $\gamma_0=\gamma_0(|E(U)|,E_{\max}(U))$
is as in Proposition~\ref{p:gaus.poly}. Altogether it follows that
	\[
	|Ny_k|
	\le N\tau \bigg(\f{e^{14}}{N}\bigg)^{11/2}
	+ \f{o_\ETA(1) e^6}{\delta' (\gamma_0)^{7/4}}\,,
	\]
and the claim follows by summing over $k\le M=N\alpha$.
\end{proof}
We now finally finish the proof of the main theorem: 

\begin{proof}[\hypertarget{proof:t.main.conclusion}{Proof of Theorem~\ref{t:main} (conclusion)}] The \hyperlink{proof:t.main.ubd}{proof of the upper bound} was given at the end of Section~\ref{s:analysis.first.mmt}, after the \hyperlink{proof:t.ubd}{proof of Theorem~\ref{t:ubd}}. The \hyperlink{proof:t.main.lbd.bdd}{proof of the lower bound in the case $\|u\|_\infty<\infty$} was given at the end of Section~\ref{s:second.mmt}, after the \hyperlink{proof:t.lbd}{proof of Theorem~\ref{t:lbd}}. It remains to prove the lower bound in the case $\|u\|_\infty=\infty$. We follow the \hyperlink{sketch:proof.t.main}{proof sketch} given at the end of Section~\ref{s:intro}. It follows from Propositions~\ref{p:poly.conc} and  \ref{p:Z.eps.approx.Z} that
	\beq\label{e:log.N.tau.Z.approx.tilde.Z}
	\P\bigg(
	\bigg|\f1N\log_{N\tau}\bigg(\f{\bZ}{2^N}\bigg)
	- \f1N\E\log_{N\tau}\bigg(
		\f{\bZ(\ETA)}{2^N}\bigg)
	\bigg| 
	\ge
	\f{(\log N)^2}{N^{1/2}}
	+ o_\ETA(1)
	\bigg) \le o_N(1)\,.
	\eeq
Given $\epsilon>0$, we can choose $\ETA$ small enough such that
the $o_\ETA(1)$ error above is at most $\epsilon$ in absolute value.
By Proposition~\ref{p:RS.U.eps} together with
Corollary~\ref{c:smoothed.free.energy},
for $0<\alpha\le\alpha'(U)$ and $\ETA$ small enough we have
	\beq\label{e:tilde.Z.close.to.RS}
	\P\bigg(
		\f1N\log \bZ(\ETA)
		\le \RS(\alpha;U)- 2\epsilon\bigg)
		\le \P\bigg(
		\f1N\log \bZ(\ETA)
		\le\RS(\alpha;U_\ETA)-\epsilon\bigg)
	\le o_N(1)\,,
	\eeq
It follows from Corollary~\ref{c:rs.estimate} that 
$\RS(\alpha;U)\ge \log 2-\tau/4$ for $0<\alpha\le\alpha(U)$, so taking $\epsilon\le\tau/8$ in the above gives
	\[
	0\le
	\f1N\bigg\{
	\E\log_{N\tau}
	\bigg(\f{\bZ(\ETA)}{2^N}\bigg)
	-\E\log
	\bigg(\f{\bZ(\ETA)}{2^N}\bigg)
	\bigg\}\le\tau\P\bigg(
	\f1N\log \f{\bZ(\ETA)}{2^N} \le -\f{\tau}{2}
	\bigg)
	\stackrel{\eqref{e:tilde.Z.close.to.RS}}{\le} o_N(1)\,.
	\]
It follows by combining with \eqref{e:log.N.tau.Z.approx.tilde.Z} and \eqref{e:tilde.Z.close.to.RS} that
	\[\P\bigg(
	\f1N\log_{N\tau}\bigg(\f{\bZ}{2^N}\bigg)
	\le \RS(\alpha;U) -\log2-4\epsilon\bigg)
	\le o_N(1)\,.\]
Since $\RS(\alpha;U)-\log2 -4\epsilon \ge -\tau/4-4\epsilon \ge -\tau$,	it follows that in fact
	\[\P\bigg(\f1N\log \bZ	\le \RS(\alpha;U) - 4\epsilon\bigg)
	\le o_N(1)\,,\]
as claimed.
\end{proof}

\appendix \section{Review of AMP for perceptron}
\label{s:AMP}

In \S\ref{ss:gauss.conditioning} and \ref{ss:resampling}
we prove Lemma~\ref{l:resampling}. In the rest of the section, we give a heuristic derivation of the state evolution recursions introduced in Definition~\ref{d:state}. We emphasize that \S\ref{ss:gauss.conditioning} and \ref{ss:resampling} are rigorous, while \S\ref{ss:amp.heuristic.one}--\ref{ss:amp.heuristic.conclusion} are not (and are intended only to provide intuition). For rigorous derivations of the asymptotics described in \S\ref{ss:amp.heuristic.one}--\ref{ss:amp.heuristic.conclusion}, we again refer the reader to
\cite{MR2810285,MR3147441}.

\subsection{Gaussian conditioning results} 
\label{ss:gauss.conditioning}
Suppose for simplicity that $t,F:\R\to\R$ are two smooth functions. Let $Z$ denote a standard gaussian random variable. Suppose we have $(q,\psi)$ such that
(cf.\ \eqref{e:fp})
	\beq\label{e:q.psi.recurs}
	\begin{pmatrix}q \\ \psi \end{pmatrix}
	= \begin{pmatrix}
	\E[ t( \psi^{1/2} Z)^2]\\
	\alpha \E[ F( q^{1/2} Z)^2 ]
	\end{pmatrix}\,.\eeq
Let $\bmag^{(0)}=\zroN\in\R^N$, $\bn^{(0)}=\zroM\in\R^M$,
$\bmag^{(1)}=q^{1/2}\oneN\in\R^N$,
$\bn^{(1)}=(\psi/\alpha)^{1/2}\oneM\in\R^M$. 
The AMP iteration in this setting is given by
 (cf.\ \eqref{e:perceptron.AMP.col}
and \eqref{e:perceptron.AMP.row})
	\begin{align}
	\label{e:general.AMP.col}
	\bmag^{(t+1)}
	&= t\bigg( \f{\bG^\st\bn^{(t)}}{N^{1/2}}
		-\beta \bmag^{(t-1)}
	\bigg) \in\R^N\,, \\
	\bn^{(t+1)}
	&= F\bigg(\f{\bG\bmag^{(t)}}{N^{1/2}}
		-\acute{\beta}\bn^{(t-1)}
		\bigg)\in\R^M \,,
	\label{e:general.AMP.row}
	\end{align}
where the Onsager coefficients are defined as
(cf.\ \eqref{e:perceptron.onsager})
	\[\begin{pmatrix}
	\beta \\ \acute{\beta}
	\end{pmatrix}
	=
	\begin{pmatrix}
	\alpha \E F'(q^{1/2}Z) \\
	\E t'(\psi^{1/2}Z)
	\end{pmatrix}\,.
	\]
A preliminary observation is the following:

\begin{lem}\label{l:G.given.row.col}
Let $\bG$ be an $M\times N$ matrix with jointly gaussian entries.
Suppose $\br$ is a unit vector in $\R^N$, while $\bc$ is a unit vector in $\R^M$. Then
$\E(\bG\,|\,\bG\br)=\bG\br\br^\st$,
$\E(\bG\,|\,\bG^\st\bc)=\bc\bc^\st\bG$,
and
	\beq\label{e:G.given.row.col}
	\E\bigg(\bG\,\bigg|\,\bG\br,\bG^\st\bc\bigg)
	= (\bG\br)\br^\st + \bc(\bG^\st\bc)^\st
		- (\bc^\st\bG\br)\bc\br^\st
	\equiv \bGam\Big(\br,\bc,\bG\br,\bG^\st\bc\Big)
	\equiv\bGam
	\,.\eeq

\begin{proof} 
Denote $\bR_a\equiv \eM_a\br^\st$ for $a\le M$,
and denote $\bC_i\equiv \bc (\eN_i)^\st$ for $i\le N$.
The equations $\bG\br=\bbx$ and
$\bG^\st\bc=\bby$ are equivalent to the linear constraints
	{\setlength{\jot}{0pt}\begin{align*}
	x_a &= \eM_a \bG \br
	= (\bG, \eM_a\br^\st)
	= (\bG,\bR_a)\,,\\
	y_i &= \bc \bG (\eN_i)^\st
	=(\bG, \bc (\eN_i)^\st)
	=(\bG,\bC_i)\,,
	\end{align*}}%
where $(\cdot,\cdot)$ above denotes the Frobenius inner product. Note then that
	\begin{align*}
	\E(\bG\,|\,\bG\br=\bbx)
	&=\sum_{a\le M} (\bG,\bR_a)\bR_a
	=\sum_{a\le M} x_a \eM_a\br^\st
	=\bbx\br^\st = \bG\br\br^\st\,,\\
	\E(\bG\,|\,\bG^\st\bc=\bby)
	&=\sum_{i\le N}(\bG,\bC_i)\bC_i
	=\sum_{i\le N} y_i \bc (\eN_i)^\st
	=\bc\bby^\st = \bc\bc^\st\bG\,.
	\end{align*}
The claim \eqref{e:G.given.row.col} follows by noting that $\bGam\br=\bbx$, $\bGam^\st\bc=\bby$, and $\bGam$ is in the span of the $(\bR_a,\bC_i)$.
\end{proof}
\end{lem}

Let $\br^{(1)},\ldots,\br^{(t)}$ be the Gram--Schmidt orthogonalization of the vectors $\bmag^{(1)},\ldots,\bmag^{(t)}$. Likewise let $\bc^{(1)},\ldots,\bc^{(t)}$ be the Gram--Schmidt orthogonalization of the vectors $\bn^{(1)},\ldots,\bn^{(t)}$. Let $\bG^{(1)}\equiv\bG$, and suppose recursively that $\bG^{(s)}$ has been defined. Let
$\bG^{(s)}\br^{(s)} = \bbx^{(s)}$, $(\bG^{(s)})^\st\bc^{(s)} =\bby^{(s)}$, 
and define (cf.\  \eqref{e:G.given.row.col})
	\beq\label{e:modified.G}
	\bG^{(s+1)}
	\equiv \bG^{(s)}
	- \bGam\bigg(\br^{(s)},\bc^{(s)},
		\bbx^{(s)},\bby^{(s)}\bigg)
	\equiv \bG^{(s)}-\bGam^{(s)}\,.
	\eeq
We also define a corresponding $\sigma$-field
	\beq\label{e:FF.star}
	\FF_\star(t)
	\equiv \sigma\bigg(
	(\bbx^{(s)}: s\le t ),
	(\bby^{(s)}: s\le t )
	\bigg)\,.
	\eeq
The next lemma records some basic facts about $\FF_\star(t)$.

\begin{lem}\label{l:FF.star}For the AMP iteration described above, 
the random variables
	\[\bigg(
	\Big(\bmag^{(s)}, \bn^{(s)}, \br^{(s)}, \bc^{(s)}
		: s\le t+1\Big),
	\Big(
	\bG\bmag^{(\ell)},\bG^\st\bn^{(\ell)},
	\bG\br^{(\ell)},\bG^\st\bc^{(\ell)},
	\bGam^{(\ell)}
		: \ell \le t
	\Big)
	\bigg)
	\]
are all measurable with respect to $\FF_\star(t)$.
\begin{proof}
Recall that the initial vectors $\bmag^{(0)}$, $\bn^{(0)}$,  $\bmag^{(1)}$, $\bn^{(1)}$ are fixed and deterministic, so they are measurable with respect to the trivial $\sigma$-field $\FF_\star(0)$. From these we can also obtain the deterministic vectors $\br^{(1)}$ and $\bc^{(1)}$.
Next we consider the $\sigma$-field $\FF_\star(1)$: it is clear that $\bGam^{(1)}$ is  $\FF_\star(1)$-measurable. Next note that
	\[
	\bbx^{(1)} = \bG\br^{(1)}
	=\f{ \bG\bmag^{(1)} }{\|\bmag^{(1)}\|}\,,\quad
	\bby^{(1)} =\bG^\st\bc^{(1)}
	=\f{\bG^\st\bn^{(1)}}{\|\bn^{(1)}\|}\,,
	\]
so we see that $ \bG\bmag^{(1)}$ and $\bG^\st\bn^{(1)}$ are measurable with respect to $\FF_\star(1)$. We can then apply the AMP iteration \eqref{e:general.AMP.col} and \eqref{e:general.AMP.row} to obtain $\bmag^{(2)}$ and $\bn^{(2)}$, so these are also 
measurable with respect to $\FF_\star(1)$. It follows by Gram--Schmidt orthogonalization that $\br^{(2)}$ and $\bc^{(1)}$ are also 
 $\FF_\star(1)$-measurable.

Now suppose inductively that the claim holds up to $\FF_\star(t-1)$, and consider the $\sigma$-field $\FF_\star(t)$. Then the matrix $\bGam^{(t)}$ is clearly $\FF_\star(t)$-measurable. Next note that \eqref{e:modified.G} implies
	\[
	\bG=\bG^{(1)}
	=\bGam^{(1)}+\bG^{(2)} 
	=\ldots
	=
	\sum_{s=1}^{t-1} \bGam^{(s)} + \bG^{(t)}\,,
	\]
where the $\bGam^{(s)}$, $s\le t-1$, are all measurable with respect to $\FF_\star(t-1)\subseteq\FF_\star(t)$. Therefore
	\[
	\bG\br^{(t)}
	=\sum_{s=1}^{t-1} \bGam^{(s)}\br^{(t)} + \bbx^{(t)}
	\]
is $\FF_\star(t)$-measurable, as is $\bG^\st\bc^{(t)}$. Recall from the Gram--Schmidt orthogonalization that
	\[
	\br^{(t)}
	=\f{ \bmag^{(t)}-\sum_{s\le t-1} (\bmag^{(t)},\br^{(s)})\br^{(s)}}
		{ \|\bmag^{(t)}-\sum_{s\le t-1} (\bmag^{(t)},\br^{(s)})\br^{(s)}\|}\,,
	\]
which we can rearrange to obtain an expression for $\bmag^{(t)}$. It follows from this that $\bG\bmag^{(t)}$ is $\FF_\star(t)$-measurable, as is 
$\bG^\st\bn^{(t)}$. We can then apply the AMP iteration 
\eqref{e:general.AMP.col} and \eqref{e:general.AMP.row} to obtain $\bmag^{(t+1)}$ and $\bn^{(t+1)}$, so these are also 
measurable with respect to $\FF_\star(t)$. Finally, it follows by Gram--Schmidt orthogonalization that $\br^{(t+1)}$ and $\bc^{(t+1)}$ are also 
 $\FF_\star(1)$-measurable. This verifies the inductive hypothesis and proves the claim.
\end{proof}
\end{lem}

\subsection{Projection and resampling}\label{ss:resampling}
In this subsection we give the
\hyperlink{proof:l.resampling}{proof of Lemma~\ref{l:resampling}}. For notational convenience, the roles of $\bG$ and $\Gprime$ through this section are switched from the main body of the paper.

\begin{dfn}[similar to Definition~\ref{d:row.col}]\label{d:appx.row.col}
Given $\FF_\star(t-1)$ as in \eqref{e:FF.star}, consider the linear subspaces
	\begin{align*}
	V_\ROW(t)
		&\equiv \spn\bigg\{
	\eM_a(\bmag^{(s)})^\st
	: 1\le a\le M, 1\le s\le t\bigg\}\,,\\
	V_\COL(t)
		&\equiv \spn\bigg\{
	\bn^{(\ell)}(\eN_i)^\st
	: 1\le i\le N, 1\le \ell\le t\bigg\}\,.
	\end{align*}
It follows from Lemma~\ref{l:FF.star} that these (random) subspaces are measurable with respect to $\FF_\star(t-1)$. Let $V_\star(t)=V_\ROW(t)+V_\COL(t)$, and let $\proj_t$ denote orthogonal projection onto $V_\star(t)$.
\end{dfn}

We remark that $V_\star(t)$ is very similar to the (random) subspace
$V_{\ROW\COL}=V_\ROW(t)+V_\COL(t-1)$ which appears in the proof of Lemma~\ref{l:resampling}. We will address the discrepancy between $V_\star(t)$ and $V_{\ROW\COL}$ in the \hyperlink{proof:l.resampling}{proof of Lemma~\ref{l:resampling}}, below.
 The following is a straightforward consequence of the preceding lemmas and the definition:

\begin{cor}\label{c:FF.star}Let $\bG$ be an $M\times N$ matrix with i.i.d.\ standard gaussian entries.
With $\FF_\star(t)$ as in \eqref{e:FF.star},
	\[
	\E\Big(\bG \,\Big|\, \FF_\star(t)\Big) 
	=\sum_{s \leq t}\bGam^{(s)}
	= \bG-\bG^{(t+1)}
	=\proj_t(\bG)\,,
	\]
where $\proj_t$ is the orthogonal projection onto
the (random) subspace $V_\star(t)$  from Definition~\ref{d:appx.row.col}.
Moreover, conditional on $\FF_\star(t-1)$, $\bG^{(t+1)}$ is distributed as a standard gaussian element of the ($\FF_\star(t-1)$-measurable) subspace $V_\star(t)^\perp$, and is independent of $\FF_\star(t)$.

\begin{proof} Note that the recursive definition~\eqref{e:modified.G} implies
	\beq\label{e:G.t.decomp}
	\bG^{(t)}
	= \bG^{(t-1)}-\bGam^{(t-1)}
	= \ldots
	= \bG - \sum_{s=1}^{t-1} \bGam^{(s)}\,.
	\eeq
By induction, conditional on $\FF_\star(t-2)$, the random matrix $\bG^{(t)}$ is distributed as a standard gaussian element of the ($\FF_\star(t-2)$-measurable) subspace $V_\star(t-1)^\perp$, and is independent of $\FF_\star(t-1)$. It follows that $\bG^{(t)}$ has jointly gaussian entries conditional on $\FF_\star(t-1)$. We also have from Lemma~\ref{l:FF.star} that the vectors $\br^{(t)}$ and $\bc^{(t)}$ are measurable with respect to $\FF_\star(t-1)$. It follows by applying Lemma~\ref{l:G.given.row.col} (conditional on $\FF_\star(t-1)$) that
	\beq\label{e:apply.conditioning.lemma}
	\E\Big(\bG^{(t)}\,\Big|\, \FF_\star(t)\Big)
	=\bGam\bigg(\br^{(t)},\bc^{(t)},
		\bbx^{(t)},\bby^{(t)}\bigg)
	=\bGam^{(t)}\,,
	\eeq
and $\bG^{(t+1)}=\bG^{(t)}-\bGam^{(t)}$ 
and $\FF_\star(t)$ are independent given $\FF_\star(t-1)$.
It follows that
	\begin{align*}
	\E\Big(\bG \,\Big|\, \FF_\star(t)\Big) 
	&\stackrel{\eqref{e:modified.G}}{=}  
	\E\bigg( \sum_{s \leq t-1}\bGam^{(s)}
		+\bG^{(t)}\,\bigg|\, \FF_\star(t)\bigg) 
	=\sum_{s \leq t-1}\bGam^{(s)}
		+ \E\bigg( \bG^{(t)}\,\bigg|\, \FF_\star(t)\bigg) \\
	&\stackrel{\eqref{e:apply.conditioning.lemma}}{=}
		\sum_{s \leq t}\bGam^{(s)}
	\stackrel{\eqref{e:G.t.decomp}}{=} \bG-\bG^{(t+1)} \,.
	\end{align*}
We also note that $\bG^{(s+1)}\br^{(s)}
	=(\bG^{(s)}-\bGam^{(s)})\br^{(s)}=\zroM\in\R^M$ by the construction of $\bGam^{(s)}$, and as a result
	\[
	\bGam^{(s+1)}\br^{(s)}
	= \bigg(
	\bG\br\br^\st
		+ \bc\bc^\st\bG
		- (\bc^\st\bG\br)
		\bc\br^\st\bigg)^{(s+1)}
		\br^{(s)}
	=\zroM \in\R^M\,,
	\]
and likewise $(\bGam^{(s+1)})^\st\bc^{(s)}=\zroN\in\R^N$.
One can then show by induction that for all $s<t$ we have $\bG^{(t)}\br^{(s)}=\zroM\in\R^M$, and likewise $(\bG^{(t)})^\st\bc^{(s)}=\zroN\in\R^N$.  This implies, for all $\ell \le t$,
	\begin{align*}
	\bigg(
	\bG-\sum_{s \leq t}\bGam^{(s)}
	\bigg) \br^{(\ell)}
	\stackrel{\eqref{e:G.t.decomp}}{=} \bG^{(t+1)}\br^{(\ell)}
	&=\zroM\in\R^M\,,\\
	\bigg(
	\bG-\sum_{s \leq t}\bGam^{(s)}
	\bigg)^\st \bc^{(\ell)}
	\stackrel{\eqref{e:G.t.decomp}}{=} (\bG^{(t+1)})^\st
		\bc^{(\ell)}
	&=\zroN\in\R^N\,.
	\end{align*}
It follows from this that $\bG^{(t+1)}$ is orthogonal to $V_\star(t)$. On the other hand
	\[
	\bG-\bG^{(t+1)}
	= \sum_{s \leq t}\bGam^{(s)}
	\]
lies in $V_\star(t)$, so we see that this the orthogonal projection of $\bG$ onto $V_\star(t)$, as claimed. We also see that $\bG^{(t+1)}$ is the orthogonal projection of $\bG^{(t)}$ onto $V_\star(t)^\perp$, which confirms the inductive hypothesis.
\end{proof}
\end{cor}

The next result is similar to Lemma~\ref{l:resampling}:

\begin{lem}\label{l:resampling.variant}
Let $\bG$ be an $M\times N$ matrix with i.i.d.\ gaussian entries, and use it to define $\FF_\star(t)$ as in \eqref{e:FF.star}. As in Definition~\ref{d:appx.row.col}, let $\proj_t$ denote the orthogonal projection onto the $\FF_\star(t-1)$-measurable subspace $V_\star(t)$. Then, for any bounded measurable function $f:\R^{M\times N}\to\R$, we have
	\begin{align}\label{e:G.cond.indep}
	\E \bigg(f(\bG^{(t+1)}) \,\bigg|\, \FF_\star(t)\bigg)
	&=\E \bigg(f(\bG^{(t+1)}) \,\bigg|\, \FF_\star(t-1)\bigg)\\
	&= \E \bigg(f\Big(
	\Gprime-\proj_t(\Gprime)
	\Big) \,\bigg|\, \FF_\star(t-1)\bigg)
	\label{e:G-G'}
	\end{align}
where $\Gprime$ is an independent copy of $\bG$.

\begin{proof} We saw in Corollary~\ref{c:FF.star} that conditional on $\FF_\star(t-1)$,
the random matrix $\bG^{(t+1)}$ is distributed as a standard gaussian element of the ($\FF_\star(t-1)$-measurable) subspace $V_\star(t)^\perp$. As a result, $\bG^{(t+1)}$ and $\FF_\star(t)$ are independent conditional on $\FF_\star(t-1)$, so the first claim \eqref{e:G.cond.indep} follows. Since $\bG$ and $\Gprime$ are independent, if we condition on $\FF_\star(t-1)$ then the random matrix $\Gprime-\proj_t(\Gprime)$ is also distributed as a standard gaussian element of $V_\star(t)^\perp$. This implies \eqref{e:G.cond.indep}.
\end{proof}
\end{lem}

\begin{rmk} We can also give a more explicit description of the projection of $\Gprime$ onto $V_\star(t)$, although it is not needed in the above proof of Lemma~\ref{l:resampling.variant}. Define $\bG^{\bullet(1)}\equiv\Gprime$, and recursively
	\beq\label{e:G.prime.recursion}
	\bG^{\bullet(t+1)}
	\equiv \bG^{\bullet(t)}
	-\bGam^{\bullet(t)}
	\equiv\tilde{\bG}^{(t)}
	-\bGam\bigg(
		\br^{(t)},\bc^{(t)},\bG^{\bullet(t)}\br^{(t)},
		(\bG^{\bullet(t)})^\st\br^{(t)}
		\bigg)\,.
	\eeq
Note that $\bGam^{\bullet(t)}$ is defined using the vectors $\br^{(t)}$ and $\bc^{(t)}$ that came from $\bG$, not $\Gprime$. As in Definition~\ref{d:appx.row.col}, we let $\proj_t$ denote the orthogonal projection onto the $\FF_\star(t-1)$-measurable subspace $V_\star(t)$. We then claim that
	\beq\label{e:proj.G.prime}
	\proj_t(\Gprime)
	= \Gprime-\bG^{\bullet(t+1)}
	\stackrel{\eqref{e:G.prime.recursion}}{=}
	\sum_{s\le t}\bGam^{\bullet(t)}\,.
	\eeq
This is very similar to the proof of Corollary~\ref{c:FF.star}, but in fact simpler because $\bG$ and $\Gprime$ are independent, which implies that $\Gprime$ is independent of the random subspace $V_\star(t)$. 
Arguing as before, we have by construction $\bG^{\bullet(s+1)}\br^{(s)}=\zroM$ and $(\bG^{\bullet(s+1)})^\st\bc^{(s)}=\zroN$. One can then show by induction that for all $s<t$ we have
$\bG^{\bullet(t)}\br^{(s)}=\zroM$
and $(\bG^{\bullet(t)})^\st\bc^{(s)}=\zroN$. This implies, for all $\ell\le t$,
	\begin{align*}
	\bigg(
	\Gprime-\sum_{s \leq t}\bGam^{\bullet(s)}
	\bigg) \br^{(\ell)}
	\stackrel{\eqref{e:G.prime.recursion}}{=}
	\bG^{\bullet(t+1)}\br^{(\ell)}
	&=\zroM\in\R^M\,,\\
	\bigg(
	\Gprime-\sum_{s \leq t}\bGam^{\bullet(s)}
	\bigg)^\st \bc^{(\ell)}
	\stackrel{\eqref{e:G.prime.recursion}}{=}
	(\bG^{\bullet(t+1)})^\st
		\bc^{(\ell)}
	&=\zroN\in\R^N\,.
	\end{align*}
It follows from this that $\bG^{(t+1)}$ is orthogonal to $V_\star(t)$.
This verifies \eqref{e:proj.G.prime}, since we see that the right-hand side of \eqref{e:proj.G.prime} lies in $V_\star(t)$. 
\end{rmk}

\begin{proof}[\hypertarget{proof:l.resampling}{Proof of Lemma~\ref{l:resampling}}] Recall that, for notational convenience, the roles of $\bG$ and $\Gprime$ in this section are switched from the statement of Lemma~\ref{l:resampling}. Thus, for the purposes of the proof, we use $\bG$ for the AMP iteration \eqref{e:general.AMP.col} and \eqref{e:general.AMP.row}, and this defines $\FF_\star(t)$ as in \eqref{e:FF.star}. We also let $\ROW$ and $\COL$ be as in Definition~\ref{d:R.C.events}, but with $\bG$ and $\Gprime$ switched.
The $\sigma$-field $\FF(t)$ from
\eqref{e:tap.condition} is very closely related to $\FF_\star(t-1)$, but is not exactly the same: indeed, we can see from the proof of Lemma~\ref{l:FF.star} that
	\[
	\FF(t)=\sigma\bigg(
	\FF_\star(t-1), \bbx^{(t)}
	\bigg)
	=\sigma\bigg(
	\FF_\star(t-1),
		\bG\bmag^{(t)},\bn^{(t+1)}
	\bigg)\,.
	\]
By a similar (but simpler) argument as in Corollary~\ref{c:FF.star}, we see that
	\[
	\E\Big(\bG\,\Big|\,\FF(t)\Big)
	=\sum_{s\le t-1}\bGam^{(s)}
		+\br^{(t)}(\bbx^{(t)})^\st
	=\proj_{\ROW\COL}(\bG)
	= \bG_{\ROW\COL}\,,
	\]
where $\proj_{\ROW\COL}$ denotes orthogonal projection onto $V_{\ROW\COL}$ as in Definition~\ref{d:row.col},
except that $V_{\ROW\COL}$ here is defined for $\bG$ rather than $\Gprime$.
 Conditional on
$\FF_\star(t-1)$, the random matrix $\bG-\proj_{\ROW\COL}(\bG)$ is distributed as a standard gaussian
element of the $\FF_\star(t-1)$-measurable
vector space $(V_{\ROW\COL})^\perp$, and is conditionally independent of $\FF(t)$. Therefore
	\[ 
	\E\Big( f(\bG)\,\Big|\,\FF(t)\Big)
	= \E\bigg[ f\Big( \bG_{\ROW\COL}
		+ \Big(\Gprime-\proj_{\ROW\COL}(\Gprime)\Big)
		 \Big)\,\bigg|\,\FF(t)\bigg]
	=\E\Big( f(\Gprime)\,\Big|\,\ROW,\COL,
		\bG_{\ROW\COL}\Big)\,.
	\]
This concludes the proof.
\end{proof}
\subsection{AMP iterates at $t=2$ and $t=3$}\label{ss:amp.heuristic.one}

Returning to the AMP iteration 
\eqref{e:general.AMP.col} and \eqref{e:general.AMP.row}
 we have (cf.\ \eqref{e:gs.H.y} and \eqref{e:gs.h.x})
	\begin{align}\nonumber
	\bmag^{(2)} 
	&\equiv t(\bH^{(2)})
	= t\bigg( \f{ \bG^\st \bn^{(1)}}{N^{1/2}} \bigg)
	= t\bigg(
	\psi^{1/2} \bG^\st \f{\bn^{(1)}}{(N\psi)^{1/2}}
	\bigg)
	=t(\psi^{1/2} \bG^\st \bc^{(1)})
	=t(\psi^{1/2} \bby^{(1)})\,,\\
	\bn^{(2)}
	&\equiv F(\bh^{(2)})
	=F\bigg( \f{\bG\bmag^{(1)}}{N^{1/2}}
	\bigg)
	=F\bigg(
	q^{1/2} \bG \f{\bmag^{(1)}}{(Nq)^{1/2}}
	\bigg)
	=F(q^{1/2} \bG\br^{(1)})
	=F(q^{1/2} \bbx^{(1)})\,.
	\label{e:step.two}
	\end{align}
It follows using \eqref{e:q.psi.recurs} that $\|\bmag^{(2)}\|^2 \simeq Nq$,
$\|\bn^{(2)}\|^2 \simeq N\psi$, and moreover
(cf.\ \eqref{e:intro.lm.gm.one})
	\begin{align}
	\nonumber\f{(\bmag^{(2)},\bmag^{(1)})}{Nq}
	& \simeq
		\bigg(\f1q\bigg)^{1/2} \E t(\psi^{1/2} Z)
	\equiv \lm_1\equiv\rho_1\,,\\
	\f{(\bn^{(2)},\bn^{(1)})}{N\psi}
	&\simeq\bigg(
	\f{\alpha}{\psi}\bigg)^{1/2} \E F(q^{1/2}Z)
	\equiv\gm_1
	\equiv\mu_1\,.
	\label{e:product.with.first.iterate}
	\end{align}
Therefore in the Gram--Schmidt orthogonalization we have
	\begin{align}\nonumber
	\br^{(2)}
	&= \f{(\bmag^{(2)})^\perp}{\|(\bmag^{(2)})^\perp\|}
	\simeq \f{\bmag^{(2)} - \lm_1\bmag^{(1)}}
		{[Nq(1-(\lm_1)^2)]^{1/2}}\,,\\
	\bc^{(2)}
	&= \f{(\bn^{(2)})^\perp}{\|(\bn^{(2)})^\perp\|}
	\simeq \f{\bn^{(2)} - \gm_1 \bn^{(1)}}
		{[N\psi(1-(\gm_1)^2)]^{1/2}}\,.
	\label{e:gs.two.r.c}
	\end{align}
We can express the $\bmag$, $\bn$ vectors in terms of the $\br$, $\bc$ vectors as
	\begin{align}\nonumber
	\f{\bmag^{(2)}}{(Nq)^{1/2}}
	&\simeq \lm_1 \br^{(1)}
		+\Big(1-(\lm_1)^2\Big)^{1/2}
		\br^{(2)}\,,\\
	\f{\bn^{(2)}}{(N\psi)^{1/2}}
	&\simeq \gm_1\bc^{(1)} 
		+\Big( 1-(\gm_1)^2\Big)^{1/2}
		\bc^{(2)}\,.
	\label{e:gs.two.m.n}
	\end{align}
At the next step of the AMP iteration we have
(cf.\ \eqref{e:step.two})
	\begin{align}\nonumber
	\bmag^{(3)}
	&\equiv t(\bH^{(3)})
	= t\bigg( \f{\bG^\st\bn^{(2)}}{N^{1/2}}
		-\beta \bmag^{(1)}\bigg)
	= t\bigg( \f{(\bG^{(2)})^\st
		(\bn^{(2)})^\perp}{N^{1/2}}
	+ \f{(\bGam^{(1)}) ^\st\bn^{(2)}}{N^{1/2}}
		-\beta \bmag^{(1)}\bigg)\,,\\
	\bn^{(3)}
	&\equiv t(\bh^{(3)})
	=F\bigg(\f{\bG\bmag^{(2)}}{N^{1/2}}
		-\acute{\beta}\bn^{(1)}\bigg)
	= F\bigg( \f{\bG^{(2)}(\bmag^{(2)})^\perp}{N^{1/2}}
	+ \f{ \bGam^{(1)} \bmag^{(2)}}{N^{1/2}}
		-\acute{\beta}\bn^{(1)}
	\bigg)\,.
	\label{e:step.three}
	\end{align}
In order to evaluate $(\bGam^{(1)})^\st \bn^{(2)}/N^{1/2}$, we calculate
	\begin{align*}
	\bigg\{
	\f{\br\br^\st \bG^\st}{N^{1/2}}
	\bigg\}^{(1)}\bn^{(2)}
	&= \f{\bmag^{(1)}}{Nq^{1/2} }
		(\bbx^{(1)},\bn^{(2)})
	= \f{\bmag^{(1)}}{Nq^{1/2} }
		(\bbx^{(1)}, F(q^{1/2}\bbx^{(1)}) )
	\simeq \beta\bmag^{(1)}\,,\\
	\bigg\{
	\f{\bG^\st\bc\bc^\st}{N^{1/2}}
	\bigg\}^{(1)}\bn^{(2)}
	&=\f{ \bby^{(1)}}{N^{1/2}}
		\f{(\bn^{(1)},\bn^{(2)})}{(N\psi)^{1/2}}
	\simeq \psi^{1/2}\gm_1\bby^{(1)}\,,\\
	\bigg\{
	\f{(\bc^\st\bG\br)\br\bc^\st }{N^{1/2}}
	\bigg\}^{(1)}\bn^{(2)}
	&= \f1{N^{1/2}}
	\f{\oneM^\st\bG\oneN}{N\alpha^{1/2}}
	\f{(\bn^{(1)},\bn^{(2)})}{N(q\psi)^{1/2}}
	\bmag^{(1)}
	= O\bigg(\f1{N^{1/2}}\bigg)\bmag^{(1)}\,.
	\end{align*}
In order to evaluate $\bGam^{(1)} \bmag^{(2)}/N^{1/2}$, we calculate
	\begin{align*}
	\bigg\{
	\f{\bG\br\br^\st}{N^{1/2}}
	\bigg\}^{(1)}
	\bmag^{(2)}
	&= \f{\bbx^{(1)} }{N^{1/2}} 
		\f{(\bmag^{(1)},\bmag^{(2)})}{(Nq)^{1/2}}
	\simeq  q^{1/2}\lm_1 
		\bbx^{(1)}\,,\\
	\bigg\{
	\f{\bc\bc^\st\bG}{N^{1/2}}
		\bigg\}^{(1)} \bmag^{(2)}
	&=\f{\bn^{(1)}}{N\psi^{1/2}}
		(\bby^{(1)},\bmag^{(2)})
	= \f{\bn^{(1)}}{N\psi^{1/2}}
		(\bby^{(1)}, t(\psi^{1/2} \bby^{(1)}) )
	\simeq \acute{\beta} \bn^{(1)}\,, \\
	\bigg\{
	\f{(\bc^\st\bG\br)\bc\br^\st}{N^{1/2}}
		\bigg\}^{(1)}
		\bmag^{(2)}
	&= \f1{N^{1/2}}
	\f{\oneM^\st \bG \oneN}{N\alpha^{1/2}}
	\f{(\bmag^{(1)},\bmag^{(2)})}{N(q\psi)^{1/2}}
	\bn^{(1)} 
	= O\bigg(\f1{N^{1/2}}\bigg) \bn^{(1)}\,.
	\end{align*}
Substituting this back into \eqref{e:step.three} gives
the decomposition (cf.\ 
\eqref{e:gs.H.y}, \eqref{e:gs.h.x}, and
\eqref{e:step.two})
	\begin{align}\nonumber
	\bmag^{(3)}
	&\equiv t(\bH^{(3)})
	\simeq t\bigg( 
	\psi^{1/2} \bigg\{ 
	\gm_1 \bby^{(1)}
	+\Big(1-(\gm_1)^2\Big)^{1/2}
	\bby^{(2)}
	\bigg\} \bigg)\,,\\
	\bn^{(3)}
	&\equiv F(\bh^{(3)})
	\simeq F\bigg(
	q^{1/2} \bigg\{
	\lm_1 \bbx^{(1)}+
	\Big(1-(\lm_1)^2\Big)^{1/2}
	\bbx^{(2)}
	\bigg\} \bigg)\,.\label{e:decomp.three}
	\end{align}
It follows that 
\eqref{e:product.with.first.iterate} continues to hold (approximately) with $\bmag^{(3)}$, $\bn^{(3)}$
 in place of $\bmag^{(2)}$, $\bn^{(2)}$.
We also see by combining \eqref{e:step.two} with \eqref{e:decomp.three} that
	\begin{align}\nonumber
	\f{(\bH^{(2)},\bH^{(3)})}{N\psi}
	\simeq \f1N
	\bigg(
	\bby^{(1)},
	\gm_1\bby^{(1)}
		+\Big(1-(\gm_1)^2\Big)\bby^{(2)}
	\bigg)
	&\simeq \gm_1 \equiv \mu_1\,,\\
	\f{(\bh^{(2)},\bh^{(3)})}{M q}
	\simeq
	\f1M \bigg(
	\bbx^{(1)},
	\lm_1\bbx^{(1)}
	+\Big(1-(\lm_1)^2\Big)\bbx^{(2)}
	\bigg)
	&\simeq \lm_1 \equiv \rho_1\,,
	\label{e:rho.mu.one}
	\end{align}
from which we obtain (cf.\ \eqref{e:def.rho.mu})
 	\begin{align}\nonumber
	\f{(\bmag^{(2)},\bmag^{(3)})}{Nq}
	\simeq
	\f1{q} \E\Bigg[ t\bigg( \psi^{1/2}\bigg\{
	\gm_1 Z
	+\Big(1-(\gm_1)^2\Big)^{1/2}
	\xi
	\bigg\} \bigg)
	t(\psi^{1/2} Z) \Bigg]
	&\equiv \rho(\gm_1)  = \rho(\mu_1) \equiv
	\rho_2\,,\\
	\f{(\bn^{(2)},\bn^{(3)})}{N\psi}
	\simeq
	\f{\alpha}{\psi} \E\Bigg[
	F\bigg(q^{1/2}\bigg\{
	\lm_1 Z+
	\Big(1-(\lm_1)^2\Big)^{1/2}
	\xi
	\bigg\}
	\bigg) F(q^{1/2}Z)
	\Bigg]
	&\equiv \mu(\lm_1) = \mu(\rho_1) \equiv
	\mu_2\,.
	\label{e:rho.mu.two}
	\end{align}
It follows that
(cf.\ \eqref{e:intro.lm.gm.rec})
 	\begin{align}\nonumber
	\bigg( \f{\bmag^{(3)}}{(Nq)^{1/2}},\br^{(2)} \bigg)
	&\simeq \bigg(
		\f{\bmag^{(3)}}{(Nq)^{1/2}}, \f{ 
		\bmag^{(2)}	-\lm_1\bmag^{(1)}}
		{[Nq(1-(\lm_1)^2)]^{1/2}} \bigg)
	\simeq \f{\rho_2 -(\lm_1)^2}
		{[1-(\lm_1)^2]^{1/2}}
	\equiv \lm_2\,, \\
	\bigg( \f{\bn^{(3)} }{(N\psi)^{1/2}},\bc^{(2)}\bigg)
	&\simeq
	\bigg( \f{\bn^{(3)} }{(N\psi)^{1/2}},
	\f{\bn^{(2)} - \gm_1 \bn^{(1)}}
		{ [N\psi(1-(\gm_1)^2)]^{1/2}
		}\bigg)
	\simeq \f{\mu_2 - (\gm_1)^2}
		{[1-(\gm_1)^2]^{1/2}}
	\equiv \gm_2\,.
	\label{e:product.with.second.iterate}
	\end{align}
Then in the Gram--Schmidt orthogonalization we have (cf.\ \eqref{e:gs.two.r.c})
	\begin{align}\nonumber
	\br^{(3)}
	&\simeq
	\f{\bmag^{(3)}
	- (Nq)^{1/2}(\lm_1 \br^{(1)}
	+ \lm_2\br^{(2)})}
	{[Nq( 1-(\lm_1)^2-(\lm_2)^2 )]^{1/2}}\\
	\bc^{(3)}
	&\simeq
	\f{\bn^{(3)} - (N\psi)^{1/2}
		( \gm_1\bc^{(1)} + \gm_2 \bc^{(2)} ) }
		{
		[N\psi(( 1-(\gm_1)^2-(\gm_2)^2 )]^{1/2}
		}\,.
	\label{e:gs.three.r.c}
	\end{align}
We can express the $\bmag$, $\bn$ vectors in terms of the $\br$, $\bc$ vectors as (cf.\
\eqref{e:gs.m.r.EXACT}, 
\eqref{e:gs.n.c.EXACT}, and \eqref{e:gs.two.m.n})
	\begin{align}\nonumber
	\f{\bmag^{(3)}}{(Nq)^{1/2}}
	&\simeq \lm_1\br^{(1)}
		+\lm_2\br^{(2)}
		+\Big( 1-(\lm_1)^2-(\lm_2)^2\Big)^{1/2}
		\br^{(3)}\,,\\
	\f{\bn^{(3)}}{(N\psi)^{1/2}}
	&\simeq \gm_1\bc^{(1)}+\gm_2\bc^{(2)}
		+\Big(1-(\gm_1)^2-(\gm_2)^2\Big)^{1/2}
		\bc^{(3)}
	\,.\label{e:gs.three.m.n}
	\end{align}

\subsection{AMP iterates at $t=4$}\label{ss:amp.heuristic.conclusion}

At the next step of the AMP iteration we have
(cf.\ \eqref{e:step.three})
	\begin{align}\nonumber
	\bmag^{(4)}
	&\equiv t(\bH^{(4)})
	=t\bigg( \f{(\bG^{(3)})^\st
		(\bn^{(3)})^\perp}{N^{1/2}}
	+ \f{(\bGam^{(2)}) ^\st\bn^{(3)}}{N^{1/2}}
	+ \f{(\bGam^{(1)}) ^\st\bn^{(3)}}{N^{1/2}}
		-\beta \bmag^{(2)}\bigg)\,,\\
	\bn^{(4)}
	&\equiv F(\bh^{(4)})
	= F\bigg( \f{\bG^{(3)}(\bmag^{(3)})^\perp}{N^{1/2}}
	+ \f{ \bGam^{(2)} \bmag^{(3)}}{N^{1/2}}
	+ \f{ \bGam^{(1)} \bmag^{(3)}}{N^{1/2}}
		-\acute{\beta}\bn^{(2)}
	\bigg)\,.\label{e:step.four}
	\end{align}
For the purposes of evaluating $(\bGam^{(s)})\bn^{(3)}$ for $s=1,2$ we calculate
	\begin{align*}
	\bigg\{\f{\br\br^\st\bG}{N^{1/2}}\bigg\}^{(1)}
	\bn^{(3)}
	&=\f{\br^{(1)}}{N^{1/2}} (\bbx^{(1)},\bn^{(3)})
	\simeq \f{\br^{(1)}}{N^{1/2}}
	N\alpha \E 
	\bigg[ 
	\bigg(\lm_1Z 
		+ \Big( 1-(\lm_1)^2\Big)^{1/2}
		\xi\bigg)F(q^{1/2}Z) 
	\bigg]
	= (Nq)^{1/2} 
	\beta \lm_1
	\br^{(1)}\,,\\
	\bigg\{\f{\bG\bc\bc^\st}{N^{1/2}}\bigg\}^{(1)}\bn^{(3)}
	&= \f{\bby^{(1)}}{N^{1/2}}
		(\bc^{(1)},\bn^{(3)})
	= \psi^{1/2}\gm_1 \bby^{(1)}\,.
	\end{align*}
Substituting back into \eqref{e:step.four} gives
(cf.\
\eqref{e:gs.H.y}, \eqref{e:gs.h.x}, 
\eqref{e:step.two}, and
 \eqref{e:decomp.three})
	\begin{align}\nonumber
	\bmag^{(4)}
	&\equiv t(\bH^{(4)})
	\simeq t\bigg(\psi^{1/2}\bigg\{
	\gm_1\bby^{(1)}+\gm_2\bby^{(2)}
	+ \Big(1-(\gm_1)^2-(\gm_2)^2\Big)^{1/2}\bby^{(3)}
	\bigg\}
	\bigg)\,,\\
	\bn^{(4)}
	&\equiv F(\bh^{(4)})
	\simeq F\bigg(
	q^{1/2}\bigg\{
	\lm_1\bbx^{(1)} + \lm_2\bbx^{(2)}
	+ \Big(1-(\lm_1)^2-(\lm_2)^2\Big)^{1/2}\bbx^{(3)}
	\bigg\} \bigg)\,.
	\label{e:decomp.four}
	\end{align}
It follows that 
\eqref{e:product.with.first.iterate} continues to hold (approximately) with $\bmag^{(4)}$, $\bn^{(4)}$ in place of $\bmag^{(2)}$, $\bn^{(2)}$. 
Likewise, \eqref{e:product.with.second.iterate}
continues to approximately hold with
$\bmag^{(4)}$, $\bn^{(4)}$ in place of $\bmag^{(3)}$, $\bn^{(3)}$. We also have  (cf.\ \eqref{e:rho.mu.one})
	\begin{align}\nonumber
	\f{(\bH^{(3)},\bH^{(4)})}{N\psi}
	\simeq
	(\gm_1)^2 + \gm_2\Big(1-(\gm_1)^2\Big)^{1/2}
	&\stackrel{\eqref{e:product.with.second.iterate}}{=}
	\mu_2\,,\\
	\f{(\bh^{(3)},\bh^{(4)})}{Mq}
	\simeq
	(\lm_1)^2 + \lm_2\Big(1-(\lm_1)^2\Big)^{1/2}
	&\stackrel{\eqref{e:product.with.second.iterate}}{=}
	\rho_2\,,
	\label{e:rho.mu.one.again}
	\end{align}
from which it we obtain (cf.\ 
\eqref{e:def.rho.mu} and 
\eqref{e:rho.mu.two})
	\beq\label{e:rho.mu.three}
	\f{(\bmag^{(4)},\bmag^{(3)})}{Nq}
	\simeq \rho(\mu_2)\equiv\rho_3\,,\quad
	\f{(\bn^{(4)},\bn^{(3)})}{N\psi}
	\simeq \mu(\rho_2)\equiv\mu_3\,.
	\eeq
It then follows that (cf.\ 
\eqref{e:intro.lm.gm.rec} and
\eqref{e:product.with.second.iterate})
	\begin{align}\nonumber
	\bigg(\f{\bmag^{(4)}}{(Nq)^{1/2}},\br^{(3)}
	\bigg)
	&\simeq\bigg(\f{\bmag^{(4)}}{(Nq)^{1/2}},
	\f{
	\bmag^{(3)}/(Nq)^{1/2}
		-\lm_1\br^{(1)}-\lm_2\br^{(2)}}{
	[1-(\lm_1)^2-(\lm_2)^2 ]^{1/2}}
		\bigg)
	\simeq\f{ \rho_3-(\lm_1)^2-(\lm_2)^2}
		{[1-(\lm_1)^2-(\lm_2)^2 ]^{1/2}}
		\equiv\lm_3\,,\\
	\bigg( \f{\bn^{(4)} }{(N\psi)^{1/2}},\bc^{(3)}\bigg)
	&\simeq\bigg( \f{\bn^{(4)} }{(N\psi)^{1/2}},
	\f{\bn^{(3)}/(N\psi)^{1/2}
		- \gm_1\bc^{(1)}- \gm_2\bc^{(2)}
		}
	{[1-(\gm_1)^2-(\gm_2)^2]^{1/2}}
	\bigg)
	\simeq\f{ \mu_3 -(\gm_1)^2-(\gm_2)^2}
	{[1-(\gm_1)^2-(\gm_2)^2]^{1/2}}
	\equiv \gm_3\,.
	\label{e:product.with.third.iterate}
	\end{align}
In summary, 
using the notation \eqref{e:sum.of.squares},
we have
(cf.\
\eqref{e:gs.H.y},
\eqref{e:gs.h.x},
\eqref{e:step.two},
 \eqref{e:decomp.three}, and \eqref{e:decomp.four})
	\begin{align}
	\label{e:amp.decomp.H.y}
	\bmag^{(t+1)}
	&\equiv t(\bH^{(t+1)})
	\simeq t\bigg(
	\psi^{1/2}
	\bigg\{
	\gm_1\bby^{(1)}+\ldots+\gm_{t-1} \bby^{(t-1)}
	+(1-\Gamma_{t-1})^{1/2}\bby^{(t)}
	\bigg\}
	\bigg)\,,\\
	\bn^{(t+1)}
	&\equiv F(\bh^{(t+1)})
	\simeq F\bigg( q^{1/2}\bigg\{
	\lm_1\bbx^{(1)}+\ldots+\lm_{t-1}\bbx^{(t-1)}
	+(1-\Lambda_{t-1})^{1/2}\bbx^{(t)}
	\bigg\} \bigg)\,.
	\label{e:amp.decomp.h.x}
	\end{align}
where the coefficients are defined recursively:
we start with $\lm_1\equiv\rho_1$ and $\gm_1\equiv\mu_1$ as in \eqref{e:product.with.first.iterate}
(cf.\ \eqref{e:intro.lm.gm.one}).
For $s\ge1$ we let 
$\rho_{s+1}\equiv\rho(\mu_s)$
and $\mu_{s+1}\equiv\mu(\rho_s)$
as in \eqref{e:rho.mu.two} and \eqref{e:rho.mu.three}
(cf.\ \eqref{e:def.rho.mu}). Then, as in \eqref{e:product.with.second.iterate} and \eqref{e:product.with.third.iterate} (cf.\ \eqref{e:intro.lm.gm.rec}), we can define recursively the constants
	\beq\label{e:lm.gm.recursion}
	\lm_s
	=\f{\rho_s-\Lambda_{s-1}}{(1-\Lambda_{s-1})^{1/2}}\,,
	\quad
	\gm_s=\f{\mu_s-\Gamma_{s-1}}
		{(1-\Gamma_{s-1})^{1/2}}\,.
	\eeq
We use these to define the matrices $\bGam$ and $\bLam$
as in \eqref{e:Gamma.matrix} and \eqref{e:Lambda.matrix}. Then \eqref{e:amp.decomp.H.y}
and \eqref{e:amp.decomp.h.x} 
can be rewritten as
\eqref{e:gs.H.y}
and \eqref{e:gs.h.x}.
The Gram--Schmidt orthogonalization 
\eqref{e:gs.two.m.n} and \eqref{e:gs.three.m.n}
then correspond (approximately) to \eqref{e:gs.m.r.EXACT}
and \eqref{e:gs.n.c.EXACT}.

\pagebreak

{\raggedright\bibliographystyle{alphaabbr}
\bibliography{prefs}}

\end{document}